\documentclass[a4paper]{article}

\usepackage{skript}

\renewcommand{\lecturer}{\href{https://olehmelnyk.xyz}{Oleh Melnyk} and \href{https://www.tu-chemnitz.de/mathematik/invimg/index.en.php}{Sebastian Neumayer}}
\renewcommand{\university}{Technische Universität Berlin \& Technische Universität Chemnitz}
\renewcommand{\lecture}{Data-driven Regularization \\ for \\ Inverse Problems}
\renewcommand{\lecSubtitle}{principles, algorithms and their convergence guarantees}
\renewcommand{\semester}{}

\renewcommand{\mythanks}{}

\renewcommand{\mysols}{%With solutions to the exercises. \\ 
}
\NewEnviron{Solution}{%
}%
\newcommand{\Reading}{\textbf{Textbooks:}\,\,}

\begin{document}
\selectlanguage{english}
\mytitlepage
%%%%%%%%%%%%%%%%%%%%%%%%%%%%%%%%%%%%%%%%%%%%%%%%%%%%%%%%%%%%%%%%%%%%%%%%%%%%%%%%%%%%%%%%%%%%%%%%%%%%%%%%%%%%%%%%%%%%%%%%%%%%%%%%%%%%%%%%%%%%%%%%%%%%%%%%%%%%%%%%%%%%%%%%%%%%%%%%%%%%%%%%%%%%%%%%
\pagenumbering{arabic}

\newgeometry{ 
            %headheight=13pt, 
            includehead,
            includefoot,
            top = 2cm,
            bottom = 2cm,
            left = 3cm,
            right = 4cm}
\edef\marginnotetextwidth{\the\textwidth} %  While some environments (e.g., those of package framed) change \textwidth, marginnote defines its own text width macro. If you change type area after \begin{document} you should add this
\setlength\headheight{25pt}

\section*{Foreword}

In many applications, we are interested in studying an object of interest $x$.
However, its properties are often not directly measurable.
Instead, we only observe $y = F(x)$, where the (known) function $F$ is referred to as the forward mapping.
To infer about $x$ from the observation $y$, we seek an (approximate) inverse mapping $x \approx F^{-1}(y)$.
This process is known as solving the \emph{inverse problem}.

In practice, the inverse mapping $F^{-1}$ may be ill-defined or discontinuous, which forces us to choose a certain solution from many possible ones based on prior assumptions on $x$.
Such assumptions about $x$ are typically incorporated into the inversion process through regularization.
Therefore, the first part of this course focuses on the classical regularization theory of inverse problems.
This part is primarily based on the notes by C.\ Clason \cite{clason2020regularization}.

In the second part, we transition to algorithmic approaches for solving inverse problems.
The foundation for this lies in the variational formulation of inverse problems, where reconstructions are obtained by minimizing specific functionals.
Here, we will explore key optimization principles, discuss their practical implementation, and examine how they contribute to robust reconstructions.
 
Today, classical methods are often modified to incorporate data-driven machine learning techniques, which improves the empirical performance.
However, this often comes at the cost of explainability or convergence guarantees.
Such properties are critical for example in medical applications, where one wants to be confident that reconstructed abnormalities are of physical nature and not phantoms resulting from the method. 

For these reasons, we discuss a number of machine learning methods for inverse problems, such as unrolling techniques, learned regularizers and plug-and-play methods, and explore their performance guarantees.
We refer to the review articles \cite{habring2024neural,mukherjee2023learned} and references therein for further reading beyond these notes.

\section{Introduction to Inverse Problems} \label{sec: inverse problems}

In this chapter, we review the theoretical background for this lecture. 
Moreover, we introduce what an inverse problem is and provide examples from practical applications.

%\datum{Lecture 1, 03.04.2024}

\Reading
\cite{EnglHankeNeubauer1996},\cite{Kirsch2021}, and \cite{ScherzerGrasmair2009}.

\subsection{Normed Spaces}
Let $\bb X$ be a $\mathbb K$-vector space.
A mapping $\norm{\cdot}\colon \bb X\to \R_+\coloneqq [0,\infty)$ is a \emphi{norm} if for all $x\in \bb X$ the following hold
\begin{enumerate}
	\item $\norm{\lambda x} = |\lambda| \norm{x}$ for all $\lambda \in \mathbb{K}$;
	\item $\norm{x+y} \leq \norm{x} + \norm{y}$ for all $y\in \bb X$;
	\item $\norm{x} = 0$ if and only if $x = 0$.
\end{enumerate}

If $\norm{\cdot}$ is a norm on $\bb X$, the pair $(\bb X,\norm{\cdot})$ is a \emphi{normed space}, and we denote this by $\norm{\cdot}_{\bb X}$.
We write $\norm{ \cdot }$ to simplify the notation unless it is necessary.
Throughout these notes, we mainly work with $\bb K = \R$ and only occasionally refer to the case $\bb K = \bb C$.

A sequence $\{x_n\}_{n\in \N}\subset \bb X$ \emphi{converges strongly} to $x\in \bb X$, denoted by $x_n\to x$, if
\begin{equation}
	\lim_{n\to\infty} \norm{x_n-x} = 0.
\end{equation}

\begin{example}\label{ex:funktan:norm}
	\begin{enumerate}
		\item On $\bb X = \R^N$, we can define the norms
		\begin{equation}
			\begin{aligned}
				\norm{x}_p &= \left(\sum_{i=1}^N |x_i|^p\right)^{1/p},\qquad1\leq p<\infty,\\
				\norm{x}_\infty &= \max_{i=1,\dots,N} |x_i|.
			\end{aligned}
		\end{equation}
		\item On $\bb X = \ell^p$ (the real-valued sequences with finite $p$-norm), we can define the norms
		\begin{equation}
			\begin{aligned}
				\norm{x}_p &= \left(\sum_{i=1}^\infty |x_i|^p\right)^{1/p},\qquad 1\leq p<\infty,\\
				\norm{x}_\infty &= \sup_{i=1,\dots,\infty} |x_i|.
			\end{aligned}
		\end{equation}
		\item Let  $\Omega \subset \R^d$ and $\mu$ be a measure.
		On $\bb X = L^p(\Omega, \mu)$ (the  $\mu$-measurable functions on $\Omega$ with finite $p$th moment) we can define the norms
		\begin{equation}
			\begin{aligned}
				\norm{x}_p &= \left(\int_\Omega |x(v)|^p\,d\mu(v) \right)^{1/p},\qquad 1\leq p<\infty,\\
				\norm{x}_\infty &= \mathop\mathrm{ess\,\sup}_{v\in \Omega} |x(v)|.
			\end{aligned}
		\end{equation}
		We write $L^p(\Omega)$ if $\mu$ is the Lebesgue measure on $\R^d$.
		
		\item For compact $\Omega \subset \R^d$, a norm on the continuous functions $\bb X = C(\Omega)$ is given by
		\begin{equation}
			\norm{x}_C = \sup_{v\in \Omega} |x(v)|.
		\end{equation}
		Similarly, a norm on the space $C^k(\Omega)$ of $k$ times continuously differentiable functions is defined by $\norm{x}_{C^k}=\sum_{j=0}^k\norm{x^{(j)}}_{C}$.
	\end{enumerate}
\end{example}

For a midpoint $x\in \bb X$ and a radius $r>0$, we define
\begin{itemize}
	\item the \emph{open ball} $U_r(x) \coloneqq \{ z\in \bb X \, : \, \norm{x-z} < r\}$;
	\item the \emph{closed ball} $B_r(x) \coloneqq \{z\in \bb X  \, : \,  \norm{x-z} \leq r \}$.
\end{itemize}
The closed ball $B_\bb X \coloneqq B_1(0)$ is also referred to as the \emph{unit ball}.
A set $U\subset \bb X$ is 
\begin{itemize}
	\item \emphi{open} if for all $x\in U$ there exists an $r>0$ such that $U_r(x)\subset U$.
	We denote by $U^\circ$ the \emphi{interior} of $U$, the largest open subset of $U$;
	\item \emphi{closed} if every convergent sequence $\{x_n\}_{n\in\N}\subset U$ satisifies $\lim_{n \to \infty} x_n \in U$.
	The union of $U$ with the set of all limit points of convergent sequences in $U$ (namely the accumulation points of $U$) is the \emph{closure} $\overline U$ of $U$. It is the smallest closed set containing $U$;
	\item \emphi{bounded} if $U \subset B_r(0)$ for some $r>0$;
	\item \emphi{compact} if every sequence $\{x_n\}_{n\in\N}\subset U$ contains a convergent subsequence $\{x_{n_k}\}_{k\in\N}$ with limit $x\in U$;
	\item \emphi{convex} if for all $x,y\in U$ and $\lambda\in [0,1]$ also $\lambda x + (1-\lambda)y\in U$.
	The definition of  a norm $\Vert \cdot \Vert$ implies that open and closed balls are convex.
	\item \emphi{dense} in $\bb X$ if for all $x\in \bb X$ there exists a sequence $\{x_n\}_{n\in\N}\subset U$ with $x_n\to x$.
	Obviously, $U$ is dense in $\overline U$.
\end{itemize}
In normed spaces, the complement of an open set is closed and vice versa.
In particular, the closed sets in the topological sense are exactly the sequentially closed ones from above.

A normed vector space $\bb X$ is \emph{complete} if every Cauchy sequence converges.
Such spaces $\bb X$ are called \emphi{Banach spaces}.
By $\overline{\bb X}$, we denote the \emph{completion} (with respect to the norm $\|\cdot\|_{\bb X}$) of $\bb X$.
All spaces in \cref{ex:funktan:norm} are complete.
Moreover, we call $\bb X$ \emphi{separable} if there exists a countable dense $U \subset \bb X$.
Finally, the unit ball $B_\bb X$ is compact \emph{if and only if} $\dim(\bb X) < \infty$, which is a fundamental result with many consequences for our investigation of inverse problems.
\subsection{Hilbert Spaces}

The mathematical analysis of inverse problems is easier in Hilbert spaces since we can rely on compactness arguments and the notions of orthogonality and bases.  
A real \emphi{Hilbert space} is a real Banach space with an \emph{inner product} $\inner{\cdot}{\cdot} \colon \bb X \times \bb X \to \R$ that satisfies
\begin{enumerate}
	\item $\inner{\alpha x+ y}{z} = \alpha\inner{x}{z} + \inner{y}{z}$ for all $x,y,z\in \bb X$ and $\alpha \in \R$;
	\item $\inner{x}{y} = \inner{y}{x}$ for all $x,y\in \bb X$;
	\item $\inner{x}{x} \geq 0$ for all $x\in \bb X$ with $\inner{x}{x}=0$ if and only if $x=0$.
\end{enumerate}
The induced norm $\norm{x}_\bb X \coloneqq  \sqrt{\inner{x}{x}_\bb X}$ satisfies the \emph{Cauchy--Schwarz inequality}
\begin{equation}\label{eq:CauchySchwarz}
	|\inner{x}{y}_\bb X| \leq \norm{x}_\bb X\norm{y}_\bb X.
\end{equation}

\begin{example}
	In \Cref{ex:funktan:norm}\, \circled{1}, \circled{2}, and \circled{3} for $p=2$, the underlying space is a Hilbert space with inner product
	\begin{enumerate}
		\item  $\displaystyle \inner{x}{y}_{\R^N} = \sum_{i=1}^N x_i y_i$,
		\item  $\displaystyle \inner{x}{y}_{\ell^2} = \sum_{i=1}^\infty x_i y_i$,
		\item $\displaystyle \inner{x}{y}_{L^2(\Omega, \mu)} = \int_\Omega x(v)y(v)\,\diff \mu(v)$.
	\end{enumerate}
\end{example}

A sequence $\{x_n\}_{n\in\N}\subset \bb X$ \emphi{converges weakly} to $x\in \bb X$, denoted by $x_n \wkto x$, if
\begin{equation}\label{eq:WeakConv}
	\inner{z}{x_n}_\bb X\to  \inner{z}{x}_\bb X \qquad\text{for all }z\in \bb X.
\end{equation}
One can show that weakly convergent sequences are always bounded.
For the case $\bb X = \R^N$, choosing $z=e_i$, the $i$th unit vector, gives the componentwise convergence implying that weak and strong convergence coincide if $\dim(\bb X) < \infty$.
By \eqref{eq:CauchySchwarz}, strong convergence always implies weak convergence.
In general, the reverse is only true if we additionally have $\norm{x_n}_\bb X \to \norm{x}_\bb X$.
We have the following important compactness result.

\begin{thm}{Banach--Alaoglu}{Banach_Alaoglu}
	Every bounded sequence in a Hilbert space has a weakly convergent subsequence.
\end{thm}

Two elements $x,y\in \bb X$ are \emphi{orthogonal} if $\inner{x}{y}_\bb X = 0$.
For a set $U\subset \bb X$, we define the  \emph{orthogonal complement} as
\begin{equation}
	U^\bot \coloneqq \bigl\{ x\in \bb X  : \inner{x}{u}_\bb X = 0 \text{ for all } u\in U \bigr\}.
\end{equation}
Since $x \mapsto \langle x, u \rangle$ is continuous, $U^\bot \subset \bb X$ is a closed subspace.
Furthermore, $\bb X^\bot = \{0\}$ and $U\subset (U^\bot)^\bot$.
If $U$ is a subspace, it holds that $\overline U=(U^\bot)^\bot$ (and hence $\{0\}^\bot = \bb X$).

Given a nonempty closed subspace $U \subset \bb X$, every $x\in \bb X$ can be represented uniquely as
\begin{equation}
	x = u +u_\bot\qquad\text{with}\qquad u\in U,\ u_\bot\in U^\bot.
\end{equation}
The  associated \emphi{orthogonal projection} $P_U$ with $x\mapsto u$ is a linear and continuous mapping.
Equivalently, we can characterize $P_U x$ as the distance minimizing element in $U$, namely that $\norm{x-P_U x}_\bb X = \min_{u\in U} \norm{x-u}_\bb X$.
By construction, it holds that $P_U^2 = P_U$ and $I - P_U = P_{U^\bot}$.
Moreover, we have  the \emph{Pythagorean} theorem $\|x\|^2_\bb X = \|u\|^2_\bb X + \|u_\perp\|^2_\bb X$. 

A set $U\subset \bb X$ with pairwise orthogonal elements is called \emph{orthogonal system}. 
If additionally $\Vert x \Vert = 1$ for all $x \in U$, then $U$ is \emph{orthonormal}.
Note that we do not require that $U$ is countable.
An orthonormal system with 
\begin{equation}
\bb X = \overline{\mathrm{span} \{U\}} \coloneqq \overline{ \left\{ \sum_{i = 1}^N \alpha_i u_i \, \colon \, N \in \N, \ \alpha_i \in \R, \ u_i \in U \right\} }  
\end{equation}
is called complete and we refer to such systems as \emphi{orthonormal basis} (ONB).
An orthonormal system $U$  is complete if and only if $\langle x, u \rangle = 0$ for all $u \in U$ implies $x=0$.
Interestingly, every Hilbert space $\bb X$ contains an ONB $U$.
A countable ONB exists if and only if $\bb X$ is separable.

Every orthonormal system $U\subset \bb X$ satisfies the \emphi{Bessel inequality}
\begin{equation}\label{eq:funktan:bessel}
	\sum_{u\in U} |\inner{x}{u}_\bb X|^2 \leq \norm{x}_\bb X^2\qquad\text{for all } x \in \bb X,
\end{equation}
where at most countably many terms $ |\inner{x}{u}_\bb X|$ are nonzero.
Hence, the summation is well-defined.
Moreover, equality in \eqref{eq:funktan:bessel} holds if $U$ is an ONB.
Then, we also have that 
\begin{equation}\label{eq:funktan:ons_reihe}
	x  = \sum_{u\in U} \inner{x}{u}_\bb X u \qquad\text{for all } x\in \bb X.
\end{equation}
Vice versa, if equality holds in \eqref{eq:funktan:bessel} and $U$ is countable, then $U$ is an ONB.

\begin{example}
	For $\bb X=L^2((0,1))$, an ONB $\{u_n\}_{n\in\bb Z}$ is given by
	\begin{equation}
		u_n(x) = \begin{cases}
			\sqrt{2}\sin(2 \pi n x), & n>0,\\
			\sqrt{2}\cos(2 \pi n x), & n<0,\\
			1, & n=0.
		\end{cases}
	\end{equation}
	This ONB is known as the \emph{Fourier} basis.
	It consists of simple wave functions and plays an important role in harmonic analysis and differential equations. 
\end{example}

\subsubsection{Exercises}
\begin{exerciseB}
	Construct a bounded sequence in some Hilbert space $\bb X$ that has no strongly convergent subsequence.
	Then, find a subsequence and that is weakly convergent. 
\end{exerciseB}

\begin{Solution}
	Consider the standard orthonormal basis \( (e_n)_{n \in \mathbb{N}} \) in the Hilbert space \( \ell^2(\mathbb{N}) \), where each \( e_n \) has a 1 in the \( n \)-th position and 0 elsewhere.
	
	The sequence \( (e_n) \) is bounded, since \( \|e_n\| = 1 \) for all \( n \). However, it has no strongly convergent subsequence, because
	\[
	\|e_n - e_m\| = \sqrt{2} \quad \text{for } n \ne m,
	\]
	so no subsequence is Cauchy in norm.
	
	Nevertheless, \( e_n \rightharpoonup 0 \) weakly. Indeed, for any \( x = (x_k) \in \ell^2(\mathbb{N}) \), we have
	\[
	\langle e_n, x \rangle = x_n \to 0 \quad \text{as } n \to \infty,
	\]
	since \( x_n \to 0 \) for square-summable sequences.
\end{Solution}

\subsection{Examples of Forward Mappings}\label{sec:ForwardMappings}

Let $(\bb X,\norm{\cdot}_\bb X)$ and $(\bb Y,\norm{\cdot}_\bb Y)$ be normed spaces, $U\subset \bb X$, and $F\colon U\to \bb Y$ be a mapping. 

\begin{defn}[Forward mapping and inverse problem] \label{defn:forward mapping}
	We refer to $F\colon U\to \bb Y$ as the \emph{forward mapping}.
	The \emph{measurements} $y \in \bb Y$ generated by some $x \in U$ are \emph{noiseless} if $F(x) = y$ and \emph{noisy} if $y \approx F(x)$. The \emph{inverse problem} consists of finding $x$ corresponding to the measurements $y$.    
\end{defn}

\begin{figure}
\begin{center}
	\subfloat[Blur from wrong focusing distance. \label{fig:blur_photography}]{
		\includegraphics[width= 200pt, height = 100pt]{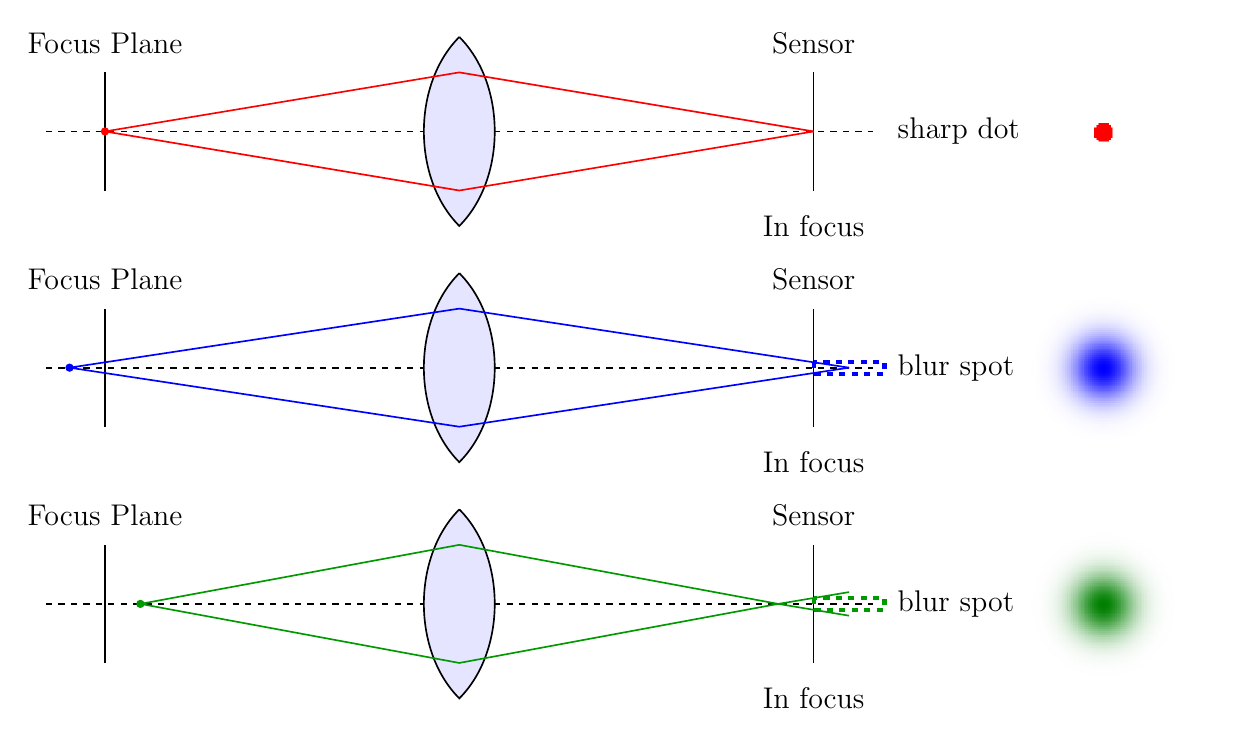}% 
	} \hfill
	\subfloat[Transmitter-receiver system.  \label{fig:blur_signal}]{
		\includegraphics[width= 175pt, height = 100pt]{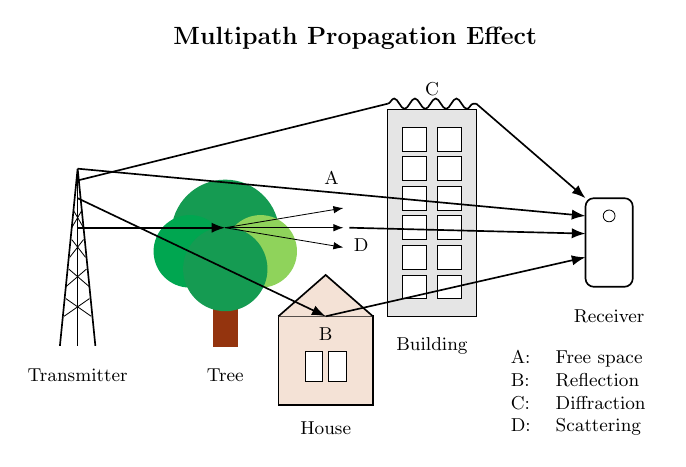}%
	}
	\caption{Examples of blurring effects in imaging and signal processing.}
\end{center}
\end{figure}

\begin{example}[Deconvolution]\label{ex: convolution}
	A fundamental problem in image and signal processing is \emph{deblurring} or \emphi{deconvolution}, which typically occurs as result of a nonideal measurement process.
	Figure~\ref{fig:blur_photography} shows examples of blur in photography due to a bad adjustment of the focus plane.
	Another instance are transmitter-receiver systems as depicted in Figure~\ref{fig:blur_signal}.
	Here, a signal is sent from an antenna in all directions towards a cellphone.
	Since the environment contains obstacles, the signal takes paths of different length.
	This results in delays and varying strengths for the receiver. 
	Other application examples are atmospheric perturbations for space telescopes, motion blur in videos, and tracer kinetics in medicine.
	      
	For a mathematical description, let $x \in L^p(\bb R)$, $1 \le p \le \infty$, represent the image and let $g \in L^1(\bb R)$ be the convolution kernel.
	Then, the blurred $x$ corresponds to the convolution 
	\begin{equation}\label{eq:conv}
	F[x](t) = x * g (t) \coloneqq \int_{\bb R} x(s) g(t-s) \diff s, \quad t \in \bb R. 
	\end{equation}
	In applications, signals and images are represented by vectors and matrices.
	Hence, \eqref{eq:conv} has to be discretized.
	Typically, either circular or noncircular models are used.
	The first one assumes that $x \in \R^N$ and $g \in \R^N$ are periodic, yielding
	\begin{equation}\label{eq: convolution}
		F(x)_k = (x *_N g)_k \coloneqq \sum_{j = 0}^{N-1} x_j g_{(k-j)_N}, \quad k = 0, \ldots, N-1, 
	\end{equation}
	where $(k)_N = k \MOD N$, for all $k \in \bb Z$.
  	The latter one appends $g$ with zeros, i.e., $g_k = 0$, $k \notin \{ 0, \ldots, N-1 \}$, leading to
  	\begin{equation}
  		F(x)_k = (x * g)_k \coloneqq \sum_{j = 0}^{N-1} x_j g_{k-j}, \quad k = -N + 1, \ldots, N-1.
  	\end{equation}
%  	We explore the differences between these definitions in the exercises.
\end{example}

\begin{figure}
	\begin{center}
		\subfloat[Illustration of a tomographic experiment. \label{fig: tomography}]{
			\includegraphics[width= 200pt, height = 100pt]{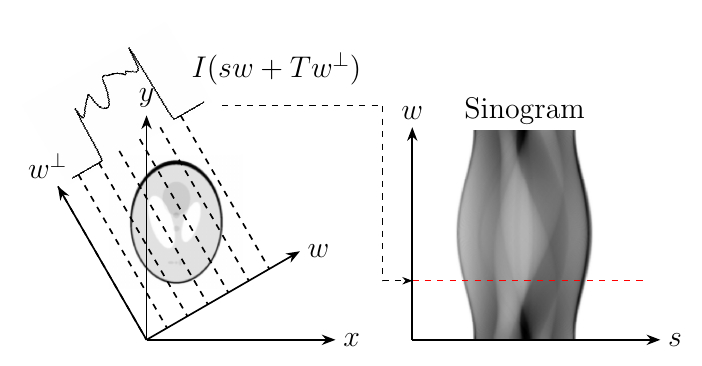}% 
		} \hfill
		\subfloat[Reconstruction examples \cite{gadelha2019shape}. \label{fig: tomography reconstruction}]{
			\includegraphics[width= 175pt, height = 100pt]{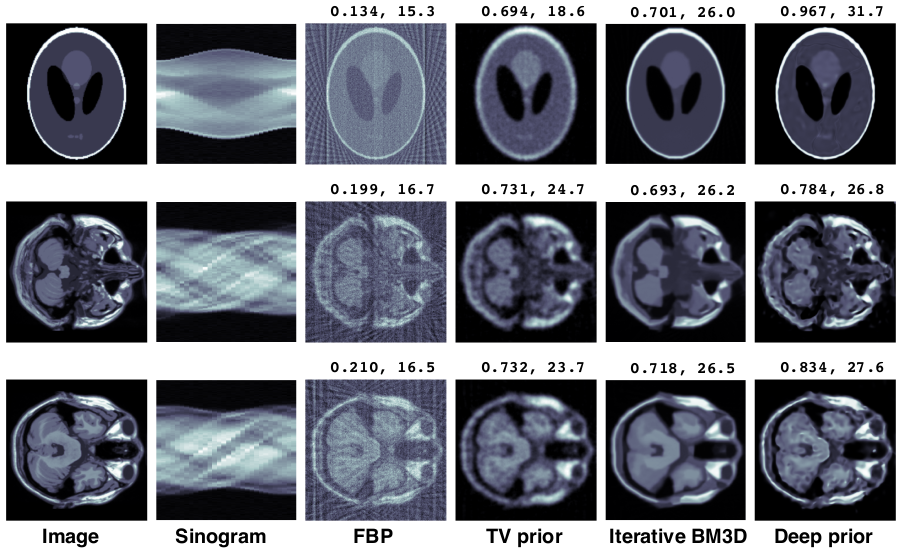}%
		}
		\caption{Principle of tomography and reconstruction examples.}
	\end{center}
\end{figure}

\begin{example}[Computerized Tomography]\label{ex: tomography}
	In X-ray tomography, we are interested in observing the density of an object without  physically penetrating it.
	A standout example is medical screening of internal organs for abnormalities such as cancer. 
	
	Let us denote the density of the object as a nonnegative function
	$x\colon D \to \bb R_+$ with bounded support $D \subset \bb R^2$.
	The object is illuminated with parallel X-ray beams from some direction.
	These can be viewed as parallel lines defined by the normal vector $w \in \bb R^2$ with $\| w \|_2 = 1$ and the shift from the origin $s \in \bb R$.
	More precisely, for each pair $(s,w)$, photons travel from the emitter at position $t=0$ to the detector at $t = T$ along the line
	\begin{equation}
	s w + t w^\perp, \quad t \in [0,T]
	\quad \text{with} \quad w^\perp \in \R^2, \ w^\perp \perp w.
	\end{equation}
	As the beam penetrates the object, it loses energy proportional to the density of the object 
	\begin{equation}
	\frac{\diff I(s w + t w^\perp)}{\diff t} = - I(s w + t w^\perp) x(s w + t w^\perp).
	\end{equation}
	By bringing the intensity to the left-hand side and integrating for $t \in [0,T]$, we obtain
	\begin{equation}
	\log I(s w + T w^\perp) - \log I(sw) = - \int_{0}^T x(s w + t w^\perp) \diff t.
	\end{equation}
	The first and the second term are the resulting and the initial intensities measured by the detector and the emitter, respectively.
	Assuming that $D$ is completely within the measured region, the integral on the right-hand side becomes
	\begin{equation}
	F[x](s,w) = \int_{\bb R} x(s w + t w^\perp) \diff t, \quad s \in \bb R, \  w \in \bb R^2, \  \| w \|_2 =1.
	\end{equation}
	This is known as the \emphi{Radon transform}.
	The measurements are performed multiple times from different directions, commonly referred to as angles, and then arranged in the form of a \emph{sinogram}, see \Cref{fig: tomography}. \Cref{fig: tomography reconstruction} shows some examples, their sinograms and reconstructions using classical, variational and learned methods.  
\end{example}

\begin{figure}
	\centering
	\begin{minipage}{.45\textwidth}
		\centering
			\includegraphics[width= 175pt, height = 100pt]{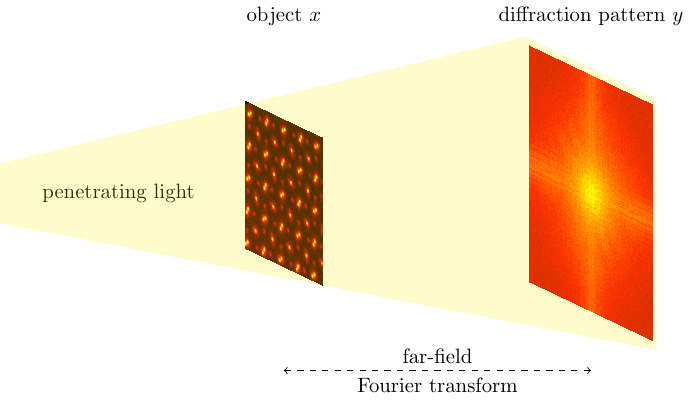}% 
		\caption{Fourier phase retrieval.}
		\label{fig: Fourier pr}
	\end{minipage}%
	\hfill
	\begin{minipage}{.45\textwidth}
		\centering
		\vspace{30pt}
		\includegraphics[width= 175pt, height = 70pt]{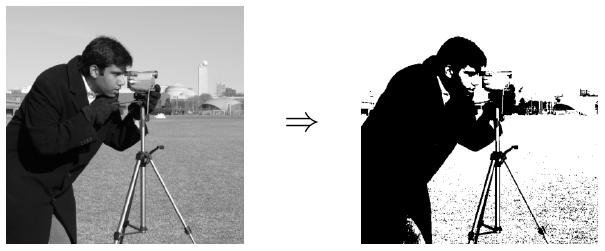}
		\caption{Black-and-white conversion.}
		\label{fig: quantization}
	\end{minipage}
\end{figure}

\begin{example}[Phase Retrieval]\label{ex: Fourier phase retrieval}
	In material science, penetrating light is used to study microscopic objects. Since the size of an object is of the same magnitude as the wavelength of light, the phenomena of diffraction occurs. 
%	Due to resolution limit of conventional optics, such as microscopes, a principle of diffraction is used instead when the  
	
	We represent a thin object of interest by a complex-valued  function $x \in L^1(\R^2)$, where $|x|$ and $x/|x|$ model absorption and scattering of light by the object. 
	The object is illuminated with plane wave illumination perpendicular to the object and the resulting diffraction pattern is captured by a detector in the far field, see \Cref{fig: Fourier pr}.
	This is described by the measurements
	\begin{equation}\label{eq: fourier phase retrieval}
		y(\omega) = F[x](\omega) = \left| \int_{\bb R^2} x(t) e^{- 2 \pi i (t_1 \omega_1 + t_2 \omega_2)}  \diff t \right|^2, \quad \omega \in \bb R^2
	\end{equation}
	with the integral being two-dimensional \emph{Fourier transform}. The problem of recovering $x$ from $y$ is known as the \emph{Fourier phase retrieval} problem.
\end{example}

\begin{example}[Heat Equation]\label{ex: heat equation}
	Consider a metal rod of length one that is isolated at both ends.
	Given its initial temperature $x\colon [0,1] \to \bb R_+$, the heat distribution  $u(s,t)$  at position $s \in [0,1]$ and time $t \ge 0$ is governed by the differential equation
	\begin{equation}
	\frac{\partial u}{\partial t} = \frac{\partial^2 u}{\partial s^2} 
	\end{equation}
	with boundary conditions
	\begin{equation}
	u(0,t) = u(1,t) = 0, \quad t \ge 0,
	\end{equation}
	and initial condition
	\begin{equation}
	u(s,0) = x(s), \quad s \in [0,1].
	\end{equation}
	Now, we can view the final state $y(s) = u(s,T)$, $s \in [0,1]$ at time $t =T$ as forward mapping.
	Then, the inverse problem is to infer the initial temperature distribution $x$ from the measurements $y$.
	Inverse problems in this area are not restricted to finding solutions of differential equations, and a common application is also parameter inference.   
\end{example}

\begin{example}[Quantization]\label{ex: quantization}
	In signal processing, a standing problem is digital-to-analog conversion.
	While the physical quantities are, for instance, elements of $\R$, the measuring device is only able to capture them up to a determined precision. An extreme example is one-bit quantization, where we represent the signal as a sequence of bits in $\{0,1\}$. 
	
	Let us consider a naive black-and-white conversion.
	In this case, $x \in [0,1]^{N \times M}$ is a gray-scale image, which is converted to the black-and-white color palette as
	\begin{equation}\label{eq: grayscale}
		y_{k,j} = F[x](k,j) = \begin{cases}
			1, & x_{k,j} \ge 1/2, \\
			0, & x_{k,j} < 1/2,
		\end{cases}  
		\quad k = 1, \ldots, N, \ j = 1, \ldots, M.
	\end{equation}     
	As can be seen in \Cref{fig: quantization}, recovering $x$ from $y$ is impossible and for this reason more practical quantization techniques were developed.
\end{example}

We call a mapping $F\colon U\to \bb Y$
\begin{itemize}
	\item \emphi{continuous} in $x\in U$ if for all $\eps>0$ there exists a $\delta >0$ with
	\begin{equation}
		\norm{F(x)-F(z)}_\bb Y \leq \eps\qquad \text{for all } z\in U \cap B_\delta(x);
	\end{equation}
	\item \emph{$L$-}\emphi{Lipschitz continuous} if there exists a \emph{Lipschitz constant} $L>0$ with
	\begin{equation}
		\norm{F(x_1)-F(x_2)}_\bb Y \leq L \norm{x_1-x_2}_\bb X \qquad\text{for all } x_1,x_2 \in U;
	\end{equation}
\end{itemize}
A mapping $F\colon\bb X\to \bb Y$ is continuous if and only if $x_n\to x$ implies $F(x_n)\to F(x)$.
If $U \subset \bb X$ is a subspace of $\bb X$ and
\begin{equation}
	F(\lambda_1 x_1 +\lambda_2 x_2) = \lambda_1 F(x_1) + \lambda_2 F(x_2) 
	\quad \text{for all } \lambda_1,\lambda_2 \in \R \text{ and } x_1,x_2 \in U,
\end{equation}
then we call $F$ \emphi{linear}.
In this case, continuity is equivalent to the existence of a constant $C>0$ such that
\begin{equation}\label{eq:LinOpCont}
	\norm{F(x)}_\bb Y \leq C\norm{x}_\bb X \qquad\text{for all }x\in \bb X.
\end{equation}
Hence, continuous linear mappings are also called \emphi{bounded}.

Examples \ref{ex: convolution}, \ref{ex: tomography}, and \ref{ex: heat equation}, feature a continuous linear $F$, in \Cref{ex: Fourier phase retrieval} $F$ is continuous but not linear, and in \Cref{ex: quantization} the $F$ is discontinuous and nonlinear.

Finally, we denote by
\begin{itemize}
	\item $\q D(F) \coloneqq U$ the \emphi{domain} of $F$;
	\item $\q N(F) \coloneqq \{ x\in U \, : \, F(x) = 0 \}$ the \enquote{\emphi{kernel}} or \enquote{\emph{null space}} of $F$;
	\item $\q R(F) \coloneqq \{ F(x)\in \bb Y \, : \, x\in U \}$ the \enquote{\emphi{range}} of $F$.
\end{itemize}

\emph{From now on, we mostly consider} bounded linear operators as forward mappings.
We denote these by $A$ and omit the parentheses around the argument, i.e., $F(x) = Ax$.
When a definition applies to the nonlinear case, we keep using $F$ instead.

\subsubsection{Exercises}

\begin{exerciseB}[Convolution theorem]\label{ex: convolution theorem}
	Let $F_N \in \bb C^{N \times N}$ be the discrete Fourier transform matrix $(F_N)_{k,j} = e^{-\frac{2 \pi i k j}{N}}$ with $i$ being the complex unit, i.e., $i^2 =-1$. Show the convolution theorem for the discrete Fourier transform,
	\begin{equation}
		F_N(f\ast g) = (F_N f)\circ (F_N g),
	\end{equation}
	where $\circ$ denotes element-wise multiplication.
\end{exerciseB}

\begin{Solution}
	We compute the discrete Fourier transform of the cyclic convolution \( f * g \) for each \( \ell \in \{0, 1, \ldots, N-1\} \):
	\begin{align*}
		(F_N(f * g))(\ell)
		&= \sum_{n=0}^{N-1} (f * g)(n) \cdot e^{-2\pi i \frac{n\ell}{N}} \\
		&= \sum_{n=0}^{N-1} \left( \sum_{m=0}^{N-1} f(m) \cdot g((n - m) \bmod N) \right) e^{-2\pi i \frac{n\ell}{N}} \\
		&= \sum_{m=0}^{N-1} f(m) \sum_{n=0}^{N-1} g((n - m) \bmod N) \cdot e^{-2\pi i \frac{n\ell}{N}} \\
		&= \sum_{m=0}^{N-1} f(m) \cdot e^{-2\pi i \frac{m\ell}{N}} \sum_{r=0}^{N-1} g(r) \cdot e^{-2\pi i \frac{r\ell}{N}} \\
		&= (F_N f)(\ell) \cdot (F_N g)(\ell),
	\end{align*}
	where we used the substitution \( r = (n - m) \bmod N \). This proves the convolution theorem:
	\begin{equation}
		F_N(f * g) = (F_N f) \circ (F_N g).
	\end{equation}
\end{Solution}

\begin{exerciseB}\label{ex: circulant matrix}
	Given a vector \( g \in \mathbb{C}^N \), the \emph{circulant matrix} \( C(g) \in \mathbb{C}^{N \times N} \) is defined such that its first column is \( g \), and each subsequent column is a cyclic shift of the previous one, namely $C(g)_{n,m} = g_{(n - m) \bmod N}$.
	
	\begin{itemize}
		\item[(i)] Verify that the cyclic convolution from \Cref{ex: convolution} can be written as matrix-vector multiplication
		\begin{equation}
			f *_N g = C(g) f = C(f) g.
		\end{equation}
		\item[(ii)] Consider the discrete derivative $\nabla f$ of a vector $f \in \mathbb \R^N$ is defined by
		\begin{equation}
		(\nabla f)_n = f_{n+1} - f_n, \qquad n=1,\ldots,N
		\end{equation}
		where $f$ is continued periodically, i.e., $f_{N+1} \coloneqq f_1$.
		Show that it can be presented as a convolution $f *_N g$, where the kernel $g$ has two nonzero entries $g_0 = -1$, $g_{N-1} =1$.
		\item[(iii)] Show that discrete convolution is associative, i.e., it satisfies
		\[
		(f *_N g) *_N h = f *_N (g *_N h) 
		\]
	\end{itemize}
\end{exerciseB}

\begin{Solution}
	(i):~By construction, 
	\begin{equation}
		(C(g) f)_n = \sum_{k =0}^{N-1} C(g)_{n,k} f_k
		= \sum_{k =0}^{N} g_{(n-k) \bmod N} f_k
		= (f * g)_n
	\end{equation}
	Changing the summation order gives $f*g = g*f$.
	
	(ii):~By definition, we get
	\begin{equation}
		(f * g)_n = \sum_{k =0}^{N} f_k g_{(n-k) \bmod N} 
		= f_{n} g_{0} + f_{n+1} g_{-1 \bmod N}
		= f_{n} g_{0} + f_{n+1} g_{N-1}
		= f_{n+1} - f_{n} = (\nabla f)_n.
	\end{equation}	
	
	(iii):~ Follows from the convolution theorem, \Cref{ex: convolution theorem}, and associativity of entrywise product by
	\begin{align*}
		F_N [(f *_N g) *_N h ] 
		& = F_N [f *_N g] \circ F_N h 
		= [F_N  f \circ F_N g] \circ F_N h \\
		& = F_N  f \circ [F_N g \circ F_N h]
		= F_N f \circ F_N [ g *_N h]
		= F_N [f *_N (g *_N h)].
	\end{align*} 
	Applying $F^{-1}_N$ from both sides gives the desired equality.
\end{Solution}

\subsection{Bounded Linear Operators}

As we have seen in Section \ref{sec:ForwardMappings}, the forward mapping $A$ of many inverse problems is bounded.
Here, we briefly review the key theoretical results from functional analysis for such operators.
As a motivation, we consider the following example.

\begin{example}\label{ex: linear operators as infinite matrices}
	Let $\bb X$ and $\bb Y$ be separable Hilbert spaces with bases $\{u_j\}_{j \in \N}$ and $\{v_k\}_{k \in \N}$, respectively.
	Let $A \colon \bb X \to \bb Y$ be a bounded linear operator. For any $x \in \bb X$, we can write
	\begin{equation}
		Ax = \sum_{k \in \N} \inner{Ax}{v_k} v_k.
	\end{equation}
	By continuity and linearity of $A$, the coefficients $(Ax)_k \coloneqq \inner{Ax}{v_k}$ are given by
	\begin{equation}
		(Ax)_k = \inner{Ax}{v_k} 
		= \Bigl\langle A \sum_{j \in \N} \inner{x}{u_j} u_j, v_k \Bigr\rangle
		= \sum_{j \in \N} \inner{x}{u_j} \inner{A u_j}{v_k}.
	\end{equation}
	By setting $A_{k,j} \coloneqq \inner{A u_j}{v_k}$ and $x_j \coloneqq \inner{x}{u_j}$, we obtain that 
	\begin{equation}
		(Ax)_k = \sum_{j \in \N} A_{k,j} x_j,
	\end{equation}
	and $A$ is essentially an infinitely long matrix.
	Hence, in separable infinite-dimensional Hilbert spaces, bounded operators can be seen as a generalization of matrices.   
\end{example}

Given normed spaces $(\bb X,\norm{\cdot}_\bb X)$ and $(\bb Y,\norm{\cdot}_\bb Y)$, we endow the vector space $\q L(\bb X,\bb Y)$ of bounded linear operators  $A \colon \bb X \to \bb Y$ with the \emphi{operator norm}
\begin{equation}
	\norm{A}_{\q L(\bb X,\bb Y)} = \sup_{x\in \bb X \, \setminus\{0\}}\frac{\norm{Ax}_\bb Y}{\norm{x}_\bb X}
	= \sup_{\norm{x}_\bb X\leq 1} \norm{Ax}_\bb Y = \sup_{\norm{x}_\bb X = 1} \norm{Ax}_\bb Y.
\end{equation}
By definition, we have that $\norm{Ax}_\bb Y \leq \norm{A}_{\q L(\bb X, \bb Y)} \norm{x}_\bb X$ for all $x\in \bb X$.
Moreover, the composition of bounded operators is bounded.
If $\bb Y$ is complete, then $\q L(\bb X,\bb Y)$ is a Banach space.
For $\bb Y=\R$, we refer to the space $\q L(\bb X,\R)$ as \emph{bounded functionals} on $\bb X$.

We call a linear operator $A \in \q L(\bb X,\bb Y)$
\begin{itemize}
	\item \emph{injective} if $\q N(A) = \{0\}$;
	\item \emph{surjective} if $\q R(A)= \bb Y$;
	\item \emph{bijective} if $A$ is injective and surjective.
\end{itemize}
For a bijective $A\in \q L(\bb X,\bb Y)$, the \emph{inverse} $A^{-1}\colon \bb Y\to \bb X$ is continuous if and only if there exists a $c>0$ with
\begin{equation}
	\label{eq:funktan:inverse_bd}
	c\norm{x}_\bb X \leq \norm{Ax}_\bb Y \qquad\text{for all } x\in \bb X.
\end{equation}
In this case, $\norm{A^{-1}}_{\q L(\bb Y,\bb X)} = c^{-1}$ holds for the largest $c$ satisfying \eqref{eq:funktan:inverse_bd}.
If $\bb X,\bb Y$ are Banach spaces, then the theorem of the continuous inverse states that $A^{-1} \in \q L(\bb Y,\bb X)$.

In practice, we frequently encounter that $A$ is injective but not surjective. 
Then, the question arises if a continuous inverse on  $\q R(A)$ exists 
However, this does not hold in general, which is a fundamental issue for infinite-dimensional problems.
\begin{thm}{Closed range}{closedrange}
	If $\bb X,\bb Y$ are Banach spaces and $A \in\q L(\bb X,\bb Y)$ is injective, then $A^{-1}\colon \q R(A)\to \bb X$ is continuous if and only if $\q R(A)$ is closed.
\end{thm}

For the rest of this section, let $\bb X$ and $\bb Y$ be Hilbert spaces.
For a bounded operator $A \in \q L(\bb X,\bb Y)$, we define the \emphi{adjoint operator} $A^*\in \q L(\bb Y,\bb X)$ via
\begin{equation}
	\inner{A^* y}{x}_\bb X = \inner{A x}{y}_\bb Y \qquad\text{for all } x\in \bb X, y \in \bb Y.
\end{equation}
We can interpret the adjoint operator as a generalization of taking the transposed of a matrix $A$.
The adjoint satisfies
\begin{enumerate}
	\item $(A^*)^* = A$;
	\item $(AB)^* = B^* A^*$,
	\item $\norm{A^*}_{\q L(\bb Y,\bb X)} = \norm{A}_{\q L(\bb X,\bb Y)}$;
	\item $\norm{A^*A}_{\q L(\bb X,\bb X)} = \norm{A}_{\q L(\bb X,\bb Y)}^2$.
\end{enumerate}
For $A\in \q L(\bb X,\bb Y)$ we have $\q R(A)^\bot = \q N(A^*)$ and hence $\q N(A^*)^\bot = \overline{\q R(A)}$.
Thus, $\q N(A)$  is closed and $A$ is injective if and only if $\overline{\q R(A^*)} = \bb X$.
If $A^* = A$, then $A$ is called \emphi{selfadjoint}.
\begin{example}
	The orthogonal projection $P_U$ is selfadjoint with $\norm{P_U}_{\q L(\bb X,\bb X)} = 1$.
\end{example}

Many results of this lecture rely on the spectral theory for compact operators.
An operator $A \in \q L(\bb X,\bb Y)$ is \emphi{compact} if the image of every bounded sequence $\{x_n\}_{n\in\N}\subset \bb X$ contains a convergent subsequence $\{A x_{n_k}\}_{k\in\N}\subset \bb Y$.
We denote such operators by $K$.

For $A \in \q L(\bb X,\bb Y)$ with $\dim(\q R(A)) < \infty$, compactness follows by the Bolzano--Weierstrass theorem.  
In particular, the \emph{identity} $I \colon \bb X \to \bb X$ is compact if \emph{and only if} $\dim(\bb X) < \infty$.
The space $\q K(\bb X,\bb Y)$ of compact operators from $\bb X$ to $\bb Y$ is a closed subspace of $\q L(\bb X,\bb Y)$ and hence it is a Banach space. 

For $A\in \q L(\bb X,\bb Y)$ and $B \in \q L(\bb Y,\bb E)$, the composition $B \circ A$ is compact if $A$ or $B$ is compact. Moreover, $A^*$ is compact if and only if $A$ is compact.

\begin{example}[Integral Operators]\label{ex:integral-operator}
	For $\bb X =L^2([0,1])$ and a kernel $k\in L^2([0,1] \times [0,1])$, we consider the operator $K\colon \bb X\to  \bb X$ defined pointwise via
	\begin{equation}
		[K x](t) = \int_0^1 k(s,t)x(s)\,\diff s \qquad\text{for almost every }t\in [0,1].
	\end{equation}
	By the Cauchy--Schwarz inequality and Fubini's Theorem we have
	\begin{align}
		\norm{K x}_{L^2([0,1])}^2
		& = \int_{0}^1 \bigg| \int_0^1 k(s,t)x(s)\,\diff s \bigg|^2  \,\diff t \\
		& \le \int_{0}^1 \int_0^1 |k(s,t)|^2 \diff s  \,\diff t \, \int_{0}^1 \int_{0}^1 |x(s)|^2 \diff s \,\diff t
		= \norm{k}_{L^2([0,1]^2)}^2 \norm{x}_{L^2([0,1])}^2.
	\end{align}
	This gives that $K \in \q L(\bb X,\bb X)$ with $\norm{K}_{\q L(\bb X,\bb X)} \leq \norm{k}_{L^2([0,1]^2)}.$
	Since $k\in L^2([0,1]^2)$ is measurable, there exists a sequence of simple functions $\{k_n\}_{n\in \N}$ with
	\begin{equation}\label{eq:SimpFunc}
		k_n(s,t) = \sum_{i,j=1}^n \alpha_{ij}\1_{E_i}(s)\1_{E_j}(t)
	\end{equation}
	such that $k_n\to k$ in $L^2([0,1]^2)$.
	In \eqref{eq:SimpFunc}, the $E_i$ form a disjoint decomposition of $[0,1]$ and $\1_E$ denotes the characteristic function of $E$, i.e, $\1_E(x) = 1$ for $x\in E$ and $\1_E(x) = 0$ otherwise.
	By linearity of the integral, the operators $K_n$ corresponding to the kernels $k_n$ satisfy
	\begin{equation}
		\norm{K_n -K}_{\q L(\bb X,\bb X)} \leq \norm{k_n-k}_{L^2([0,1]^2)} \to 0,
	\end{equation}
	i.e., $K_n\to K$.
	By 
	\begin{equation}
		[K_n x](t) = \int_0^1 k_n(s,t)x(s)\,\diff s = \sum_{j=1}^{n} \left(\sum_{i=1}^n \alpha_{ij}  \int_{E_i} x(s)\,\diff s\right)\1_{E_j}(t),
	\end{equation}
	their range is $\mathrm{span}(\{\1_{E_j} \, \colon \, j \in n\})$, which is finite dimensional.
	Consequently, $K$ is the limit of the compact operators $\{K_n\}_{n\in\N}$ and therefore compact.
	Using the definition of the inner product and Fubini's Theorem, we compute 
	\begin{align}
		\inner{K x}{y} 
		&= \int_0^1 Kx (t) y(t) \, \diff t
		= \int_0^1 \int_0^1 k(s,t) x (s) \, \diff s y(t)  \, \diff t \\
		& =  \int_0^1 \int_0^1 k(s,t) y(t) \, \diff t x (s) \, \diff s
		= \langle \int_0^1 k(s,t) y(t) \, \diff t, x  \rangle.
	\end{align}
	Hence, the adjoint $K^*\in \q L(\bb X,\bb X)$ is
	\begin{equation}\label{eq: int op adjoint}
		[K^*y](s) = \int_0^1 k(s,t)y(t)\,\diff t \qquad\text{for almost every }s\in [0,1].
	\end{equation}
	Thus, an integral operator is selfadjoint if and only if the kernel is symmetric, i.e., $k(s,t) = k(t,s)$ for almost every $s,t\in [0,1]$.
\end{example}

Many practically relevant operators such as solution operators to (partial) differential equations or convolutions are integral operators.
In particular, Examples \ref{ex: convolution} and \ref{ex: tomography} are instances of integral operators.
Similarly to Example \ref{ex:integral-operator}, it can be shown that these are compact, see also Exercise \ref{ex:IntOpGen}.
Next, we state a spectral theorem in Hilbert spaces.

\begin{thm}{Spectral theorem}{spektral}
	Let $\bb X$ be a Hilbert space and $K\in \q K(\bb X,\bb X)$ be selfadjoint. Then, there exists a ONB $\{u_n\}_{n\in \N}$ of $\overline{\q R(K)}$ and a null sequence $\{\lambda_n\}_{n\in\N}\subset \R \setminus \{0\}$ with
	\begin{equation}
		Kx = \sum_{n\in \N} \lambda_n \inner{x}{u_n}_\bb X u_n \qquad \text{for all }x \in \bb X.
	\end{equation}
	Both $\{u_n\}_{n\in \N}$ and $\{\lambda_n\}_{n\in\N}$ can be finitely terminating. 
\end{thm}
\begin{remark}
	In  \Cref{th:spektral}, we do not require $\bb X$ to be separable.
	In particular, no countable basis might exist for $\q N(K)$, which is the eigenspace to the eigenvalue $\lambda=0$.
	Hence, it is excluded in the theorem.
\end{remark}
For $x=u_n$, we obtain $K u_n = \lambda_n u_n$, so that $u_n$ is an eigenvector corresponding to the eigenvalue $\lambda_n$.
By convention, the eigenvalues are sorted by decreasing magnitude, i.e.,
\begin{equation}
	|\lambda_1|\geq |\lambda_2|\geq \cdots > 0.
\end{equation}
In particular, $\norm{K}_{\q L(\bb X,\bb X)} = |\lambda_1|$.
If  a selfadjoint  operator $K \in \q K(\bb X, \bb X)$ satisfies 
\begin{equation}
	\inner{K x}{x}_\bb X \ge 0 \quad \text{for all } x \in \bb X,
\end{equation}
it is called \emphi{positive semidefinite}.
For $K \in \q K(\bb X, \bb X)$, positive semidefinitness is equivalent to the eigenvalues being nonnegative, $\lambda_n \ge 0$, $n \ge 0$.
For any $K \in \q K(\bb X, \bb Y)$, the operators $K^* K \in \q K(\bb X, \bb X)$ and $K K^* \in \q K(\bb Y, \bb Y)$ are selfadjoint and positive semidefinite.
We will use this fact later to establish a singular value decomposition of compact operators.

\subsubsection{Exercises}
\begin{exerciseB}
	Show that a compact operator $ K \in \mathcal K(\bb X, \bb Y)$ maps weakly convergent sequences in $\bb X$ to strongly convergent ones in $\bb Y$.
\end{exerciseB}

\begin{exerciseB}
	Compute the adjoint of the continuous convolution operator $F \colon L^2(\R) \to L^2(\R)$ from  \Cref{ex: convolution}.
	When is it self adjoint?
\end{exerciseB}

\begin{Solution}
	Let $g \in L^1(\R)$, and define the convolution operator
	\begin{equation}
	(A_g f)(x) = (f * g)(x) = \int_{\R} f(y) h(x - y) \, dy.
	\end{equation}
	To compute the adjoint $A_g^*$, we consider
	\begin{align}
		\inner{A_g f}{h} 
		&= \int_{\R} (f * g)(x) \, h(x) \, \diff x 
		= \int_{\R} \int_{\R} f(y) g(x - y) \, \diff y h(x) \, \diff x  \notag \\
		&= \int_{\R} f(y) \left( \int_{\R} g(x - y) h(x) \, \diff x \right) \diff y 
		= \int_{\R} f(y) (h * \tilde g)(y) \, \diff y
		= \inner{f}{h * \tilde g},
	\end{align}
	where $\tilde g(x) \coloneqq \overline{g(-x)}$. Thus,
	\begin{equation}
	A_g^* = A_{\tilde{g}}, \quad \text{and} \quad A_g \text{ is self-adjoint if and only if } g(x) = \overline{g(-x)}.
	\end{equation}
\end{Solution}

\begin{exerciseB}\label{ex:IntOpGen}
	Here, we generalize \Cref{ex:integral-operator}. 
	Let $\Omega, \Omega' \subset \R^d$ be open and $ k \in L^2(\Omega \times \Omega') $. Let 
	$K \colon L^2(\Omega) \to L^2(\Omega') $ be the associated Hilbert-Schmidt operator
	\begin{equation}
		[Kx](t) = \int_\Omega k(s,t)x(s) \diff s \qquad \text{for almost every } t \in \Omega'. 
	\end{equation}
	Show that $K$ is compact.
%		\item[(b)] Conclude that the singular values of $ K$ satisfy $\sigma_n = \mathcal{O}(n^{-\frac{1}{2}})$.
%	\end{itemize}
%	
%	\noindent
%	\textbf{Hint:} \textit{Given an ONB $ U $ and $ V $ of $ L^2(\Omega) $ and $ L^2(\Omega')$, respectively, the set}
%	\begin{equation}
%		J = \{ u \otimes v \mid u \in U,\, v \in V \}, \quad \text{with } [f \otimes g](s,t) = f(s)g(t) \text{ a.e.}
%	\end{equation}
%	\textit{is an ONB of $ L^2(\Omega \times \Omega')$.
%		Also note the $L^2$ spaces are separable, hence any orthonormal system is at most countable.}
\end{exerciseB}

\begin{Solution}
	We will prove that $K$ is compact by showing it can be approximated in operator norm by finite-rank operators.
	
	\textbf{Step 1:} First, let us verify that $K$ is a bounded operator. For any $x \in L^2(\Omega)$, we have
	\begin{align}
		\|Kx\|^2_{L^2(\Omega')} 
		&= \int_{\Omega'} |[Kx](t)|^2 \, \diff t 
		= \int_{\Omega'} \left|\int_\Omega k(s,t)x(s) \, \diff s\right|^2 \, \diff t.
	\end{align}
	The Cauchy-Schwarz inequality yields
	\begin{align}
		\left|\int_\Omega k(s,t)x(s) \, \diff s\right|^2 &\leq \left(\int_\Omega |k(s,t)|^2 \, \diff s\right) \left(\int_\Omega |x(s)|^2 \, \diff s\right)
		= \left(\int_\Omega |k(s,t)|^2 \, \diff s\right) \|x\|^2_{L^2(\Omega)}.
	\end{align}
	Therefore,
	\begin{align}
		\|Kx\|^2_{L^2(\Omega')} &\leq \int_{\Omega'} \left(\int_\Omega |k(s,t)|^2 \, \diff s\right) \|x\|^2_{L^2(\Omega)} \, \diff t \notag \\
		& = \|x\|^2_{L^2(\Omega)} \int_{\Omega'} \int_\Omega |k(s,t)|^2 \, \diff s \, \diff t
		= \|k\|^2_{L^2(\Omega \times \Omega')} \cdot \|x\|^2_{L^2(\Omega)}.
	\end{align}
	Thus, $\|K\| \leq \|k\|_{L^2(\Omega \times \Omega')} < \infty$, showing that $K$ is bounded.
	
	\textbf{Step 2:} Now, we show that $K$ can be approximated by finite-rank operators.
	
	Since $k \in L^2(\Omega \times \Omega')$, we can find a sequence of simple functions $\{k_n\}_{n=1}^{\infty}$ such that $k_n \to k$ in $L^2(\Omega \times \Omega')$. Each simple function $k_n$ can be written as
	\begin{equation}
		k_n(s,t) = \sum_{i=1}^{n} a_i^{(n)} \cdot \chi_{A_i^{(n)}}(s) \cdot \chi_{B_i^{(n)}}(t),
	\end{equation}
	where $a_i^{(n)}$ are constants, $\{A_i\}_{i = 1}^{n}$ and $\{B_i\}_{i = 1}^{n}$ are partitions of $\Omega$ and $\Omega'$, and $\1_{E}$ denotes a characteristic functions of a measurable set $E$.
	
	Let $K_n$ be the integral operator with kernel $k_n$. We rewrite it as
	\begin{align}
		[K_n x](t) 
		&= \int_\Omega k_n(s,t)x(s) \, \diff s
		= \sum_{i=1}^{n} a_i^{(n)} \cdot \left(\int_\Omega \1_{A_i^{(n)}}(s)x(s) \, \diff s\right) \cdot \1_{B_i^{(n)}}(t) \notag \\
		&= \sum_{i=1}^{n} a_i^{(n)} \cdot \langle x, \1_{A_i^{(n)}} \rangle_{L^2(\Omega)} \cdot \1_{B_i^{(n)}}(t).
	\end{align}
	This shows that $K_n$ is a finite-rank operator for each $n$, as its range is contained in the finite-dimensional space spanned by $\{\1_{B_i^{(n)}}\}_{i=1}^{n}$.
	
	\textbf{Step 3:} We show that $K_n \to K$ in operator norm.
	
	For any $x \in L^2(\Omega)$ with $\|x\|_{L^2(\Omega)} \leq 1$ we obtain
	\begin{align}
		\|(K-K_n)x\|^2_{L^2(\Omega')} &= \int_{\Omega'} \left|\int_\Omega (k(s,t)-k_n(s,t))x(s) \, \diff s\right|^2 \, \diff t.
	\end{align}
	Using Cauchy-Schwarz as before we get
	\begin{align}
		\|(K-K_n)x\|^2_{L^2(\Omega')} &\leq \|k-k_n\|^2_{L^2(\Omega \times \Omega')} \cdot \|x\|^2_{L^2(\Omega)}
		\leq \|k-k_n\|^2_{L^2(\Omega \times \Omega')}.
	\end{align}
	Hence,
	\begin{equation}
		\|K-K_n\| \leq \|k-k_n\|_{L^2(\Omega \times \Omega')}.
	\end{equation}
	Since $k_n \to k$ in $L^2(\Omega \times \Omega')$, we have $\|K-K_n\| \to 0$ as $n \to \infty$.
	
	\textbf{Conclusion:} We have shown that $K$ is a bounded linear operator that can be approximated in operator norm by finite-rank operators. Hence $K$ is compact.
\end{Solution}

\section{Ill-posed Inverse Problems}
For an operator $A \in \q L(\bb X,\bb Y)$, we call the equation $Ax=y$ \emphi{well-posed} (in the sense of \href{http://mathshistory.st-andrews.ac.uk/Biographies/Hadamard.html}{Hadamard}) if 
\begin{enumerate}
	\item for all $y\in \bb Y$ there exists an $x\in \bb X$ with $Ax = y$;
	\item the solution $x$ is unique, i.e., $z\neq x$ implies $Az\neq y$;
	\item $x$ depends continuously on $y$, i.e., for $\{x_n\}_{n\in\N}$ with $Ax_n\to y$ we have $x_n\to x$.
\end{enumerate}
In \circled{3}, the slightly technical formulation ensures that the measurements $y_n \coloneqq A x_n$ have pre-images.
 If one of these conditions is violated, the equation is \emphi{ill-posed}.
Note that the conditions are not completely independent.
 For example, if $\bb X$ and $\bb Y$ are Banach spaces and $A\in \q L(\bb X,\bb Y)$ satisfies \circled{1} and \circled{2}, then $A^{-1}\in \q L(\bb Y,\bb X)$ by \Cref{th:closedrange} and thus \circled{3} holds.
 
Usually, a violation of \circled{1} occurs due to undersampling, namely if we have less data than unknowns.
Then, we resort to \emph{least-squares solutions}.
Further, if \circled{2} is violated and multiple solutions exist, we simply choose the one with the smallest norm.
This leads to the Moore-Penrose inverse as introduced in \Cref{sec:MoorePenrose}.
For finite-dimensional $\bb X$, we thereby get a continuous inverse and the problem of well-posedness appears to be resolved (although an efficient implementation may take significant effort).
The situation is completely different in infinite dimensions as the following example shows.
\begin{example}\label{ex:illposed_diff}
	Given $y\in C^1([0,1])$ we want to find its derivative $y'$.
	The associated forward mapping $A \in \q L(\bb X,\bb Y)$ with $A y^\prime = y$ is an integration operator, for which we need to specify the norms on its domain and range.
	To guarantee integrability, a natural choice is $\bb X = C([0,1])$.
	Further, let us assume that the data is noisy, i.e., we are only given
	\begin{equation}
		\tilde y = y + \eta.
	\end{equation}
	To ensure $\circled{1}$, we would need to require $\eta \in C^1([0,1])$ (even though noise is usually irregular).
	We ignore this problem for now and equip $\bb Y$ with the $C([0,1])$ norm.
	However, \circled{3} is violated in this case:
	Let  $\{\delta_n\}_{n\in\N}$ with $\delta_n\to 0$ and choose $k\in \N$ arbitrary.
	Then, we can define data $\tilde y_n \coloneqq y +\eta_n$ with noise
	\begin{equation}
		\eta_n\coloneqq  \delta_n\sin\bigl(\tfrac{k}{\delta_n} \cdot \bigr) \in C^1([0,1]).
	\end{equation}
	For the data $\tilde y_n $, we have
	\begin{equation}
		\norm{\tilde y_n-y}_{C} = \norm{\eta_n}_{C}=\delta_n\to 0
	\end{equation}
	and
	\begin{equation}
		\tilde x_n(t)\coloneqq \tilde y_n'(t) = y'(t) + k\cos\bigl(\tfrac{k}{\delta_n} t \bigr).
	\end{equation}
	Thus, it holds that $A \tilde x_n \to y$, but for $x \coloneqq y'$ we have
	\begin{equation}
		\norm{\tilde x_n - x}_{C} = \norm{\eta'_n}_{C} = k \qquad\text{for all }n\in \N,
	\end{equation}
	namely that the inverse is discontinuous.
	
	In contrast, the problem is well-posed for $\bb Y = C^1([0,1])$, since then $\norm{\eta_n}_{C^1}\to 0$ implies that $\norm{x-x_n}_{C} \leq \norm{\eta_n}_{C^1}\to 0$.
	However, taking $\bb Y = C^1([0,1])$ implies that besides $y$ also $y'$ is measured -- so we are no longer considering an inverse problem.
\end{example}

\begin{remark}
	Given these problems, one might wonder why inverse problems in infinite-dimensional spaces should be considered at all, especially since any practical implementation is inherently finite-dimensional.
	However, this perspective overlooks an important point.
	Discretizing an inverse problem already introduces a form of regularization, and as the discretization becomes finer, the operator norm of $A^{-1}$ often grows (sometimes even very fast).
	Investigating the inverse problem in an infinite-dimensional setting therefore provides crucial insight into its stability and the behavior of numerical implementations.
\end{remark}

\subsection{Moore-Penrose  Inverse}\label{sec:MoorePenrose}
Let $\bb X$ and $\bb Y$ be Hilbert spaces and consider the equation $Ax=y$ with $A \in \q L(\bb X, \bb Y)$.
A solution only exists for $y\in\q R(A)$.
For $y\notin\q R(A)$, we instead take a $x\in \bb X$ that minimize the distance $\norm{Ax-y}_{\bb Y}$.
If $\q N(A)\neq \{0\}$, infinitely many such $x$ exist and we take the one with the minimal norm.
By doing so, we can resolve the first two sources of ill-posedness.
\begin{defn}
	An element $x^\dag\in \bb X$ is called
	\begin{enumerate}
		\item \emphi{least squares solution} of $Ax=y$ if
		\begin{equation}\label{eq:LeastSquaresDef}
			\norm{Ax^\dag-y}_{\bb Y} = \min_{z\in \bb X}\norm{Az-y}_{\bb Y};
		\end{equation}
		\item \emphi{minimum norm solution} of $Ax=y$ if additionally
		\begin{equation}
			\norm{x^\dag}_{\bb X} = \min \bigl\{\norm{z}_{\bb X} :  \text{$z$ is least squares solution of $Ax=y$} \bigr\}.
		\end{equation}
	\end{enumerate}
\end{defn}
For bijective $A$, we have that $x=A^{-1}y$ is the minimum norm solution.
Throughout this section, we will answer the question for which $y\in \bb Y$ a minimum norm solution exists.
Moreover, we will establish the connection to the \emphi{Moore--Penrose inverse} (or \emph{pseudoinverse}) as defined below.
This operator generalizes the usual inverse for non-invertible $A$.
\begin{thm}{}{pseudoinverse}
	Let $A\in \q L(\bb X,\bb Y)$ and consider the restriction of $A$ to $\q N(A)^{\bot}$, namely
	\begin{equation}
		\tilde A \coloneqq A|_{\q N(A)^{\bot}}\colon \q N(A)^\bot \to \q R(A).
	\end{equation}
	Then, the \emph{Moore--Penrose inverse} $A^\dag = \tilde A^{-1} P_{\overline{\q R(A)}}$ is the unique linear extension of $\tilde A^{-1}$ with
	\begin{align}\label{eq:RestrictedOperator}
		\q D(A^\dag)&=\q R(A) \oplus \q R(A)^\bot,\\
		\q N(A^\dag)&=\q R(A)^\bot.
	\end{align}
\end{thm}
\begin{pr}
	The operator $\tilde A \colon \q N(A)^\bot \to \q R(A)$ is bijective.
	Further, any $y\in \q D(A^\dag)$ can be uniquely decomposed as $y=y_1+y_2$ with $y_1 \in \q R(A)$ and $y_2\in \q R(A)^\bot$.
	Since we require $\q N(A^\dag) = \q R(A)^\bot$, we get that $A^\dag$ is uniquely defined by
	\begin{equation}\label{eq:inverse:proj1}
		A^\dag y = A^\dag y_1 + A^\dag y_2 = A^\dag y_ 1= \tilde A^{-1} y_1.
	\end{equation}
	Finally, we need to take the closure to ensure the existence of the projection $P_{\overline{\q R(A)}}$.
\end{pr}
If $A$ is invertible, we indeed have $A^\dag = A^{-1}$ and $\q D(A^\dag) = \bb Y$.
However, if $\q R(A)$ is not closed, then $D(A^\dag) \neq \bb Y$ and we will show later in \Cref{th:inverse:cont} that $A^\dag$ is discontinuous.
Now, we prove some properties of the Moore--Penrose inverse.
\begin{lemma}\label{lem:inverse}
	The pseudoinverse $A^\dag$ satisfies $\q R(A^\dag) = \q N(A)^\bot$ and the \emph{Moore--Penrose equations}
	\begin{enumerate}
		\item $AA^\dag A = A$,
		\item $A^\dag A A^\dag = A^\dag$,
		\item $A^\dag A = I-P_{\q N}$,
		\item $AA^\dag = (P_{\overline{\q R}})|_{\q D(A^\dag)}$,
	\end{enumerate}
	where $P_{\q N}$ and $P_{\overline{\q R}}$ denote the orthogonal projections on $\q N(A)$ and $\overline{\q R(A)}$, respectively.
\end{lemma}
\begin{pr}
	Since $A^\dag = \tilde A^{-1} P_{\overline{\q R}}$ and $\tilde A \colon \q N(A)^\bot \to \q R(A)$ is bijective, we get $\q R(A^\dag) = \q N(A)^\bot$.
	
	\circled{4} For $y\in\q D(A^\dag)$, we have from $A^\dag = \tilde A^{-1} P_{\overline{\q R}}$ and $\tilde A^{-1} P_{\overline{\q R}}y\in \q N(A)^\bot$  that
	\begin{equation}
		AA^\dag y = A\tilde A^{-1} P_{\overline{\q R}}y = \tilde A\tilde A^{-1} P_{\overline{\q R}}y = P_{\overline{\q R}}y.
	\end{equation}
	
	\circled{3} It holds that $A^\dag Ax = \tilde A^{-1} P_{\overline{\q R}} A x = \tilde A^{-1} Ax$ for all $x\in \bb X$ and hence that
	\begin{equation}
		A^\dag A x = \tilde A^{-1} A\left(P_{\q N} x + (I - P_{\q N})x\right)
		= \tilde A^{-1} A P_{\q N} x + \tilde A^{-1} \tilde A (I - P_{\q N}) x = (I - P_{\q N})x.
	\end{equation}
	
	\circled{2} Using \circled{4} %and \eqref{eq:inverse:mpgl1} 
	yields for $y\in \q D(A^\dag)$ that
	\begin{equation}
		A^\dag y
		= A^\dag P_{\overline{\q R}}y + A^\dag \underbrace{( I - P_{\overline{\q R}} )y}_{\in \q R(A)^\perp}
		= A^\dag P_{\overline{\q R}}y + 0
		= A^\dag AA^\dag y.
	\end{equation}
	
	\circled{1}: From \circled{4} it follows for all $x\in \bb X$ that
	\begin{equation}
		AA^\dag Ax = P_{\overline{\q R}} Ax = Ax.
	\end{equation}
	This concludes the proof.
\end{pr}
In Exercise \ref{ex:ReverseStatementInverse}, we will prove that the Moore--Penrose equations are an equivalent characterization of $A^\dag$.
Next, we show that $A^\dagger$ assigns the minimum norm solution.
\begin{thm}{}{inverse:minnorm}
	For any $y\in \q D(A^\dag)$, the equation $Ax=y$ admits
	\begin{enumerate}
		\item least squares solutions.
		These coincide with the solutions of
		\begin{equation}\label{eq:inverse:ausgleichung}
			Ax = P_{\overline{\q R}} y;
		\end{equation}
		\item the unique minimum norm solution $x^\dag = A^\dag y \in  \q N(A)^\bot$.
		Moreover, the set of all least squares solutions is given by $x^\dag + \q N(A)$.
	\end{enumerate}
\end{thm}
\begin{pr}
	\circled{1} Since $P_{\overline{\q R}}y \in \q R(A)$ for $y\in\q D(A^\dag)$, equation \eqref{eq:inverse:ausgleichung} admits at least one solution.
	By definition of $P_{\overline{\q R}}$, we further have
	\begin{equation}\label{eq:EqualityProjection}
		\norm{P_{\overline{\q R}}y-y}_{\bb Y} \leq \min_{w\in{\q R(A)}}\norm{w-y}_{\bb Y} = \min_{x\in \bb X}\norm{Ax-y}_{\bb Y} \leq \norm{P_{\overline{\q R}}y-y}_{\bb Y}.
	\end{equation}
	Thus, any solution $z\in \bb X$ of equation \eqref{eq:inverse:ausgleichung} satisfies
	$\norm{Az-y}_{\bb Y} = \norm{P_{\overline{\q R}}y-y}_{\bb Y}  \leq \norm{Ax-y}_{\bb Y}$ for all $x\in \bb X$, namely that $z$ is a least squares solutions of $Ax=y$.
	On the other hand, by the Pythagorean theorem and  \eqref{eq:EqualityProjection} any least squares solution $z\in \bb X$ satisfies
	\begin{align}
	\norm{Az-y}_{\bb Y}^2 
	& = \| \underbrace{Az - P_{\overline{\q R}}y}_{\in \q R(A)} + \underbrace{P_{\overline{\q R}}y - y}_{\in \q R(A)^\perp} \|^2_{\bb Y}
	= \| Az - P_{\overline{\q R}}y \|_{\bb Y}^2 + \| P_{\overline{\q R}}y - y \|_{\bb Y}^2 \notag\\
	&= \| Az - P_{\overline{\q R}}y \|_{\bb Y}^2 + \| Az-y \|_{\bb Y}^2.
	\end{align}
	Hence, such a $z$ is a solution of equation \eqref{eq:inverse:ausgleichung}. 
	
	\circled{2} By \Cref{lem:inverse}\,\circled{4}, we have that $x^\dag = A^\dag y$ solves \eqref{eq:inverse:ausgleichung}, namely that $x^\dag$ is a least squares solution.
	Any other solution of \eqref{eq:inverse:ausgleichung} can be uniquely decomposed as $x= x^\dag + x_0$ with $x_0\in\q N(A)$.
	Since $x^\dag \in  \q N(A)^\bot$, we get that
	\begin{equation}
		\norm{x}_{\bb X}^2 = \norm{x^\dag +x_0}_{\bb X}^2 = \norm{x^\dag }_{\bb X}^2 + \norm{x_0}_{\bb X}^2 \geq \norm{x^\dag }^2_{\bb X}.
	\end{equation}
	Hence, $x^\dag$ is the unique minimum norm solution.
\end{pr}

We now provide an useful alternative characterization of $A^\dagger$.
\begin{cor}\label{thm:inverse:normalen}
	Let $y\in \q D(A^\dag)$. Then $x\in \bb X$ is a least squares solution of $Ax =y$ if and only if $x$ satisfies the \emphi{normal equation}
	\begin{equation}
		\label{eq:inverse:normal}
		A^* A x = A^*y.
	\end{equation}
\end{cor}
\begin{pr}
	By \Cref{th:inverse:minnorm}\,\circled{1}, $x\in \bb X$ is a least squares solution if and only if $Ax = P_{\overline{\q R}} y$.
	The latter is equivalent to $Ax-y\in \q R(A)^\bot = \q N(A^*)$, namely that $A^*(Ax-y) = 0$.
\end{pr}
The minimum norm solution of $Ax=y$ and \eqref{eq:inverse:normal} therefore coincide, which implies
\begin{equation}\label{eq:inverse:normal_pseudo}
	A^\dag y = (A^*A)^\dag A^*y.
\end{equation}

So far, we did not investigate $\q D(A^\dag) = \q R(A)\oplus \q R(A)^\bot$.
Since $\q R(A)^\bot$ is closed, we have that $\overline{\q D(A^\dag)} = \overline{\q R(A)} \oplus \q R(A)^\bot = \bb Y$.
Thus, $\q R(A)$ is closed if and only if $\q D(A^\dag)=\bb Y$.
\begin{thm}{}{inverse:cont}
	Let $A\in \q L(\bb X,\bb Y)$. Then $A^\dag \colon \q D(A^\dag)\to \bb X$ is bounded if and only if $\q R(A)$ is closed.
\end{thm}
\begin{pr}
	$\Leftarrow$: Let $\q R(A)$ be closed.
	Then  $\tilde A^{-1} \colon \q R(A) \to \q N(A)^\bot $ from \Cref{th:pseudoinverse} is bounded by \Cref{th:closedrange}.
	 Thus, we get that  $A^\dag = \tilde A^{-1} P_{\q R(A)}$  is bounded.
	
	$\Rightarrow$:
	For $y\in\overline{\q R(A)}$ and $\{y_n\}_{n\in\N} \subset \q R(A)$ with $y_n\to y$ it holds by \Cref{lem:inverse}\,\circled{4} that
	\begin{equation}
		y = P_{\overline{\q R}}y = \lim_{n\to\infty} P_{\overline{\q R}}y_n = \lim_{n\to\infty}A{A^\dag}y_n.
	\end{equation}
	Since $A^\dag$ is bounded, $\{A^\dag y_n\}_{n \in \N}$ is a Cauchy sequence by
	\begin{equation}
		\| A^\dagger y_n - A^\dagger y_m\|_{\bb X} \le \| A^\dagger \|_{\q L(\bb Y, \bb X)} \|y_n - y_m\|_{\bb Y}.
	\end{equation}
	Thus, it converges to some $x \in \bb X$, and we have that 
	\begin{equation}
		y = \lim_{n\to\infty} A A^\dag y_n = A \lim_{n\to\infty} A^\dag y_n = A x \in \q R(A),
	\end{equation}
	which in turn implies $\overline{\q R(A)}=\q R(A)$.
\end{pr}

Therefore, if $\q R(A)$ is not closed, we call $Ax=y$ \emph{ill-posed in the sense of \href{https://en.wikipedia.org/wiki/Zuhair_Nashed}{Nashed}}.
As we show below, this occurs for most $K\in\q K(\bb X,\bb Y)$, and condition \circled{3} for well-posedness in the sense of Hadamard is violated for their pseudoinverse $K^\dagger$.
\begin{cor}\label{thm:inverse:discont}
	Let $K\in\q K(\bb X,\bb Y)$. If $\dim(\q R(K)) = \infty$, then $K^\dag$ is unbounded.
\end{cor}
\begin{pr}
	Assume that $K^\dag$ is bounded.
	Then,  \Cref{th:inverse:cont} implies that $\q R(K)$ is closed.
	Thus, the restriction $\tilde K \colon \q N(K)^\bot \to \q R(K)$ admits a bounded inverse $\tilde K^{-1}$ by \Cref{th:closedrange}.
	Since $K \tilde K^{-1} y = y$ for all  $y\in \q R(K)$ and $K\circ\tilde K^{-1}$ is compact, we obtain that $I\colon \q R(K)\to \q R(K)$ is compact.
	This is only possible if $\dim(\q R(K))<\infty$, which contradicts the assumption.
\end{pr}

\subsubsection{Exercises}
\begin{exerciseB}\label{ex:ReverseStatementInverse}
	Let $A \in \q L(\bb X, \bb Y)$ and consider a linear operator
	\begin{equation}
	B\colon \mathcal{R}(A) \oplus \mathcal{R}(A)^\perp \rightarrow \mathcal{N}(A)^\perp
	\end{equation}
	that satisfies
	\begin{equation}
	BA = P_{\mathcal{N}(A)^\perp} \quad \text{and} \quad AB = P_{\overline{\mathcal{R}(A)}},
	\end{equation}
	Show that $B = A^\dagger$.
	Thus, the properties \circled{1}-\circled{4} in \Cref{lem:inverse}  uniquely characterize $A^\dagger$.
\end{exerciseB}

\begin{Solution}
	{\emph{Step 1:}} We show first that by the two properties $A^\dagger$ is uniquely defined on $\q R(A)$. Let $y\in \q R(A)$. Then there exists a unique $u \in \q N(A)^\perp$ such that $y=Au$. Then 
	\begin{equation}
	A^\dagger y = A^\dagger A u = P_{\q N(A)^\perp}u = BAu = By.
	\end{equation}
	{\emph{Step 2:}} We now show  that by the two properties $A^\dagger$ is uniquely defined on $\q R(A)^\perp $. If $y\in \q R(A)^\perp$ then 
	\begin{equation}
	A A^\dagger y = P_{\q R(A)} y = 0=ABy.
	\end{equation}
	Consequently $B f\in \q N(A)$. But $B f$ also lies  in $ \q N(A)^\perp$ by definition. Hence, $A^\dagger y = 0 = By$.
	
	{\emph{Step 3:}} Steps 1-2 completely determine $A^\dagger$ on the orthogonal decomposition $\q R(A) \oplus \q R(A)^\perp$. So $A^\dagger$ is uniquely defined. 
\end{Solution}

\begin{exerciseB}
	Using the result of \Cref{ex:ReverseStatementInverse}, proof the explicit formulas for the Moore--Penrose inverse of $ A \in \R^{M \times N} $ given below.
	
	\begin{itemize}
		\item[(i)] If $ \rank(A) = N $ and $M \geq N$, the Moore--Penrose inverse is given by $A^\dagger = (A^\tT A)^{-1} A^\tT$.
		
		\item[(ii)] If $ \rank(A) = M $ and $M \leq N$, the Moore--Penrose inverse is given by $A^\dagger = A^\tT (A A^\tT)^{-1}$.
	\end{itemize}
	These formulas are very useful for practical computations.
\end{exerciseB}

\begin{Solution}
	\begin{itemize}
		\item[(i)] By \Cref{thm:inverse:normalen}, the minimal norm solution of $Ax = y$ is the minimal norm solution of $A^\tT A x = A^\tT y$, that is 
		\begin{equation}
		A^\dagger y = (A^\tT A)^\dagger A^\tT. 
		\end{equation}  
		Since $\rank(A) = N$, $A^\tT A$ is invertible and $(A^\tT A)^\dagger = (A^\tT A)^{-1}$ giving the assertion.		
		
		\item[(ii)] Assume that $A \in \R^{M \times N} $ has full row rank, that is, $\rank(A) = m$ and $m \leq n$. Then $A A^\tT$ is invertible.
		We define
		\begin{equation}
		B \coloneqq A^\tT (A A^\tT)^{-1}.
		\end{equation}
		We now verify that $B$ satisfies
		\begin{enumerate}
			\item $B \colon \mathcal{R}(A) \oplus \mathcal{R}(A)^\perp \rightarrow \mathcal{N}(A)^\perp$,
			\item $B A = P_{\mathcal{N}(A)^\perp} $,
			\item $AB = \left. P_{\overline{\mathcal{R}(A)}} \right|_{\mathcal{D}(A^\dagger)} $.
		\end{enumerate}
		
		Since $A$ has full row rank, $\q R(A) = \R^M$, so $\mathcal{R}(A) \oplus \mathcal{R}(A)^\perp= \R^M$ and $P_{\overline{\q R(A)}} = P_{\q R(A)} = I_M$. We show that $\q R(B) = \q N(A)^\perp$. For any $x = By \in \q R(B)$ and $u \in \q N(A)$ we have
		\begin{equation}
		\inner{x}{u} 
		= \inner{A^\tT (A A^\tT)^{-1} y}{u}
		= \inner{(A A^\tT)^{-1} y}{A u}
		= 0.
		\end{equation}
		Hence, $x \in \q N(A)^{\perp}$ and $\q R(B) \subseteq \q N(A)^{\perp}$.
		On the other hand, for $x \in \q N(A)^{\perp} = \q R(A^\tT)$ using bijectivity of $A A^\tT$, we can find $y,z \in \R^M$ such that $x = A^\tT y$ and $y = (A A^\tT)^{-1}z$. Then, $x = A^\tT (A A^\tT)^{-1} z = B z$ and $\q N(A)^{\perp} \subseteq \q R(B)$.
		
		Next, we have
		\begin{equation}
		A B = A A^{\tT} (A A^{\tT})^{-1} = I_M
		\end{equation}
		yielding the second point. 
		
		Define
		\begin{equation}
		B A = A^\tT (A A^\tT)^{-1} A =: P,
		\end{equation}
		and we show that $P$ is a projection onto $\q N(A)^\perp$. Since
		\begin{equation}
		P^2 = A^\tT (A A^\tT)^{-1} A A^\tT (A A^\tT)^{-1} A = A^\tT (A A^\tT)^{-1} A = P,
		\end{equation}
		and
		\begin{equation}
		P^\tT = (A^\tT (A A^\tT)^{-1} A)^\tT = A^\tT (A A^\tT)^{-1} A = P,
		\end{equation}
		we see that $P$ is symmetric and idempotent. Hence $P$ is an orthogonal projection. 
		Moreover, for any $x \in \q N(A)$, we have $P x = B A x = B 0 = 0$ and for all $x \in \q N(A)^\perp = \q R(A^\tT)$, there exists $y \in \bb R^M$ such that $x = A^\tT y$. Then,
		\begin{equation}
		P x = A^\tT (A A^\tT)^{-1} A A^\tT y = A^\tT y = x.
		\end{equation}
		Hence, $B A = P = P_{\q N(A)^\perp}$.
		Thus, by \Cref{ex:ReverseStatementInverse}, we get that $B = A^\dagger$.
	\end{itemize}
\end{Solution}

\begin{exerciseB}\label{ex: some pseudoinverses}
	Let $M, M \in \N$ with $M \geq N \geq 2$. Compute $A^\dagger$ for the following matrices:
	\begin{itemize}
		\item[(i)] $A = (1, 1, \dots, 1) \in \R^{1 \times N}$;
		\item[(ii)] $A = \diag(a_1, \dots, a_n) \in \R^{N \times N}$ with $a_j \in \R$ for $j \in \{1, \dots, N\}$;
		\item[(iii)] $A \in \mathbb{R}^{M \times N}$ with $A^\tT A = I_N$.
	\end{itemize}
\end{exerciseB}

\begin{Solution}
	{\bf  Solution (i):} The matrix $A$ has full row-rank. Hence, the Moore-Penrose Inverse is given by
	\begin{equation}
	A^\dagger = A^\tT(AA^\tT)^{-1} = A^\tT (N)^{-1} = \frac{1}{N} A^\tT = \frac{1}{N} \begin{pmatrix}1 \\ \vdots\\1\end{pmatrix} \in \R^{N\times 1}.
	\end{equation}
	
	{\bf  Solution (ii):} 
	Set $\q J \coloneqq \{ j : a_j \neq 0 \}$. The minimizers of 
	\begin{equation}
		\|Ax -y\|^2_2 = \sum_{j = 1}^N (a_j x_j - y_j)^2 = \sum_{j \in \q J} (a_j x_j - y_j)^2
	\end{equation}
	are given by $x_j = a_j^{-1} y_j$, $j \in \q J$ and $x_j$ arbitrary for $j \notin \q J$. 
	Then, the minimum norm solution is given by $x\dagger = Dy$.
	
	{\bf  Solution (iii):} Since the matrix $A \in \R^{M \times N}$ fulfills $A^\tT A=I_N$, it has full row-rank, i.e., $\rank(A)=N$. Thus, the Moore-Penrose Inverse is given by
	\begin{equation}
	A^\dagger = (A^\tT A)^{-1} A^\tT = A^\tT.
	\end{equation}
\end{Solution}

\begin{exerciseB}
	Let $a, b \in \mathbb{R}$ with $a < b$.
	Compute $A^\dagger$ for $A\colon L^2([a,b]) \to \mathbb{R}$ with $Au = \int_a^b u(x) \diff x$.
\end{exerciseB}

\begin{Solution}
	For, $c\in\mathbb{R}$, we set $B: \R \to \textrm{span}\{\1\}$ as 
	\begin{equation}
	B c = \frac{c}{b-a}  \1_{[a,b]}.
	\end{equation}
	One can easily check that $AB = I_\mathbb{R} = P_{\overline{\q R(A)}}$ and $B A = \textrm{span}\{\1\} = \mathcal{N}(A)^\perp$. 
	This shows that $B=A^\dagger$ is the Moore-Penrose Inverse by \Cref{ex:ReverseStatementInverse}.
\end{Solution}

\subsection{Singular Value Decomposition}

In this section, we compute the Moore--Penrose inverse $K^\dagger$ of a compact operator $K\in \q K(\bb X,\bb Y)$ based on the singular value decomposition, which is introduced below.
\begin{thm}{Singular Value Decomposition (SVD)}{svd}
	For every $K\in\q K(\bb X,\bb Y)$, there exists a \emphi{singular system} $\{(\sigma_n,u_n,v_n)\}_{n\in\N}$ (possibly finitely terminating) consisting of
	\begin{enumerate}
		\item a null sequence $\{\sigma_n\}_{n\in\N}$ with $\sigma_1\geq \sigma_2\geq \dots >0$;
		\item an ONB $\{u_n\}_{n\in\N}\subset \bb Y$ of $\overline{\q R(K)}$;
		\item an ONB $\{v_n\}_{n\in\N}\subset \bb X$ of $\overline{\q R(K^*)}$,
	\end{enumerate}
	with the following properties.
	First, we have that
	\begin{equation}\label{eq:sing_vec}
		K v_n = \sigma_n u_n\quad\text{and}\quad K^* u_n = \sigma_n v_n \qquad\text{for all }n\in\N.
	\end{equation}
	Additionally, we have the \emphi{singular value decomposition} (SVD)
	\begin{equation}\label{eq:svd}
		Kx = \sum_{n\in\N} \sigma_n \inner{x}{v_n}_{\bb X} u_n \qquad\text{for all }x\in \bb X.
	\end{equation}
\end{thm}
\begin{pr}
	As $K^*K\in \q K(\bb X, \bb X)$ is selfadjoint, by \Cref{th:spektral} we can find a null sequence $\{\lambda_n\}_{n\in\N}\subset \R\setminus \{0\}$ and an ONB $\{v_n\}_{n\in\N}\subset \bb X$ of $\overline{\q R(K^*K)}$ with
	\begin{equation}
		K^*Kx = \sum_{n\in\N} \lambda_n \inner{x}{v_n}_{\bb X} v_n \qquad \text{for all }x\in \bb X.
	\end{equation}
	
	Next, we show $\overline{\q R(K^* K)} = \overline{\q R(K^*)}$.   
	By construction $\overline{\q R(K^* K)} \subseteq \overline{\q R(K^*)}$.
	On the other hand, since $\bb Y = \q N(K^*) \oplus \overline{\q R(K)}$, we get that $K^* y = K^* P_{ \overline{\q R(K)}} y$.
	Moreover, for any $y \in \bb Y$ there exists $\{x_n\}_{n \in \bb N} \subset \bb X$ with $\lim_{n \to \infty} K x_n  = P_{ \overline{\q R(K)}} y$. Consequently, $K^* y = \lim_{n \to \infty} K^* K x_n$ and $\q R(K^*) \subseteq \overline{\q R(K^* K)}$.
	This yields the desired equality.
	
	Since $\lambda_n = \inner{\lambda_n v_n}{v_n}_{\bb X} =\inner{K^*K v_n}{v_n}_{\bb X} = \norm{Kv_n}_{\bb Y}^2 >0$, we can now define
	\begin{equation}
		\sigma_n \coloneqq \sqrt{\lambda_n} >0 \qquad\text{ and }\qquad u_n \coloneqq \sigma_n^{-1}{Kv_n} \in \q R(K) \subset \bb Y.
	\end{equation}
	The system $\{u_n\}_{n\in\N} \subset \bb Y$ is orthonormal since
	\begin{equation}
		\inner{u_i}{u_j}_{\bb Y} = \frac1{\sigma_i\sigma_j}\inner{Kv_i}{Kv_j}_{\bb Y} = \frac1{\sigma_i\sigma_j} \inner{K^*Kv_i}{v_j}_{\bb X} = \frac{\lambda_i}{\sigma_i\sigma_j}\inner{v_i}{v_j}_{\bb X} = \delta_{i j}.
	\end{equation}
	Further, it holds	$K^*u_n = \sigma_n^{-1} K^*K v_n = \sigma_n v_n$ and $Kv_n = \sigma_n u_n$, namely that \eqref{eq:sing_vec} is satisfied.
	Hence, it remains to show that \eqref{eq:svd} holds and that $\{u_n\}_{n\in\N}$ is an ONB of $\overline{\q R(K)}$.
	
	Since $\{v_n\}_{n\in\N}$ is an ONB of $\q N(K)^\bot$, we get that
	\begin{equation}\label{eq:PreSingDeco}
		Kx = K P_{\q N^\bot} x = K \biggl(\sum_{n\in\N} \inner{x}{v_n}_{\bb X} v_n\biggr) \qquad\text{for all }x\in \bb X.
	\end{equation}
	From equation \eqref{eq:PreSingDeco}, we obtain the SVD \eqref{eq:svd} by interchanging $K$ and the summation.
	We justify this frequently applied step once in detail.
	First, we set $x_N \coloneqq \sum_{n=1}^N \inner{x}{v_n}_{\bb X} v_n$ for any $x\in \bb X$ and $N\in\N$.
	Then, we have $x_N \to P_{\q N^\bot} x$ as $N\to \infty$ and hence by continuity of $K$ also
	\begin{equation}\label{eq:svd:operator}
			Kx = \lim_{N\to\infty} Kx_N = \lim_{N\to\infty} \sum_{n=1}^N \inner{x}{v_n}_{\bb X} Kv_n = \sum_{n\in\N} \inner{x}{v_n}_{\bb X} Kv_n.
	\end{equation}
	Thus, using the relation \eqref{eq:sing_vec}, we get
	\begin{equation}\label{eq:SVD_proof}
		Kx = \sum_{n\in\N} \inner{x}{v_n}_{\bb X} Kv_n = \sum_{n\in\N} \inner{x}{v_n}_{\bb X} \sigma_n u_n = \sum_{n\in\N} \inner{Kx}{u_n}_{\bb X} u_n.
	\end{equation}
	The second equality is \eqref{eq:svd}.
	
	To show that $\{u_n\}_{n \in \bb N}$ is an ONB of $\overline{\q R(K)}$, we consider $\q U = \overline{\{ \sum_{n =1}^N \alpha_n u_n  : \alpha_n \in \bb R, \ N \in \bb N \}}$.
	Since $u_n \in \q R(K)$ by construction, we have $\q U \subset \overline{\q R(K)}$.
	Furthermore, equation \eqref{eq:SVD_proof} implies that $\q R(K) \subseteq \q U$.
	Thus,  $\q U = \overline{\q R(K)}$ and $\{u_n\}_{n \in \bb N}$ is an ONB of $\overline{\q R(K)}$. 
\end{pr}
From our construction, we directly obtain that $\{(\sigma_n,v_n,u_n)\}_{n\in\N}$ is a singular system of $K^*$ since
\begin{equation}
	\label{eq:svd_adj}
	K^* y = \sum_{n\in\N} \inner{y}{u_n}_{\bb Y} K^* u_n = \sum_{n\in\N} \sigma_n \inner{y}{u_n}_{\bb Y} v_n \qquad\text{for all }y\in \bb Y.
\end{equation}

\begin{remark}\label{rem:svd matrix}
For a matrix $K \in \R^{M \times N}$ with $r = \rank(K)$, our SVD construction can be written in the form $K = U \Sigma V^\tT$, where 
\begin{itemize}
	\item $U \in \R^{M \times r}$ and $V \in \R^{N \times r}$ are \emphi{semi-orthogonal}, namely $U^\tT U = V^\tT V= I_r$;
	\item $\Sigma \in \mathbb{R}^{r \times r}$ is a diagonal matrix with $\Sigma_{i,i} = \sigma_i > 0$.
\end{itemize}
By using an orthonormal basis completion in $\R^M$ and $\R^N$, respectively, we can also write $K = \tilde U \tilde \Sigma \tilde V^\tT$ with 
\begin{equation}
\tilde U = \begin{bmatrix}
	U & U_\perp 
\end{bmatrix}  
\in \R^{M \times M},
\quad 
\tilde V = \begin{bmatrix}
	V & V_\perp 
\end{bmatrix} 
\in \R^{N \times N},
\quad \text{and} \quad
\tilde \Sigma = \begin{bmatrix}
	\Sigma & 0 \\
	0 & 0 
\end{bmatrix}
\in \R^{M \times N}.
\end{equation}
In this case $\tilde U$ and $\tilde V$ are \emph{orthogonal} matrices.
The SVD extends to the complex case.
In particular, for $K \in \bb C^{M \times N}$ the SVD is $K = U \Sigma V^H$ with $V^H$ denoting the conjugate transpose.
\end{remark}

Now, we characterize all $y\in \q D(K^\dag)$ and the Moore--Penrose inverse $K^\dag$ itself using the SVD.
\begin{thm}{}{inverse:picard}
	Let $K\in\q K(\bb X,\bb Y)$ with singular system $\{(\sigma_n,u_n,v_n)\}_{n\in\N}$ and $y\in \bb Y$. Then $y\in \q D(K^\dag)$ if and only if the \emphi{Picard condition}
	\begin{equation}
		\label{eq:inverse:picard}
		\sum_{n\in\N} \sigma_n^{-2} |\inner{y}{u_n}_{\bb Y}|^2 <\infty
	\end{equation}
	is satisfied.
	In this case, it holds 
	\begin{equation}\label{eq:inverse:picard_pseudo}
		K^\dag y = \sum_{n\in\N} \sigma_n^{-1} \inner{y}{u_n}_{\bb Y} v_n.
	\end{equation}
\end{thm}
\begin{pr}
	$\Rightarrow$: For any $y\in \q D(K^\dag)$, there exists $x\in \bb X$ with $Kx= P_{\overline{\q R(K)}} y$.
	From \eqref{eq:sing_vec} and $\{u_n\}_{n \in \N}$ being an ONB of $\overline{\q R(K)}$, we get for all $n\in\N$ that
	\begin{equation}\label{eq:coefficientsY}
		\inner{y}{u_n}_{\bb Y} = \inner{P_{\overline{\q R(K)}} y}{u_n}_{\bb Y} = \inner{Kx}{u_n}_{\bb Y} = \inner{x}{K^*u_n}_{\bb X} = \sigma_n \inner{x}{v_n}_{\bb X}.
	\end{equation}
	Inserting \eqref{eq:coefficientsY} into the Bessel inequality \eqref{eq:funktan:bessel} implies
	\begin{equation}
		\sum_{n\in\N} \sigma_n^{-2} |\inner{y}{u_n}_{\bb Y}|^2 = \sum_{n\in\N} |\inner{x}{v_n}_{\bb X}|^2 \leq \norm{x}_{\bb X}^2 < \infty.
	\end{equation}
	
	$\Leftarrow$: Let $y\in \bb Y$ satisfy \eqref{eq:inverse:picard}.
	Then $\{\sum_{n=1}^N \sigma_n^{-2}|\inner{y}{u_n}_{\bb Y}|^2\}_{N\in\N}$ is a Cauchy sequence.
	Since $\{v_n\}_{n\in\N}$ is orthonormal, we can define a Cauchy sequence $\{x_N\}_{N\in\N}$ via
	\begin{equation}
		x_N\coloneqq \sum_{n=1}^N \sigma_n^{-1} \inner{y}{u_n}_{\bb Y}v_n \in \overline{\q R(K^*)}= \q N(K)^\bot.
	\end{equation}
	By the closedness of $\q N(K)^\bot$, we get that $x \coloneqq  \lim_{N \to \infty} x_N  \in  \q N(K)^\bot$.
	Interchanging $K$ and the summation (similar as for \eqref{eq:svd:operator}), we get that
	\begin{equation}\label{eq:ProjProb}
		Kx = \sum_{n\in\N} \sigma_n^{-1} \inner{y}{u_n}_{\bb Y}Kv_n = \sum_{n\in\N} \inner{y}{u_n}_{\bb Y}u_n = P_{\overline{\q R}} y.
	\end{equation}
	This implies that $y\in \q D(K^\dag)$.
	
	Finally, \eqref{eq:ProjProb}, $x \in  \q N(K)^\bot$, and \Cref{th:inverse:minnorm} imply that $x = K^\dagger y$ and \eqref{eq:inverse:picard_pseudo} holds.
\end{pr}

\Cref{th:inverse:picard} states that $y \in \q D(A^\dagger)$ if and only if $\{\inner{y}{u_n}_{\bb Y}\}_{n \in \N}$ decays fast enough compared to $\{\sigma_n\}_{n \in \N}$.
In particular, we get that $K^\dag$  can only be continuous if $\{\sigma_n\}_{n \in \N}$ is not a null sequence.
Moreover, \eqref{eq:inverse:picard_pseudo} shows how the errors propagate in $x^\dag$.
If $y^\delta = y+\delta u_n$ for some $\delta>0$ and $n\in\N$, then
\begin{equation}\label{eq: SVD error scaling}
	\norm{K^\dag y^\delta -K^\dag y}_{\bb X} = \delta \norm{K^\dag u_n}_{\bb X} = \sigma_n^{-1}\delta \to \infty \quad\text{as }n\to \infty.
\end{equation}
If the $\sigma_n$ decay rapidly, high error amplification is more likely.
Thus, one distinguishes
\begin{itemize}
	\item \emphi{moderately ill-posed} problems, for which there exist $c,r>0$ with $\sigma_n \geq c n^{-r}$ for all $n\in\N$, i.e., $\sigma_n$ decays at most polynomially;
	\item otherwise, the problem is  \emphi{severely ill-posed}.
	If there exist  $c,r>0$ with $\sigma_n \leq c e^{-n^r}$ for all $n\in\N$, the problem is called \emph{exponentially ill-posed}.
\end{itemize}
If $\dim(\q R(K))< \infty$, then $\{\sigma_n\}_{n\in\N}$ terminates finitely and $\Vert K^\dag \Vert = \frac{1}{\sigma_\mathrm{min}}$.
Then, one typically uses the \emphi{condition number} $\kappa(K) = \frac{\sigma_\mathrm{max}}{\sigma_\mathrm{min}}$ to classify the hardness of the inverse problem.

The SVD is a valuable analytical tool.
Unfortunately, there is no generic approach or recipe for its computation.
As an example, we consider a differentiation operator.
\begin{example}\label{ex:integration_svd}
	Let $\bb X=L^2([0,1])$ and $K\in \q K(\bb X,\bb X)$ with
	\begin{equation}
		[Kx](t) = \int_0^1 k(s,t) x(s) \diff s, \qquad\text{where}\qquad
		k(s,t) =\begin{cases} 1 & \text{if }s\leq t,\\0 & \text{else}.\end{cases}
	\end{equation}
	Then, for $y \in C^1([0,1])$ with $y(0)=0$, we have
	\begin{equation}
		[Ky' ](t) = \int_0^t y' (s) \diff s = y(t) - y(0) = y(t),
	\end{equation}
	namely that $y'$ solves the equation $Kx=y$.
	In general, a solution may not exist.
	
	Now, we want to construct the SVD and then $K^\dag$ based on the spectral decomposition of $K^*K$. 
	Using the definition of $k$ and \eqref{eq: int op adjoint}, we obtain the adjoint 
	\begin{equation}
		[K^*y](t) = \int_0^1 k(t,s) y(s) \diff s = \int_t^1 y(s) \diff s.
	\end{equation}
	Thus, any eigenvector $v\in L^2([0,1])$ of $K^*K$ with eigenvalue $\lambda>0$ must satisfy
	\begin{equation}
		\label{eq:svd_diff_ev}
		\lambda v(t) = \int_t^1\int_0^s v(r)\diff r\diff s.
	\end{equation}
	From \eqref{eq:svd_diff_ev}, we can derive multiple necessary conditions.
	Inserting $t=1$ into \eqref{eq:svd_diff_ev} gives $v(1) = 0$.
	Next, differentiating \eqref{eq:svd_diff_ev} with respect to $t$ yields
	\begin{equation}\label{eq:DiffOnce}
		\lambda v'(t) = \frac{\diff }{\diff t}\left( -\int_1^t\int_0^s v(r)\diff r\diff s\right) = -\int_0^t v(r)\diff r.
	\end{equation}
	Inserting $t=0$ into \eqref{eq:DiffOnce} leads to $v'(0) = 0$.
	Finally, differentiating \eqref{eq:DiffOnce} a second time with respect to $t$ yields that any eigenvector $v$ solves the harmonic oscillator equation
	\begin{equation}\label{eq:HarmonicOscillator}
		\lambda v''(t) + v(t) = 0.
	\end{equation}
	Using the boundary condition  $v(1) = 0$ and $v'(0) = 0$, we can determine the constants $c_1,c_2$ in the general solution of \eqref{eq:HarmonicOscillator}, which is given by
	\begin{equation}\label{eq:DiffEqSol}
		v(t) = c_1 \sin(\sigma^{-1} t) + c_2 \cos(\sigma^{-1} t), \qquad \sigma \coloneqq \sqrt{\lambda}.
	\end{equation}
	From $v'(0) = 0$, we get $c_1 = 0$.
	Since $c_2 = 0$ would lead to the eigenvector $v=0$, we conclude from $v(1) = 0$ that $\cos(\sigma^{-1})=0$.
	Since the latter can only hold if $\sigma_n = \frac{2}{(2n-1)\pi}$ for some $n \in \N$, any (normalized) eigenvector of $K^*K$ must be of the form 
	\begin{equation}\label{eq:Eigenfunctions}
		v_n(t) = \sqrt{2}\cos\left((n-\tfrac12)\pi\,t\right),\qquad n\in\N.
	\end{equation}
	By plugging the $v_n$ into \eqref{eq:svd_diff_ev}, we verify that they are indeed eigenvectors for $\lambda_n = \sigma_n^2$.
	Since the $\lambda_n$ are pairwise different, we further have that $\{v_n\}_{n \in \N}$  is an orthonormal system.
	Moreover, one can show that $\{v_n\}_{n \in \N}$ is an ONB of $L^2([0,1])$, namely that no other eigenvectors can exist.
	
	Thus, we can proceed as in the proof of \Cref{th:svd} by defining the left singular vectors
	\begin{equation}
		u_n \coloneqq \sigma_n^{-1} Kv_n = (n-\tfrac12)\pi\int_0^t \sqrt{2}\cos\left((n-\tfrac12)\pi\,s\right)\diff s = \sqrt{2} \sin\left((n-\tfrac12)\pi\,t\right), 
	\end{equation}
	which leads to the singular system $\{(\sigma_n,u_n,v_n)\}_{n\in\N}$ of $K$.
	Since $\sigma_n =\mathcal{O}(\frac1n)$, we have that the inverse problem $Kx=y$ is moderately ill-posed.
	Moreover, for $y\in L^2([0,1])$ satisfying the Picard condition
	\begin{equation}\label{eq:PicardSeries}
		\sum_{n\in\N} \tfrac14{(2n-1)^2\pi^2} |\inner{y}{u_n}_{L^2}|^2 <\infty,
	\end{equation}
	\Cref{th:inverse:picard} implies that 
	\begin{equation}\label{eq:MoorePenroseDiff}
		K^\dag y =
		\sum_{n\in\N} \left(n-\tfrac12\right)\pi \inner{y}{u_n}_{L^2} v_n = \sum_{n\in\N} \inner{y}{u_n}_{L^2} u_n'.
	\end{equation}
	If $y \in  C^1([0,1])$, then the second series in \eqref{eq:MoorePenroseDiff} converges even uniformly.
	In this case, we can exchange summation and differentiation to obtain $K^\dag y = y'$ as expected.
\end{example}

\subsubsection{Exercises}

\begin{exerciseB}
	Let $A \in \R^{M \times N}$, $I \in \R^{N \times N}$ be the identity matrix, $\delta > 0$, and
		\begin{equation}
		K = \begin{bmatrix}
			A \\
			\sqrt{\delta} I
		\end{bmatrix} \in \R^{(M+N) \times N}.
		\end{equation}
		Show that the singular values of $K$ satisfy $\sigma_j(K) \geq \sqrt{\delta}, \quad j = 1, \dots, N$.
\end{exerciseB}
\begin{Solution}
	The singular values of $K$ are the square roots of the eigenvalues of $K^\tT K$. We have
	\begin{align}
		K^\tT K = A^\tT A + \delta I.
	\end{align}
	Note that $K^\tT K$ and $A^\tT A$ are symmetric positive semidefinite matrices and have non-negative eigenvalues. Then, for all eigenvectors $x_j \in \R^N$ of $K^* K$ with corresponding eigenvalues $\sigma_j^2(K)$ we have
	\begin{equation}
	\sigma_j^2(K) = \sigma_j^2(K) \inner{x_j}{x_j} 
	= \inner{K^\tT K x_j}{x_j} 
	= \inner{A^\tT A x_j}{x_j} + \delta \| x_j \|^2 \ge \delta.
	\end{equation}
\end{Solution}

\begin{exerciseB}[SVD of Convolution]\label{ex: svd convolution}
	Let $\bb X = \bb Y = \R^N$ and let $g \in \R^N$ be a convolution kernel.
	Using that the scaled discrete Fourier transform $\frac{1}{\sqrt{N}} F_N$ is unitary, find the SVD of the convolution matrix $C(g)$. Furthermore, find the pseudoinverse $C(g)^\dagger$.
	\textbf{Hint:} Use \Cref{ex: circulant matrix} and \Cref{ex: convolution theorem}.
\end{exerciseB}
\begin{Solution}
	By \Cref{ex: circulant matrix} and \Cref{ex: convolution theorem}, we can rewrite 
	\begin{equation}
	C(g) f = f *_N g = F_N^{-1} \diag(F_N g) F_N f =  (\tfrac{1}{\sqrt{N}} F_N)^{-1} \diag(F_N g) ((\tfrac{1}{\sqrt{N}} F_N)^{-1})^{H} f.
	\end{equation}
	Let $S$ be a diagonal matrix with entries 
	\begin{equation}
	S_{j,j} = \begin{cases}
		(F_N g)_j / |(F_N g)_j|, & (F_N g)_j \neq 0, \\
		1, & (F_N g)_j =0,
		\end{cases}
	\end{equation}
	$j = 0, \ldots, N-1$. Then, the SVD of $C(g)$ is given by $U = (\tfrac{1}{\sqrt{N}} F_N)^{-1}$, $\Sigma = \diag(|F_N g|)$ and $V = U S^H$.
	
	Consequently,
	\begin{equation}
	f = C(g)^\dagger y = F_N^{-1} \diag(F_N g)^\dagger F_N y,
	\end{equation}
	where $\diag(F_N g)^\dagger$ is as in \Cref{ex: some pseudoinverses}.
\end{Solution}

\begin{exerciseB}[Eckart--Young Theorem]
	Suppose $ A, B \in \mathbb{R}^{m \times n}$ and $ \mathrm{rank}(B) = k \leq \mathrm{rank}(A) = r$.
	Show that for the the truncated SVD $ A_k $ (all singular values $\sigma_n$ with $n \geq k$ are set to zero) it holds 
	\begin{equation}
	\|A - B\|_2 \geq \|A - A_k\|_2 = \sigma_{k+1}.
	\end{equation}
	That is, $A_k$ is the \emph{best rank-$k$ approximation} to $ A $ with respect to the $\ell_2$-norm.
\end{exerciseB}
\begin{Solution}
	We begin by clarifying the notation. Let $A \in \mathbb{R}^{m \times n}$ with $\operatorname{rank}(A) = r$, and let $B \in \mathbb{R}^{m \times n}$ be any matrix with $\operatorname{rank}(B) = k \leq r$. Here, $A_k$ denotes the rank-$k$ approximation of $A$ obtained by keeping only the $k$ largest singular values, and $\|\cdot\|$ refers to the $L^2$ operator norm (spectral norm).
	
	By \cref{rem:svd matrix}, we can rewrite $A$ as
	\begin{equation}
	A = U\Sigma V^\tT,
	\end{equation}
	where:
	\begin{itemize}
		\item $U \in \R^{M \times M}$ is an orthogonal matrix with columns $u_1, u_2, \ldots, u_M$
		\item $V \in \R^{N \times N}$ is an orthogonal matrix with columns $v_1, v_2, \ldots, v_N$
		\item $\Sigma \in \R^{M \times N}$ is a diagonal matrix with singular values $\sigma_1 \geq \sigma_2 \geq \cdots \geq \sigma_r > 0$ and zeros elsewhere
	\end{itemize}
	
	The matrix $A_k$ is given by
	\begin{equation}
	A_k = U\Sigma_k V^\tT,
	\end{equation}
	where $\Sigma_k$ is $\Sigma$ with only the $k$ largest singular values retained and the rest set to zero. Now, let us consider the difference $A - A_k$ 
	\begin{align}
		A - A_k 
		= U\Sigma V^\tT - U\Sigma_k V^\tT
		& = U(\Sigma - \Sigma_k)V^\tT
	\end{align}
	Since $U$ and $V$ are orthogonal matrices, they preserve the norm. Therefore, we get
	\begin{align}
		\|A - A_k\|_2 &= \|U(\Sigma - \Sigma_k)V^\tT\|_2 
		= \|\Sigma - \Sigma_k\|_2 
		= \sigma_{k+1}.
	\end{align}
	This establishes the second part of the statement.
	
	For the inequality, we define the following subspace: $\q R(V_k) = \operatorname{span}\{v_1, v_2, \ldots, v_k\}$, the span of the first $k$ right singular vectors of $A$.
	Since $\rank(B) = k$, we have $\dim(\q N(B)) = N-k$. By a dimensionality argument, there must exist a unit-length vector 
	\begin{equation}
	w\in \q N(B) \cap \q R(V_{k+1}),
	\end{equation} 
	since otherwise $\q N(B)$ and $\q R(V_{k+1})$ would together span a subspace of dimension $N-k + k+1 = N+1 > N$, which contradicts $\q N(B),\q R(B) \subset \R^N$. We compute
	\begin{align}
		\|A -B\|^2_2 
		& \geq \|(A-B)w\|_2^2 
		= \|Aw\|_2^2 
		= \|\sum_{i=1}^{k+1} \sigma_i u_i v_i^\tT w \|_2^2 
		= \sum_{i=1}^{k+1} \sigma_i^2 (v_i^\tT w)^2 
		\geq \sigma_{k+1}^2 \sum_{i=1}^{k+1} (v_i^\tT w)^2 
		= \sigma_{k+1}^2.
	\end{align}
	Hence, we showed that
	\begin{equation}
	\|A -B\|_2 \geq  \sigma_{k+1} = \|A-A_k\|_2.
	\end{equation}
\end{Solution}

\begin{exerciseB} \label{ex:SVD conv L2}
	Let $\Omega \coloneq [0,1]^{2}$, $k \in L^{2}(\Omega)$ and $\tilde{k}$ be the periodic extension of $k$ with
	\begin{equation}
	\tilde{k}(z)=\left\{\begin{array}{ll}
		k(z) & z \in \Omega \\
		k(\bmod (z, 1)) & z \in \R^{2} \backslash \Omega
	\end{array},\right.
	\end{equation}
	where $\bmod (z, 1)=\left(\bmod \left(z_{1}, 1\right), \bmod \left(z_{2}, 1\right)\right)$.
	
	\begin{itemize}
		\item[(a)]  Compute the SVD of the convolution  $K \colon L^{2}(\Omega) \rightarrow$ $L^{2}(\Omega)$ with
		\begin{equation}
		[K u](x)\coloneq \int_{\Omega} \tilde{k}(x-y) u(y) \diff y .
		\end{equation}
		
		\textbf{Hint:} \textit{The set $\{\varphi_{m, n}\}_{n,m \in \mathbb Z}$ with $\varphi_{m, n}\left(x_{1}, x_{2}\right)=$ $\exp \left(-2 \pi i\left(m x_{1}+n x_{2}\right)\right)$ is an ONB of $L^{2}(\Omega)$.}
		\item [(b)] Investigate if the associated deconvolution problem is ill-posed for the specific choices
		\begin{itemize}
			\item[i)] $k\left(x_{1}, x_{2}\right)=\frac{1}{h^{2}} \1_{\left[-\frac{1}{2}, \frac{1}{2}\right]}\left(\frac{x_{1}-1 / 2}{h}\right) \1_{\left[-\frac{1}{2}, \frac{1}{2}\right]}\left(\frac{x_{2}-1 / 2}{h}\right)$ for $0<h<1$;
			\item[ii)] $k\left(x_{1}, x_{2}\right)=\varphi\left(x_{1}\right) \varphi\left(x_{2}\right)$ with $\varphi(x)\coloneq \left\{\begin{array}{ll}\exp \left(-\frac{1}{1 / 4-(x-1 / 2)^{2}}\right) & x \in [0,1] \\ 0 & \text { else }\end{array}\right.$.
		\end{itemize}
	\end{itemize}
\end{exerciseB}

\begin{Solution}
	{\bf (a):} 
	We are given the convolution operator $K \colon L^{2}(\Omega) \rightarrow L^{2}(\Omega)$ defined as:
	\begin{equation}
	[K u](x) \coloneqq \int_{\Omega} \tilde{k}(x-y) u(y) \, \diff y .
	\end{equation}
	We compute $K$ on the orthonormal basis functions 
	\begin{equation}
	\varphi_{m, n}(x) = \exp(-2 \pi i (m x_1 + n x_2)).
	\end{equation}
	We have
	\begin{align} 
		[K \varphi_{m, n}](x) &= \int_{\Omega} \tilde{k}(x-y) \exp(-2 \pi i (m y_1 + n y_2)) \, \diff y.
	\end{align}
	Changing the variable to $z = x - y$ with $y = x - z$ and $\diff y = \diff z$ gives 
	\begin{align} 
		[K \varphi_{m, n}](x) &= \int_{\Omega} \tilde{k}(z) \exp(-2 \pi i (m (x_1 - z_1) + n (x_2 - z_2))) \, \diff z \\ 
		&= \exp(-2 \pi i (m x_1 + n x_2)) \int_{\Omega} \tilde{k}(z) \exp(2 \pi i (m z_1 + n z_2)) \, \diff z 
	\end{align}
	The integral is the $(m, n)$-th Fourier coefficient of $\tilde k$ given by 
	\begin{equation}
	\hat{k}_{m, n} = \int_{\Omega} k(z) \exp(2 \pi i (m z_1 + n z_2)) \, \diff z.
	\end{equation}
	Thus, we have
	\begin{equation}
	[K \varphi_{m, n}](x) = \hat{k}_{m, n} \varphi_{m, n}(x)
	\end{equation}
	This shows that $\varphi_{m, n}(x)$ are the {eigenfunctions} of the operator $K$, and the Fourier coefficients $\hat{k}_{m, n}$ are the corresponding {eigenvalues}.
	
	To find the SVD, we need to relate this eigenvalue equation ($K \varphi_{m, n} = \hat{k}_{m, n} \varphi_{m, n}$) to the SVD equation ($K v_i = \sigma_i u_i$).
	We identify the orthonormal eigenfunctions $\varphi_{m, n}$ as the {right singular functions} ($v_{m, n}$).
	So, $K v_{m, n} = \hat{k}_{m, n} v_{m, n}$.
	We must equate this to $K v_{m, n} = \sigma_{m, n} u_{m, n}$.
	By definition, singular values $\sigma_{m, n}$ must be real and non-negative. Since $\hat{k}_{m, n}$ can be complex, we take the $L^2(\Omega)$ norm of both sides, leading to 
	$\| K v_{m, n} \|_2 = \| \hat{k}_{m, n} v_{m, n} \|_2 = |\hat{k}_{m, n}| \| v_{m, n} \|_2 = |\hat{k}_{m, n}|$.
	Also, $\| K v_{m, n} \|_2 = \| \sigma_{m, n} u_{m, n} \|_2 = \sigma_{m, n} \| u_{m, n} \|_2 = \sigma_{m, n}$.
	Therefore, we deduce that $\sigma_{m, n} = |\hat{k}_{m, n}|$. We naturally ignore $m,n,$ with $|\hat{k}_{m, n}| =0 $ corresponding to the zero singular values. 
	With this, we can find $u_{m, n}$ such that $\sigma_{m, n} u_{m, n} = \hat{k}_{m, n} v_{m, n}$.
	
	Thus, the SVD of $K$ is given by:
	\begin{itemize}
		\item {right singular functions} $v_{m, n}(x) = \varphi_{m, n}(x) = \exp(-2 \pi i (m x_1 + n x_2))$;
		\item {singular values} $\sigma_{m, n} = |\hat{k}_{m, n}|$, and
		\item {left singular functions} $u_{m, n}(x) = \frac{\hat{k}_{m, n}}{|\hat{k}_{m, n}|} \varphi_{m, n}(x)$ (for $\hat{k}_{m, n} \ne 0$). 
	\end{itemize}

	{\bf (b):} 
	To show the deconvolution problem is {ill-posed}, it suffices to show the singular values decay to zero sufficiently fast. We investigate the decay of $\sigma_{m, n} = |\hat{k}_{m, n}|$ for the given kernels.
	
	\begin{enumerate}[label=\roman*)]
		\item  
		The kernel
		\begin{equation}
			k\left(x_{1}, x_{2}\right)=\frac{1}{h^{2}} \1_{\left[-\frac{1}{2}, \frac{1}{2}\right]}\left(\frac{x_{1}-1 / 2}{h}\right) \1_{\left[-\frac{1}{2}, \frac{1}{2}\right]}\left(\frac{x_{2}-1 / 2}{h}\right).
		\end{equation}
		represents an averaging filter. The integral for $\hat{k}_{m, n}$ is separable:
		\begin{equation}
			\hat{k}_{m, n} = \left( \frac{1}{h} \int_{1/2-h/2}^{1/2+h/2} \exp(2 \pi i m x_1) \, \diff x_1 \right) \left( \frac{1}{h} \int_{1/2-h/2}^{1/2+h/2} \exp(2 \pi i n x_2) \, \diff x_2 \right) \eqqcolon I_m I_n.
		\end{equation}
		Let's calculate $I_m$. If $m=0$, then $I_0 = \frac{1}{h} \int_{1/2-h/2}^{1/2+h/2} 1 \, \diff x = \frac{1}{h} \cdot h = 1$.
		If $m \ne 0$, then
		\begin{align} 
			I_m &= \frac{1}{h} \left[ \frac{\exp(2 \pi i m x)}{2 \pi i m} \right]_{1/2-h/2}^{1/2+h/2} \\ 
			&= \frac{1}{2 \pi i m h} (\exp(2 \pi i m (1/2+h/2)) - \exp(2 \pi i m (1/2-h/2))) \\
			&= \frac{1}{2 \pi i m h} (\exp(\pi i m + \pi i m h) - \exp(\pi i m - \pi i m h)) \\
			&= \frac{\exp(\pi i m)}{2 \pi i m h} (\exp(\pi i m h) - \exp(-\pi i m h)) \\
			&= \frac{(-1)^m}{2 \pi i m h} (2 i \sin(\pi m h)) \quad [\text{Using } \sin\theta = (e^{i\theta}-e^{-i\theta})/2i] \\
			&= (-1)^m \frac{\sin(\pi m h)}{\pi m h} = (-1)^m \sinc(m h) \quad [\text{Using } \sinc(x) = \sin(\pi x)/(\pi x)]
		\end{align}
		So, $\hat{k}_{m, n} = (-1)^{m+n} \sinc(m h) \sinc(n h)$.
		The singular values are $\sigma_{m, n} = |\sinc(m h) \sinc(n h)|$.
		We have $|\sinc(m h)| = \frac{|\sin(\pi m h)|}{|\pi m h|} \le \frac{1}{|\pi m h|}$.
		Thus, the singular values decay at a rate close to $O(\frac{1}{mn})$ as $m, n \to \infty$ and the problem is almost \emph{moderately ill-posed} as the decay rate is polynomial (ignoring sine).
		
		\item  We consider the kernel
		\begin{equation} 
		k\left(x_{1}, x_{2}\right)
		=\varphi\left(x_{1}\right) \varphi\left(x_{2}\right),
		\end{equation}
		with $\varphi(x)$ being a $C_0^\infty$ (infinitely differentiable with compact support) function on $[0,1]$.
		We have $\hat{k}_{m, n} = \hat{\varphi}_m \hat{\varphi}_n$. We need to understand the decay of $\hat{\varphi}_m$.
		A key result in Fourier analysis is that the smoother a function, the faster its Fourier coefficients decay. Since $\varphi(x)$ is infinitely differentiable and zero at the boundaries (0 and 1), we can integrate by parts repeatedly
		\begin{align} 
			\hat{\varphi}_m &= \int_0^1 \varphi(x) \exp(2 \pi i m x) \, \diff x \\
			&= \left[ \varphi(x) \frac{\exp(2 \pi i m x)}{2 \pi i m} \right]_0^1 - \int_0^1 \varphi'(x) \frac{\exp(2 \pi i m x)}{2 \pi i m} \, \diff x \\
			&= 0 - \frac{1}{2 \pi i m} \int_0^1 \varphi'(x) \exp(2 \pi i m x) \, \diff x \quad (\text{Since } \varphi(0)=\varphi(1)=0) \\
			&= \left( - \frac{1}{2 \pi i m} \right) \widehat{\varphi'}_m
		\end{align}
		Since $\varphi'(x)$ is also smooth and zero at the boundaries, we can repeat this $k$ times
		\begin{align} \label{eq: diff rule fourier}
			\hat{\varphi}_m = \left( - \frac{1}{2 \pi i m} \right)^k \widehat{\varphi^{(k)}}_m 
		\end{align}
		Since $\varphi^{(k)}(x)$ exists and is bounded for any $k$, its Fourier coefficient $\hat{\varphi^{(k)}}_m$ is bounded. This means that $|\hat{\varphi}_m|$ decays faster than $1/|m|^k$ for \emph{any} integer $k > 0$, i.e., \emph{super-polynomial decay}.
		Therefore, the singular values $\sigma_{m, n} = |\hat{\varphi}_m \hat{\varphi}_n|$ decay  faster  than $1/(mn)^k$ for any $k$.
		Thus the problem is \emph{severely ill-posed}. 
	\end{enumerate}
\end{Solution}

\section{Regularization Theory for Inverse Problems}
So far, we established the existence of a unique minimum norm solution $x^\dag = A^\dag y$ for any $y\in\q D(A^\dag)$.
In applications, we can often only access \enquote{noisy data} $y^\delta\in B_\delta(y)$ (with \emphi{noise level} $\delta>0$) instead of the \enquote{exact data} $y$.
Unfortunately, if $A^\dag$ is unbounded, we cannot control the mismatch between  $A^\dag y^\delta$ and $A^\dag y$.
To resolve this issue, we construct a family of approximations $\{x^\delta_\alpha\}_{\alpha>0}$ for $A^\dag y^\delta$, where each $x^\delta_\alpha$ depends Lipschitz continuously on $\delta$.
In a second step, we aim to guarantee $x^\delta_{\alpha(\delta)}\to x^\dag$ as $\delta \to 0$ with an appropriate \emph{parameter choice}.
In particular, we want to approximate the pseudoinverse for vanishing noise.
Such approximations are called \emph{regularization methods}.
To develop our results, we require another fundamental result from functional analysis.
\begin{thm}{Banach--Steinhaus}{banach-steinhaus}
	Let $\bb X$, $\bb Y$ be Banach spaces, and $\{A_\lambda\}_{\lambda \in \Lambda}\subset \q L(\bb X,\bb Y)$ be pointwise bounded, i.e., for all $x\in \bb X$ there exists an $M_x>0$ with $\sup_{\lambda \in \Lambda}\|A_\lambda x\|_{\bb Y} \leq M_x$.
	Then it holds that
	\begin{equation}
		\sup_{\lambda \in \Lambda} \| A_\lambda \|_{ \q L(\bb X,\bb Y)} < \infty.
	\end{equation}
\end{thm}
\begin{cor}\label{cor:banach-steinhaus}
	Let $\bb X,\bb Y$ be Banach spaces and $\{A_n\}_{n\in\N}\subset  \q L(\bb X,\bb Y)$. Then the following are equivalent:
	\begin{enumerate}
		\item $\{A_n\}_{n\in\N}$ converges uniformly on compact subsets of $\bb X$;
		\item $\{A_n\}_{n\in\N}$ converges pointwise, i.e., $A_nx\to Ax$ (strongly in $\bb Y$) for all $x\in \bb X$;
		\item $\{A_n\}_{n\in\N}$ converges pointwise on a dense subset $U\subset \bb X$ and
		\begin{equation}
			\sup_{n\in \N} \|A_n\|_{ \q L(\bb X,\bb Y)} < \infty.
		\end{equation}
	\end{enumerate}
	Moreover, if $A_n$ converges pointwise to $A \colon \bb X\to \bb Y$, then $A$ is bounded.
\end{cor}

\subsection{Regularization and Parameter Choice}\label{sec:Regularization}
For Hilbert spaces $\bb X$, $\bb Y$ and $A\in\q L(\bb X,\bb Y)$, we introduce \emph{regularization operators} as continuous approximations of the (potentially unbounded) pseudoinverse $A^\dag$.
\begin{defn}\label{def:regularisierung}
	Let $A\in\q L(\bb X,\bb Y)$.
	A family $\{R_\alpha\}_{\alpha>0}$ is a \emphi{regularization} (of $A^\dag$) if
	\begin{enumerate}
		\item $R_\alpha\in\q L(\bb Y,\bb X)$ for all $\alpha>0$;
		\item $R_\alpha y\to A^\dag y$ as $\alpha\to 0$ and all $y\in\q D(A^\dag)$.
	\end{enumerate}
\end{defn}
Outside of $\q D(A^\dag)$, these approximations can diverge.
\begin{thm}{}{regularisierung:glm}
	Let $A\in\q L(\bb X,\bb Y)$,  $A^\dag$ be unbounded, and $\{R_\alpha\}_{\alpha>0}\subset\q L(\bb Y,\bb X)$ be a regularization.
	Then, there exists $y\in \bb Y$ and  $\{\alpha_n\}_{n\in\N}$ with $\alpha_n \to 0$ such that $\norm{R_{\alpha_n} y}_{\bb X}\to \infty$.
\end{thm}
\begin{pr}
	If no such $y\in \bb Y$ exists, then $\{R_\alpha\}_{\alpha>0}\subset\q L(\bb Y,\bb X)$ is pointwise bounded and thus $\sup_\alpha \|R_\alpha\|_{\q L(\bb Y,\bb X)} < \infty$ by \Cref{th:banach-steinhaus}.
	Since $R_\alpha(y)\to A^\dag(y)$ for $ y \in \q D(A^\dag)$ and $\q D(A^\dag)$ is dense in $\bb Y$, \Cref{cor:banach-steinhaus} then yields the contradiction that $A^\dag$ is bounded.
\end{pr}
Under one additional assumption, $R_{\alpha}y$ even diverges for \emph{all} $y\notin\q D(A^\dag)$.
\begin{thm}{}{regularisierung:div}
	Let $A\in\q L(\bb X,\bb Y)$,  $A^\dag$ be unbounded, and $\{R_\alpha\}_{\alpha>0}\subset\q L(\bb Y,\bb X)$ be a regularization.
	If $\{R_\alpha\}_{\alpha>0}$ is bounded in the data space norm
	\begin{equation}\label{eq:regularisierung:bounded}
		\sup_{\alpha>0} \norm{AR_\alpha}_{\q L(\bb Y,\bb Y)}<\infty,
	\end{equation}
	then $\norm{R_\alpha y}_{\bb X}\to \infty$ as $\alpha\to0$ for every $y\notin\q D(A^\dag)$.
\end{thm}
\begin{pr}
	Let $y\notin \q D(A^\dag)$ and assume that there is $\{\alpha_n\}_{n\in\N}$ with $\alpha_n \to 0$ such that $\{R_{\alpha_n}y\}_{n\in\N}$ is bounded.
	Then, we can find a subsequence $\{\alpha_{n_k}\}_{k\in\N}$ such that $x_k\coloneqq R_{\alpha_{n_k}}y$ satisfies $x_k\wkto x\in \bb X$.
	For all $z \in \bb Y$, we have  
	\begin{equation}
			\inner{z}{Ax_k}_\bb Y 
			= \inner{A^*z}{x_k}_{\bb X} 
			\to \inner{A^*z}{x}_{\bb X} 
			= \inner{z}{A x}_\bb Y, 
			\quad \text{as} \quad k \to \infty, 
	\end{equation} 
	which yields $Ax_k\wkto Ax$ and is also known as weak continuity of bounded operators.
	
	At the same time, $A \in \q L(\bb X, \bb Y)$ and $R_\alpha\to A^\dag$ pointwise on $\q D(A^\dag)$ imply that $AR_\alpha \to AA^\dag $ pointwise on $\q D(A^\dag)$.
	 Thus, \Cref{lem:inverse}\,\circled{4} , assumption \eqref{eq:regularisierung:bounded} and \Cref{cor:banach-steinhaus} imply that $AR_{\alpha_n} \to P_{\overline{\q R}}$ pointwise.
	Hence, $Ax_k \to  P_{\overline{\q R}}y$, which combined with $Ax_k\wkto Ax$ implies that $Ax =  P_{\overline{\q R}}y$.
	This  contradicts our assumption $y \notin D(A^\dag)$.
\end{pr}

In general, we cannot assume $B_\delta(y) \subset \q D(A^\dag)$.
Thus, we have to consider for $y^\delta \in B_\delta(y)$ the \emphi{regularization error}
\begin{equation}
	\label{eq:regularisierung:splitting}
		\norm{R_{\alpha}y^\delta - A^\dag y}_{\bb X} \leq  \norm{R_{\alpha}y^\delta - R_{\alpha}y}_{\bb X} +  \norm{R_{\alpha}y - A^\dag y}_{\bb X} \leq \delta \norm{R_{\alpha}}_{\q L(\bb Y,\bb X)} + \norm{R_{\alpha}y - A^\dag y}_{\bb X}.
\end{equation}
The first summand is called \emph{data error}.
By \Cref{th:regularisierung:glm}, this term might diverge for $\alpha\to 0$ as long as $\delta>0$ does not vanish.
The second term is called \emph{approximation error}.
Due to \Cref{def:regularisierung} \circled{2} , this term converges to zero.
In summary, we should choose $\alpha$ depending on $\delta$ and $y^\delta$ such that \eqref{eq:regularisierung:splitting} vanishes as $\delta\to 0$.
\begin{defn}\label{def:parameterwahl}
	 We differentiate 3 types of \emphi{parameter choice rules} $\alpha\colon \R^+ \times \bb Y \to\R^+$ with $(\delta,y^\delta)\mapsto \alpha(\delta,y^\delta)$, namely
	\begin{enumerate}
		\item \emph{a priori choice rules} that only depend on $\delta$;
		\item \emph{a posteriori choice rules} that depend on $\delta$ and $y^\delta$;
		\item \emph{heuristic choice rules} that only depend on $y^\delta$.
	\end{enumerate}
	
	If $\{R_\alpha\}_{\alpha>0}$ is a regularization of $A^\dag$ with parameter choice rule $\alpha$, the pair $(R_\alpha,\alpha)$ is called a (convergent) \emphi{regularization method} if
	\begin{equation}\label{eq:regularisierung:paramwahl}
		\lim_{\delta\to 0}\sup_{y^\delta\in B_\delta(y)}\norm{R_{\alpha(\delta,y^\delta)}y^\delta - A^\dag y}_{\bb X} = 0 \qquad\text{for all }y\in \q D(A^\dag).
	\end{equation}
\end{defn}

Below, we characterize a priori rules that lead to regularization methods.
\begin{thm}{}{regularisierung:apriori_char}
	Let $A^\dag$ be unbounded, $\{R_\alpha\}_{\alpha>0}$ be a regularization, and $\alpha\colon \R^+\to\R^+$ be an a priori parameter choice rule.
	Then $(R_\alpha,\alpha)$ is a regularization method if and only if
	\begin{enumerate}
		\item $\displaystyle\lim_{\delta\to 0} \alpha(\delta) = 0$,
		\item $\displaystyle \lim_{\delta\to 0} \delta\norm{R_{\alpha(\delta)}}_{\q L(\bb Y,\bb X)} = 0$.
	\end{enumerate}
\end{thm}
\begin{pr}
	$\Leftarrow$: Using the decomposition \eqref{eq:regularisierung:splitting}, we get
	\begin{equation}
		\norm{R_{\alpha(\delta)}y^\delta - A^\dag y}_{\bb X} \leq \delta \norm{R_{\alpha(\delta)}}_{\q L(\bb Y,\bb X)} + \norm{R_{\alpha(\delta)}y - A^\dag y}_{\bb X}.
	\end{equation}
	As $\delta\to 0$, the first term vanishes by \circled{2}, while the second one vanishes due \circled{1} and the pointwise convergence of $\{R_\alpha\}_{\alpha>0}$.
	
	$\Rightarrow$:
	If \circled{1} is violated, there exists $\{\delta_n\}_{n \in \N}$ with $\alpha(\delta_n) \to c > 0$ as $n \to \infty$.
	For $y^{\delta_n} = y$, we then have
	\begin{equation}
	\| R_{\alpha(\delta_n)} y - A^\dagger y \|_\bb Y \to \| R_{c} y - A^\dagger y \|_\bb Y, \quad \text{as} \quad n \to \infty.
	\end{equation}
	Since $A^\dagger$ is unbounded whereas $R_c$ is bounded, there exists 
	$y \in \q D(A^\dagger)$ such that $R_c y \neq A^\dagger y$. 
	This contradicts the assumption that $(R_\alpha,\alpha)$ is a regularization method.
	
	If \circled{2} is violated, we have  $\limsup_{\delta \to 0} \delta \norm{R_{\alpha(\delta)}}_{\q L(\bb Y,\bb X)}>0$.
	Thus, we can find $\{\delta_n\}_{n\in\N} $ with $\delta_n \to 0$ and $\{z_n\}_{n\in\N}\subset \bb Y$ with $\norm{z_n}_{\bb Y}=1$ such that $\delta_n\norm{R_{\alpha(\delta_n)}z_n}_{\bb X}\geq  \varepsilon$ for some $ \varepsilon>0$.
	For any $y\in\q D(A^\dag)$, it holds with $y^{\delta_n}\coloneqq y+\delta_n z_n \in B_{\delta_n}(y)$ that
	\begin{equation}
	\limsup_{n \to \infty} \Vert R_{\alpha(\delta_n)} y^{\delta_n} - A^\dag y \Vert_{\bb X} \geq \limsup_{n \to \infty} \delta_n \Vert R_{\alpha(\delta_n)} z_n \Vert_{\bb X} -  \Vert R_{\alpha(\delta_n)} y - A^\dag y \Vert_{\bb X} \geq \varepsilon.
	\end{equation}
	Hence, the property \eqref{eq:regularisierung:paramwahl} of regularization methods is violated.
\end{pr}
Assumption \circled{2} states that $\alpha$ must not vanish too fast compared to $\delta$.
Thus, a priori rules usually take the form $\alpha(\delta) = \delta^r$ with $r\in(0,1)$ depending on the regularization $\{R_\alpha\}_{\alpha>0}$ and the inverse problem.

\begin{remark}
	An optimal parameter choice makes use of both $\delta$ and $y^\delta$ (an a posteriori rule).
	For this, one typically aims to choose $\alpha=\alpha(\delta)$ as large as possible while still satisfying $\|A R_\alpha y^\delta - y^\delta\| \leq \tau \delta$,
	for some constant $\tau>1$ that is independent of $\delta$ and $y^\delta$.
	However, this only works if $\overline{\q R(A)}=\mathbb{Y}$.
	This approach is known in the literature as the \emph{Morozov discrepancy principle} and it typically yields convergent regularization methods.
	
	Moreover, no convergent heuristic parameter choice rule can exist if $A^\dag$ is unbounded.
	The latter is the case, e.g., for any compact operator with $\dim(\q R(K)) = \infty$.
	This result is commonly known as \emph{Bakushinskii veto} in the literature.
	Still, heuristic rules often work in practice, for example, if the noisy data satisfies $y^\delta \notin \q D(A^\dag)$.
\end{remark}
\subsection{Spectral Regularization}\label{chap:spektral}
Next, we want to derive the so-called spectral regularization methods.
These involve functions of compact selfadjont operators.
Let $\phi\colon[0,\infty)\to\R$ be piecewise continuous and locally bounded. 
Based on the singular system $\{(\sigma_n,u_n,v_n)\}_{n\in\N}$ of $K\in \q K(\bb X,\bb Y)$, we define $\phi(K^*K)\colon \bb X\to \bb X$ by applying $\phi$ to the eigenvalues, namely
\begin{equation}
	\label{eq:functional}
	\phi(K^*K)x = \sum_{n\in\N} \phi(\sigma_n^2) \inner{x}{v_n}_{\bb X} v_n+ \phi(0)P_{\q N(K)} x \qquad \text{for all }x\in \bb X.
\end{equation}
Due to the Bessel inequality \eqref{eq:funktan:bessel} and $\sigma_1 = \norm{K}_{\q L(\bb X,\bb Y)}$, we get that
\begin{equation}\label{eq:functional_continuous}
	\norm{\phi(K^*K) x}_{\bb X} \leq \Bigl(\sup_{n\in\N} |\phi(\sigma_n^2)| + \vert \phi(0) \vert \Bigr) \Vert x \Vert_{\bb X}  \leq 2 \Vert x \Vert_{\bb X}  \sup_{\lambda\in [0,\norm{K}_{\q L(\bb X,\bb Y)}^2]} |\phi(\lambda)|<\infty.
\end{equation}
Further, one can show that $\phi(K^*K) x$  is independent of the chosen singular system.
Hence, $\phi(K^*K) x$ is well defined.
Moreover, we have that $\phi(K^*K)\in \q L(\bb X,\bb X)$.

Now, we can construct a regularization $\{R_\alpha\}_{\alpha>0}$ according to \Cref{def:regularisierung} based on the SVD and \eqref{eq:inverse:normal_pseudo}.
The latter states for $y\in\q D(K^\dag)$ that
\begin{equation}
	K^\dag y = (K^*K)^\dag K^* y.
\end{equation}
We observe that $\{(\sigma_n^2,v_n)\}_{n\in\N}$ is an eigen system of $K^*K$ and that $K^* y \in \overline{\q R(K^*)} = \q N(K)^\bot$, which implies $P_{\q N(K)} x K^* y  = 0$.
Thus,  \Cref{th:inverse:picard} and \eqref{eq:functional} yield for $\phi(\lambda) = \lambda^{-1}$ that
\begin{equation}\label{eq:FormulaReg}
		(K^*K)^\dag K^* y 
		= \sum_{n\in\N} \sigma_n^{-2} \inner{K^*y}{ v_n}_{\bb X} v_n 
		= \varphi(K^*K) K^* y .
\end{equation}
This becomes unstable for ill-posed problem where $\sigma_n \to 0$.
To get a bounded approximation of $K^\dag$, we replace the unbounded function $\phi$ by a family of bounded approximations $\{\phi_\alpha\}_{\alpha>0}$.
Below, we denote $\kappa\coloneqq  \norm{K}_{\q L(\bb X, \bb Y)}^2 = \norm{K^*K}_{\q L(\bb X, \bb X)}$ for brevity.
\begin{defn}\label{def:RegFilter}
	Let $\{\phi_\alpha\}_{\alpha>0}$ be a family of piecewise continuous and bounded functions $\phi_\alpha\colon [0,\kappa] \to \R$. If
	\begin{enumerate}
		\item $\displaystyle \lim_{\alpha\to 0}\phi_\alpha(\lambda) = \frac1\lambda$  for all $\lambda\in (0,\kappa]$ and
		\item $\displaystyle \lambda|\phi_\alpha(\lambda)|\leq C_\phi$ for some $C_\phi>0$ and all $\lambda\in (0,\kappa]$, $\alpha>0$,
	\end{enumerate}
	then $\{\phi_\alpha\}_{\alpha>0}$ is a \emph{(regularizing)} \emphi{filter} for $K\in\q K(\bb X, \bb Y)$.
\end{defn}
If \circled{1} and \circled{2} hold independent of $\kappa>0$, then $\{\phi_\alpha\}_{\alpha>0}$ is a filter for \emph{all} $K\in\q K(\bb X, \bb Y)$.
As discussed, the idea is to now define the regularization $R_\alpha \coloneqq  \phi_\alpha(K^*K)K^*$, which according to \eqref{eq:FormulaReg} satisfies
\begin{equation}
	\label{eq:spektral:regularisierung}
		R_\alpha = \sum_{n\in\N} \phi_\alpha(\sigma_n^2)\sigma_n \inner{y}{u_n}_{\bb Y} v_n.
\end{equation}
Before we establish its regularizing properties, we provide several examples.
\begin{example}[Filters]\label{ex:spektral}
	\begin{enumerate}
		\item The \emphi{truncated SVD} is induced by the bounded and piecewise continuous function
		\begin{equation}
			\label{eq:spektral:cutoff}
			\phi_\alpha(\lambda) = \begin{cases}\frac1\lambda & \text{if }\lambda\geq \alpha,\\ 0 & \text{else.}
			\end{cases}
		\end{equation}
		Indeed, $\phi_\alpha$ satisfies $\lim_{\alpha \to 0} \phi_\alpha(\lambda) = \frac1\lambda$  and $\lambda|\phi_\alpha(\lambda)|\leq 1$ for  all $\lambda, \alpha>0$.
		It holds
		\begin{equation} \label{eq:spektral:cutoff_R}
			R_\alpha y = \sum_{n\in\N}\phi_\alpha(\sigma_n^2)\sigma_n \inner{y}{u_n}_{\bb Y} v_n = \sum_{\sigma_n\geq \sqrt\alpha}\frac1{\sigma_n} \inner{y}{u_n}_{\bb Y} v_n,
		\end{equation}
		which also explains the name.
		
		\item The \emphi{Tikhonov regularization} corresponds to the bounded and continuous function
		\begin{equation}
			\label{eq:spektral:tikhonov}
			\phi_\alpha(\lambda) = \frac{1}{\lambda+\alpha}.
		\end{equation}
		Again, $\lim_{\alpha \to 0} \phi_\alpha(\lambda) = \frac1\lambda$  and $\lambda|\phi_\alpha(\lambda)|\leq 1$ for  all $\lambda,\alpha>0$.
		The regularization reads
		\begin{equation}
			R_\alpha y = \sum_{n\in\N}\frac{\sigma_n}{\sigma_n^2+\alpha} \inner{y}{u_n}_{\bb Y} v_n.
		\end{equation}
		As we will show in \Cref{chap:tikhonov}, $R_\alpha y $ can be computed without the singular system.
		
		\item The \emphi{Landweber regularization} is induced by the bounded and continuous function
		\begin{equation}
			\label{eq:spektral:landweber}
			\phi_\alpha(\lambda) = \frac{1-(1-\omega\lambda)^{1/\alpha}}{\lambda}, \qquad \omega>0.
		\end{equation}
		One can show that $\phi_\alpha$ defines a filter if $\omega$ is small enough.
		Moreover, a characterization without the singular system is again possible.
	\end{enumerate}
\end{example}

To show that $R_\alpha\coloneqq \phi_\alpha(K^*K)K^*$ defines a regularization if $\{\phi_\alpha\}_{\alpha>0}$ is a filter, we require three fundamental lemmas.

\begin{lemma}\label{lem:spektral:beschraenkt_K}
	Any filter $\{\phi_\alpha\}_{\alpha>0}$ for $K\in\q K(\bb X, \bb Y)$ satisfies
	\begin{equation}
		\norm{KR_\alpha}_{\q L(Y,Y)} \leq \sup_{n\in\N} |\phi_\alpha(\sigma_n^2)|\sigma_n^2 \leq C_\phi \qquad\text{for all }\alpha>0.
	\end{equation}
\end{lemma}
\begin{pr}
	For $y\in \bb Y$ and $\alpha>0$, we push $K$ through \eqref{eq:spektral:regularisierung} (as detailed for \eqref{eq:svd:operator}) to obtain
	\begin{equation}\label{eq:spektral:beschraenkt_K}
			KR_\alpha y = \sum_{n\in\N} \phi_\alpha(\sigma_n^2)\sigma_n \inner{y}{u_n}_{\bb Y} Kv_n =  \sum_{n\in\N} \phi_\alpha(\sigma_n^2)\sigma_n^2 \inner{y}{u_n}_{\bb Y} u_n.
	\end{equation}
	Combining the Bessel inequality \eqref{eq:funktan:bessel} and equation \eqref{eq:spektral:beschraenkt_K} yields
	\begin{equation}
			\norm{KR_\alpha y }_{\bb Y}^2 \leq
			\sup_{n\in\N} |\phi_\alpha(\sigma_n^2)\sigma_n^2|^2 \sum_{n\in\N} |\inner{y}{u_n}_{\bb Y}|^2\leq  \sup_{n\in\N} |\phi_\alpha(\sigma_n^2)\sigma_n^2|^2 \norm{y}_{\bb Y}^2.
	\end{equation}
	Since $0<\sigma_n^2 \leq \kappa$, the claim now follows from the boundedness condition \circled{2} of filters.
\end{pr}
\begin{lemma}\label{lem:spektral:beschraenkt}
	Let $\{\phi_\alpha\}_{\alpha>0}$ be a filter for $K\in\q K(\bb X, \bb Y)$.
	Then $R_\alpha\in \q L(\bb Y, \bb X)$ for all $\alpha>0$ and
	\begin{equation}
		\norm{R_\alpha}_{\q L(\bb Y,\bb X)} \leq \sqrt{C_\phi} \sup_{\lambda\in(0,\kappa]} \sqrt{|\phi_\alpha(\lambda)|}.
	\end{equation}
\end{lemma}
\begin{pr}
	By \eqref{eq:spektral:regularisierung}, we  
	get for all $y\in \bb Y$ and $\alpha>0$ that
	\begin{equation}
	\norm{R_\alpha y}_{\bb X}^2 = \sum_{n\in\N} \bigl(\phi_\alpha(\sigma_n^2) \sigma_n \inner{y}{u_n}_{\bb Y} \bigr)^2.
	\end{equation}
	Using the Bessel inequality \eqref{eq:funktan:bessel} and the boundedness condition \circled{2} of filters, we obtain
		\begin{equation}
				\norm{R_\alpha y}_{\bb X}^2 \leq \sup_{n\in\N} |\phi_\alpha(\sigma_n^2) \sigma_n|^2 \ \norm{y}_{\bb Y}^2 \leq \sup_{n\in\N} |\phi_\alpha(\sigma_n^2)| \ C_\phi \norm{y}_{\bb Y}^2.
	\end{equation}
	Taking the square root yields the claim.
\end{pr}
\begin{lemma}[Fundamental Lemma of Spectral Regularization]\label{lem:spektral:verfahrensfehler}
	Let $\{\phi_\alpha\}_{\alpha>0}$ be a filter for $K\in\q K(\bb X, \bb Y)$, and $x^\dag = K^\dag y$ with $y\in \q D(K^\dag)$.
	Then
	\begin{equation}\label{eq:spektral:fehler}
		K^\dag y - R_\alpha y = \sum_{n\in\N} r_\alpha(\sigma_n^2) \inner{x^\dag}{v_n}_{\bb X} v_n,
	\end{equation}
	where $r_\alpha(\lambda) \coloneqq   1-\lambda\phi_\alpha(\lambda)$ satisfies for all  $\lambda\in(0,\kappa]$ and $\alpha>0$ that
	\begin{align}
	\lim_{\alpha\to 0}r_\alpha(\lambda) = 0 \qquad \text{and} \qquad |r_\alpha(\lambda)|\leq 1+C_\phi.
	\end{align}
\end{lemma}
\begin{pr}
	Since $K^*Kx^\dag = K^*y$ by \eqref{eq:inverse:normal}, we can rewrite \eqref{eq:FormulaReg} as
	\begin{equation}
		R_\alpha y 
		= \sum_{n \in \N} \phi_\alpha(\sigma_n^2) \inner{K^*K x^\dag}{v_n} v_n 
		= \sum_{n \in \N} \phi_\alpha(\sigma_n^2) \sigma_n^2 \inner{x^\dag}{v_n} v_n.
	\end{equation}
	Furthermore, we can decompose $x^\dag = K^\dag y \in \overline{\q R(K^*)}$ as $x^\dag = \sum_{n \in \N} \inner{x^\dag}{v_n}_{\bb X} v_n$.
	Thus, the definition of $r_\alpha$ immediately yields that
	\begin{equation}
		K^\dag y - R_\alpha y 
		= \sum_{n \in \N} (1 - \phi_\alpha(\sigma_n^2) \sigma_n^2) \inner{x^\dag}{v_n} v_n.
		= \sum_{n\in\N} r_\alpha(\sigma_n^2) \inner{x^\dag}{v_n}_{\bb X} v_n.
	\end{equation}
	The rest follows from the properties \circled{1} and \circled{2} of  filters.
\end{pr}

Now, we establish the regularization property of $\{R_\alpha\}_{\alpha>0}$.
\begin{thm}{}{spektral:konvergenz}
	If $\{\phi_\alpha\}_{\alpha>0}$ is a filter for $K$,
	then $\{R_\alpha\}_{\alpha>0} \subset  \q L(\bb Y, \bb X)$ defines a regularization, i.e.,
	\begin{equation}
		\lim_{\alpha\to 0} R_\alpha y = K^\dag y\qquad\text{for all } y\in\q D(K^\dag).
	\end{equation}
	If $K^\dag$ is unbounded, then $\lim_{\alpha\to 0} \norm{R_\alpha y}_{\bb X} = \infty$ for every $y\notin\q D(K^\dag)$.
\end{thm}
\begin{pr}
	For $x^\dag= K^\dag y$ with $y\in \q D(K^\dag)$, \Cref{lem:spektral:verfahrensfehler} implies that
	\begin{equation}\label{eq:EstPointwise}
		\norm{K^\dag y-R_\alpha y}_{\bb X}^2 =\sum_{n\in\N} |r_\alpha(\sigma_n^2)|^2 |\inner{x^\dag}{v_n}_{\bb X}|^2.
	\end{equation}
	To show that \eqref{eq:EstPointwise} converges to zero, let $\eps>0$ be arbitrary.
	With the Bessel inequality \eqref{eq:funktan:bessel}, we obtain an $N\in \N$ with
	\begin{equation}
		\sum_{n=N+1}^\infty | \inner{x^\dag}{v_n}_{\bb X} |^2 < \frac{\eps^2}{2(1+C_\phi)^2}.
	\end{equation}
Due to $\lim_{\alpha\to 0} \sup_{\lambda \in \{\sigma_1^2,\dots,\sigma_N^2\}}r_\alpha(\lambda) = 0$, there exists $\alpha_0>0$ with
	\begin{equation}
		|r_\alpha(\sigma_n^2)|^2 < \frac{\eps^2}{2\norm{x^\dag}_{\bb X}^2} \qquad\text{for all } n\leq N \text{ and }\alpha<\alpha_0.
	\end{equation}
	Thus, \Cref{lem:spektral:verfahrensfehler} implies for all $\alpha<\alpha_0$ that
	\begin{align}
			\norm{K^\dag y-R_\alpha y}_{\bb X}^2 
			&= \sum_{n=1}^N |r_\alpha(\sigma_n^2)|^2 |\inner{x^\dag }{v_n}_{\bb X}|^2 
			+ \sum_{n=N+1}^\infty |r_\alpha(\sigma_n^2)|^2 |\inner{x^\dag }{v_n}_{\bb X}|^2\notag\\
			&\leq  \frac{\eps^2}{2\norm{x^\dag}_{\bb X}^2}  \sum_{n=1}^N  |\inner{x^\dag}{v_n}_{\bb X}|^2 +  (1+C_\phi)^2 \frac{\eps^2}{2(1+C_\phi)^2} \leq \eps^2.
	\end{align}
	Thus, $\norm{K^\dag y-R_\alpha y}_{\bb X} \to 0$ as $\alpha\to 0$.
	The divergence for $y\notin\q D(K^\dag)$ follows from \Cref{th:regularisierung:div} and \Cref{lem:spektral:beschraenkt_K}.
\end{pr}

In particular, all the filters from \Cref{ex:spektral} define regularizations.
To turn these into regularization methods $(R_\alpha,\alpha)$, we can use an a priori parameter choice $\alpha \colon \R^+ \to \R^+$.
By \Cref{th:regularisierung:apriori_char}, the latter has to satisfy $\alpha(\delta)\to 0$ and $\delta\norm{R_\alpha}_{\q L(\bb Y,\bb X)}\to 0$ as $\delta\to 0$.
Using \Cref{lem:spektral:beschraenkt}, we can thus establish a condition on $\alpha$.
Below, we provide two examples.
\begin{example}[Truncated SVD]\label{ex:RegTruncSVD}
	For the filter $\phi_\alpha$ in \eqref{eq:spektral:cutoff}, we get
	\begin{equation}
		\norm{R_\alpha}_{\q L(\bb Y,\bb X)} \leq \sqrt{C_\phi} \sup_{n\in\N} \sqrt{|\phi_\alpha(\sigma_n^2)|} = \frac1{\sqrt\alpha}.
	\end{equation}
	Thus, $\alpha(\delta)\to 0$ and $\frac {\delta^2}{\alpha(\delta)} \to 0$ is a sufficient condition.
	For the rule $\alpha(\delta) \coloneqq   \delta$, we get that $R_{\alpha(\delta)} y = \sum_{\sigma_n\geq \sqrt\delta}\frac1{\sigma_n} \langle y, u_n \rangle_{\bb Y} v_n$.
\end{example}

\begin{example}[Tikhonov]\label{ex:RegTikhonov}
	For the filter in \eqref{eq:spektral:tikhonov}, we again get
	\begin{equation}
		\norm{R_\alpha}_{\q L(\bb Y,\bb X)} \leq \frac1{\sqrt{\alpha}}.
	\end{equation}
	Thus, as for \Cref{ex:RegTruncSVD}, we get that the a priori rule $\alpha(\delta) = \delta$ leads to a regularization method with $R_{\alpha(\delta)} y = \sum_{n\in\N}\frac{\sigma_n}{\sigma_n^2+\delta} \inner{y}{u_n}_{\bb Y} v_n$.
\end{example}
\subsection{Convergence Rates}

So far, we have only answered the question if the regularization error \eqref{eq:regularisierung:splitting} vanishes.
In particular, we have established that for a regularization method $(R_\alpha,\alpha)$ and $y\in \q D(A^\dag)$, the \emphi{worst-case error}
\begin{equation}
	\label{eq:rates:worstcase}
	\q E(y,\delta)\coloneqq \sup_{y^\delta\in B_\delta(y)} \norm{R_{\alpha(\delta,y^\delta)} y^\delta - A^\dag y}_{\bb X}.
\end{equation}
vanishes as $\delta \to 0$ by \eqref{eq:regularisierung:paramwahl}.
From a practical perspective, a central goal is to obtain estimates of the form
\begin{equation}\label{eq:errorEst}
	\norm{R_{\alpha(\delta,y^\delta)} y^\delta - A^\dag y}_{\bb X} \leq \psi(\delta)
\end{equation}
for some $\psi\colon \R^+\to\R^+$ (ideally with $\lim_{\delta\to 0}\psi(\delta) =0$).
Then, we can quantitatively control the error depending on the noise level $\delta$.

Unfortunately, as we show below, any such function $\psi$ has to depend on the data $y$.
\begin{thm}{}{rates:counter}
	Let $(R_\alpha,\alpha)$ be a regularization method.
	If there exists an increasing $\psi\colon \R^+\to\R^+$ with $\lim_{\delta\to 0}\psi(\delta) =0$ and
	$\sup_{y\in \q D(A^\dag)\cap B_{\bb Y}} \q E(y,\delta) \leq \psi(\delta)$, then $A^\dag$ is bounded.
\end{thm}
\begin{pr}
	Let $\tilde y, y \in \q D(A^\dag) \cap B_{\bb Y}$.
	Then, we have that with $\delta = \Vert \tilde y - y \Vert_{\bb Y}$ that 
	\begin{equation}
			\norm{A^\dag \tilde y - A^\dag y}_{\bb X} \leq \norm{A^\dag \tilde y - R_{\alpha(\delta, \tilde y)} \tilde y}_{\bb X} + \norm{R_{\alpha(\delta,\tilde y)} \tilde y - A^\dag y}_{\bb X} \leq \q E(\tilde y,\delta) + \q E(y,\delta) \leq 2\psi(\delta).\!
	\end{equation}
	Hence, $A^\dag \colon \q D(A^\dag) \cap B_{\bb Y} \to \bb X$ is uniformly continuous.
	Thus, we can extend it to a bounded operator $A^\dag \colon B_{\bb Y} \to \bb X$ and the claim follows by linearity.
\end{pr}
Consequently, we require additional assumptions on the exact data $y$ (or equivalently $x^\dag = A^\dag y$) if we want to establish error estimates.
From the proof of \Cref{th:rates:counter}, we conclude that these are closely tied to the uniform continuity of $A^\dag$ on closed sets.
Inspired by this observation, given $\q M\subset \q R(A^\dagger)$ and $\delta>0$, we restrict $A^\dag$ to $A(\q M) \cap B_{\delta}(0)$ and compute the radius of its range.
This leads to the \emphi{modulus of conditional continuity} 
\begin{equation}\label{eq:ModulusContinuity}
	\eps(\q M,\delta) \coloneqq  \sup\bigl\{\norm{x}_{\bb X} \, : \, x\in\q M,\ \norm{Ax}_{\bb Y}\leq \delta \bigr\},
\end{equation}
which we can use to lower-bound the worst-case error \eqref{eq:UpperBoundError}.
\begin{thm}{}{lower:bound}
	Let $(R_\alpha,\alpha)$ be a regularization method.
	Then for all $\delta>0$ and $\q M\subset  \q R(A^\dagger)$ it holds 
	\begin{equation}
		\sup_{y\in\q D(A^\dag): A^\dag y\in \q M} \q E(y,\delta) \geq \eps(\q M,\delta).
	\end{equation}
\end{thm}
\begin{pr}
	Let $x\in\q M$ with $\norm{Ax}_{\bb Y} \leq \delta$.
	From $x\in \q R(A^\dagger)$ and \Cref{lem:inverse}\,\circled{3}, we get that
	\begin{equation}
		\norm{x}_{\bb X} =  \norm{A^\dag Ax}_{\bb X} = \norm{A^\dag Ax - R_{\alpha(\delta,0)} 0}_{\bb X} \leq \q E(Ax,\delta).
	\end{equation}
	By inserting this into \eqref{eq:ModulusContinuity}, we further have
	\begin{equation}
		\eps(\q M,\delta) \leq  \sup_{x\in\q M:  \norm{Ax}_{\bb Y}\leq \delta} \q E(Ax,\delta) \leq  \sup_{x\in\q M} \q E(Ax,\delta) \leq \sup_{y: A^\dag y\in\q M} \q E(y,\delta).
	\end{equation}
\end{pr}
In the following, we aim to find sets $\mathcal M \subset \bb X$ (i.e., ground truth data) for which we can control the error decay.
To this end, we  introduce for $s>0$ the operator $|K|^s \coloneqq \phi(K^*K)$ with $\phi(t) =t^{s/2}$, which according to \eqref{eq:functional} can be written using a singular system $\{(\sigma_n,u_n,v_n)\}_{n\in\N}$ as
\begin{equation}\label{eq:functional_sqrt}
	|K|^sx = \sum_{n\in\N} \sigma^s_n \inner{x}{v_n}_{\bb X} v_n \qquad\text{for all }x\in \bb X.
\end{equation}
For selfadjoint and positive definite $K$, it holds that $K = \vert K \vert$.
As we show in \Cref{ex:functional_range}, the operator $|K|^s$ has the following properties.
\begin{lemma}\label{lem:functional_range}
	Let $K\in \q K(\bb X,\bb Y)$. Then
	\begin{enumerate}
		\item $|K|^{r+s} = |K|^r \circ |K|^s$ for all $r,s\geq 0$;
		\item $|K|^r$ is selfadjoint for all $r\geq 0$;
		\item $\norm{|K|x}_{\bb X} = \norm{Kx}_{\bb Y}$ for all $x\in \bb X$;
		\item $\q R(|K|) = \q R(K^*)$.
	\end{enumerate}
\end{lemma}

Based on \eqref{eq:functional_sqrt}, we now consider the sets 
\begin{equation}
	X_{\nu,\rho} = \bigl\{|K|^\nu w \in \bb X : \norm{w}_{\bb X}\leq \rho\bigr\} \subset\q R(|K|^\nu) \subset \q R(K^*).
\end{equation}
By definition of $|K|^\nu w$, it holds that $X_{\nu,\rho}\subset \q R(A^\dagger)$.
\begin{thm}{}{rates:limit}
	Let $K\in\q K(\bb X,\bb Y)$ and $\nu,\rho>0$.
	Then, for all $\delta>0$, it holds that
	\begin{equation}\label{eq:UpperBoundError}
		\eps(X_{\nu,\rho},\delta)\leq \delta^{\frac\nu{\nu+1}} \rho^{\frac1{\nu+1}}.
	\end{equation}
\end{thm}
\begin{pr}
	For $x\in X_{\nu,\rho}$, there exists $w\in \bb X$ with $x=|K|^\nu w$ and $\norm{w}_{\bb X}\leq \rho$.
	For any $\nu>0$, we have by definition of $|K|^\nu$ that
	\begin{equation}\label{eq:interpolation1}
		\norm{|K|^\nu w}_{\bb X}^2 = \sum_{n\in\N} \sigma_n^{2\nu} |\inner{w}{v_n}_{\bb X}|^2.
	\end{equation}
	For $r>\nu$ and sequences $\{a_n\}_{n \in \N}$ and $\{b_n\}_{n \in \N}$, the Hölder inequality implies
	\begin{equation}\label{eq:HölderEq}
		\sum_{n\in\N} a_n b_n \leq \biggl(\sum_{n\in\N} a_n^\frac{r}{\nu}\biggr)^{\frac{\nu}{r}}\biggl(\sum_{n\in\N} b_n^\frac{r}{r-\nu}\biggr)^{1-\frac{\nu}{r}}.
	\end{equation}
	Inserting the estimate \eqref{eq:HölderEq} with
	\begin{equation}
		a_n \coloneqq \sigma_n^{2\nu} |\inner{w}{v_n}_{\bb X}|^{2\frac{\nu}{r}},\qquad
		b_n \coloneqq |\inner{w}{v_n}_{\bb X}|^{2-2\frac{\nu}{r}},
	\end{equation}
	into \eqref{eq:interpolation1} and applying Bessel inequality \eqref{eq:funktan:bessel} leads to
	\begin{equation}
		\norm{|K|^\nu w}_{\bb X}^2 \leq
		\biggl(\sum_{n\in\N} \sigma_n^{2r} |\inner{w}{v_n}_{\bb X}|^{2}\biggr)^{\frac{\nu}{r}}
		\biggl(\sum_{n\in\N} |\inner{w}{v_n}_{\bb X}|^{2}\biggr)^{1-\frac{\nu}{r}} \leq \norm{|K|^rw}_{\bb X}^{2\frac{\nu}{r}} \norm{w}_{\bb X}^{2(1-\frac{\nu}{r})}.
	\end{equation}
	By taking the square root, we get the \emph{interpolation inequality}
	\begin{equation}
		\label{eq:interpolation}
		\norm{|K|^\nu w}_{\bb X} \leq \norm{|K|^{r}w}_{\bb X}^{\frac{\nu}{r}} \norm{w}_{\bb X}^{1-\frac{\nu}{r}}.
	\end{equation}
	
	The application of \eqref{eq:interpolation} with $r =\nu + 1$ together with \Cref{lem:functional_range} leads to
	\begin{equation}
		\begin{aligned}
			\norm{x}_{\bb X} &= \norm{|K|^\nu w}_{\bb X} \leq \norm{|K|^{\nu+1}w}_{\bb X}^{\frac{\nu}{\nu+1}} \norm{w}_{\bb X}^{\frac{1}{\nu+1}} = \norm{\vert K \vert |K|^{\nu}w}_{\bb X}^{\frac{\nu}{\nu+1}} \norm{w}_{\bb X}^{\frac{1}{\nu+1}} \\
			&= \norm{\vert K \vert x}_{\bb X}^{\frac{\nu}{\nu+1}} \norm{w}_{\bb X}^{\frac{1}{\nu+1}} = \norm{Kx}_{\bb Y}^{\frac{\nu}{\nu+1}} \norm{w}_{\bb X}^{\frac{1}{\nu+1}}
			\leq \norm{Kx}_{\bb Y}^{\frac{\nu}{\nu+1}} \rho^{\frac{1}{\nu+1}}.
		\end{aligned}
	\end{equation}
	The claim follows by taking the supremum over all $x\in X_{\nu,\rho}$ for which $\norm{Kx}_{\bb Y}\leq \delta$.
\end{pr}

Actually, we can show that the upper bound \eqref{eq:UpperBoundError} is asymptotically sharp.
\begin{thm}{}{}
	Let $K\in\q K(\bb X,\bb Y)$ and $\nu,\rho>0$.
	Then there exists a sequence $\{\delta_n\}_{n\in\N}$ with $\delta_n \to 0$ such that $\eps(X_{\nu,\rho},\delta_n) = \delta_n^{\frac\nu{\nu+1}} \rho^{\frac1{\nu+1}}$.
\end{thm}
\begin{pr}
	We define $\delta_n \coloneqq  \rho \sigma_n^{\nu+1} \to 0$ and $x_n\coloneqq |K|^\nu (\rho v_n) \in X_{\nu,\rho}$.
	Since $\sigma_n^\nu = (\rho^{-1}\delta_n)^{\frac\nu{\nu+1}}$, it follows that
	\begin{equation}
		x_n = \rho |K|^\nu v_n  = \rho \sigma_n^\nu v_n =  \delta_n^{\frac\nu{\nu+1}}\rho^{\frac{1}{\nu+1}}v_n.
	\end{equation}
	Thus, we get $\norm{x_n}_{\bb X} =  \delta_n^{\frac\nu{\nu+1}}\rho^{\frac{1}{\nu+1}}$.
	Similarly, we can show that
	\begin{equation}
		K^*K x_n = \delta_n^{\frac\nu{\nu+1}}\rho^{\frac{1}{\nu+1}} \sigma_n^2 v_n = \delta_n^{\frac{\nu+2}{\nu+1}} \rho^{-\frac{1}{\nu+1}} v_n,
	\end{equation}
	which in turn leads to
	\begin{equation}
		\norm{Kx_n}_{\bb Y}^2 = \inner{K^*Kx_n}{x_n}_{\bb X} = \delta_n^{2}.
	\end{equation}
	Hence, the vectors $x_n \in X_{\nu,\rho}$ satisfy $\norm{Kx_n}_{\bb Y} \leq \delta_n$ for all $n\in\N$.
	Therefore, it holds
	\begin{equation}\label{eq:LowerBound}
		\eps(X_{\nu,\rho},\delta_n) = \sup_{x\in X_{\nu,\rho},\ \norm{Kx}_{\bb Y}\leq \delta_n}\norm{x}_{\bb X} \geq \norm{x_n}_{\bb X} = \delta_n^{\frac\nu{\nu+1}} \rho^{\frac1{\nu+1}}.
	\end{equation}
	Due to the matching upper bound in \Cref{th:rates:limit}, we must have equality in \eqref{eq:LowerBound}.
\end{pr}

Hence, for $K\in\q K(\bb X,\bb Y)$ with $\dim(\q R(K)) = \infty$, the worst-case error \eqref{eq:rates:worstcase} can at best decay as $\delta_n^{\frac\nu{\nu+1}} \rho^{\frac1{\nu+1}}$ for $\delta\to 0$.
Since ${\frac\nu{\nu+1}} < 1$ (which corresponds to Hölder type continuity), this rate is always slower than the data error, which decays as $\delta$.

Based on these findings, we say that a regularization method $(R_\alpha,\alpha)$ is \emphi{optimal} (for $\nu$ and $\rho$) if
\begin{equation}
	\sup_{x^\dag \in X_{\nu,\rho}} \q E(Kx^\dag,\delta) = 	\sup_{y\in\q D(K^\dag): K^\dag y\in X_{\nu,\rho}} \q E(y,\delta) = \delta^{\frac\nu{\nu+1}} \rho^{\frac1{\nu+1}},
\end{equation}
and \emphi{order optimal} (for $\nu$ and $\rho$) if there is merely $c >1$ with
\begin{equation}\label{eq:regularisierung:ordnung}
	\sup_{x^\dag \in X_{\nu,\rho}} \q E(Kx^\dag,\delta) = \sup_{y\in\q D(K^\dag): K^\dag y\in X_{\nu,\rho}} \q E(y,\delta) \leq c \delta^{\frac\nu{\nu+1}} \rho^{\frac1{\nu+1}}.
\end{equation}
To obtain asymptotic \emph{convergence rates}, we instead set $X_\nu = \q R(|K|^\nu)$ and say that $(R_\alpha,\alpha)$ is order optimal (for $\nu$) if there exists $c(x^\dag) \geq 1$ such that
\begin{equation}
	 \q E(Kx^\dag,\delta) \leq c(x^\dag)  \delta^{\frac\nu{\nu+1}} \qquad \text{for all } x^\dag \in X_{\nu}.
\end{equation}

If $x^\dag = K^\dag y \in X_{\nu,\rho}$ for some $\nu,\rho>0$, there indeed exists an a priori parameter choice rule that leads to an order optimal regularization for the spectral methods from \Cref{chap:spektral}.
\begin{thm}{\cite[Corollary 4.4]{EnglHankeNeubauer1996}}{spektral:apriori}
	Let $y\in\q D(K^\dag)$ and $x^\dag = K^\dag y\in X_{\nu,\rho}$ for some $\nu,\rho>0$. If $\alpha(\delta)$ is an a priori rule with
	\begin{equation}\label{eq:spektral:apriori}
		c \Bigl(\frac\delta\rho\Bigr)^{\frac2{\nu+1}} \leq \alpha(\delta) \leq C  \Bigl(\frac\delta\rho\Bigr)^{\frac2{\nu+1}} \qquad\text{for } C>c>0
	\end{equation}
	and the filter $\{\phi_\alpha\}_{\alpha>0}$ satisfies for some $C_r, C_\nu>0$ the conditions
	\begin{align}
		\sup_{\lambda\in(0,\kappa]} |\phi_\alpha(\lambda)| &= C_r \alpha^{-1},
		\label{eq:spektral:ordnung_phi}\\
		\sup_{\lambda\in(0,\kappa]} \lambda^{\nu/2}|r_\alpha(\lambda)| &\leq C_\nu \alpha^{\nu/2},\label{eq:spektral:ordnung_omega}
	\end{align}
	then $(R_\alpha,\alpha)$ is a (for this $\nu$ and all $\rho$) order optimal regularization method with
	\begin{equation}
		\sup_{y^\delta\in B_\delta(y)} \norm{R_{\alpha(\delta)} y^\delta - x^\dag}_{\bb X} \q  \leq (C_r c^{-1/2}+ C_\nu C^{\nu/2}) \delta^{\frac\nu{\nu+1}} \rho^{\frac1{\nu+1}}.
	\end{equation}
	The maximal $\nu_0$ such that \eqref{eq:spektral:ordnung_omega} holds for all $\nu \in (0, \nu_0]$ is the \emphi{qualification} of $\phi_\alpha$.
\end{thm}

\begin{example}[Truncated SVD]
	\label{ex:spectral:tsvd}
	From \eqref{eq:spektral:cutoff} we get that $\sup_{\lambda\in(0,\kappa]} |\phi_\alpha(\lambda)| \leq \alpha^{-1}$.
	Hence, the filter $\phi_\alpha$ satisfies \eqref{eq:spektral:ordnung_phi} with $C_r = 1$.
	Further, for all $\nu>0$ and $\lambda\in(0,\kappa]$, we have that
	\begin{equation}
		\lambda^{\nu/2} |r_\alpha(\lambda)| = \lambda^{\nu/2} |1-\lambda \phi_\alpha(\lambda)| =
		\begin{cases}
			0 & \text{if }\lambda\geq \alpha,\\
			\lambda^{\nu/2} & \text{if }\lambda<\alpha.
		\end{cases}
	\end{equation}
	Thus, it holds $\sup_{\lambda\in(0,\kappa]} \lambda^{\nu/2} |r_\alpha(\lambda)| \leq \alpha^{\nu/2}$  for all $\alpha\in(0,\kappa]$, and \eqref{eq:spektral:ordnung_omega} is  satisfied for all $\nu>0$ with $C_\nu = 1$.
	Hence, the truncated SVD is order optimal for any $\nu>0$.
\end{example}

\begin{example}[Tikhonov]
	\label{ex:spectral:tikhonov}
	From \eqref{eq:spektral:tikhonov}, we get that $\sup_{\lambda\in(0,\kappa]} |\phi_\alpha(\lambda)| \leq \alpha^{-1}$.
	Hence, the filter $\phi_\alpha$ satisfies \eqref{eq:spektral:ordnung_phi} with $C_r = 1$.
	Furthermore, it holds that
	\begin{equation}
		r_\alpha(\lambda) = 1-\lambda \phi_\alpha(\lambda) = \frac{\alpha}{\lambda + \alpha}\leq 1
		\qquad\text{for all }\alpha >0, \lambda\in(0,\kappa].
	\end{equation}
	Thus, regarding condition \eqref{eq:spektral:ordnung_omega}, we have to find  $\nu$ and a constant $C_{\nu}$ with
	\begin{equation}
		\sup_{\lambda\in(0,\kappa]}\lambda^{\nu/2} |r_\alpha(\lambda)| = \sup_{\lambda\in(0,\kappa]} \frac{\lambda^{\nu/2}\alpha}{\lambda+\alpha} \leq C_\nu\alpha^{\nu/2}.
	\end{equation}
	To this end, we have to find the maximum of  $\lambda \mapsto h_{\alpha,\nu}(\lambda) \coloneqq \frac{\lambda^{\nu/2}\alpha}{\lambda+\alpha}$.
	For its derivative, we get 
	\begin{equation}
		h_{\alpha,\nu}'(\lambda) = \frac{\alpha\frac\nu2\lambda^{\nu/2-1}(\lambda + \alpha) -  \alpha\lambda^{\nu/2}}{(\lambda+\alpha)^2} =  \frac{\alpha\lambda^{\nu/2-1}}{(\lambda+\alpha)^2}\left(\tfrac\nu2 \alpha + \left(\tfrac\nu2 -1 \right)\lambda\right).
	\end{equation}
	For $\nu\in(0,2)$, the root of $h_{\alpha,\nu}'(\lambda)$ is $\lambda^*= \frac{\alpha\nu}{2-\nu}$.
	Since $h_{\alpha,\nu}''(\lambda^*)<0$, this implies that
	\begin{equation}
		\sup_{\lambda\in(0,\kappa]} h_{\alpha,\nu}(\lambda) = \frac{\alpha\left(\alpha\nu(2-\nu)^{-1}\right)^{\nu/2}}{\alpha+\alpha\nu(2-\nu)^{-1}} \leq \left(\nu(2-\nu)^{-1}\right)^{\nu/2} \alpha^{\nu/2}
	\end{equation}
	as desired.
	For $\nu=2$, we observe that $h_{\alpha,\nu}$ is increasing.
	Then, the maximum is attained in $\kappa$ and
	\begin{equation}
		\sup_{\lambda\in(0,\kappa]} h_{\alpha,\nu}(\lambda) = \frac{\alpha\kappa^{\nu/2}}{\kappa + \alpha} \leq \kappa^{\nu/2-1} \alpha.
	\end{equation}
	Thus, Tikhonov regularization has (at least) qualification $\nu_0=2$.
	For $\alpha \sim \delta^{2/3}$ and some $C(x^\dag) > 0$, \Cref{th:spektral:apriori} now implies that
	\begin{equation}
		\norm{x_{\alpha(\delta)}^\delta-x^\dag}_{\bb X} \leq C(x^\dag) \delta^{\frac{2}{3}} \qquad \text{for all }x^\dag\in \q R(K^*K)\text{ and } y^\delta\in B_\delta(Kx^\dag).
	\end{equation}
\end{example}

To conclude this section, we provide more intuition on the assumption $x^\dag \in X_{\nu,\rho}$.
The latter is commonly referred to as \emphi{source condition}. 
As we have already seen, these are a prerequisite for establishing convergence rates.
Further, we call any $w\in \bb X$ for which $|K|^\nu w = x^\dag$ a \emph{source representer}.

\begin{example}
	For the integral operator $K$ from \Cref{ex:integration_svd}, the condition $x^\dag \in X_{2,\rho}$ implies the existence of some $w \in L^2([0,1])$ with $x^\dag =K^*Kw = \int_t^1 \int_0^s w(r)\,dr\,ds$.
	Consequently, we have that  $x^\dag$ is twice (weakly) differentiable with second derivative $w$.
	Moreover, it has to hold that $\Vert w \Vert_{L^2} \leq \rho$.
	Thus, the source condition $x^\dag \in X_{2,\rho}$ is equivalent to a Sobolev type regularity requirement for $x^\dag$.
\end{example}

The interpretation from the previous example can be extended. 
As we show below, $x^\dag \in X_\nu$ implies that the decay of $\inner{K x^\dag}{u_n}_{\bb Y}$  is faster than required by the Picard condition in \Cref{th:inverse:picard}.
Thus, we can see source conditions as an abstract regularity requirement.
\begin{lemma}\label{lem:regularisierung:quell}
	Let $K\in\q K(\bb X,\bb Y)$ and $y \in \q D(K^\dag)$.
	Then we have that $x^\dag = K^\dag y\in X_\nu$ if and only if
	\begin{equation}\label{eq:regularisierung:quell}
		\sum_{n\in\N} \sigma_n^{-2-2\nu} |\inner{y}{u_n}_{\bb Y}|^2 < \infty.
	\end{equation}
\end{lemma}
\begin{pr}
	$\Rightarrow$: For any $K^\dag y  \in X_\nu$, there exists $w\in \bb X$ with $K^\dag y = |K|^\nu w$.
	From the series representation \eqref{eq:inverse:picard_pseudo} for $K^\dag y$ and with the definition of $|K|^\nu$ in \eqref{eq:functional_sqrt}, we get that
	\begin{equation}
		\sum_{n\in\N} \sigma_n^{-1} \inner{y}{u_n}_{\bb Y} v_n  = \sum_{n\in\N} \sigma_n^{\nu} \inner{w}{v_n}_{\bb X} v_n.
	\end{equation}
	Since $\{v_n\}_{n \in \N}$ is orthonormal, the coefficients in both series must coincide, namely
	\begin{equation}\label{eq:regularisierung:quell1}
		\sigma_n^{-1} \inner{y}{u_n}_{\bb Y} =  \sigma_n^{\nu} \inner{w}{v_n}_{\bb X} \qquad\text{for all }n\in\N.
	\end{equation}
	By the Bessel inequality \eqref{eq:funktan:bessel} it holds $\sum_{n\in\N} |\inner{w}{v_n}_{\bb X}|^2 < \infty$.
	Inserting \eqref{eq:regularisierung:quell1} yields \eqref{eq:regularisierung:quell}.
	
	$\Leftarrow$: We can define $w = \sum_{n\in\N} \sigma_n^{-1-\nu} \inner{y}{u_n}_{\bb Y} v_n \in \bb X$.
	Then, we get that
	\begin{equation}
		|K|^\nu w =  \sum_{n\in\N}  \sigma_n^{\nu} \sigma_n^{-1-\nu} \inner{y}{u_n}_{\bb Y} v_n =   \sum_{n\in\N}  \sigma_n^{-1} \inner{y}{u_n}_{\bb Y} v_n= K^\dagger y.
	\end{equation}
\end{pr}

\subsubsection{Exercises}
\begin{exerciseB}\label{ex:functional_range}
	Establish the proof of \Cref{lem:functional_range}.
\end{exerciseB}
\begin{Solution}
		\circled{1}: Since $\{v_n\}_{n\in\N}$ is orthonormal, we get
		\begin{align}
			|K|^r ( |K|^s x) &= \sum_{n\in\N} \sigma_n^{r} \inner{|K|^s x}{v_n}_{\bb X} v_n =
			\sum_{n\in\N} \sigma_n^{r} \biggl \langle \sum_{m\in\N}\sigma_m^s\inner{x}{v_m}_{\bb X} v_m, v_n \biggr \rangle_{\bb X}v_n \notag\\
			=&\sum_{n\in\N} \sigma_n^{r} \left(\sigma_n^s\inner{x}{v_n}_{\bb X}\right) v_n = \sum_{n\in\N} \sigma_n^{r+s} \inner{x}{v_n}_{\bb X} v_n 
			= |K|^{r+s}x.
		\end{align}
		
		\circled{2}: For $x,z\in \bb X$ and $r\geq 0$, the defining properties of an inner product yield
		\begin{equation}
			\inner{|K|^r x}{z }_{\bb X} = \sum_{n\in\N} \sigma_n^r \inner{x}{v_n}_{\bb X}\inner{v_n}{z}_{\bb X} = \inner{x}{|K|^rz}_{\bb X}.
		\end{equation}
		
		\circled{3}:  This claim follows from \circled{1}, \circled{2}, and
		\begin{equation}
			\norm{|K|x}_{\bb X}^2 = \inner{|K|x}{|K|x}_{\bb X} = \inner{|K|^2 x }{x}_{\bb X} = \inner{K^*Kx}{x}_{\bb X} = \inner{Kx}{Kx}_{\bb Y} = \norm{Kx}_{\bb Y}^2.
		\end{equation}
		
		\circled{4}:
		Given a singular system $\{(\sigma_n,u_n,v_n)\}_{n\in\N}$ of $K$, the triples $\{(\sigma_n,v_n,u_n)\}_{n\in\N}$ and $\{(\sigma_n,v_n,v_n)\}_{n\in\N}$ are singular systems of $K^*$ and $|K|$, respectively.
		Now, for $x \in \bb X$, we define $y = \sum_{n\in\N}  \inner{x}{v_n}_{\bb X} u_n \in \bb Y$.
		Since
		\begin{equation}
			K^*y = \sum_{n\in\N} \sigma_n \inner{y}{u_n}_{\bb X} v_n = \sum_{n\in\N} \sigma_n \inner{x}{v_n}_{\bb X} v_n = |K| x,
		\end{equation}
		the range of $K^*$ and $|K|$ coincide.
\end{Solution}
\begin{exerciseB}[Regularization of Derivative Operator]
	Let $\bb X = L^2([0,1])$ and $K \colon \bb X \to \bb X$ with $Ku[s] = \int_0^s u(x) \diff x$.
	As approximation to $K^\dagger$, we study for $\alpha\in (0,1/2)$ the differential quotients $R_\alpha \colon \bb X \to \bb X$ with
	\begin{equation}
		(R_\alpha f)(x) \coloneq \frac{1}{\alpha}
		\begin{cases}
			f(x+\alpha) - f(x) & x \in \left[0, \tfrac{1-\alpha}{2}\right) \\
			f(x + \tfrac{\alpha}{2}) - f(x - \tfrac{\alpha}{2}) & x \in \left[\tfrac{1-\alpha}{2}, \tfrac{1+\alpha}{2}\right) \\
			f(x) - f(x - \alpha) & x \in \left[\tfrac{1+\alpha}{2}, 1\right].
		\end{cases}
	\end{equation}
	In the following, $H^2([0,1])$ denotes the Hilbert space
	\begin{equation}
		H^2([0,1]) = \left\{ f \in L^2([0,1]) \mid f'', f' \in L^2([0,1]) \right\}.
	\end{equation}
	\begin{enumerate}
		\item Assume that $f \in H^2([0,1]) \cap  \mathcal{D}(K^\dagger)$ and let $\|f - f^\delta\|_{L^2([0,1])} \leq \delta$.
		Show that
		\begin{equation}
			\|K^\dagger f - R_\alpha f^\delta\|_{L^2([0,1])} \leq \frac{\sqrt{6}}{\alpha} \delta + \frac{\sqrt{17}}{4} \alpha \|f''\|_{L^2([0,1])} .
		\end{equation}
		
		\item Determine an a-priori parameter choice rule such that $R_\alpha \colon L^2([0,1]) \to L^2([0,1])$ becomes a convergent regularization method.
		
		\item Discretize $R_\alpha$ by evaluating $R_\alpha$ at $2N$ discrete points $x_k \coloneqq \frac{(k - 1)\alpha}{2}$, $k \in \{1, \ldots, 2N\}$, for $\alpha = \frac{1}{N-1}$ and $N \in \N \setminus \{1\}$.
		Implement the associated $\tilde R_{\alpha,N} \colon \R^{2N} \to \R^{2N}$.
		
		\item Test your implementation for $\alpha = 2^{-k}$, $k \in \{2, 4, \ldots, 8\}$ and the functions
		\begin{enumerate}
			\item $f(x) = \cos(\pi x)$ for $x \in [0,1]$;
			\item $f(x) = 
			\begin{cases}
				0 & x \in [0, \tfrac{1}{3}[ \\
				x - \tfrac{1}{3} & x \in [\tfrac{1}{3}, \tfrac{2}{3}[ \\
				\tfrac{1}{3} & x \in [\tfrac{2}{3}, 1]
			\end{cases}$.
		\end{enumerate}
		Plot the maximum error $\|\tilde{R}_{\alpha,N} f -f^\prime \|_\infty$ in dependence of $\alpha$.
	\end{enumerate}
\end{exerciseB}
\begin{Solution}
	\circled{1}\,:  By the triangle inequality we have 
	\begin{align}
		& \|K^\dagger f - R_\alpha f^\delta\|_{L^2([0,1])} \leq \|K^\dagger f - R_\alpha f \|_{L^2([0,1])}  + \|R_\alpha f - R_\alpha f^\delta\|_{L^2([0,1])} \\
		&  =  \|K^\dagger f - R_\alpha f \|_{L^2([0,1])}  + \|R_\alpha(f-f^\delta)\|_{L^2([0,1])}   \leq  \|K^\dagger f - R_\alpha f \|_{L^2([0,1])}  + \delta\|R_\alpha\|_{\q L(\bb X,\bb X)}. 
	\end{align}
	We need to show
	$
	\|K^\dagger f - R_\alpha f \|_{L^2([0,1])}  + \delta\|R_\alpha\|  \leq \frac{\sqrt{6}}{\alpha} \delta + \frac{\sqrt{17}}{4} \alpha \|f^{\prime\prime}\|_{L^2([0,1])}^2.
	$
	
	\underline{Bounding $\|R_\alpha \|_{\q L(\bb X,\bb X)}$.} We bound 
	\begin{align}
		\alpha^2 \|R_\alpha f\|_{L^2([0,1])}^2 = & \int_{0}^{\frac{1-\alpha}{2}} \left(f\left(x+\alpha\right)  - f\left(x\right)\right)^2 dx + \int_{\frac{1-\alpha}{2}}^{\frac{1+\alpha}{2}} \left(f\left(x+\frac{\alpha}{2}\right) - f\left(x-\frac{\alpha}{2}\right)\right)^2 dx \\
		& + \int_{\frac{1+\alpha}{2}}^1 \left(f\left(x\right) - f\left(x-\alpha\right)\right)^2 dx \\
		&\leq 4\|f\|_{L^2([0,1])}^2  + 4\|f\|_{L^2([0,1])}^2  + 4\|f\|_{L^2([0,1])}^2  = 12\|f\|_{L^2([0,1])}^2 .
	\end{align}
	Thus, we have $\|R_\alpha \|_{\q L(\bb X,\bb X)} \leq \frac{\sqrt{12}}{\alpha}.$  
	
	\underline{Bounding \( \|K^\dagger f - R_\alpha f \|_{L^2([0,1])}\): }
	
	We use Taylor's theorem to bound $\|K^\dagger f - R_\alpha f \|_{L^2([0,1])}$.
	In particular, we first bound the integral in the left region $x\in [0,\frac{1-\alpha}{2})$. Here, we have
	\begin{equation}
	R_\alpha f(x) = \frac{f(x+\alpha) - f(x)}{\alpha}.
	\end{equation}
	By Taylor's theorem there exists $\xi(x)\in[x,x+\alpha)$ such that
	\begin{equation}
	f(x+\alpha) - f(x)  = \alpha f^\prime(x) + \frac{\alpha^2}{2} f^{\prime\prime}(\xi(x)).
	\end{equation}
	Hence, we obtain
	\begin{equation}
	R_\alpha f(x) - f^\prime(x) 
	= \frac{f(x + \alpha) - f(x)}{\alpha} - f^\prime(x) 
	= \frac{\alpha}{2} f^{\prime\prime}(\xi(x)).
	\end{equation}
	Computing the integral of the square over  $x\in [0,\frac{1-\alpha}{2})$ gives
	\begin{equation}
	E_L = \int_{0}^{\frac{1-\alpha}{2}} |R_\alpha f(x) - f^\prime(x)|^2 dx \leq \int_{0}^{\frac{1-\alpha}{2}} \frac{\alpha^2}{4} |f^{\prime\prime}(\xi(x))|^2 dt \leq  \frac{\alpha^2}{4} \int_0^1 |f^{\prime\prime}(t)| dt =  \frac{\alpha^2}{4} \|f^{\prime\prime}\|_{L^2([0,1])}^2.
	\end{equation}
	Analogously, we obtain the bound for the integral in the right region $x\in [\frac{1+\alpha}{2}, 1)$
	\begin{equation}
	E_R = \int_{\frac{1+\alpha}{2}}^1 |R_\alpha f(x) - f^\prime(x)|^2 dx \leq  \frac{\alpha^2}{4} \|f^{\prime\prime}\|_{L^2([0,1])}^2.
	\end{equation}
	For the the bound of the integral in the middle region $x\in [\frac{1-\alpha}{2}, \frac{1+\alpha}{2})$ we apply Taylor's theorem twice around $x$,
	\begin{align}
		& f\left(x+\frac{\alpha}{2}\right) = f\left(x\right)  + \frac{\alpha}{2} f^\prime\left(x\right) + \frac{1}{2}\left(\frac{\alpha}{2}\right)^2 f^{\prime\prime}\left(\xi_+\left(x\right)\right),\\
		& f\left(x-\frac{\alpha}{2}\right) = f\left(x\right)  - \frac{\alpha}{2} f^\prime\left(x\right) - \frac{1}{2}\left(\frac{\alpha}{2}\right)^2 f^{\prime\prime}\left(\xi_-\left(x\right)\right).
	\end{align}
	Then, we have
	\begin{align}
		R_\alpha f(x) = f^\prime(x) + \frac{\alpha}{8} \left [ f^{\prime\prime}\left(\xi_+\right) + f^{\prime\prime}\left(\xi_-\right)\right].
	\end{align}
	Moreover, we get
	\begin{align}
		&E_M = \int_{\frac{1-\alpha}{2}}^{\frac{1+\alpha}{2}} |R_\alpha f(x) - f^\prime(x)|^2 dx \leq \int_{\frac{1-\alpha}{2}}^{\frac{1+\alpha}{2}}  \frac{\alpha^2}{64 }| f^{\prime\prime}\left(\xi_+\right) + f^{\prime\prime}\left(\xi_-\right)|^2 dx  \\
		&\leq \frac{\alpha^2}{64 }  2 \|f^{\prime\prime}\|_{L^2([0,1])}^2 =  \frac{\alpha^2}{32 }   \|f^{\prime\prime}\|_{L^2([0,1])}^2.
	\end{align}
	Now, we get the desired bound
	\begin{align}
		&\|K^\dagger f - R_\alpha f \|_{L^2([0,1])} =  \|K^\dagger f - f^\prime \|_{L^2([0,1])} = \sqrt{E_L + E_M + E_R} \\
		&\leq \sqrt{\frac{\alpha^2}{4} \|f^{\prime\prime}\|_{L^2([0,1])}^2 + \frac{\alpha^2}{32} \|f^{\prime\prime}\|_{L^2([0,1])}^2 +  \frac{\alpha^2}{4} \|f^{\prime\prime}\|_{L^2([0,1])}^2 }
		= \alpha \frac{\sqrt{17}}{4\sqrt{2}} \|f^{\prime\prime}\|_{L^2([0,1])}. %\leq  \frac{\sqrt{17}}{4} \alpha c.
	\end{align}
	
	\circled{2}\,:  We set the a priori choice rule 
	\[
	\alpha(\delta, f^\delta) \coloneqq \sqrt{\delta}.
	\] 
	We immediately see that \(\{R_\alpha\}_\alpha\) with the above parameter choice rule \(\alpha\) is a convergent regularization method by applying the result from part \circled{1}
	\begin{align}
		&\lim_{\delta\to 0} \sup_{f^\delta\in B_\delta(f)} \left\| R_{\alpha\left(\delta,f^\delta\right)}f^\delta - K^\dagger f\right\|_{L^2([0,1])} \leq \lim_{\delta\to 0} \frac{\sqrt{6}}{\alpha(\delta, f^\delta)} \delta + \frac{\sqrt{17}}{4} \alpha(\delta, f^\delta) \|f^{\prime\prime}\|_{L^2([0,1])} \\
		&=  \lim_{\delta\to 0} \sqrt{17\delta} + \frac{\sqrt{17\delta}}{4} \|f^{\prime\prime}\|_{L^2([0,1])} = 0.
	\end{align}
	
	\circled{3} \, + \circled{4}\,:  Check Jupyter Notebook (\url{https://github.com/Sniper2k/inverse_problems_ai_exercises/blob/main/Ex3_3_12.ipynb}) for programming solution.
\end{Solution}

\begin{exerciseB}
	Let $\Omega = [0,1)$ and $k \in L^2(\Omega)$.
	We consider the periodic convolution $K\in \q L(L^2(\Omega), L^2(\Omega))$ with
	\begin{equation}
		[Kf](x) = \int_\Omega f(y)k(x - y) \diff y \quad \text{for almost every } x \in \Omega,
	\end{equation}
	By \cref{ex:SVD conv L2}, $K$ is compact and that $K$ possesses the right singular functions $\{e_n\}_{n \in \bb Z}$ with $e_n = e^{2\pi i nx}$.
	The corresponding singular values and left singular function are $\hat{k}_n = \left|\int_\Omega k(x)\overline{e_n} \diff x\right| = |\langle k, e_n \rangle|$ and $u_n = \frac{\langle k, e_n \rangle}{|\langle k, e_n \rangle|} e_n$.
	
	We say $f \in H_\#^\ell(\Omega)$ if $f \in L^2(\Omega)$ is $\ell$ times weakly differentiable (on $\R$) with derivatives in $L^2(\Omega)$. Recall that $f$ is weakly differentiable in $L^2([a,b])$ with weak derivative $g$ if
	\begin{equation}\label{eq: weak derivative}
		\int_{a}^b f(t) \varphi^\prime(t) \diff t = \int_{a}^b g(t) \varphi(t) \diff t
	\end{equation} 
	for all $\varphi \in C^{\infty}([a,b])$ with $\varphi(a) = \varphi(b) = 0$.
	Note that the weak differentiation $\partial^\ell \colon H_\#^\ell(\Omega) \subset L^2(\Omega) \to L^2(\Omega)$ is closed.
	
	\begin{enumerate}
		\item Show that $f \in H_\#^\ell(\Omega)$ if and only if $\sum_{n \in \bb Z} |n|^{2\ell} |\langle f, e_n \rangle|^2 < \infty$.
		
		\item Assume $\hat{k}_0 = 1$ and $c|n|^{-2} \leq |\hat{k}_n| \leq C|n|^{-2}$ for constants $c, C > 0$ and all $n \in \bb Z \setminus \{0\}$.
		Show that $K' \colon f \mapsto \partial K f$ is well-defined on $L^2(\Omega)$, with closed range and $\ker(K') = \operatorname{span}\{1\}$.
		Furthermore, derive $X_\ell = \q R(|K^\prime|^\ell) = H_\#^{1+\ell}(\Omega)$ for $\ell \in \N$.
	\end{enumerate}
\end{exerciseB}

\begin{Solution}
	\circled{1}\,: As derived in \eqref{eq: diff rule fourier}, the Fourier coefficients of the weak derivative satisfy
	$\hat{f}_n = \left( - \frac{1}{2 \pi i n} \right)^\ell \widehat{f^{(\ell)}}_n$
	$\Rightarrow$: Given that $\{e_n\}_{n \in \bb Z}$ is the ONB of $L^2(\Omega)$, for all $f \in H_\#^\ell(\Omega)$ by Bessel's inequality \eqref{eq:funktan:bessel} we get
	\begin{equation}
		\infty > \|f^{(\ell)} \|_{L^2(\Omega)}
		= \sum_{n \in \bb Z} |\widehat{f^{(\ell)}}_n|^2 
		= |2 \pi|^{2\ell} \sum_{n \in \bb Z} |n|^{2\ell} |\hat{f}_n|^2
		= |2 \pi|^{2\ell} \sum_{n \in \bb Z} |n|^{2\ell} |\inner{f}{e_n}|^2.
	\end{equation}
	$\Leftarrow$: We show it for $\ell = 1$ and the case $\ell > 1$ is analogous. Consider $z_N = \sum_{|n|\le N} (-2 \pi i n ) \inner{f}{e_n} e_n$, $N \in \N$. By assumption on the series, this sequence has a limit in $L^2(\Omega)$, which we denote by $z$. Then, there exists a subsequence $\{N(m)\}_{m \in \N}$ such that both $z_{N(m)}(x) \to z(x)$ and $\sum_{|n| \le N(m)} \inner{f,e_n} e_n(x) \to f(x)$ as $m \to \infty$ for almost every $x \in \Omega$. Now, we verify \eqref{eq: weak derivative}. Let $\varphi \in C^{\infty}(\Omega)$ with $\varphi(0) = \varphi(1) = 0$. Then,
	\begin{equation}
		\int_\Omega z(t) \varphi^{\prime}(t) \diff t
		= \int_\Omega \lim_{m \to \infty } z_{N(m)}(t) \varphi^{\prime}(t) \diff t
		= \lim_{m \to \infty } \int_\Omega  z_{N(m)}(t) \varphi^{\prime}(t) \diff t,
	\end{equation} 
	where the last equality follows by Lebesque's dominated convergence, since $z_{N(m)}(t) \varphi^{\prime}(t) \to z(t) \varphi^{\prime}(t)$ and we can construct the integrable majorant $g(t) =  \|\varphi^{\prime}\|_{C} \sum_{n \in \N} |2 \pi n \inner{f}{e_n} | \in L^1(\Omega)$.
	Substituting $z_{N(m)}$ into the integral gives
	\begin{align}
		& \int_\Omega z_{N(m)}(t) \varphi^{\prime}(t) \diff t
		= - \sum_{|n|\le N(m)} \inner{f}{e_n} (2 \pi i n ) \int_\Omega e^{2 \pi i n t} \varphi^{\prime}(t) \diff t\\
		& \qquad =  - \sum_{|n| \le  N(m)} \inner{f}{e_n} \left[ e^{2 \pi i n t} \varphi(t) \bigg|_0^1 - \int_\Omega e^{2 \pi i n t} \varphi(t) \diff t \right] 
		= \int_\Omega \sum_{|n| \le  N(m)} \inner{f}{e_n}  e_n(t) \varphi(t) \diff t. \notag
	\end{align} 
	Again, by Lebesque dominated convergence (with majorant $g(t) = \|f\|_{L^{2}(\Omega)} \|\varphi\|_{C} \in L^1(\Omega)$), we can bring the limit inside the integral,
	\begin{equation}
		\int_\Omega z(t) \varphi^{\prime}(t) \diff t
		= \int_\Omega  \lim_{m \to \infty } \sum_{|n| \le  N(m)} \inner{f}{e_n}  e_n(t) \varphi(t) \diff t
		= \int_\Omega f(t) \varphi(t) \diff t.
	\end{equation} 
	Hence, $z$ is the weak derivative of $f$. 
	
	\circled{2}\,: By construction $K f = \sum_{n \in \bb Z} \hat k_n \inner{f,u_n} e_n \in L^2(\Omega)$. Since $c n^{-2} \le |\hat k_n|$, we have 
	\begin{equation}
		\sum_{n \in \bb Z} |n|^{2} |\hat k_n \inner{f}{u_n}|^2 
		\le c^{-1} \sum_{n \in \bb Z} |n|^{-2} |\inner{f}{u_n}|^2 
		\le c^{-1} \sum_{n \in \bb Z} |\inner{f}{u_n}|^2 
		= c^{-1} \|f\|_{L^2(\Omega)} < \infty.
	\end{equation}
	Hence, $Kx \in H_\#^1(\Omega)$ and its weak derivative exists.
	Therefore, $K^{\prime}$ is well-defined and we can show that its actions is given by
	\begin{equation}
		K^{\prime} f = \sum_{n \in \bb Z} (2 \pi i n) \hat k_n \inner{f}{u_n} e_n.
	\end{equation}
	Note that the singular system of $K^{\prime}$ is $\{ ( 2 \pi |n| \hat k_n, -\sign(n) i u_n, e_n) \}_{n \in \bb Z}$. Therefore, by \Cref{lem:regularisierung:quell}, we have $f \in X_\ell = \q R(|K^\prime|^\ell) $ if and only if $\sum_{n \in \bb Z} |2 \pi n i \hat k_n|^{- 2 - 2\ell} |\inner{f}{u_n}|^2 < \infty$. Given the bounds $c |n|^{-2} \le |k_n| \le C |n|^{-2}$, the series is convergent if and only if $\sum_{n \in \bb Z} |n|^{2 + 2\ell} |\inner{f}{u_n}|^2 < \infty$, which by \circled{1} is equivalent to $f \in H_\#^{1+\ell}(\Omega)$.
\end{Solution}
\subsection{Tikhonov Regularization}\label{chap:tikhonov}

Due to its pioneering role for the development of variational regularization, we study the Tikhonov regularization $R_\alpha = \phi_\alpha(K^*K)K^*$ with the filter $\phi_\alpha = \frac{1}{\lambda+\alpha}$ from \Cref{ex:spektral}  in more detail.
First, we establish the promised alternative characterization.

\begin{lemma}\label{lem:tikhonov:normalen}
	Let $K\in\q K(\bb X, \bb Y)$ , $y\in \bb Y$ and $\alpha>0$.
	Then it holds $ x_\alpha =  R_\alpha y$ if and only if
	\begin{equation}
		\label{eq:tikhonov:normalen}
		(K^*K + \alpha I)x_\alpha = K^*y.
	\end{equation}
\end{lemma}
\begin{pr}
	$\Rightarrow$: Using a singular system $\{(\sigma_n,u_n,v_n)\}_{n\in\N}$ of $K$ and the filter definition \eqref{eq:spektral:tikhonov}, we get
	\begin{align}
		\alpha x_\alpha  &= \sum_{n\in\N} \alpha \frac{\sigma_n}{\sigma_n^2+\alpha} \inner{y}{u_n}_{\bb Y} v_n,\notag
		\\
		K^*Kx_\alpha &= \sum_{n\in\N} \frac{\sigma_n}{\sigma_n^2+\alpha} \inner{y}{u_n}_{\bb Y} K^*Kv_n =\sum_{n\in\N} \sigma_n^2 \frac{\sigma_n}{\sigma_n^2+\alpha} \inner{y}{u_n}_{\bb Y} v_n.
	\end{align}
	Adding both equations yields the claim
	\begin{equation}
		(K^*K + \alpha I)x_\alpha =  \sum_{n\in\N} \sigma_n  \inner{y}{u_n}_{\bb Y} v_n = K^* y.
	\end{equation}
	
	$\Leftarrow$: Let $x\in \bb X$ be a solution of \eqref{eq:tikhonov:normalen}. Inserting its basis decomposition
	\begin{equation}\label{eq:tikhonov:normal1}
		x = \sum_{n\in\N} \inner{x}{v_n}_{\bb X}v_n + P_{\q N(K)}x
	\end{equation}
	into both sides of \eqref{eq:tikhonov:normalen} gives
	\begin{equation}
		\sum_{n\in\N} (\sigma_n^2+\alpha)  \inner{x}{v_n}_{\bb X} v_n + \alpha  P_{\q N(K)}x =  \sum_{n\in\N} \sigma_n\inner{y}{u_n}_{\bb Y} v_n.
	\end{equation}
	Since $\{v_n\}_{n\in\N} \subset \q N(K)^\bot$ is orthonormal, we must have $P_{\q N (K)}x = 0$ and
	\begin{equation}\label{eq:RepresentationScalar}
		\inner{x}{v_n}_{\bb X} = \frac{\sigma_n}{\sigma_n^2+\alpha} \inner{y}{u_n}_{\bb Y} \qquad\text{for all }n\in \N.
	\end{equation}
	Inserting \eqref{eq:RepresentationScalar} into \eqref{eq:tikhonov:normal1} implies $x = \sum_{n\in\N} \inner{x}{v_n}_{\bb X}v_n = \sum_{n\in\N} \frac{\sigma_n}{\sigma_n^2+\alpha} \inner{y}{u_n}_{\bb Y} v_n$, namely that $x_\alpha$ is the only solution of \eqref{eq:tikhonov:normalen}.
\end{pr}

Thus, instead of computing the SVD of $K$, it suffices to solve a \emph{well-posed} linear equation (with a positive definite operator).
For the latter, we can rely on efficient iterative solvers such as the conjugate gradient algorithm.

In \Cref{ex:spectral:tikhonov}, we have shown that a priori rules of the form $\alpha \sim \delta^{2/3}$  lead to the error decay $\mathcal O(\delta^{2/3})$.
There, we could not prove that this rate is optimal.
Based on \Cref{lem:tikhonov:normalen}, we now show that indeed no a priori rule $\alpha$ with
$\norm{x_{\alpha(\delta)}^\delta-x^\dag}_{\bb X} \delta^{-2/3} \to 0$ can exist.
\begin{thm}{}{}
	Let $K\in\q K(\bb X,\bb Y)$ with $\dim(\q R(K)) = \infty$ and $y\in\q D(K^\dag)$.
	If  $\lim_{\delta\to 0}\alpha(\delta) = 0$ and
	\begin{equation}\label{eq:tikhonov:sat_ass}
		\lim_{\delta\to 0}\sup_{y^\delta\in B_\delta(y)}\norm{R_{\alpha(\delta)}y^\delta-x^\dag}_{\bb X}\delta^{-\frac23} = 0,
	\end{equation}
	then $x^\dag =0$.
\end{thm}
\begin{pr}
	Assume that $x^\dag\neq 0$.
	Using a singular system $\{(\sigma_n,u_n,v_n)\}_{n\in\N}$ of $K$, we define
	\begin{equation}
		\delta_n \coloneq  \sigma_n^3\quad\text{and}\quad y_n \coloneq  y + \delta_n u_n \in B_{\delta_n}(y),\qquad n\in\N.
	\end{equation}
	Here, we require that $\dim(\q R(K)) = \infty$ to have $\delta_n>0$.
	Using the shorthands $\alpha_n \coloneq  \alpha(\delta_n)$ and $x_{\alpha_n} = R_{\alpha_n}y$, we have that
	\begin{align}
			R_{\alpha_n} y_n - x^\dag &= R_{\alpha_n} (y_n - y) +  x_{\alpha_n} - x^\dag  \notag \\
			& = \sum_{m\in\N}\frac{\sigma_m}{\sigma_m^2+\alpha_n}\inner{\delta_n u_n}{u_m}_{\bb Y} v_m + x_{\alpha_n} - x^\dag \notag \\
			&= \frac{\delta_n\sigma_n}{\sigma_n^2+\alpha_n} v_n + (x_{\alpha_n} - x^\dag). 
	\end{align}
	Now, using assumption \eqref{eq:tikhonov:sat_ass} twice for $y^\delta = y_n$ and $y^\delta = y$ yields
	\begin{equation}\label{eq:DecayRate}
		\frac{\sigma_n\delta_n^{1/3}}{\sigma_n^2+\alpha_n} \leq \norm{R_{\alpha_n} y_n - x^\dag}_{\bb X}\delta_n^{-2/3} + \norm{x_{\alpha_n} - x^\dag}_{\bb X}\delta_n^{-2/3}\to 0 \qquad\text{as }n\to\infty.
	\end{equation}	
	At the same time, the left-hand side of \eqref{eq:DecayRate} can be written as
	\begin{equation}
		\frac{\sigma_n\delta_n^{1/3}}{\sigma_n^2+\alpha_n} = \frac{\delta_n^{2/3}}{\delta_n^{2/3} + \alpha_n} = \frac{1}{1+\alpha_n\delta_n^{-2/3}}.
	\end{equation}
	For a contradiction with \eqref{eq:DecayRate}, we have to show that $\lim_{\delta \to 0} \alpha(\delta)\delta^{-2/3} = 0$.	
	With the characterization \eqref{eq:tikhonov:normalen} of $R_{\alpha(\delta)}y^\delta$ and $K^*K x^\dag = K^*y$, we get
	\begin{equation}
		\left(K^*K+\alpha(\delta) I\right)\bigl(R_{\alpha(\delta)}y^\delta-x^\dag\bigr) = K^*y^\delta - K^* y -\alpha(\delta) x^\dag.
	\end{equation}
	Rearranging the terms such that $\alpha(\delta) x^\dag$ is isolated and using that  $y^\delta \in B_\delta(y)$ implies
	\begin{equation}\label{eq:EstAlphaDagger}
		|\alpha(\delta)| \norm{x^\dag}_{\bb X} \leq \norm{K}_{\q L(\bb Y,\bb X)} \delta + \bigl(\alpha(\delta) + \norm{K}_{\q L(\bb Y,\bb X)}^2 \bigr)\norm{R_{\alpha(\delta)}y^\delta-x^\dag}_{\bb X}.
	\end{equation}
	Multiplying \eqref{eq:EstAlphaDagger} with $\delta^{-2/3}$ and dividing by $\norm{x^\dag}_{\bb X}$ leads to
	\begin{equation}\label{eq:DecayRate2}
		|\alpha(\delta)|\delta^{-2/3}\leq \norm{x^\dag}_{\bb X}^{-1} \left(\norm{K}_{\q L(\bb Y,\bb X)}\delta^{\frac13} + (\alpha(\delta)+\norm{K}_{\q L(\bb Y,\bb X)}^2)\norm{R_{\alpha(\delta)}y^\delta-x^\dag}_{\bb X}\delta^{-\frac23}\right).
	\end{equation}
	Now, $\lim_{\delta \to 0} \alpha(\delta)\delta^{-2/3} = 0$ follows by using again assumption \eqref{eq:tikhonov:sat_ass}.
\end{pr}

Looking at equation \eqref{eq:tikhonov:normalen}, one might suspect that the Tikhonov regularization minimizes some functional.
This is indeed the case and for a generalization we refer to \Cref{ex: quadratic regularization}.
\begin{thm}{}{tikhonov:funktional}
	Let $y\in\bb Y$ and $\alpha>0$.
	Then, the unique minimizer of the \emphi{Tikhonov functional}
	\begin{equation}\label{eq:tikhonov:funktional}
		J_\alpha(x)\coloneq \frac12\norm{Kx-y}_{\bb Y}^2 + \frac\alpha2 \norm{x}_{\bb X}^2
	\end{equation}
	is given by $x_\alpha = R_\alpha y$ .
\end{thm}
\begin{pr}
	$\Rightarrow$:
	Using \eqref{eq:tikhonov:normalen}, we obtain for all $x\in \bb X$ that
	\begin{align}
			2(J_\alpha(x)-J_\alpha(x_\alpha)) 
			&= \| K x - y \|_{\bb Y}^2 + \alpha \| x  \|_{\bb X}^2 
			- \| K x_\alpha - y \|_{\bb Y}^2 - \alpha \| x_\alpha \|_{\bb X}^2 \notag \\
			&= \norm{Kx-Kx_\alpha}_{\bb Y}^2 + \alpha \norm{x-x_\alpha}_{\bb X}^2
			+ 2\bigl \langle K^*(Kx_\alpha -y)+\alpha x_\alpha, x-x_\alpha \bigr\rangle_{\bb X}\notag \\
			&= \norm{Kx-Kx_\alpha}_{\bb Y}^2 + \alpha \norm{x-x_\alpha}_{\bb X}^2 \geq 0. \notag
	\end{align}
	Hence, $x_\alpha$ minimizes $J_\alpha$.
	
	$\Leftarrow$: If $J_\alpha(x)-J_\alpha(\bar x)\geq 0$ for all $x\in \bb X$, we get for $x = \bar x + tz$ with $t>0$ and $z\in \bb X$ that
	\begin{equation}\label{eq:ScalarTik}
		0\leq 2\bigl(J_\alpha(\bar x + tz)-J_\alpha(\bar x)\bigr) = t^2\norm{Kz}_{\bb Y}^2 + t^2\alpha \norm{z}_{\bb X}^2 + 2t\bigl\langle K^*(K\bar x -y)  + \alpha\bar x,z\bigr\rangle_{\bb X}.
	\end{equation}
	Dividing \eqref{eq:ScalarTik} by $t>0$ and taking the limit $t\to 0$ gives
	\begin{equation}\label{eq:ScalarTik2}
		\inner{K^*(K\bar x -y)  + \alpha\bar x}{z}_{\bb X} \geq 0.
	\end{equation}
	Since $z\in \bb X$ was arbitrary, \eqref{eq:ScalarTik2} can only hold if $K^*K \bar x+\alpha \bar x =K^*y$. As $x_\alpha$ is the unique solution of \eqref{eq:tikhonov:normalen}, we obtain $\bar x = x_\alpha$. Hence, $x_\alpha$ is the unique minimizer of \eqref{eq:tikhonov:funktional}.
\end{pr}

 Instead of a least squares solution $x^\dag$, we look for minimizers of $\norm{Kx-y}_{\bb Y}^2 + \norm{x}_{\bb X}^2$ with trade-off parameter $\alpha$.\footnote{Tikhonov regularization was introduced in this form by \href{http://www-history.mcs.st-andrews.ac.uk/Biographies/Tikhonov.html}{Andre\u{\i} Nikolaevich Tikhonov} \cite{tikhonov1963}.}
For small noise levels $\delta$, we can put more emphasize on minimizing $\norm{Kx-y}_{\bb Y}^2$ (i.e., smaller $\alpha$ can be chosen).
On the other hand, the \emph{penalty} $\norm{x}_{\bb X}^2$  ensures the stability for large $\delta$.

\begin{remark}
	Based on \Cref{th:tikhonov:funktional}, we can establish a new interpretation for the source condition $x^\dag\in X_1=\q R(K^*)$.
	By construction, $x_\alpha^\delta$ also minimizes
	\begin{equation}\label{eq:tikhonov:funkt_ref}
		\frac{1}{2\alpha}\norm{Kx-y^\delta}_{\bb Y}^2 + \frac12\norm{x}_{\bb X}^2.
	\end{equation}
	Passing to the limit $\alpha\to 0$ and $\delta \to 0$  in \eqref{eq:tikhonov:funkt_ref}, we formally get the limiting problem
	\begin{equation}
		\label{eq:tikhonov:funkt_const}
		\min_{x\in \bb X,\ Kx=y}  \frac12\norm{x}_{\bb X}^2.
	\end{equation}
	Using a so-called Lagrange multiplier $p\in \bb Y$, we can rewrite \eqref{eq:tikhonov:funkt_const} as 
	\begin{equation}
		\min_{x\in \bb X}\max_{p\in \bb Y} \frac12\norm{x}_{\bb X}^2 - \inner{p}{Kx-y}_{\bb Y}.
	\end{equation}
	Computing the partial derivatives with respect to $x$ and $p$ leads to the optimality conditions $\bar x = K^*\bar p$ and $K\bar x = y$.
	For $y\in \q R(K)$, the minimizer of \eqref{eq:tikhonov:funkt_const} is $x^\dag$, i.e., $\bar x = x^\dag$.
	The existence of $\bar p$ with $x^\dag = K^*\bar p$ is thus equivalent to $x^\dag \in \q R(K^*)$.
\end{remark}

\subsubsection{Exercises}

\begin{exerciseB}
	Let $N=1024$ and consider the signal $f \in \bb C^{N}$ with 
	\begin{equation}
		f_k = \sin\bigl(\tfrac{2 \pi k}{N}\bigr) + 0.5 \sin\bigl(\tfrac{6 \pi k}{N}\bigr)  - 0.2 \sin\bigl(\tfrac{40 \pi k}{N}\bigr) + 0.3 \sin\bigl(\tfrac{10 \pi k}{N}\bigr)
	\end{equation}
	for  $k = 0, \ldots, N-1$.
	Further, we set $g \in \bb C^{N}$ to be the Gaussian kernel
	\begin{equation}
		g_k = \tilde g_{(k - N/2) \bmod N} 
		\quad \text{with} \quad 
		\tilde g_k = \exp\Bigl( - \frac{(k - \frac{N-1}{2})^2}{N} \Bigr), \quad k = 0, \ldots, N-1.
	\end{equation}
	
	\begin{enumerate}
		\item[(i)] Implement a Python function to compute $f * g$.
		Then, plot $f$, $g$ and $f * g$. 
		
		\textbf{Hint:} Use \texttt{numpy.fft.fft} and \texttt{numpy.fft.ifft} to evaluate $F_N$ and $F_N^{-1}$.
		\item[(ii)] Using \Cref{ex: svd convolution}, implement a function \texttt{pseudoinv} that computes $f^\dag = C(g)^\dag y$ for given $g$ and $y$. 
		Does it successfully reconstruct $f$ when applied to $y = f *_N g$? 
		\item[(iii)] To understand why \texttt{pseudoinv} fails, plot $|F_N g|$ in $\log_{10}$ scale. What do you observe? 
		Recall that we can only store numbers up to a certain precision.
		In this light, how does  the above plot explain the failure of \texttt{pseudoinv}? 
		\item[(iv)] Implement the truncated SVD and Tikhonov regularization using  \Cref{ex: svd convolution}.
		Then, apply them to $y = f *_N g$ with $\alpha \in \{ 2, 10^{-6}\}$.
		Can you recover $f$ now? 
		\item[(v)] Let $e \in \bb C^{N}$ with $e_k = \tfrac{0.1}{\sqrt N}\exp(60 \pi i k/N)$, where $i$ is the complex unit $i^2 = -1$.
		Plot the output of the implemented regularizations for $y = f *_N g + e$.
	\end{enumerate}
\end{exerciseB}

\begin{Solution}
	See Jupyter notebook \url{https://github.com/Sniper2k/inverse_problems_ai_exercises/blob/main/Ex3_4_5.ipynb}.
\end{Solution}

\begin{exerciseB}\label{ex: quadratic regularization}
	Let $\bb X = \R^N$, $\bb Y = \R^M$, $\bb E = \R^K$, $y \in \bb Y$, $w \in \bb E$, $A \in \R^{M \times N}$ and $B \in \R^{K \times N}$.
	Consider the functional 
	\begin{equation}
	J_\alpha(x) = \| Ax - y \|^2_2 + \alpha \| B x - w\|^2_2
	\end{equation}
	and show that it can be rewritten as $J_\alpha(x) = \| Kx - v\|^2_2$, where 
	\begin{equation}
	K = \begin{bmatrix}
		A \\ \sqrt{\alpha} B
		\end{bmatrix} \in \R^{(M+K) \times N},
	\quad \text{and} \quad 
	v = \begin{bmatrix}
		y \\ \sqrt{\alpha} w 
	\end{bmatrix} \in \R^{M + K}.
	\end{equation}
	From this, conclude that the minimizers of $J_\alpha$ are given by
	\begin{equation}
	(A^\tT A + \alpha B^\tT B) x = A^\tT y + \alpha B^\tT w.
	\end{equation}
\end{exerciseB}

\begin{Solution}
	We expand $\| Kx - v\|^2_2$ as
	\begin{align}
		\| Kx - v \|_2^2 &= \sum_{i = 1}^{M+K} (Kx - v)_i^2
		= \sum_{i = 1}^{M} (Kx - v)_i^2 + \sum_{i = M+1}^{M+K} (Kx - v)_i^2 \\
		& = \sum_{i = 1}^{M} (Ax - y)_i^2 + \sum_{i = 1}^{K} (\sqrt{\alpha}Bx - \sqrt{\alpha}w)_i^2
		= \|Ax - y\|_2^2 + \alpha \| Bx - w\|_2^2.
	\end{align} 
	Consequently, the minimizers of $J_\alpha$ are the least squares solutions of $Kx = v$, which by \cref{thm:inverse:normalen} are also the solution of normal equation
	\[
	K^\tT K x = K^\tT v. 
	\]
	Expanding both sides gives the desired formula.
\end{Solution}

\section{Optimization for Inverse Problems} \label{sec: classical algorithms}

\Reading
 \cite{bauschke2011convex}, \cite{bonnans2013perturbation} and \cite{ScherzerGrasmair2009}.

In \Cref{th:tikhonov:funktional}, we identified the Tikhonov regularized solution to an  inverse problem with the minimizer of the associated Tikhonov functional \eqref{eq:ScalarTik2}.
This idea can be generalized by considering the generic minimization problem
\begin{equation}\label{eq:data_fidelity_regularizer}
	\argmin_{x\in\bb X}J_\alpha(x),\qquad J_\alpha(x)=  D(F(x),y)+ \alpha R(x),\quad \alpha>0.
\end{equation}
Its first component is the \emphi{data-fidelity} $D(F(x),y)$, for which $D\colon  \bb Y \times \bb Y \to \R$ measures the distance between $F(x)$ and the data $y$.
Since  we assume $y \approx F(\hat x)$, any \enquote{good} reconstruction $x$ should satisfy $D(F(x),y) \approx 0$.
A common choice for $D$ (under a Gaussian noise model) is $ D(x,y)=\frac12\|x-y\|_{\bb Y}^2$, which also appears in the Tikhonov functional.
From a statistical perspective, the chosen data-fidelity $D$ should fit to the underlying noise model.

\begin{figure}[H]
	\centering
	\begin{subfigure}[t]{.45\textwidth}
		\includegraphics[width=\textwidth]{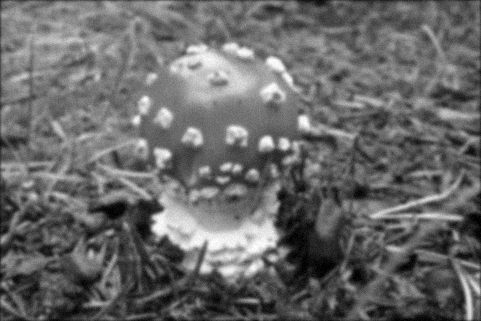}
		\caption*{Blurred image.}
	\end{subfigure}
	\begin{subfigure}[t]{.45\textwidth}
		\includegraphics[width=\textwidth]{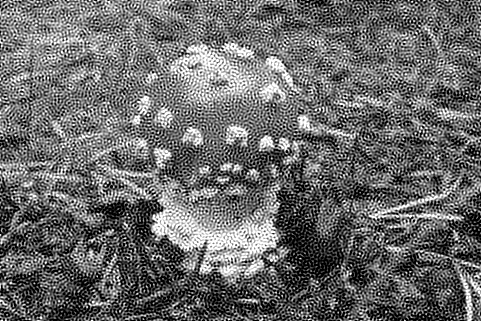}
		\caption*{Least squares solution.}
	\end{subfigure}
	\begin{subfigure}[t]{.45\textwidth}
		\includegraphics[width=\textwidth]{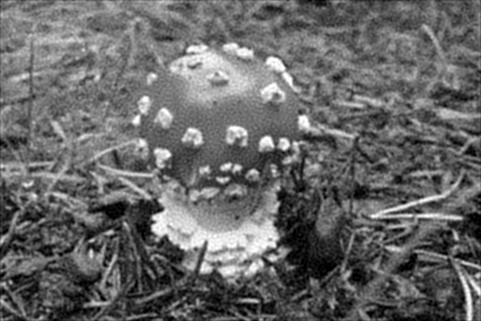}
		\caption*{Tikhonov regularization.}
	\end{subfigure}
	\begin{subfigure}[t]{.45\textwidth}
		\includegraphics[width=\textwidth]{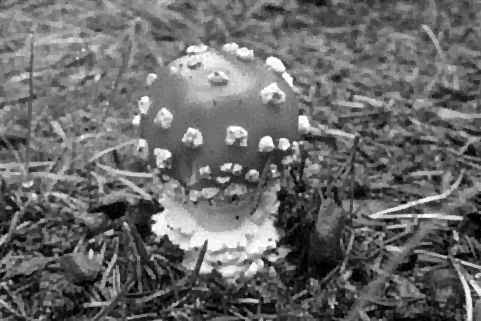}
		\caption*{Total variation regularization.}
	\end{subfigure}
	\caption{Deblurring by solving the variational problem \eqref{eq:data_fidelity_regularizer} for different regularizers.}
	\label{fig:deblurring_reconstructions}
\end{figure}

Secondly, we have the  \emphi{regularizer} $R(x)$.
It incorporates prior knowledge about the image $x$.
More precisely, $R \colon \bb X \to \R_+$ should be small if $x$ is a \enquote{natural} image and large otherwise.
The parameter $\alpha>0$ balances between the two.
In \Cref{fig:deblurring_reconstructions}, we show how the minimization of \eqref{eq:data_fidelity_regularizer} behaves for Tikhonov and total variation regularization. 
In comparison, the least squares solution without any regularization contains many artifacts.

\subsection{Minimization of Functions}

First, we review definitions and existence results from optimization theory.
Here, we continue to work with Hilbert space $\bb X$.
For any $U \subset \bb X$ and $f \colon U \to \R$, we can define an extension $\bar f \colon \bb X \to \overline{\R}$ with $\overline{\R} \coloneqq \R \cup \{\infty\}$ by setting $\bar f=f$ on $U$ and $\bar f(x) = \infty$ for $x\notin U$.
The set $\q D(f) \coloneqq \{x \in \bb X: f(x) < \infty\}$ is the \emphi{domain} of $f$.
This construction does not change the minimizers of $f$ and simplifies the notation.

\begin{defn}[Local / global minimum / minimizer]
	Let $f \colon \bb X\to \overline{\R}$.
	Then, $\hat{x} \in \q D(f)$ is a \emphi{local minimizer} of $f$ if there exist $r >0$ and an open ball $U_r(\hat{x})$ around $\hat{x}$ such that $f(\hat{x}) \le f(x)$ for all $x \in U_r(\hat{x})$.    
	Further, $\hat{x} \in U$ is a \emphi{global minimizer} of $f$ if $f(\hat{x}) \le f(x)$ for all $x \in \bb X$.
	In this case, we call $f(\hat{x})$ a \emph{local/global minimum} of $f$.
	Moreover, we denote the set of global minimizers by
	\begin{equation}
		\argmin(f)
		\coloneqq \big\{ \hat{x} \in \q D(f): f(\hat{x}) \le f(x) \ \forall x \in \bb X \big\}.
	\end{equation}
\end{defn}
The set of global minimizers may be empty, such as for the function
\begin{equation}\label{eq:ExampleLSC}
f_1(x) = x 
\quad \text{or} \quad
f_2(x) = \begin{cases}
	-x, & x < 0, \\
	1 + x, & x \ge 0.
\end{cases}
\end{equation}
The key requirements for the existence of global minimizers are given in the next definitions.

\begin{defn}\label{defn:lsc}
A function $f\colon \bb X\to \overline{\R}$ is \emphi{weakly lower-semicontinuous (lsc) at} $x_0 \in \q D(f)$ if for all sequences $\{x_n\}_{n\ge 0} \subset \bb X$ with $x_n \wkto x _0$ as $n \to \infty$ we have
\begin{equation}
f(x_0) \le \liminf_{n \to \infty} f(x_n).
\end{equation}
If $f$ is weakly lsc at all $x_0 \in \q D(f)$, then $f$ is \emphi{weakly lower-semicontinuous}.
\end{defn}

\begin{figure}[tbp]
	\centering
	\hspace{-2mm}
	\subfloat{
		\includegraphics[width=0.3\linewidth]{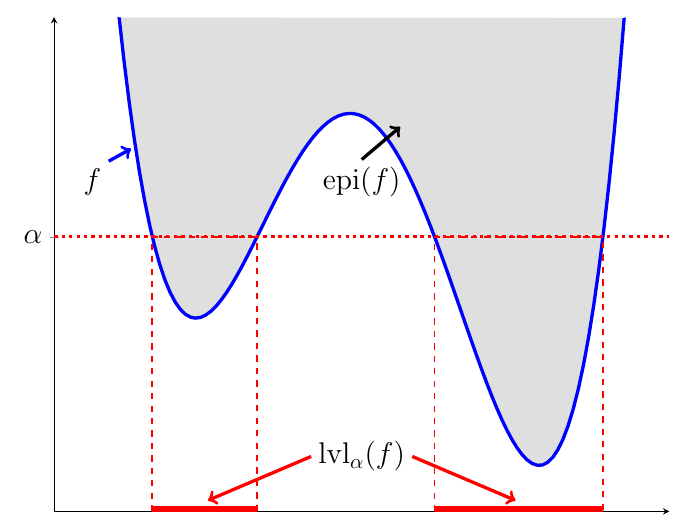}
		\vspace{2mm}
	}
	\hspace{-2mm}
	\subfloat{
		\includegraphics[width=0.3\linewidth]{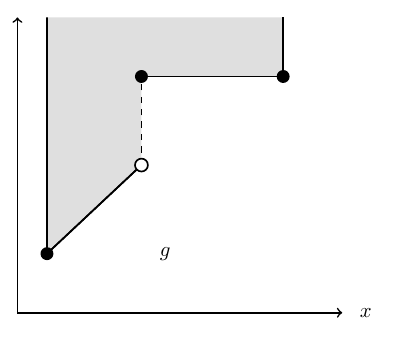}
	}
	\hspace{-2mm}
	\subfloat{
		\includegraphics[width=0.3\linewidth]{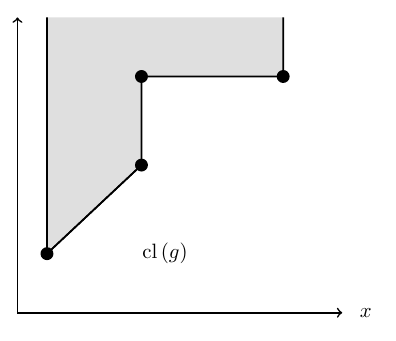}
		\label{fig:closure2}
	}
	\caption{\textbf{Left:} Epigraph $\epi(f)$ and level set $\lev_\alpha(f)$ of a function $f$; \textbf{Center}: A non lsc function $g$; \textbf{Right}: a lsc modification (closure) $\cl(g)$. The gray area illustrates the epigraph, which is only closed on the right and not the center plot.}
	\label{fig:closure}
\end{figure}

\begin{figure}[H]
	\centering

	\label{fig:closure}
\end{figure}

As the name suggests, weakly continuous functions, i.e., $f(x_n) \to f(x)$ as $x_n \wkto x$, are weakly lsc. 

\begin{example}
	\begin{itemize}
		\item[\circled{1}] The bounded operators $A \in \q L(\bb X, \bb Y)$ are weakly continuous.
		\item[\circled{2}] The norm $\norm{\cdot}_{\bb X}$ is weakly lsc, but not weakly continuous. 
	\end{itemize} 
\end{example}

However, the word weak in \Cref{defn:lsc} can be misleading: it refers to the fact that the property holds with respect to the weak topology.
Hence, it is a \emph{stronger} property than lower semi-continuity with respect to the norm topology.
Weak lower semicontinuity is equivalent to the geometric property that the \emphi{epigraph} $\mathrm{epi}(f)=\{ (x, a) \in \bb X \times \R \, : \, f(x) \le a\}$ of $f$ is weakly closed.
For convex $f$, this condition is equivalent to $\mathrm{epi}(f)$ being closed.
Therefore, (weakly) lsc functions are sometimes called closed. 

\begin{defn} \label{defn:coercive}
	A function $f \colon \bb X \to \overline{\R}$ is \emphi{coercive}, if $f(x_n) \to \infty$ as $\| x_n \|_\bb X \to \infty$.
\end{defn}
For $\alpha \in \R$, the \emphi{level set} $\lev_{\alpha }(f)$ is given by
\begin{equation}
	\lev_{\alpha }(f) \coloneqq \{ x \in \bb X \, : \, f(x) \le \alpha \}.
\end{equation}

The coercivity of $f$ is equivalent to the level sets $\lev_\alpha (f)$ being bounded for all $\alpha  \in \R$, see Exercise \ref{lemma:coerciveLvlbd}.
Furthermore, level sets of weakly lsc functions are weakly closed. 
In \eqref{eq:ExampleLSC}, $f_1$ is not coercive and $f_2$ is not (weakly) lsc. 
\begin{defn}\label{defn:convexity}
	A function $f \colon \bb X \to \overline{\R}$ is \emphii{convex}\index{convex function} if
	\begin{equation} \label{eq:convex}
		f(t x + (1 - t) z)
		\le t f(x) + (1 - t) f(z)
		\qquad \forall t \in [0, 1], \ \forall x, z \in \q D(f);
	\end{equation}
	\emphii{strictly convex}\index{strictly convex function} if
	\begin{equation}
		f(t x + (1 - t) z)
		< t f(x) + (1 - t) f(z)
		\qquad \forall t \in (0, 1), \ \forall x, z \in \q D(f), \ x \ne z;
	\end{equation}
	and \emphii{$\mu$-strongly convex}\index{strongly convex function} (with \emph{strong convexity constant} $\mu > 0$) if 
	\begin{equation}\label{eq:strong_cvx}
		f(t x + (1 - t) z)
		\le t f(x) + (1 - t) f(z) - \frac{\mu}{2} t (1 - t) \| x -  z \|^2_\bb X
		\ \forall t \in [0, 1], \ \forall x, z \in \q D(f).
	\end{equation}
\end{defn}
\begin{figure}[H]
	\centering
	\includegraphics[width=.4\textwidth, align=c]{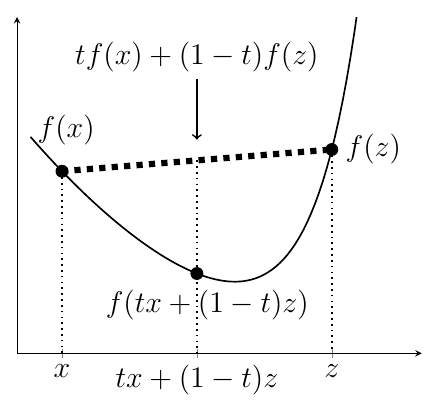} \hspace*{2cm}
	\includegraphics[width=.4\textwidth, align=c]{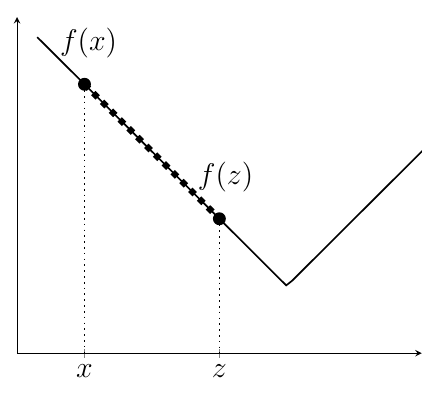}
	\caption{\textbf{Left}: Strictly convex $f$; \textbf{Right}: Convex, but not strictly convex $f$.}
	\label{fig:strictlyConvex}
\end{figure}
\begin{example}\label{example:convexity}
	\begin{itemize}
		\item[\circled{1}] For $\bb X = \R^N$ and $1 \leq p \leq \infty$, we have that $\Vert \cdot \Vert_p$ is convex.
		\item[\circled{2}] Linear functionals and their absolute value are convex, but not strictly convex.
		\item[\circled{3}] The function $\frac{1}{2} \| \cdot - y \|_{\bb X}^2$ with $y \in \bb X$ is 1-strongly convex.
	\end{itemize}
\end{example}

Strong convexity implies both coercivity and strict convexity, and strictly convex functions are convex. By setting $\mu = 0$ in \eqref{eq:strong_cvx}, we get back \eqref{eq:convex}, namely that $0$-strong convexity and convexity coincide.

Now, we prove existence of minimizers based on the properties from above.

\begin{thm}{Existence and uniqueness of global minima}{minima}
Let $f \colon \bb X \to \overline{\R}$ and $\q D(f) \ne \emptyset $.
	\begin{enumerate}
		\item 
		If $f$ is weakly lsc and coercive, then $\argmin(f) \ne \emptyset$.
		Furthermore, $\argmin(f)$ is weakly closed and bounded.
		\item
		Every local minimizer of a convex function $f$ is a global one.
		\item If $f$ is convex, $\argmin(f)$ is convex (but possibly empty).
		\item
		If $f$ is strictly convex, then \emph{at most} one minimizer exists.
	\end{enumerate}
\end{thm}
\begin{pr}
	\begin{enumerate}
		\item  Since $\q D(f) \ne \emptyset $, there exists $x_0\in \q D(f)$.
		For $a \coloneqq f(x_0)\in \R$, we choose $\{x_n\}_{n \in \N} \subset \q D(f)$ with $f(x_n) \leq a$ and $\lim_{n\to\infty} f(x_n)=\inf_{x\in \bb X} f(x)$.
		Since $f$ is coercive, the set $\lev_a(f)$ is bounded and therefore $\{x_n\}_{n \in \N} \subset  \lev_a(f)$ is bounded.
		By the Banach--Alaoglu \Cref{th:Banach_Alaoglu}, there exists a weakly convergent subsequence $\{x_{n_k}\}_{k \in \N}$ with limit $x \in \lev_a(f)$. 
		Since $f$ is weakly lsc, this implies
		\begin{equation}
			\inf_{x\in \bb X} f(x)\leq f(\hat x)\leq \liminf_{k\to\infty} f(x_{n_k})=\inf_{x\in \bb X} f(x).
		\end{equation}
		Hence, $f(\hat x)=\inf_{x\in \bb X} f(x)$ and $\hat x\in\argmin(f)$.
		The set $\argmin(f)=\lev_{f(\hat x)}(f)$ is bounded due to coercivity of $f$, and weakly closed since $f$ is weakly lsc.
		
		\item 
		Let $\hat x \in \q D(f)$ be a local minimizer, which is not a global one.
		Then, there exists $\tilde x \in \q D(f)$ with $f(\tilde x)<f(\hat x)$.
		In particular, we have for any $t\in (0,1)$ that
		\begin{equation}
			f((1-t)\hat x+t\tilde x)\leq (1-t) f(\hat x)+t f(\tilde x)<f(\hat x).
		\end{equation}
		Since $(1-t)\hat x+t\tilde x\to \hat x$ as $t\to 0$, we obtain that for any $\delta>0$ there exists $x\in B_\delta(\hat x)$ with
		$f(x)<f(\hat x)$.
		This contradicts the assumption that $\hat x$ is a local minimizer.
		
		\item  Let $\hat x_1,\hat x_2\in \argmin(f)$ and $t\in[0,1]$. Then, we have that
		\begin{equation}
			\inf_{x\in \bb X}f(x)\leq f((1-t)\hat x_1+t\hat x_2)\leq (1-t)f(\hat x_1)+t f(\hat x_2)=\inf_{x\in \bb X}f(x).
		\end{equation}
		Thus, $f((1-t)\hat x_1+t\hat x_2)=\inf_{x\in \bb X}f(x)$, which implies $(1-t)\hat x_1+t\hat x_2\in\argmin(f)$.
		
		\item
		Assume that $\hat x_1\neq \hat x_2 \in \argmin(f)$.
		By strict convexity, we obtain the contradiction
		\begin{equation}\label{eq:ConvexHess}
		f\bigl(\tfrac12 \hat x_1+\tfrac12 \hat x_2\bigr)<\tfrac12(f(\hat x_1)+f(\hat x_2))=\inf_{x\in \bb X}f(x).
		\end{equation}
		This concludes the proof.
	\end{enumerate}
\end{pr}

\begin{remark}
	Part \circled{1} of the proof follows the so-called \emph{direct method in the calculus of variations} or \textsc{Tonelli}'s direct method
	for proving the existence of minimizers of $f$.
\end{remark}
\begin{example}[Tikhonov Functional]
	The Tikhonov functional $\frac12\norm{Kx-y}_{\bb Y}^2 + \frac\alpha2 \norm{x}_{\bb X}^2$ from \Cref{th:tikhonov:funktional} is strongly convex and weakly lsc. 
	Hence, \Cref{th:minima} provides a new justification for the existence of a unique minimizer.
\end{example}

In analogy to the differentiability in finite-dimensional spaces, the derivative in normed spaces is based on directional derivatives.
The definition is named in honor of the French mathematician \textsc{René Eugène Gateaux}.

\begin{defn}[\textsc{Gateaux} derivative] \label{defn:GateauxDerivative}
	Let $x_0\in \bb X$ and $U \subset \bb X$ be a neighborhood of $x_0$.
	For $f\colon U\to \bb E$, the \emphi{one-sided directional derivative} of $f$ at $x_0$ into direction $v\in \bb X$ is defined (if it exists) by
	\begin{equation} \label{eq:GateauxLimit}
		\Diff f[x_0](v)\coloneqq \lim_{t \searrow 0} \frac{f(x_0 + t v) - f(x_0)}{t} \in \bb E.
	\end{equation}
	
	If the limit in \eqref{eq:GateauxLimit} exists for all $v$ and if $\Diff f[x_0] \in \q L(\bb X,\bb E)$, we call $\Diff f[x_0]$ the \emphii{\textsc{Gateaux} derivative}\index{Gateaux derivative@\textsc{Gateaux} derivative} of $f$ at $x_0$ and $f$ \emph{\textsc{Gateaux} differentiable at $x_0$}.
	
	If $f\colon U\to \bb E$ is \textsc{Gateaux} differentiable at $x_0$ for all $x_0 \in U^\circ$, then $f$ is differentiable (on $U$), and $\Diff f\colon U^\circ\to \q L(\bb X,\bb E)$ with $x\mapsto \Diff f[x]$ is the \emph{\textsc{Gateaux} differential}.
	A \textsc{Gateaux} differentiable function $f\colon U\to \bb E$ is \emph{twice \textsc{Gateaux} differentiable} (at $x_0\in U$) if the \textsc{Gateaux} differential $\Diff f\colon U\to \q L(\bb X,\bb E)$ is again \textsc{Gateaux} differentiable (at $x_0\in U$).
\end{defn}
Any \textsc{Gateaux} differentiable $f\colon\R^N\to\R^M$ is partially differentiable and the \textsc{Gateaux} derivative $\Diff f[x_0]$ is given by the Jacobian matrix $J_f[x_0] \in \R^{M \times N}$ of $f$ at $x_0$ via the matrix-vector multiplication $\Diff f[x_0](v)=J_f[x_0] v$.
Consequently, the \textsc{Gateaux} derivative coincides with the partial derivative in Euclidean spaces.
Vice versa, not every partially differentiable function is \textsc{Gateaux} differentiable as \Cref{ex: Gateaux differentiability} shows.

We require first- and second-order \textsc{Gateaux} derivatives of functions $f \colon \bb X \to \R$.
In this case, the \textsc{Gateaux} derivative $\Diff f[x_0]$ can be identified with an element in $\bb X$.
We call the latter \emphii{\textsc{Gateaux} gradient}\index{Gateaux gradient@\textsc{Gateaux} gradient}, denoted by $\nabla f(x_0) \in \bb X$, and recall that for all $v\in \bb X$ it holds 
\begin{equation}
	\Diff f[x_0](v)=\inner{\nabla f(x_0)}{v}.
\end{equation}
In contrast to  $\Diff f[x_0]$, the gradient $\nabla f(x_0)$ depends on the inner product.
Moreover, we can identify the second \textsc{Gateaux} derivative (which is a bilinear form in $\q L(\bb X,\q L(\bb X,\R))$ with an operator $\nabla^2 f(x_0)\in \q L(\bb X,\bb X)$ that we call  \emphii{\textsc{Gateaux} Hessian}\index{Gateaux Hessian@\textsc{Gateaux} Hessian}.
In particular, we have for all $u, v\in \bb X$ that 
\begin{equation}
	\Diff ^2f[x_0](u,v) = \inner{\nabla^2 f(x_0) u}{v}.
\end{equation}

We can use derivatives to characterize convexity.  
If $\q D(f)$ is open and $f$ is Gateaux differentiable on $\q D(f)$, then $\mu$-strong convexity is equivalent to 
\begin{equation}\label{eq: strong convexity}
	f(x) \ge f(z) + \inner{\nabla f(z)}{x-z} + \frac{\mu}{2} \| x - z\|^2_{\bb X}
	\quad \text{for all} \quad
	x,z \in \q D(f).
\end{equation}
%where $\nabla f$ denotes the Gateaux gradient of $f$.
If $f$ is twice Gateaux differentiable  on $\q D(f)$ with continuous Gateaux Hessian $\nabla^2 f$, then $f$ is $\mu$-strongly convex if and only if
\begin{equation}\label{eq: convexity grad}
	\inner{\nabla^2 f(x) z}{z} \ge \mu \|z\|^2_\bb X
	\quad \text{for all} \quad
	x \in \q D(f), z \in \bb X.
\end{equation}
A similar and useful characterization of $C^2$ functions with $L$-Lipschitz continuous gradient can be found in  \Cref{ex:Lipschitz Hessian}.

If $\q D(f)$ is open, \emph{Fermat's rule} states that if $f$ is Gateaux differentiable at a local minimum $x$, then $\nabla f(x) = 0$.
The points $x \in \q D(f)$ where $\nabla f(x) = 0$ are called \emphi{critical points}. 

\begin{thm}{Second order sufficient conditions}{second_order_sufficient}
Let $\q D(f)$ be open and $f$ be twice continuously Gateaux differentiable at a critical point $x \in \q D(f)$. 
If there exists $\alpha>0$ with 
\begin{enumerate}
	\item $\inner{\nabla^2 f(x) z}{z} \geq \alpha \norm{z}^2_{\bb X}$ for all $z \in \bb X$, then $x$ is a local minimum;
	\item $\inner{\nabla^2 f(x) z}{z} \leq  - \alpha \norm{z}^2_{\bb X}$ for all $z \in \bb X$, then $x$ is a local maximum.
\end{enumerate} 
\end{thm}
Beyond local extrema, the set of critical points includes \emphi{saddle points}, which are neither local minima nor local maxima.
For convex $f \colon \bb X \to \overline{\R}$, every critical point is a global minimum.
To see this for Gateaux differentiable $f$, let $x_0\in \bb X$ with $\nabla f(x_0)= 0$.
For any $x\in \q D(f)$, we have by \eqref{eq: convexity grad} that $0= \inner{\nabla f(x_0)}{x-x_0} \leq f(x)-f(x_0)$, i.e., $f(x_0)\leq f(x)$.
Consequently, we get $x_0\in \argmin(f)$.

\subsubsection{Exercises}

\begin{exerciseB}
	Find $f \colon \R \to \R$ that is
	\begin{enumerate}
		\item[(i)] lower bounded and continuous with $\argmin(f) = \emptyset$.
		\item[(ii)] lower bounded, continuous but not coercive, convex with $\argmin(f) = \{0\}$.
	\end{enumerate}
\end{exerciseB}

\begin{Solution}
	(i): The exponential function $\exp \colon \R \to \R$, $x \mapsto e^x$ is strictly convex since $\exp''(x) = \exp(x) > 0$ for all $x \in \R$, but it has no minimizers: it is strictly increasing and $\exp(x) \searrow 0$ for $x \searrow - \infty$.
	
	(ii): We can take the function
	\begin{equation}
		f(x) = \begin{cases}
			x^2, & |x| \le 1, \\
			1, & |x| > 1.
			\end{cases}
	\end{equation}
\end{Solution}

\begin{exerciseB}[Coercivity and level-boundedness] \label{lemma:coerciveLvlbd}
	Show that $f \colon \bb X \to \overline{\R}$ is coercive if and only if it is level-bounded.
\end{exerciseB}

\begin{Solution}
	$\Rightarrow$: Assume $f$ is coercive. We want to show that for any $\alpha \in \R$, the sublevel set $\lev_\alpha(f) = \{x \in \bb X : f(x) \le \alpha\}$ is bounded.
	
	Let's assume, for the sake of contradiction, that $\lev_\alpha(f)$ is not bounded for some $\alpha \in \R$.
	If $\lev_\alpha(f)$ is not bounded, then there exists a sequence $\{x_n\}_{n\in\N}$ in $\lev_\alpha(f)$ such that $\| x_n \| \to \infty$ as $n \to \infty$.
	Since $x_n \in \lev_\alpha(f)$ for all $n$, we have $f(x_n) \le \alpha$ for all $n$.
	
	However, since $f$ is coercive and $\|x_n\| \to \infty$, we must have $f(x_n) \to \infty$.
	This means that for any $M > \alpha$, there exists an $n_0 \in \N$ such that for all $n \ge n_0$, $f(x_n) > M > \alpha$.
	This contradicts the fact that $f(x_n) \le \alpha$ for all $n$.
	Therefore, our assumption that $\lev_\alpha(f)$ is not bounded must be false.
	Hence, $\lev_\alpha(f)$ is bounded for all $\alpha \in \R$, which means $f$ is level-bounded.
	
	$\Leftarrow$: Assume $f$ is level-bounded. We want to show that $f(x) \to \infty$ as $\| x \| \to \infty$.
	This means we need to show that for any $M \in \R$, there exists an $R > 0$ such that for all $x \in \bb X$ with $\| x \| > R$, we have $f(x) > M$.
	
	Let $M \in \R$ be arbitrary.
	Since $f$ is level-bounded, the level set $\lev_M(f)$ is bounded.
	Because $\lev_M(f)$ is bounded, there exists an $R > 0$ such that for all $x \in \lev_M(f)$, $\|x\| \le R$.
	
	Then for any $x \in \bb X$ such that $\|x\| > R$,  we must have $f(x) > M$.
	Since $M$ was arbitrary, this shows that $f(x) \to \infty$ as $\|x\| \to \infty$.
	Hence, $f$ is coercive.
\end{Solution}

\begin{exerciseB}[One level set suffices] \label{exercise:oneisenough}
	Let $\bb X = \R^N$.
	Show that a convex and lsc function $f \colon \bb X\to \overline{\R}$ is \emph{coercive} if and only if \emph{there exists} $a \in \R$ such that $\lev_a(f)$ is non-empty and bounded.
\end{exerciseB}

\begin{Solution}
	$\Rightarrow$:
	Follows by \cref{lemma:coerciveLvlbd}.

	$\Leftarrow$:
	For $b \le a$ we have $\lev_b(f) \subset \lev_a(f)$.
	Thus, consider $b > a$ and assume that the level set $\lev_b(f)$ is unbounded.
	Denote by $f_* \coloneqq \inf_{x \in \bb X} f(x)$. There exists $\{y_n\}_{n \in \N} \subset \lev_a(f)$ such that $f(y_n) \to f_*$. Since $\lev_a(f)$ is bounded and closed, we can find $y \in \lev_a(f)$ such that $y_{n(m)} \to y$ as $m \to \infty$ and by $f$ being lsc, we get $f(y) \le \liminf_{m \to \infty} f(y_{n(m)}) = f_*$. Hence, $f(y) = f_*$ and the minimum is attained. Moreover, by boundedness of $\lev_a(f)$ there exists $R>0$ such that $\|z\|_\bb X \le R$ for all $z \in \lev_a(f)$.
	
	By assumption, $\lev_b(f)$ is unbounded and we can find $\{x_n\}_{n \in \N} \subset \lev_b(f)$ such that $\|x_n\|_{\bb X} \to \infty$ as $n \to \infty$. Without loss of generality assume that $\|x_n\|_{\bb X} \ge 6R$. Consider $z_n \coloneqq (1 - \frac{3R}{\|x_n\|_{\bb X}}) y + \frac{3R}{\|x_n\|_{\bb X}} x_n$. This sequence admits
	\begin{equation}
		1.5 R \le 3R \frac{\|x_n\|_{\bb X}}{\|x_n\|_{\bb X}} - \left(1 + \frac{3R}{\|x_n\|_{\bb X}}\right) \|y\|_{\bb X} \le 
		\|z\|_\bb X 
		\le \left(1 + \frac{3R}{\|x_n\|_{\bb X}}\right) \|y\|_{\bb X} + 3R \frac{\|x_n\|_{\bb X}}{\|x_n\|_{\bb X}} \le 4.5 R.
	\end{equation}
	Hence, it is bounded and $z_n \notin \lev_a(f)$. By convexity, we get 
	\begin{equation}
		f(z_n) \le \left(1 - \frac{3R}{\|x_n\|_{\bb X}} \right) f(y) + \frac{3R}{\|x_n\|_{\bb X}} f(x_n) 
		\le \left( 1 - \frac{3R}{\|x_n\|_{\bb X}} \right) f_* + \frac{3R}{\|x_n\|_{\bb X}} b.
	\end{equation}
	Rearranging gives $f(z_n) - f_* \le \frac{3R}{\|x_n\|_{\bb X}}(b - f_*) \to 0$ as $n \to \infty$. Since $\{z_n\}_{n \in \N}$ is bounded, there exists a subsequence $\{ z_{m(n)} \}_{m \in \N}$ such that $z_{m(n)} \to z$ as $m \to \infty$. The limit point satisfies $\|z\|_{\bb X} \ge 1.5R$, i.e., $z \notin \lev_a(f)$. On the other hand 
	\begin{equation}
		f_* \le f(z) \le \liminf_{m \to \infty} f(z_{m(n)}) 
		\le \limsup_{m \to \infty} f(z_{m(n)}) \le f_*,
	\end{equation}   
	and $f(z) = f_* \in \lev_a(f)$ must hold giving a contradiction.
\end{Solution}

\begin{exerciseB}\label{ex: Gateaux differentiability}
Show that the function
\begin{equation}
	f \colon \R^2 \to \R, \qquad
	(x, y) \mapsto \begin{cases}
		\frac{x^2 y}{x^2 + y^2} & \text{for } (x, y) \ne (0,0),\\
		0                       & \text{else}
	\end{cases}
\end{equation}
is partially differentiable.
Further, show that it is not Gateaux differentiable in $0$. 
\end{exerciseB}

\begin{Solution}
	Partial derivatives in $x$- and $y$ direction are given by
	\begin{eqnarray*}
		f_x(x, y) 
		= \frac{2xy^3}{(x^2+y^2)^2}
		\qquad \text{and} \qquad
		f_y(x, y)
		= \frac{x^2(x^2-y^2)}{(x^2+y^2)^2}
	\end{eqnarray*}
	for $(x, y) \ne (0, 0)$ and $f_x(0,0) = f_y(0,0) = 0$.
	Note that the $f_x$ is not continuous in $x_0 \coloneqq (0,0)$ (choose $(x_n,y_n)=(\frac{1}{n},\frac{1}{n}) \xrightarrow{n \to \infty} (0,0)$).
	The directional derivative in $x_0$ is not a linear function since
	\begin{eqnarray*}
		\lim_{t\to 0}\frac{f(x_0+t h)-f(x_0)}{t} 
		= \lim_{t\searrow 0}\frac{t^2 h_1^2 t h_2}{t (t^2 h_1^2 + t^2 h_2^2)}
		=\lim_{t \searrow 0}\frac{h_1^2 h_2}{h_1^2+h_2^2}
		=\frac{h_1^2h_2}{h_1^2+h_2^2}
	\end{eqnarray*}
	and the last term is not linear in $h$.
	Thus, $f$ is not Gateaux differentiable in $x_0$.
\end{Solution}

\begin{exerciseB}[Monotonicity of $\nabla f$]\label{ex: strong convexity first order}
Show that for $\mu$-strongly convex and Gateaux differentiable $f \colon \bb X \to \R$ it holds
\begin{equation}
\| \nabla f(x) - \nabla f(x) \|_{\bb X} \ge \mu \| x- z\|_\bb X \quad \text{for all} \quad x,z \in \bb X.
\end{equation}
\end{exerciseB}

\subsection{Variational Regularization}\label{sec: variational regularization} 

Returning to \eqref{eq:data_fidelity_regularizer}, we consider the data fidelity $D(K x,y) = \|Kx - y\|^p_{\bb Y}$, $p \geq 1$, with $K \in \q K(\bb X,\bb Y)$.
We have that $x \mapsto D(K x,y)$ is convex and weakly lsc, however, when $\q N(K) \neq \{0\}$, it is not coercive or strictly convex.
Then, $R \colon \bb X \to \overline{\R}$ has to ensure that \Cref{th:minima} is applicable.
Throughout, we write $J_{\alpha, y}$ to express the dependence on $y \in \bb Y$.

\begin{cor}[Existence]\label{cor:existence variational}
	Let $D(K x,y) = \|Kx - y\|^p_\bb Y$, $p \geq 1$.
	Additionally, let $R \colon \bb X \to \overline{\R}$ be weakly lsc and coercive with $\q D(R) \neq \emptyset$.
	Then, for every $ y\in \bb Y$, we have $\argmin(J_{\alpha, y}) \ne \emptyset$ for \eqref{eq:data_fidelity_regularizer}.
\end{cor}

\begin{cor}[Uniqueness]\label{cor:uniqueness variational}
	Let the assumptions of \Cref{cor:existence variational} hold.
	Additionally, let $R\colon \bb X \to \overline{\R}$ be strictly convex.
	Then, there exists a unique minimizer of the functional $J_{\alpha,y}$ in \eqref{eq:data_fidelity_regularizer}.
\end{cor}

If \Cref{cor:uniqueness variational} does not apply, we might need to choose among possible reconstructions.
Hence, we generalize the concept of minimal norm solutions, which correspond to $R=\Vert \cdot \Vert$.
\begin{defn}
	An element $\hat x \in \bb X$ is a \emphi{$R$-minimizing solution} for measurements $y \in \q R(K)$ if $K \hat x = y$ and $R(\hat x) = \min\{R(x) : x \in \bb X,\ Kx = y\}$.
\end{defn}
\begin{cor}
	Let the assumptions of \Cref{cor:existence variational} hold.
	If there exists $x\in \q D(R)$ with $Kx=y$, then there exists a $R$-minimizing solution.
\end{cor}
Compared to \Cref{sec:Regularization}, we now have a \emph{nonlinear} regularization method.
Hence, a new analysis becomes necessary.
For this, we require the following auxiliary result.
\begin{lemma}\label{lem:EstData}
	Let $p \geq 1$ and $J_{\alpha, y}$ as in \eqref{eq:data_fidelity_regularizer}.
	For all $x \in \q D(R)$ and $y_1, y_2 \in \bb Y$, we have that
	\begin{equation}
			J_{\alpha, y_1}(x) \leq 2^{p-1} 	J_{\alpha, y_2}(x)+2^{p-1}\left\|y_{1}-y_{2}\right\|_{\bb Y}^{p}.
	\end{equation}
\end{lemma}
\begin{pr}
	We first show that $(a + b)^p \le 2^{p-1}(a^p + b^p)$ for $a,b \ge 0$. If $p=1$, the inequality is trivial. For $p>1$, let $q$ be its Hölder dual so that $\frac{1}{p} + \frac{1}{q} = 1$. Then, by Hölder's inequality  
	\begin{equation}\label{eq:HoeldEst}
		(a + b)^p 
		= (1 \cdot a + 1 \cdot b)^p 
		\le (1^q + 1^q)^{p/q} (a^p + b^p)^{p/p} 
		= 2^{p-1} (a^p + b^p). 
	\end{equation}
	With \eqref{eq:HoeldEst}, we obtain the estimate as follows
	\begin{align}
		J_{\alpha, y_1}(x) 
		& = \|Kx-y_1\|_{\bb Y}^{p}+\alpha R (x) \leq 2^{p-1}\bigl(\|Kx-y_2\|_{\bb Y}^{p}+\|y_1 - y_2\|_{\bb Y}^{p}\bigr)+\alpha R(x) \notag\\
		& \leq 2^{p-1}\left(\left\|Kx-y_2\right\|_{\bb Y}^{p}+\alpha R (x)\right)+2^{p-1} \|y_1-y_2\|_{\bb Y}^p\notag\\ 
		&=2^{p-1}J_{\alpha, y_2}(x) + 2^{p-1} \|y_1-y_2 \|_{\bb Y}^p.
	\end{align}
	This concludes the proof.
\end{pr}
As counterpart to the continuity property in \Cref{def:regularisierung}, we get the following.
\begin{thm}{Stability}{stability variational}
	Let the assumptions of \Cref{cor:existence variational} hold.
	If $\{y_{n}\}_{n \in \N}$ converges to $y \in \bb Y$, then every sequence $\{x_{n}\}_{n \in \N}$ with $x_n \in \argmin (J_{\alpha, y_n})$ has weak accumulation points.
	Every such accumulation point minimizes the functional $J_{\alpha, y}$.
\end{thm}
\begin{pr}
	Let $\bar{x} \in \q D(R)$.
	Since $x_n \in \argmin (J_{\alpha, y_n})$, we have $J_{\alpha, y_n} (x_n) \leq J_{\alpha, y_n}(\bar{x}) < \infty$.
	By applying \Cref{lem:EstData} twice, it follows that
	\begin{align}
		J_{\alpha, y}(x_n) 
		& \leq 2^{p-1} J_{\alpha, y_n}(x_n)+2^{p-1}\|y_n-y\|_{\bb Y}^p\notag\\
		& \leq 2^{p-1} J_{\alpha, y_n}(\bar{x})+2^{p-1}\|y_n-y\|_{\bb Y}^p 
		\leq 4^{p-1} 	J_{\alpha, y}(\bar{x})+4^p \|y_n-y\|_{\bb Y}^p.
	\end{align}
	Since $\lim_{n \to \infty} y_n = y$, there exists $n_{0} \in \mathbb{N}$ such that
	\begin{equation}
		M \coloneqq 4^{p-1} J_{\alpha, y} (\bar{x})+1 \geq J_{\alpha, y} (x_n), \quad n \geq n_{0}.
	\end{equation}
	Thus, $x_{n} \in \lev_{M}(J_{\alpha, y})$ for $n \geq n_{0}$.
	Due to coercivity, $\lev_{M}(J_{\alpha, y})$ is bounded and by the Banach--Alaoglu \Cref{th:Banach_Alaoglu} at least one weak accumulation point of $\{x_{n}\}_{n \in \N}$ exists.
	
	Now, let $x_{n_k} \wkto x \in \bb X$.
	Since $K \in \q K(\bb X,\bb Y)$, it follows that $K x_{n_k} \wkto Kx$.
	Hence, we also have $K x_{n_k}  -y_{n_k}  \wkto Kx-y$.
	Because $\|\cdot\|_{\bb Y}$ and $R$ are weakly lsc, it follows that
	\begin{align}
		& \|Kx - y \|_{\bb Y}^p \leq \liminf_{k \to \infty} \| K x_{n_k}  -y_{n_k}  \|_{\bb Y}^{p}, \label{eq:liminf1}\\
		& R(x) \leq \liminf_{k \to \infty}R (x_{n_k} ). \label{eq:liminf2}
	\end{align}
	With \eqref{eq:liminf1} and \eqref{eq:liminf2} we get
	\begin{align}
		 J_{\alpha, y} (x) & = \|Kx - y \|_{\bb Y}^p +\alpha R(x) \leq \liminf_{k \to \infty} \| K x_{n_k}  -y_{n_k}  \|_{\bb Y}^{p}+\alpha \liminf_{k \to \infty}R (x_{n_k} ) \notag \\
		&  \leq \liminf_{k \to \infty} \| K x_{n_k}  -y_{n_k}  \|_{\bb Y}^{p}+\alpha R (x_{n_k} )
		= \liminf_{k \to \infty} J_{\alpha, y_{n_k} } (x_{n_k} ). 
	\end{align}
	Since $x_{n_k}  \in \argmin (J_{\alpha, y_{n_k} })$ and for every $\tilde x$ the functional $y \mapsto J_{\alpha, y}(\tilde x)$ is continuous, we conclude that 
	\begin{equation}
		J_{\alpha, y} (x) 
		\leq \liminf_{k \to \infty} J_{\alpha, y_{n_k} } (\tilde x) 
		= J_{\alpha, y} (\tilde x).
	\end{equation}
	Hence, $x$ is a minimizer of $ J_{\alpha, y}$. 
\end{pr}
\begin{remark}\label{rem:ContinuityOptim}
	In the setting of \Cref{cor:uniqueness variational}, we automatically get weak convergence of the sequence $\{x_{n}\}_{n \in \N}$ to the unique minimizer of $J_{\alpha, y}$.
\end{remark}
The following counterpart to \Cref{th:regularisierung:apriori_char} involves the weak topology, and we cannot take the supremum over $B_\delta(y)$.
Instead, we have to fix a measurement sequence.
\begin{thm}{Convergence}{convergence variational}
	Let the assumptions of \Cref{cor:existence variational} hold and assume $\{ x \in \q D(R): Kx = y\} \neq \emptyset$.
	Moreover, let $\alpha \colon (0,\infty) \to (0,\infty)$ satisfy $\lim_{\delta \to 0 }\delta^p/ \alpha(\delta) = 0$ and $\lim_{\delta \to 0 } \alpha(\delta) = 0$.
	
	For  $\delta_n \to 0$ and $\{y_n\}_{n \in \N}$ with $\|y-y_n\|_{\bb Y} \leq \delta_n$, we obtain with $\alpha_n\coloneqq \alpha\left(\delta_n\right)$ that every weak accumulation point of $\{x_n\}_{n \in \N}$ with $x_n \in \argmin(J_{\alpha_n, y_n})$  is a $R$-minimizing solution for $y$.
	If the $R$-minimizing solution $\hat x$ for $y$ is unique, then $x_n \wkto \hat{x}$.
\end{thm}
\begin{pr}
	Let $\hat x$ be a $R$-minimizing solution for $y$.
	By definition of $x_n$ and $\hat x$ it holds that
	\begin{equation}\label{eq:EstMinNorm}
	J_{\alpha_n, y_n}(x_n ) = \|K x_n -y_n\|_{\bb Y}^{p}+\alpha_n R(x_n) \leq \|K \hat{x} - y_n\|_{\bb Y}^{p}+\alpha_n R (\hat{x}) \leq \delta_n^{p}+\alpha_n R (\hat{x}). 
	\end{equation}
	By assumption, the right-hand side of \eqref{eq:EstMinNorm} converges to 0.
	Hence, we obtain $\|K x_n -y_n\|_{\bb Y} \rightarrow 0$, which implies $\|K x_n -y\|_{\bb Y} \rightarrow 0$.
	Moreover, we conclude from \eqref{eq:EstMinNorm} that
	\begin{equation}\label{eq:BoundR}
		\limsup_{n \to \infty} R \left(x_n\right) \leq 	\limsup_{n \to \infty} \frac{\delta_n^{p}}{\alpha_n} + R (\hat{x}) = R (\hat{x}).
	\end{equation}
	With $\alpha^{+}\coloneqq \max \{\alpha_n: n \in \N\ \}$, we then get
	\begin{equation}
		\limsup_{n \to \infty} \left(\|K x_n -y\|_{\bb Y}^{p}+\alpha^{+} R (x_n)\right) \leq \alpha^{+} R(\hat{x}) <  \alpha^{+} R(\hat{x}) + 1 \eqcolon M<\infty. 
	\end{equation}
	Hence, there exists $n_0 \in \N$ such that $x_n \in \lev_{M}(J_{\alpha^+, y})$ for $n \geq n_0$, and $\{x_{n}\}_{n \in \N}$  has at least one weak accumulation point.
	
	Let $x_{n_k} \wkto x \in \q D(R)$.
	Since $K\in \q K(\bb X,\bb Y)$, it holds $K x_{n_k} \wkto K x$ and therefore $Kx = y$.
	From the lower semi-continuity of $R$, \eqref{eq:BoundR}, and the definition of $\hat{x}$, it follows that
	\begin{equation}
		R(x) \leq \liminf_{k \to \infty} R(x_{n_k}) \leq \limsup_{k \to \infty} R(x_{n_k}) \leq R(\hat{x}).
	\end{equation}
	Hence, $x$ is an $R$-minimizing solution.
	If the $R$-minimizing solution $\hat{x}$ is unique, then the subsequence criterion implies convergence of the whole sequence $\{x_n\}_{n \in \N}$.
\end{pr}

The regularizer $R$ in \eqref{eq:data_fidelity_regularizer} should promote "natural" solutions.
Depending on $R$, certain properties are dominant in the solution.
\begin{example}\label{ex: smoothness}
	Let $\bb X = \R^N$, $\bb Y = \R^M$ and $A \in \R^{M \times N}$.
	We consider a discrete counterpart of the $H^1$ norm (which is a Sobolev norm), which leads to
	\begin{equation}\label{eq:JH1}
		J_\alpha(x) = \frac{1}{2} \|Ax  - y \|_2^2 + \frac{\alpha}{2} \| \nabla x \|^2_2,
	\end{equation}
	where $\nabla$ is the forward difference $(\nabla x)_k = x_{k+1} - x_k$ for $k =1, \ldots, N-1$.
	This is a discrete analog of the derivative.
	By \Cref{ex: quadratic regularization}, the minimizers of \eqref{eq:JH1} satisfy $ (A^T A + \alpha \Delta) x = A^T y$.
	This approach has a strong smoothing effect due to the Laplacian $\Delta = \nabla^T \nabla$.
\end{example}

\begin{figure}
	\begin{center}
		\subfloat[$\|\cdot\|_1$-minimizing solution. \label{fig:ell1_good}]{
			\includegraphics[width= 120pt, height = 120pt]{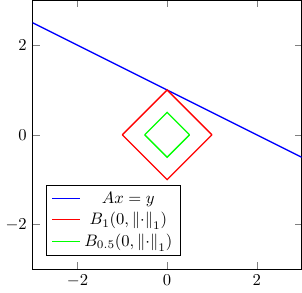}% 
		} \hfill
		\subfloat[$\|\cdot\|_1$-minimizing solutions.
		\label{fig:ell1_bad}]{
			\includegraphics[width= 120pt, height = 120pt]{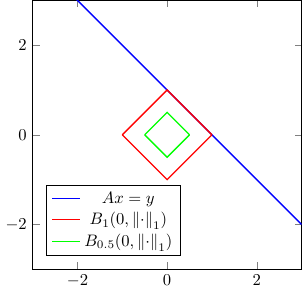}%
		}\hfill
		\subfloat[$\|\cdot\|_{3/4}$-minimizing solutions.
		\label{fig:ellq}]{
			\includegraphics[width= 120pt, height = 120pt]{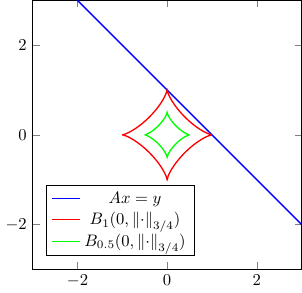}%
		}
		\caption{Sparsity promoting $R$-minimizing solutions.}
	\end{center}
\end{figure}

\begin{example}\label{ex: LASSO}
	Let $\bb X = \R^N$, $\bb Y = \R^M$ and $A \in \R^{M \times N}$. We consider the \emph{least absolute shrinkage and selection}  functional (LASSO) 
	\begin{equation}
		J_\alpha(x) = \frac{1}{2} \|Ax  -y \|_2^2 + \alpha \| x \|_1.
	\end{equation}
	The associated $\|\cdot\|_1$-minimizing solutions $\argmin_{x \in \R^N : Ax = y} \| x \|_1$ are depicted in \Cref{fig:ell1_good}.
	In most cases, they are linear combination of a few extreme point of the $\|\cdot\|_1$-ball, namely the unit vectors $e_k$.
	As $J_\alpha$ balances between data fidelity and the $\|\cdot\|_1$-norm, we expect \emph{sparse} solutions, i.e., mostly containing zero entries.  	
\end{example} 

\begin{example}\label{ex: sparsity promoting}
The idea of promoting zeros can be extended to other representations.
A common choice is $R(x) = \alpha \| B x \|_1$, where $B x$ represents features of $x$ that should be sparse.
For instance, in the case of \emphi{total variation} regularization, we choose $B = \nabla$ as in \Cref{ex: smoothness}.
The sparsity of $\nabla x$ gives that $x$ is piecewise constant with only a few jumps.
Another common choice for $B$ are wavelet transforms.   
\end{example}

In the \Cref{ex: smoothness} and \Cref{ex: sparsity promoting}, we chose regularization of the form $\| B x \|_p$, $p \ge 1$. If $\q N(A) \cap \q N(B) \neq \{ 0\}$, any element $x$ of the intersection admits $\q J_\alpha(x) = 0$, making the global minimum not unique. For instance, consider $Ax = x *_N k$ with $\sum_{n = 1}^N k_n = 0$. Then, for any $c \in \R$, $x = c\1_N = c\,(1,\ldots,1)^\tT $, we have $Ax = 0$. Therefore, taking regularizer $B = \nabla$ with $\q N(B) = \{ \1_N \}$ will not lead to the unique solutions in \Cref{ex: sparsity promoting}. 

We also see in  \Cref{fig:ell1_bad} that the (convex) $\| \cdot \|_1$ regularizer may lead to multiple solutions even with $B = I$. Working with nonconvex $R$ as depicted in \Cref{fig:ellq} can resolve this issue.
Prominent examples of nonconvex $R$ are 
\begin{itemize}
	\item quasinorms $\|x\|_q = (\sum_{k=1}^N |x_k|^q )^{1/q}$ with $0 < q < 1$;
	\item log-sum penalty $\sum_{k=1}^N \log(|x_k| + \varepsilon)$ with $\varepsilon \ge 0$;
	\item (weighted) differences of norms $\|x\|_1 - \beta \|x\|_2$ with $\beta >0$.
\end{itemize}
It is important to keep in mind that uniqueness results no longer apply and nonconvex optimization is in general harder (since local minimizers may exist).
A study on the performance of quasinorm regularizes in presence of noise can be found in \cite{zheng2017does}.

\subsection{Minimization Algorithms}\label{sec:Algs}
Now that we established that minimizers of \eqref{eq:data_fidelity_regularizer} exist, we focus on finding them. 
Since $\dim(\bb X) < \infty$ for applications, we restrict ourselves to $\bb X = \bb R^N$ equipped with an arbitrary inner product.
To emphasize the latter, we keep the notation $\Vert \cdot \Vert_{\bb X}$.
In all our examples, the inner product is the Euclidean one.
Depending on the regularity of \eqref{eq:data_fidelity_regularizer}, we can deploy different minimization techniques, which we derive for a generic $f \colon \R^N \to \overline{\R}$ instead of $J_\alpha$.

\subsubsection{Gradient Descent}

For $f \in C^1(\R^N )$, we can deploy gradient descent, which moves the iterates in the steepest descend direction represented by $-\nabla f$.
This direction \emph{depends} on the inner product in $\bb X$. 

\begin{algorithm}[H]
    \KwIn{$f \in C^1(\R^N )$ with $\argmin(f) \ne \emptyset$, initialization $x_0 \in \R^N $, step size $\nu > 0$.}    
    \For{\color{black} $k = 0, 1, \ldots$}{
        $x_{k + 1} = x_k-\nu \nabla f(x_k)$\;
    }
    \caption{Gradient Descent}
    \label{algo:gradient_descent}
\end{algorithm}

Any fixed point $x_*$ found by \Cref{algo:gradient_descent} satisfies
\begin{equation}
	x_* = x_* - \nu \nabla f(x_*).
\end{equation}
This is equivalent to $x_*$ being a critical point, namely $\nabla f(x_*) = 0$.
To establish the convergence of \Cref{algo:gradient_descent}, we require the fundamental \emph{descent lemma}.

\begin{lemma}[Descent Lemma]\label{l: descent lemma}
	Let $f \in C^1(\R^N )$ with $L$-Lipschitz continuous gradient and $x, h \in \R^N$.
	Then, it holds
	\begin{equation}\label{eq:DesLem}
	f(x + h) \le f(x) + \inner{\nabla f(x)}{h} + \frac{L}{2} \|h\|_{\bb X}^2.
	\end{equation} 
\end{lemma}

\begin{pr}
	For $x, h \in \R^N$, a first-order Taylor expansion gives
	\begin{align}
		f(x + h) &= f(x) + \int_0^1 \inner{\nabla f(x + th)}{h} \diff t \notag\\
		&= f(x) + \inner{\nabla f(x)}{h} + \int_0^1 \inner{\nabla f(x + th) - \nabla f(x)}{h} \diff t.
	\end{align}
	Applying the Cauchy-Schwarz inequality and then the Lipschitz continuity of $\nabla f$ yields
	\begin{align}
		f(x + h) & \le f(x) + \inner{\nabla f(x)}{h} +
		\int_0^1 \| \nabla f(x + th) - \nabla f(x) \|_{\bb X} \| h \|_{\bb X} \diff t  \notag\\
		& \le f(x) + \inner{\nabla f(x)}{h} + L \int_0^1 t \| h \|_{\bb X}^2
		\diff t  = f(x) + \inner{\nabla f(x)}{h} + \frac{L}{2} \|h\|_{\bb X}^2.
	\end{align}
	This concludes the proof.
\end{pr}

Now, we show that \Cref{algo:gradient_descent} returns almost stationary points, namely $x$ for which $ \| \nabla f(x) \|_{\bb X}$ is small.
Here, we  call a saddle point $x \in \R^N$ a \emphi{strict saddle} if for some $z \in \R^N$ it holds $\inner{\nabla^2 f(x) z}{z}_{\bb X} < 0$.
Put differently, this means that $z$ is a descent direction.

\begin{thm}{Convergence of Gradient Descent}{grad descent nonconvex}
	Let $f \in C^1(\R^N )$ be bounded from below with $L$-Lipschitz continuous gradient.
	If $0<\nu \leq 1/L$, then $\{x_k\}_{k \in \N}$ generated by \Cref{algo:gradient_descent} satisfies for any $K \ge 0$ that
	\begin{equation}\label{eq:OptimalityGD}
		\min_{k = 0, \ldots, K-1} \| \nabla f(x_k) \|_{\bb X}^2 
		\le 2\frac{f(x_0) - \inf_{x \in \R^N}f(x)}{\nu K}.
	\end{equation}
	If $\lim_{k \to \infty} x_k = x_*$, then $\nabla f(x_*) = 0$. 
	Moreover, if $f \in C^2(\R^N )$, then $\{x_k\}_{k \in \N}$ generated with $\nu < 1/L$ does not converge to a strict saddle for almost every $x_0 \in \R^N$. 
\end{thm} 

\begin{pr}
	We apply \Cref{l: descent lemma} with $x = x_{k}$ and $h = - \nu \nabla f(x_k)$ to obtain
	\begin{equation}\label{eq: descent key}
	f(x_{k+1}) - f(x_k) \le 
	%- \nu \|\nabla f(x_k)\|_2^2 + \frac{L \nu^2}{2} \|\nabla f(x_k)\|_2^2 = 
	- \nu \Bigl(1 - \frac{L \nu}{2}\Bigr) \|\nabla f(x_k)\|_{\bb X}^2
	\le - \frac{\nu}{2} \|\nabla f(x_k)\|_{\bb X}^2,
	\end{equation}
	where we used that $\nu \le 1/L$ for the last inequality. 
	Consequently, we have
	\begin{equation}
		\sum_{k=0}^{K-1} \|\nabla f(x_k)\|_{\bb X}^2
		\le \frac{2}{\nu} \sum_{k = 0}^{K-1}[f(x_k) - f(x_{k+1})]
		\le 2\frac{f(x_0) - \inf_{x \in \R^N}f(x)}{\nu}.
	\end{equation}
	As $f$ is bounded from below, we get $\sum_{k=0}^{\infty} \|\nabla f(x_k)\|_{\bb X}^2 < \infty$ and thus $\lim_{k \to \infty} \|\nabla f(x_k)\|_{\bb X} = 0$.
	If $\lim_{k\to \infty} x_k = x_*$, this readily implies that $ \|\nabla f(x_*)\|_{\bb X} = \lim_{k \to \infty} \|\nabla f(x_k)\|_{\bb X} = 0$.
	Finally, we use that the minimum is less or equal than the average, so that
	\begin{equation}
		\min_{k = 0, \ldots, K-1} \| \nabla f(x_k) \|_{\bb X}^2 
		\le \frac{1}{K} \sum_{k=0}^{K-1} \|\nabla f(x_k)\|_{\bb X}^2
		\le 2\frac{f(x_0) - \inf_{x \in \bb X}f(x)}{\nu K}.
	\end{equation}
	The last part of the theorem is given in \cite[Thm.\ 4.1]{LeeSimJor2016}.
\end{pr}

The convergence rate $1/K$ in \eqref{eq:OptimalityGD} is optimal for general nonconvex $f$.
It is often referred to as a \emphi{sublinear} rate, since in log scale the decay is slower than a straight line.

\begin{remark}[Convergence of Iterates]\label{rem: kl}
	In general, we cannot expect convergence of $\{x_k\}_{k \in \N}$ generated by \Cref{algo:gradient_descent} without further assumptions on $f$ even if $\{x_k\}_{k \in \N}$ is bounded.
	A well-known counterexample is the Mexican hat $f \colon \R^2 \to \R$ with
	\begin{equation}\label{eq:MexHat}
		f\bigl(r \cos(\theta),r \sin(\theta)\bigr) = \begin{cases}
			\exp^{-\frac{1}{1 - r^2}}\Bigl(1 - \frac{4r^4}{4r^4+ (1 - r^2)^4} \sin(\theta-\frac{1}{1-r^2})\Bigr) & \text{ if } r < 1,\\
				0 & \text{ if } r \geq 1,\\
		\end{cases}
	\end{equation}
	\begin{figure}
		\subfloat{
		\includegraphics[width = 0.45\textwidth]{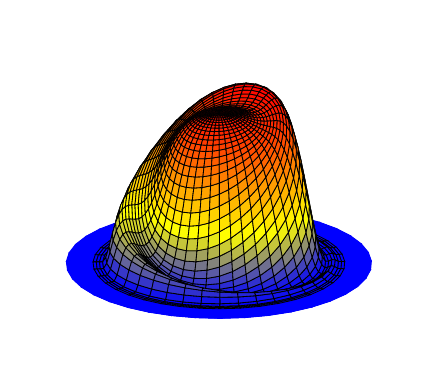}
		}
		\hfill
		\subfloat{
			\includegraphics[width = 0.45\textwidth]{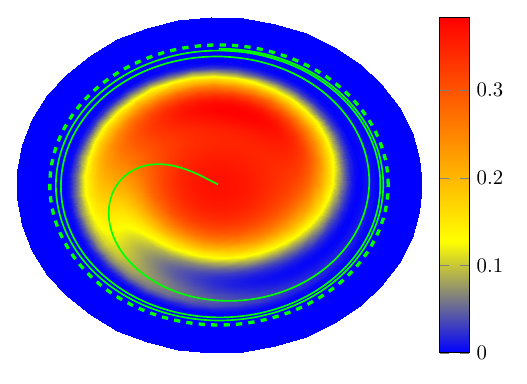}
		}
		\caption{\textbf{Left:} Plot of $f$; \textbf{Right:} Infinite length trajectory of a gradient flow $\dot x(t)=-\nabla f(x(t))$ in green with the set of limit points $r = 1$ in dashed green.}
		\label{fig: mexican hat}
		% See absil paper
	\end{figure}
	for which we can create a gradient descent trajectory $\{x_k\}_{k \in \N}$ with infinite length, namely $\sum_{k=1}^\infty \Vert x_k - x_{k-1} \Vert_{\bb X} = \infty$, see \Cref{ex:MexHat}.
	The problem with \eqref{eq:MexHat} is that it lacks sharpness around its minima.
	This has led to the formulation of the Kurdyka--{\L}ojasiewicz (KL) condition, under which bounded trajectories $\{x_k\}_{k \in \N}$ converge \cite{AbsMahAnd2005,AttBolSva2013}. 
	We discuss the KL condition and how to verify it in Section \ref{sec: KL}.
	In particular, we will see that most $f$ in applications satisfy the KL condition.
\end{remark}

\begin{remark}[Convex Case]\label{rem:ObjectiveConvex}
	If $f$ is convex and $x_* \in \argmin(f)$, we can strengthen \Cref{th:grad descent nonconvex}.
	Using \eqref{eq: strong convexity} with $\mu=0$ and \eqref{eq: descent key}, we can show (see \Cref{ex:ConvObj}) that
	\begin{equation}\label{eq:EstConvex}
		f(x_{k+1})-f(x_*) \leq \frac{1}{2\nu}\bigl(\Vert x_{k} - x_*\Vert_{\bb X}^2 - \Vert x_{k+1} - x_*\Vert_{\bb X}^2 \bigr).
	\end{equation}
	Since $\{f(x_{k})-f(x_*)\}_{k \in \N}$ is decreasing, summing up over the first $K$ iterates leads to
	\begin{equation}\label{eq:ConvexObjectiveRate}
		f(x_{K})-f(x_*) \leq \frac{1}{K}\sum_{k=1}^K f(x_{k})-f(x_*) \leq \frac{\Vert x_{0} - x_*\Vert_{\bb X}^2}{2\nu K}.
	\end{equation}
	Hence, we get that also the objective values $\{f(x_{k})\}_{k\in \N}$ converge with the rate $1/K$.
\end{remark}

For strongly convex $f$, an improved \emph{linear} convergence rate applies, and we obtain convergences of the iterates without further assumptions.

\begin{thm}{Gradient Descent for Strongly Convex $f$}{grad descent linear}
Let $f \in C^1(\R^N)$ be $\mu$-strongly convex with $L$-Lipschitz continuous gradient.
If $0<\nu \le 1/L$, the sequence $\{x_k\}_{k \in \N}$ generated by \Cref{algo:gradient_descent} converges to the global minimizer $x_*$  of $f$. Furthermore, for $K \ge 0$ we have
\begin{equation}\label{eq:ConvSpeedStronglyConvex}
	\| x_{K} - x_* \|_{\bb X}^2 \le \left(1 - \nu \mu \right)^K \| x_{0} - x_* \|_{\bb X}^2.
\end{equation}
\end{thm}

\begin{pr}
First, we perform the expansion
\begin{align}
	\| x_{k+1} - x_*\|_{\bb X}^2 &= \| x_{k} - x_* - \nu \nabla f(x_k) \|_{\bb X}^2\notag\\ 
	&= \| x_{k} - x_* \|_{\bb X}^2 - 2 \nu \inner{\nabla f(x_k)}{x_{k} - x_*} + \nu^2 \|\nabla f(x_k)\|_{\bb X}^2.\label{eq:dist_opt}
\end{align}
Due to the strong convexity property \eqref{eq: strong convexity}, we have
\begin{equation}\label{eq:EstInnerProd}
	- \inner{\nabla f(x_k)}{x_{k} - x_*} \le -(f(x_k) - f(x_*)) - \frac{\mu}{2} \| x_k - x_*\|_{\bb X}^2.
\end{equation}
Plugging \eqref{eq:EstInnerProd} into \eqref{eq:dist_opt}, we obtain
\begin{equation}
	\| x_{k+1} - x_*\|_{\bb X}^2 \le (1 - \nu \mu) \| x_{k} - x_* \|_{\bb X}^2
	- 2 \nu (f(x_k) - f(x_*)) + \nu^2 \|\nabla f(x_k)\|_{\bb X}^2.
\end{equation}
Since the assumptions of \Cref{th:grad descent nonconvex} hold, we can use \eqref{eq: descent key} to get
\begin{equation}
	\nu^2 \|\nabla f(x_k)\|_{\bb X}^2 
	\le 2 \nu (f(x_k) - f(x_{k+1}))
	\le 2 \nu (f(x_k) - f(x_*)),
\end{equation}
where the last inequality is due to $x_* \in \argmin(f)$.
Consequently,
\begin{equation}
	\| x_{k+1} - x_*\|_{\bb X}^2 \le (1 - \nu \mu) \| x_{k} - x_* \|_{\bb X}^2
\end{equation}
with $0 \le 1 - \nu \mu < 1$ due to $\mu \leq L$.
Now, a recursion gives the desired result.
\end{pr}

\begin{remark}
	The rate in \Cref{th:grad descent linear} is not optimal in the sense that we can get rid of the squares on both sides of \eqref{eq:ConvSpeedStronglyConvex}.
	We will discuss this in \Cref{ex: stongly convex via Banach}.
\end{remark}

Without strong convexity, we cannot improve the sublinear rate.
Instead, Nesterov \cite{Ne83} introduced an acceleration technique, see \Cref{algo:accelerated_gradient_descent}, which improves the rate for convex $f$.
Below, we present a second version of the method \cite{Ne04,Ts08} that is slightly easier to handle.

\begin{algorithm}[H]
	\KwIn{$f \in C^1(\R^N)$ with $\argmin(f) \ne \emptyset$, initialization $x_0 \in \R^N$,  step size $\nu > 0$.}    
	Set $z_{0} = x_0$, and $\theta_k = \frac{2}{k+2}$ for $k \ge 0$\;
	\For{\color{black} $k = 0, 1, \ldots$}{
		$u_k = (1 - \theta_k) x_k + \theta_k z_k $\;
		$x_{k + 1} = u_k - \nu \nabla f(u_k)$\;
		$z_{k+1} = x_k + \frac{1}{\theta_k}(x_{k+1} - x_k)$\;
	}
	\caption{Accelerated Gradient Descent}
	\label{algo:accelerated_gradient_descent}
\end{algorithm}

In contrast to \eqref{eq:ConvexObjectiveRate}, we now obtain the rate $1/K^2$ for the objective values $f(x_K)$, as shown for a more general setting in \cite{beck2009fast,Ts08}. 
\begin{thm}{Convergence of accelerated gradient descent}{gradient accelerated}
Let $f \in C^1(\R^N)$ be convex, coercive with $L$-Lipschitz continuous gradient. 
Then, with $0 < \nu \le 1/L$, the sequence $\{x_k\}_{k \in \N}$ generated by \Cref{algo:accelerated_gradient_descent} converges to a point in $\q X \coloneqq \argmin (f)$.
Furthermore, for $K \ge 0$ we have
\begin{equation}
	f(x_K) - \inf_{x \in \R^N}f(x) \le \frac{2 L \min_{x_* \in \q X} \| x_0 - x_* \|_{\bb X}^2}{(K+1)^2}.
\end{equation}
\end{thm}
For nonconvex $f$, acceleration often improves the performance empirically.
However, the theoretically proven rates remain $1/K$, see \cite{DanDvuGas2022} for an overview.
We conclude our discussion of gradient descent with two implementation aspects.

\begin{remark}[Line Search]\label{rem:LineSearch}
	Usually, we do not know the Lipschitz constant $L$, which is required in \Cref{algo:gradient_descent} and \Cref{algo:accelerated_gradient_descent} to select the step size $\nu$.
	By revisiting the proof of \Cref{th:grad descent nonconvex}, we see that $L$ is only required to ensure the sufficient decrease in \eqref{eq: descent key}. 
	Instead, we can select step sizes $\nu_k$ adaptively such that \eqref{eq: descent key} holds.
	In a backtracking line search, we initialize $\nu_k = \beta \nu_{k-1}$, $\beta > 1$, and keep dividing by some constant until it holds
	\begin{align}
		f(x_{k+1}) &\leq  f(x_k) -  \frac{\nu_k}{2} \|\nabla f(x_k)\|_{\bb X}^2.\label{eq:LineSearch}
	\end{align}
	This requires to compute several proposals $x_{k+1}$ and the corresponding $f(x_{k+1})$.
	Through this, the algorithm adapts to the local behavior of $f$.
	If $\nabla f$ is Lipschitz continuous, a feasible $\nu_k$ exists, and we perform the update in \Cref{algo:gradient_descent} or \Cref{algo:accelerated_gradient_descent} with the found $\nu_k$.
	In order for the proof of \Cref{th:grad descent nonconvex} to remain true, we need to enforce a lower bound, i.e., $\inf_{k \in \N} \nu_k>0$.
	An alternative to setting $\nu_k = \beta\nu_{k-1}$ are \textsc{Barzilei-Borwein} step sizes \cite{BB88} as initialization, which are related to Quasi-\textsc{Newton} approaches, see \cite{FZZ08} for an overview.
	We refer to \cite{GS2011,Ne07} for more details.
\end{remark}

\begin{remark}[Convergence Criteria]
	In practice,  we do not run \Cref{algo:gradient_descent} and \Cref{algo:accelerated_gradient_descent} until we reach the maximal number of steps. 
	Instead, we terminate early if $\Vert x_{k+1} - x_k\Vert_{\bb X} / \Vert x_k\Vert_{\bb X} \leq \eps$ for some $\eps >0$. 
	For single precision, $\eps = 10^{-4}$ is reasonable.
\end{remark}

\subsubsection{Exercises}

\begin{exerciseB}\label{ex:ConvObj}
	Prove the estimate \eqref{eq:EstConvex} in \Cref{rem:ObjectiveConvex}.
\end{exerciseB}

\begin{Solution}
	Using \Cref{l: descent lemma} and update rule $x_{k+1} = x_k - \nu \nabla f(x_k)$, we obtain
	\begin{align*}
		f(x_{k+1}) - f(x_k) 
		& \le - \tfrac{\nu}{2} \| \nabla f (x_k)\|^2
		= - \tfrac{1}{2\nu} \| x_{k+1} - x_k\|^2
		= - \tfrac{1}{2\nu} \| x_{k+1} - x_* + x_* - x_k\|^2.
	\end{align*}
	Opening the square norm gives
	\begin{align*}
		f(x_{k+1}) - f(x_k) 
		& \le - \tfrac{1}{2\nu} \left[ \| x_{k+1} - x_* \|^2 + \| x_* - x_k\|^2 + 2 \inner{x_{k+1} - x_* \pm x_k}{x_* - x_k} \right] \\
		& = - \tfrac{1}{2\nu} \left[ \| x_{k+1} - x_* \|^2 + \| x_* - x_k\|^2 - 2 \| x_* - x_k\|^2 + 2 \inner{x_{k+1} - x_k}{x_* - x_k} \right] \\
		& = - \tfrac{1}{2\nu} \left[ \| x_{k+1} - x_* \|^2 - \| x_* - x_k\|^2 \right] + \inner{\nabla f(x_k)}{x_* - x_k}. 
	\end{align*}
	By convexity of $f$ and \eqref{eq: strong convexity}, we get
	\[
	f(x_*) \ge f(x_k) + \inner{\nabla f(x_k)}{x_* - x_k},
	\]
	which leads to
	\begin{align*}
	f(x_{k+1}) - f(x_k) 
	& \le - \tfrac{1}{2\nu} \left[ \| x_{k+1} - x_* \|^2 + \| x_* - x_k\|^2 + 2 \inner{x_{k+1} - x_* \pm x_k}{x_* - x_k} \right] \\
	& = - \tfrac{1}{2\nu} \left[ \| x_{k+1} - x_* \|^2 + \| x_* - x_k\|^2 - 2 \| x_* - x_k\|^2 + 2 \inner{x_{k+1} - x_k}{x_* - x_k} \right] \\
	& = - \tfrac{1}{2\nu} \left[ \| x_{k+1} - x_* \|^2 - \| x_* - x_k\|^2 \right] + f(x_*) - f(x_k)
	\end{align*}
	Simplifying both sides gives the desired bound.
	
	\textbf{Remark:} Convexity of $f$ is a sufficient assumption. We only use convexity along the line connecting $x_k$ and $x_*$, which can be found in the literature under the name restricted convexity.
\end{Solution}

\begin{exerciseB}[Gradient Descent and Conditioning]\label{ex: gd and condition number}
	Let $\bb X = \bb R^N$, $\bb Y = \bb R^M$, $y \in \bb Y$ and $A \in \bb R^{M \times N}$. Consider $f(x) = \tfrac{1}{2} \| A x - y \|^2_{\bb Y}$. 
	\begin{enumerate}
		\item[(i)] Compute $\nabla f$ and $\nabla^2 f$.
		Use these to reinterpret the normal equation (\Cref{thm:inverse:normalen}).
		\item[(ii)] Implement a function that performs $K$ iterations of gradient descent (Algorithm 1) with initialization $x_0 \in \bb R^N$.
		It should return the iterate $x_K$, the (relative) measurement errors $\{ f(x_k) / f(0) \}_{k = 0}^K$, and the relative errors $\{\|x_k - x_*\| / \| x_* \| \}_{k =0}^K$ for a given $x_*$.
		\item[(iii)] Let $N$ be even.
		Consider the deconvolution problem with the convolution kernel $g^\alpha \in \R^N$, $0 < \alpha < 1$, which is implicitly given via
		\begin{equation}
			(F_N g^{\alpha})_k = 
			\begin{cases}
				1 - (1 - \alpha)\frac{k}{N/2 - 1}, & k = 0, \ldots, N/2 - 1 \\
				\alpha + (1 - \alpha)\frac{k - N/2}{N/2 - 1} & k = N/2, \ldots, N -1.
			\end{cases}
		\end{equation}
		  
		For the convolution matrix $C(g^\alpha) \in \R^{N \times N}$ (see \Cref{ex: circulant matrix}), find the spectral norm $\|C(g^\alpha)\|_2$ and the condition number $\kappa(C(g^\alpha)) = \frac{\sigma_\mathrm{max}(C(g^\alpha))}{\sigma_\mathrm{min}(C(g^\alpha))}$.  
		
		\item[(iv)]  Let $N = 200$ and $x_* \in \bb R^N$ with entries
		\begin{equation}\label{eq:signalExample}
			(x_*)_k =  \sin\bigl(\tfrac{2 \pi k}{N}\bigr) + 0.5 \sin\bigl(\tfrac{6 \pi k}{N}\bigr)  - 0.2 \sin\bigl(\tfrac{40 \pi k}{N}\bigr) + 0.3 \sin\bigl(\tfrac{10 \pi k}{N}\bigr),
		\end{equation} 
		$k = 0, \ldots, N-1$.
		As associated (noisy) data, we consider $y = C(g^\alpha) x_* + 0.01\q N(0,I_N)$, where $\q N(0,I_N)$ is a standard Gaussian random variable.
		For $\alpha = 10^{- k/2}$, $k = 1, \ldots, 5 $, perform $K =10000$ iterations of gradient descent with initialization $x_0 \sim \q N(0,I_N)$.
		Plot the resulting errors (in $\log_{10}$ scale). What do you observe?    
		
		\item[(v)] Implement accelerated gradient descent (\Cref{algo:accelerated_gradient_descent}) and perform the reconstructions in (iv) again.
		Does the convergence behavior improve? 
	\end{enumerate} 
\end{exerciseB}

\begin{Solution}
	(i): We compute the gradient by definition of Gateaux directional derivatives. Let $x, v \in \R^N$.
	\begin{align}
		\lim_{t \searrow \infty} \frac{f(x + tv) - f(x)}{t}
		& = \lim_{t \searrow \infty} \frac{ \|A x - y + t Av\|^2 - \| Ax - y \|^2}{2t} \\
		& = \lim_{t \searrow \infty} \frac{ \|A x - y\|^2 + 2 \inner{Ax -y}{t Av} + \|t Av\|^2 - \| Ax - y \|}{2t} \\
		& = \lim_{t \searrow \infty} \inner{A^\tT (Ax -y) }{v} + \frac{t}{2}\|Av\|^2 
		= \inner{A^\tT (Ax -y) }{v}
	\end{align}
	This is a linear functional in $v$ and, thus, Gateaux gradient is given by $\nabla f(x) = A^\tT (Ax -y)$.
	Analogously, we get
	\[
	\lim_{t \searrow \infty} \frac{\nabla f(x + tv) - \nabla f(x)}{t} = \lim_{t \searrow \infty} \frac{t A^\tT A v}{t} = A^\tT A v
	\]
	and the Gateuax Hessian is given by $\nabla^2 f(x) = A^\tT A$. 
	
	The formula for the gradient and the Fermat's rule give
	\[
	A^\tT (Ax -y) = 0, \text{ which is equivalent to } A^\tT Ax = A^\tT y.
	\] 
	
	(ii) - (v): See Jupyter notebook \url{https://github.com/Sniper2k/inverse_problems_ai_exercises/blob/main/Ex4_3_11.ipynb}.
\end{Solution}

\begin{exerciseB}[Tikhonov Regularization]\label{eq:ImplementTikhonov}
	Using the results of \Cref{ex: quadratic regularization}, implement Tikhonov regularization with $\alpha=5 \cdot 10^{-4}$ in two ways: Based on gradient descent with $K = 10^5$ iterations and convergence threshold $\eps = 10^{-8}$, and using the SVD.
	For a Gaussian matrix $A \in \R^{M \times N}$ with $A_{k,j} \sim M^{-1/2} \q N(0,1)$ and data as in \Cref{ex: gd and condition number}, we consider the  inverse problem $Ax = y$.
	Perform the reconstruction for $N =2^n$ and $M =3 \cdot 2^n$ with $n =2,\ldots,11$.
	Compare the runtime of both approaches.
	Does gradient descent admit a nice convergence behavior? 
\end{exerciseB}

\begin{Solution}
	See Jupyter Notebook \url{https://github.com/Sniper2k/inverse_problems_ai_exercises/blob/main/Ex4_3_12.ipynb}.
\end{Solution}

\begin{exerciseB}\label{ex:Lipschitz Hessian}
	Let $f \in C^2(\R^N)$.
	Then,  the Lipschitz continuity of $\nabla f$ is equivalent to 
	\begin{equation}\label{eq:EquivLip}
		| \inner{\nabla^2 f(x) h}{h} | \le L \|h\|_{\bb X}^2 \quad \text{for all } x,h \in \R^N. 
	\end{equation}
	Show that in \Cref{th:grad descent nonconvex} 
	we can relax \eqref{eq:EquivLip} to the upper bound 
	\begin{equation}
		\inner{\nabla^2 f(x) h}{h} \le \tilde L \|h\|_{\bb X}^2, \quad \text{for all } x,h \in \R^N 
	\end{equation}
	for some $\tilde L>0$. Further, argue that if $f$ is convex, both inequalities are equivalent.
\end{exerciseB}
\begin{Solution}
	The proof of \Cref{th:grad descent nonconvex} is a consequence of the descent lemma (\Cref{l: descent lemma}), where we use that $ f\in C^1(\R^N) $ has $L$-Lipschitz continuous gradient. Thus, to relax the Lipschitz assumption in \Cref{th:grad descent nonconvex}, we need to show that the descent Lemma holds with the  weaker condition
	\[
	\inner{\nabla^2 f(x)h}{h} \leq \tilde L \|h\|_\bb X^2.
	\]
	By second-order Taylor theorem, we get
	\begin{align*}
		f(x+h) &= f(x) + \inner{\nabla f(x)}{h} + \int_0^1 (1- t)\inner{\nabla^2 f(x + th) h}{h} \, dt \\
		&\leq f(x) + \inner{\nabla f(x)}{h} + \int_0^1 (1-t) \tilde L \|h\|^2 \, dt \\
		&= f(x) + \inner{\nabla f(x)}{h} + \frac{1}{2} \tilde L \|h\|^2.
	\end{align*}
	Thus, we get the descent lemma under our weaker bound on the Hessian and \Cref{th:grad descent nonconvex} continues to hold. 
	
	Note that since $\nabla^2 f$ is continuous, $\nabla f$ is Lipschitz, but it could be $\tilde L < L$.
	If $f$ is convex, then \eqref{eq: convexity grad} gives $\inner{\nabla^2 f(x)z}{z} \ge 0 \ge -\tilde  L \|z\|^2$ for all $z \in \R^N$, so that $L = \tilde L$ and both conditions coincide. 
\end{Solution}

\begin{exerciseB}[Phase retrieval]
	Let $\bb X = \R^N$, $\bb Y = \R^M$, $A \in \R^{M \times N}$ and consider the phase retrieval problem $|Ax|^2 = y$ from \Cref{ex: Fourier phase retrieval}.
	For the reconstruction, we consider the (smoothed) functional 
	\begin{equation}\label{eq:RecFunctionalPhase}
		f(x) = \frac{1}{2}\sum_{j = 1}^M | \sqrt{|(Ax)_k|^2 + \varepsilon} - \sqrt{y + \varepsilon} |^2
	\end{equation}
	with $\varepsilon>0$.
	\begin{itemize} 
		\item[(i)] Show that $x$ and $-x$ yield the same measurements and that $f(x) = f(-x)$. This is known as the global phase ambiguity in phase retrieval.
		\item[(ii)] Prove that $\nabla f(x)$ and $\nabla^2 f(x)$ are given by
		\begin{align}
			\nabla f(x) &= A^\tT \diag\biggl( 1 - \frac{\sqrt{y_i +\varepsilon}}{(| (Ax)_i|^2 + \varepsilon)^{1/2}}\biggr) A x,\\
			\nabla^2 f(x) &= A^\tT \diag\biggl(1 - \frac{\varepsilon \sqrt{y_i + \varepsilon} }{(|(Ax)_i|^2 + \varepsilon )^{3/2}}\biggr) A.
		\end{align}
		\item[(iii)] Assume that $A$ is injective and that $y_k \neq 0$ for all $k = 1,\ldots, M$.
		Show that $0 \in \R^N$ is a local maximum of $f$.
		\item[(iv)] Show that the Lipschitz constant of $f$ is at least $C/ \sqrt{\varepsilon}$ for some $C >0$.
	    \item[(v)] Show that \Cref{algo:gradient_descent} with step sizes $\nu \le 1 / \| A \|^2_2$ converges to a critical point. 
	\end{itemize}
\end{exerciseB}
\begin{Solution}
	(i): We have $|A(-x)|^2 = |Ax|^2$ and the same applies to $f$. 
	
	(ii): We compute the partial derivatives of $f$, 
	\begin{align*}
		&\frac{\partial f}{\partial x_k} (x) =  \frac{1}{2} \sum_{i =1}^M \frac{\partial}{\partial x_k}  \left( \sqrt{|(Ax)_i|^2 + \varepsilon} - \sqrt{y_i+\varepsilon} \right)^2 \\
		& = \frac{1}{2} \sum_{i=1}^M 2 \left( \sqrt{|(Ax)_i|^2 + \varepsilon} - \sqrt{y_i+\varepsilon} \right) \frac{1}{2 \sqrt{|(Ax)_i|^2 + \varepsilon}} 2 \left(Ax\right)_i A_{i,k} \\
		& = \sum_{i=1}^M A_{k,i}^\tT \left( 1 - \frac{\sqrt{y_i + \varepsilon}}{\sqrt{\left| \left(Ax\right)_i\right|^2 + \varepsilon}}\right) \left(Ax\right)_i  \\
	\end{align*}
	In the matrix form, the gradient of $f$ is
	\begin{equation}\label{eq: phase retrieval gradient}
	\nabla f(x) = A^\tT \diag\left( 1 - \frac{\sqrt{y_i +\varepsilon}}{\sqrt{\left| \left(Ax\right)_i\right|^2 + \varepsilon}}\right) A x.
	\end{equation}
	We continue computing the second-order partial derivates to obtain the Hessian
	\begin{align*}
		&\frac{\partial f}{\partial x_j \partial x_k}(x)
		= \sum_{i=1}^M - \frac{ \sqrt{y_i +\varepsilon }}{ (-2) \left( \left| \left( Ax \right)_i \right|^2 + \varepsilon \right)^{3/2}} 2 \left(Ax\right)_i \cdot \left(Ax\right)_i A_{i,j} A_{k,i}^\tT 
		+ \left(1 - \frac{\sqrt{y_i + \varepsilon } }{ \sqrt{ \left|\left(Ax\right)_i \right|^2 +\varepsilon} } \right) A_{i,j} A_{k,i}^\tT \\
		& = \sum_{i=1}^M  A_{k,i}^\tT   \left[ 1 - \frac{\sqrt{y_i+\varepsilon} \left(\left|\left(Ax\right)_i\right|^2 + \varepsilon\right) }{ \left( \left|\left(Ax\right)_i\right|^2 + \varepsilon \right)^{3/2}} + \frac{\sqrt{y_i +\varepsilon} \left|\left(Ax\right)_i\right|^2}{ \left(\left|\left(Ax\right)_i\right|^2 +\varepsilon\right)^{3/2}}\right] A_{i,j} \\
		& = \sum_{i=1}^M A_{k,i}^\tT \left( 1 - \frac{\varepsilon \sqrt{y_i + \varepsilon }  }{\left(\left|\left(Ax\right)_i\right|^2 +\varepsilon\right)^{3/2}}\right) A_{i,j}.
	\end{align*}
	In the matrix form, the Hessian of $f$ is 
	\begin{equation}\label{eq: phase retrieval hessian}
	\nabla^2 f(x) = A^\tT \diag\left(1 - \frac{\varepsilon \sqrt{y_i + \varepsilon} }{\left(\left|\left(Ax\right)_i\right|^2 + \varepsilon \right)^{3/2}}\right) A.
	\end{equation}
	
	(iii): From \eqref{eq: phase retrieval gradient}, we get 
	\[
	\nabla f(0) = A^\tT \diag\left( 1 - \frac{\sqrt{y_i +\varepsilon}}{\sqrt{\left| \left(Ax\right)_i\right|^2 + \varepsilon}}\right) A 0 = 0
	\]
	and from \eqref{eq: phase retrieval hessian}, we obtain for all $h \in \R^N$ that
	\begin{align}
	h^\tT \nabla^2 f(0) h 
	& = \sum_{i=1}^M \left[1 - \frac{\varepsilon \sqrt{y_i + \varepsilon}}{\left(\varepsilon\right)^{3/2}} \right] (A h)_i^2
%	= 1 - \frac{\varepsilon \sqrt{y_k + \varepsilon}}{\varepsilon^{3/2}} \|h\|^2
	= - \sum_{i=1}^M \frac{\sqrt{y_i + \varepsilon} - \sqrt{\varepsilon}}{\sqrt{\varepsilon}} (A h)_i^2 \\
	& \le - \min_{i = 1, \ldots, M} \frac{\sqrt{y_i + \varepsilon} - \sqrt{\varepsilon}}{\sqrt{\varepsilon}} \| A h\|^2
	\le - \min_{i = 1, \ldots, M} \frac{\sqrt{y_i + \varepsilon} - \sqrt{\varepsilon}}{\sqrt{\varepsilon}} \sigma_N^2(A) \|h\|^2.
	\end{align}
	Consequently, by \Cref{th:second_order_sufficient} is a local maximum. 
	
	(iv): From the above bound, we get that $L \ge \min_{k = 1, \ldots, M} \frac{\sqrt{y_k + \varepsilon} - \sqrt{\varepsilon}}{\sqrt{\varepsilon}} \sigma_N^2(A)$.
	
	(v): To show the gradient descent algorithm converges to a critical point we apply \Cref{th:grad descent nonconvex} with the weaker assumption from \Cref{ex:Lipschitz Hessian}. 
	The function $f$ is a non-negative as a sum of squares, hence it is bounded from below by zero.
	We show the Hessian bound from the previous exercise. For all $h \in \R^N$ we have
	\begin{align*}
		&h^\tT \nabla^2 f(x) h = h^\tT A^\tT \diag\left(1 - \frac{ \varepsilon \sqrt{y_i + \varepsilon} }{\left(\left|\left(Ax\right)_i\right|^2 + \varepsilon \right)^{3/2}}\right) Ah \\
		& = \sum_{i=1}^N \left[1 - \underbrace{\frac{\varepsilon \sqrt{y_i + \varepsilon}}{\left( | (Ax)_i|^2 + \varepsilon\right)^{3/2}}}_{\ge 0}\right] (A h)_i^2
		\leq \sum_{i=1}^N (Ah)_i^2 
		= \|Ah\|^2_2 
		\leq \|A\|^2_{2} \|h\|^2_2.
	\end{align*}
	Thus, \Cref{ex:Lipschitz Hessian} gives the desired step size. 
\end{Solution}

\begin{exerciseB}
	Let $N=50$ and $M \in \{N, 2N, 4N, 6N\}$.
	For a Gaussian matrix $A \in \R^{M \times N}$ with $A_{k,j} \sim M^{-1/2} \q N(0,1)$ and the signal $x_* \in \R^{N}$ from \eqref{eq:signalExample}, compute the associated measurements $y_j = \max\{|(Ax_*)_j|^2 + 0.01\q N(0,1),0\}$.
	Then, consider the inverse problem $|Ax|^2  = y$.
	
	\begin{enumerate}
		\item[(i)] Implement $K$ steps of gradient descent with step size $\nu$ and initial guess $x_0$ for the objective \eqref{eq:RecFunctionalPhase} with $\varepsilon = 10^{-4}$.
		Your function should return $x_K$, the relative measurement errors $\{ f(x_k) / f(0) \}_{k=0}^K$, and the  relative errors $\{ \frac{\min \{ \| x_k - x_* \|_2, \| x_k + x_* \|_2\}}{\| x_*\|_2} \}_{k=0}^K$.
		
		\item[(ii)] Run your algorithm for $K = 1000$ with initialization $x_0 \sim \q N(0,I_N)$. Plot the reconstructions. What do you observe? Plot the error curves (in $\log_{10}$ scale) and explain your observations.
		
		\item[(iii)] Now, change the initialization to $\tilde x_0 = 0.5 x_0 + 0.5 x_*$.
		Does it improve the performance of the method? 
		Try to justify your observation.
	\end{enumerate}
\end{exerciseB}

\begin{Solution}
See Jupyter notebook \url{https://github.com/Sniper2k/inverse_problems_ai_exercises/blob/main/Ex4_3_15.ipynb}.
\end{Solution}

\begin{exerciseB}[Simulation for \Cref{rem: kl}]\label{ex:MexHat}
	Plot $(x_k, y_k) = (r_k \cos( \theta_k), r_k \sin(\theta_k))$, $k \ge 2$, with $\theta_k = k$ and $r_k = \sqrt{(\theta_k - 1)/\theta_k}$.
	What do you observe?
	Note that the $x_k$ lie on the continuous gradient descent trajectory associated with the $f$ in \eqref{eq:MexHat}.
	Now, show that $\Vert \nabla f (x_k, y_k)\Vert_2 \leq 2\exp(-k)$ in Cartesian coordinates.
\end{exerciseB}

\begin{Solution}
	Since our function is given in polar coordinates, it is convenient to compute the derivatives with respect to polar coordinates. To transit back to the Cartesian coordinates, we need to establish conversion formula. By $x = r \cos \theta$ and $y = r \sin(\theta)$, we have
	\begin{align*}
	\frac{\partial x}{\partial r} = \cos \theta = \frac{x}{r}, \quad
	\frac{\partial x}{\partial \theta} = - r \sin \theta = - y, \quad
	\frac{\partial y}{\partial r} = \sin \theta = \frac{y}{r}, \quad 
	\frac{\partial y}{\partial \theta} = r \cos \theta = x.
	\end{align*}
	Then, the chain rule gives
	\begin{align}
		r \frac{\partial f}{\partial r} 
		& = r \frac{\partial f}{\partial x} \frac{\partial x}{\partial r} 
		+ r \frac{\partial f}{\partial y} \frac{\partial y}{\partial r}
		= r \cos \theta \frac{\partial f}{\partial x}  + r \sin \theta \frac{\partial f}{\partial y} 
		= x \frac{\partial f}{\partial x} + y \frac{\partial f}{\partial y} \label{eq: hat tech 1}\\ 
		\frac{\partial f}{\partial \theta} 
		& = \frac{\partial f}{\partial x} \frac{\partial x}{\partial \theta} 
		+ \frac{\partial f}{\partial y} \frac{\partial y}{\partial \theta}
		= - y \frac{\partial f}{\partial x} + x \frac{\partial f}{\partial y} \label{eq: hat tech 2}
	\end{align} 
	Multiplying \eqref{eq: hat tech 1} with $x$ and subtracting \eqref{eq: hat tech 2} multiplied by $y$ gives
	\[
	x^2 \frac{\partial f}{\partial x} + y^2 \frac{\partial f}{\partial x} = x r \frac{\partial f}{\partial r} - y \frac{\partial f}{\partial \theta}.
	\]
	Similarly, multiplying \eqref{eq: hat tech 1} with $y$ and adding to \eqref{eq: hat tech 2} multiplied by $x$ yields
	\[
	x^2 \frac{\partial f}{\partial y} + y^2 \frac{\partial f}{\partial y} = y r \frac{\partial f}{\partial r} + x \frac{\partial f}{\partial \theta}.
	\]
	Combining both equalities with $x^2 + y^2 = r^2$, we get
	\begin{equation}\label{eq: polar coord derivatives}
		\frac{\partial f}{\partial x} 
		= \frac{x}{r} \frac{\partial f}{\partial r} - \frac{y}{r^2} \frac{\partial f}{\partial \theta}, 
		\quad \text{and} \quad 
		\frac{\partial f}{\partial y} 
		= \frac{y}{r} \frac{\partial f}{\partial r} + \frac{x}{r^2} \frac{\partial f}{\partial \theta}.
	\end{equation}
	Then, the norm of the gradient is given by
	\begin{align}
		\| \nabla f \|^2 
		&= \left( \frac{\partial f}{\partial x} \right)^2 + \left( \frac{\partial f}{\partial y} \right)^2
		= \left( \frac{x}{r} \frac{\partial f}{\partial r} - \frac{y}{r^2} \frac{\partial f}{\partial \theta} \right)^2 + \left( \frac{y}{r} \frac{\partial f}{\partial r} + \frac{x}{r^2} \frac{\partial f}{\partial \theta} \right)^2 \notag \\
		& = \frac{x^2}{r^2} \left( \frac{\partial f}{\partial r} \right)^2 - 2 \frac{x y}{r^3} \frac{\partial f}{\partial r} \cdot \frac{\partial f}{\partial \theta} + \frac{y^2}{r^4} \left( \frac{\partial f}{\partial \theta} \right)^2
		+ \frac{y^2}{r^2} \left( \frac{\partial f}{\partial r} \right)^2 + 2 \frac{x y}{r^3} \frac{\partial f}{\partial r} \cdot \frac{\partial f}{\partial \theta} + \frac{x^2}{r^4} \left( \frac{\partial f}{\partial \theta} \right)^2  \notag \\
		& = \frac{x^2 + y^2}{r^2} \left( \frac{\partial f}{\partial r} \right)^2 + \frac{x^2 + y^2}{r^4} \left( \frac{\partial f}{\partial \theta} \right)^2
		= \left( \frac{\partial f}{\partial r} \right)^2 + \left( \frac{1}{r} \frac{\partial f}{\partial \theta} \right)^2. \label{eq: polar norm} 
	\end{align}
	
	Now, we turn to computation of the derivatives of $f$. Let us consider two helper functions
	\[
	g(r) \coloneqq -\frac{1}{1 - r^2}, \quad 
	h(r) \coloneqq \frac{4r^4}{4r^4 + (1 - r^2)^4}.
	\]
	With them, we get by chain rule
	\begin{align*}
	\frac{\partial f}{\partial r} 
	& = e^{g(r)} g'(r) \left[1 - h(r) \sin(\theta + g(r)) \right] 
	- e^{g(r)} \left[ h'(r) \sin(\theta + g(r)) + h(r) \cos(\theta + g(r)) g'(r) \right] \\
	\frac{\partial f}{\partial \theta} 
	& = - e^{g(r)} h(r) \cos(\theta + g(r))
	\end{align*}
	We only want to evaluate the gradients at points defined by $\theta_k = k$, $r_k = \sqrt{(\theta_k - 1) / \theta_k}$. For such points it holds
	\[
	\theta_k + g(r_k) 
	= \theta_k -\frac{1}{1 - r_k^2}
	= \theta_k -\frac{1}{1 - \frac{\theta_k - 1}{\theta_k}}
	= \theta_k - \theta_k 
	= 0,
	\]
	so that $\sin(\theta_k + g(r_k)) = 0$ and $\cos(\theta_k + g(r_k)) = 1$. Therefore, the formula for the partial derivative simplifies to 
	\begin{align*}
		\frac{\partial f}{\partial r}(r_k,\theta_k) 
		& = e^{g(r_k)} g'(r_k)  
		- e^{g(r_k)} h(r_k) g'(r_k)
		= e^{g(r_k)} g'(r_k) [1 - h(r_k)], \\
		\frac{\partial f}{\partial \theta} (r_k,\theta_k)
		& = - e^{g(r_k)} h(r_k).
	\end{align*} 
	With derivative
	\[
	g'(r) = - \left(\frac{1}{1 - r^2} \right)' 
	= - \left( - \frac{1}{(1-r^2)^2} \right) (-2r) 
	= - \frac{2 r^2}{(1-r^2)^2}
	\]
	we obtain
	\begin{align*}
		\frac{\partial f}{\partial r}(r_k,\theta_k) 
		& = - e^{- \frac{1}{1 - r_k^2}} \frac{2 r^2_k}{(1-r^2_k)^2} \left[1 - \frac{4r^4_k}{4r^4_k + (1 - r^2_k)^4} \right]
		= - e^{- \frac{1}{1 - r_k^2}} \frac{2 r^2_k (1 - r^2_k)^2}{4r^4_k + (1 - r^2_k)^4},
		\\
		\frac{\partial f}{\partial \theta} (r_k,\theta_k)
		& = - e^{- \frac{1}{1 - r_k^2}} \frac{4r^4_k}{4r^4_k + (1 - r^2_k)^4}
	\end{align*} 
	Substituting the above derivatives into \eqref{eq: polar norm} gives
	\begin{align*}
	\| \nabla f(x_k,y_k) \|^2 
	&= \left( e^{- \frac{1}{1 - r_k^2}} \frac{2 r^2_k (1 - r^2_k)^2}{4r^4_k + (1 - r^2_k)^4}\right)^2
	+ \left( \frac{1}{r_k} e^{- \frac{1}{1 - r_k^2}} \frac{4r^4_k}{4r^4_k + (1 - r^2_k)^4} \right)^2 \\
	& = e^{- \frac{2}{1 - r_k^2}} \frac{4 r^4_k [ (1- r_k^2)^4 + 4 r_k^2] }{(4r^4_k + (1 - r^2_k)^4)^2}.
	\end{align*}
	With $r_{k} = \sqrt{(\theta_k - 1) / \theta_k} = \sqrt{(k-1)/k}$ it gives
	\begin{align*}
		\| \nabla f(x_k,y_k) \|^2 
		&=  4 e^{- 2 k} \frac{\frac{(k-1)^2}{k^2} [ k^{-4} + 4 \frac{(k-1)}{k}] }{(4 \frac{(k-1)^2}{k^2} + k^{-4})^2}
		= 4 e^{- 2 k} \frac{ k^2 (k-1)^2 [ 1 + 4 (k-1) k^3 ] }{(4 (k-1)^2 k^2 + 1)^2}.
	\end{align*}
	For $k \ge 2$, we have $k \le 2(k-1)$,
	\[
	\frac{k^2 (k-1)^2}{4 (k-1)^2 k^2 + 1} \le \frac{1}{4}
	\quad \text{and} \quad 
	\frac{1 + 4 (k-1) k^3}{4 (k-1)^2 k^2 + 1} \le \frac{1 + 8 (k-1)^2 k^2}{4 (k-1)^2 k^2 + 1} \le 2 .
	\]
	Consequently,
	\[
	\| \nabla f(x_k,y_k) \|^2 \le 4 e^{- 2 k} \cdot \tfrac{1}{4} \cdot 2 = 2 e^{- 2 k} \le 4 e^{- 2 k}. 
	\]
\end{Solution}

\subsubsection{Forward-Backward Splitting}
As motivation, we consider the LASSO objective from \Cref{ex: LASSO}, namely
\begin{equation}\label{eq:LASSO}
	f(x)=\frac12\| A x -y\|_{\bb Y}^2 + \alpha \|x\|_1.
\end{equation}
Since $h(x) = \alpha \|x\|_1$ is not differentiable, we cannot use gradient descent.
Instead, we split $f$ as $f = g + h$ with $g=\frac12\| A x -y\|_{\bb Y}^2$, and perform a gradient step for $g$ and a \emph{proximal} step for $h$.
This results in the  \emph{forward-backward splitting} (also know as \emph{proximal gradient algorithm}) for composite objectives $f=g+h$, see \Cref{algo:FBS}.

\begin{algorithm}[H]
	\KwIn{$g \in C^1(\R^N)$ such that $\nabla g$ is $L$-\textsc{Lipschitz} continuous, a coercive and lsc $h\colon \R^N \to \overline{\R}$ with $\q D(h) \neq \emptyset$, initialization $x_0 \in \R^N$, step size $0 < \nu < 1/L$.}
	\For{\color{black}$k = 0, 1, \ldots$}{
		$x_{k + 1} \in \prox_{\nu h}\left(x_{k} - \nu \nabla g(x_{k})\right)$\;
	}
	\caption{Forward-Backward Splitting Algorithm (FBS)}
	\label{algo:FBS}
\end{algorithm}

Before we analyze \Cref{algo:FBS}, we have to properly introduce proximal maps.
\begin{defn}[Proximal mapping]
	Let $\nu > 0$ and $f\colon \R^N \to \overline{\R}$.
	Then, the (set-valued) \emphi{proximal mapping} $\prox_{\nu f}\colon \R^N \rightrightarrows \R^N$ is defined by
	\begin{equation}\label{eq: proximal mapping}
		\prox_{\nu f}(x)
		\coloneqq \argmin_{z \in \R^N} \bigl\{f(z) + \tfrac{1}{2 \nu} \| z - x\|_{\bb X}^2\bigr\}.
	\end{equation}
\end{defn}

Essentially, $\prox_{\nu f}(x)$ finds the smallest value of $f$ close to the anchor point $x$. 
To gain a better intuition on proximal mappings, we consider the following examples.
More proximal mappings can be found at \url{http://proximity-operator.net/}.

\begin{exampleB}[Proximal Mappings]\label{exampleB:proxofl1norm}
	\begin{itemize}
		\item For $h = 0$ and any $\nu>0$, we have that $\prox_{\nu h}(x) = x$.
		Hence,  \Cref{algo:gradient_descent} is a special case of \Cref{algo:FBS} for the choice $h=0$.
		\item Let $C \subseteq \bb \R^N$ be nonempty, closed and convex.
		Consider $\iota_C \colon \R^N \to \overline{\R}$ with $\iota_C(z) = 0$ for $z \in C$ and $\iota_C(z)=\infty$ else. 
		Then, $\prox_{\nu \iota_C}(x) = \argmin_{z \in C} \|z - x\|_2$ is the orthogonal projection onto $C$.
		Thus, $\prox_{\nu f}$ generalizes the orthogonal projection.
		\item
		For $f \colon \R \to \R$ with $y \mapsto | y |$ we have that $S_{\nu}(x) \coloneqq \prox_{\nu f}(x)$ is given by
		\begin{equation}\label{eq: soft thresholding}
		\begin{aligned}
			S_{\nu}(x) = \sign(x) \max \{ |x| - \nu, 0 \} = \begin{cases}
				x - \nu,    & \text{for } x > \nu, \\
				0,              & \text{for } | x | \le \nu, \\
				x + \nu,    & \text{for } x < -\nu,
			\end{cases}
		\end{aligned}
	\end{equation}
		which is the \emphi{soft-thresholding operator} with \emph{threshold} $\nu$ depicted in \Cref{fig:soft_thresholding}.
		For the multivariate analog $\| \cdot \|_1$ of $\vert \cdot \vert$, it holds that
		\begin{equation}
			\prox_{\nu \| \cdot \|_1}(x)
			= \argmin_{z \in \R^N} \sum_{k = 1}^{N} \left(\tfrac{1}{2} (z_k - x_k)^2 + \nu | z_k | \right),
		\end{equation}
		so that $\prox_{\nu \| \cdot \|_1}(x)$ is the componentwise application of $S_{\nu}$.
		
		\begin{figure}
			\begin{center}
				\subfloat[Soft-thresholding operator with $\nu = 1/2$. \label{fig:soft_thresholding}]{
					\includegraphics[width= 180pt, height = 180pt]{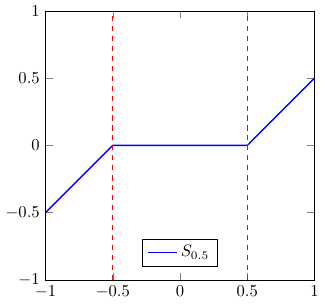}% 
				} \hfill
				\subfloat[Hard-thresholding operator with $\nu = 1/2$.
				\label{fig:hard_thresholding}]{
					\includegraphics[width= 180pt, height = 180pt]{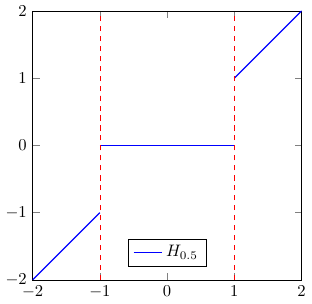}%
				}
				\caption{Examples of proximal mappings.}
			\end{center}
		\end{figure}
		
		\item For $f\colon \R \to \R$ with $f(0) = 0$ and $f(z)= 1$ else, the objective in $\prox_{\nu f}(x)$ reads
		\begin{equation}\label{eq:ObjHard}
		\frac{1}{2\nu} |z - x|^2 + f(z)
		= 
		\begin{cases}
			\frac{1}{2\nu} |z - x|^2 + 1, & z \neq 0,\\
			\frac{1}{2\nu} |x|^2, & z = 0.
		\end{cases}
		\end{equation}
		Therefore, as long as $\frac{1}{2\nu} |x|^2 \le 1$, the minimum is attained at $z = 0$. 
		On the other hand, for $\frac{1}{2\nu}|x|^2 \ge 1$, taking $z = x$ achieves the minimal value $1$ in \eqref{eq:ObjHard}.
		Consequently, we get
		\begin{equation}
			H_{\nu}(x) \coloneqq
			\prox_{\nu f}(x) = 
			\begin{cases}
				0, & |x| < \sqrt{2\nu}, \\
				\{0,x\}, & |x| = \sqrt{2\nu},\\
				x, & |x| > \sqrt{2\nu},
			\end{cases} 
		\end{equation}
		which is the \emphi{hard thresholding operator} with \emph{threshold} $\nu$, see also  \Cref{fig:hard_thresholding}. Similarly to $S_{\nu}$, we can consider the multidimensional extension $\|x\|_0 = \sum_{k=1}^N f(x_k)$.
		Again, $\prox_{\nu \|x\|_0}$ is the componentwise application of $H_{\nu}$. 
		
		\item If $f\colon \R^N \to \R$ with $f(z)=\Vert A z - y\Vert_2^2$ and $\nu>0$, then we get that $\prox_{\nu f}(x)$ is the solution of the linear equation $(\frac{1}{\nu}I + A^TA) z = A^Ty+\frac1\nu x$, see \Cref{ex: quadratic regularization}.
		This system can be solved efficiently if fast evaluation routines for $A^T A z$ exist (such as for deblurring).
		Then, implicit steps become attractive.
	\end{itemize}
\end{exampleB}

\begin{remarkB}[Gradient Flows: Explicit vs Implicit]\label{rem:GradFlow}
	We can view both gradient and proximal updates as time discretization for the gradient flow
	\begin{equation}\label{eq:GradFlow}
		\dot x(t)=-\nabla f(x(t)).
	\end{equation}
	The explicit \textsc{Euler} discretization of \eqref{eq:GradFlow} reads
	\begin{equation}
		x_{k+1} - x_k = - \nu \nabla f(x_k)
	\end{equation}
	and corresponds to \cref{algo:gradient_descent}.
	On the other hand, the proximal  update $x_{k+1} = \prox_{\nu f}(x_k)$ satisfies $0 = \nabla f(x_{k+1}) + \frac{1}{\nu} (x_{k+1} - x_{k})$ by Fermat's rule.
	Rearranging the terms gives
	\begin{equation}\label{eq: implicit}
		x_{k+1} - x_k = - \nu \nabla f(x_{k+1}),
	\end{equation}
	which is the implicit \textsc{Euler} scheme for the gradient flow \eqref{eq:GradFlow}.
	Similarly as for stiff ODEs, implicit steps are sometimes preferable in optimization even though more expensive.
\end{remarkB}

\begin{remark}
	In the following, we talk about critical points without specifying what this means for non-differentiable $f$.
	For the convex case, their definition involves the subdifferential $\partial f$ as introduced in the course \emph{Convex Analysis}. For all $x_0 \in \R^N$, $\partial f(x_0)$ is given by
	\begin{equation}\label{eq: subdifferential}
		\partial f(x_0) = \{ z \in \R^N : f(x) \ge f(x_0) + \inner{z}{x - x_0} \quad \text{for all} \ x \in \R^N \}, 
	\end{equation}
	i.e, the collection of all directions $z$ such that the line $f(x_0) + \inner{z}{x - x_0}$ stays under the graph of $f$.
	This inequality can be seen as a generalization of \eqref{eq: strong convexity}.
	For the general case, even more complicated constructions with the (limiting) Fr\'echet subdifferential are required \cite{AttBolSva2013}.
	These are far beyond the scope of this course.
	For our intuition, it suffices to keep in mind that these notions generalize differentiability and hence also our definition of critical points.
	More precisely, $\hat x$ is a critical point of $f \colon \R^N \to \overline{\R}$ if $0 \in \partial f(\hat x)$.
\end{remark}
Now, we summarize a few important properties of the proximal mapping, see \Cref{ex:ProxProp}.

\begin{thm}{Properties of Proximal Mappings}{proxprop}
	Let $f\colon \R^N \to  \overline{\R}$ with $\q D(f) \neq \emptyset$ be lsc.
	We have the following:
	\begin{enumerate}
	\item If $f$ is convex or bounded from below, then for any $x \in \R^N$ and $\nu>0$ the function $z \mapsto f(z) + \frac{1}{2 \nu} \| z - x\|_{\bb X}^2$ is coercive, lsc and it holds that $\prox_{\nu f}(x) \neq \emptyset$.
	\item An element $\hat{x} \in \q D(f)$ is a critical point of $f$ if and only if  $\hat x \in \prox_{\nu f}(\hat x)$.
	
	\item If $f$ is convex, then \smash{$\prox_{\nu f}(x)$} is single-valued and we can write $\hat{x}=\prox_{\nu f}(x)$.
	In this case, \smash{$\hat{x} = \prox_{\nu f}(x)$} if and only if the \emphi{variational inequality}
	\begin{equation}
		\tfrac{1}{\nu} \lsk x - \hat{x}, z - \hat{x} \rsk + f(\hat{x}) \le f(z) \qquad \text{holds for all $z \in \R^N$.}
	\end{equation}
	
	\end{enumerate}
\end{thm}
Now, we can prove convergence of \cref{algo:FBS}.

\begin{thm}{Convergence of FBS}{FBS}
	Assume that $g\in C^1(\R^N)$ has $L$-Lipschitz continuous gradient and that $h\colon \R^N \to \overline{\R}$ with $\q D(h) \neq \emptyset$ is lsc.
	If $g$ and $h$ are bounded from below, $x_{0} \in \R^N$, and $0<\nu < 1/L$, the sequence $\{x_k\}_{k \in \N}$  generated by \Cref{algo:FBS} satisfies
    $\lim_{k \to \infty}\| x_{k+1} - x_k \|_{\bb X} =0$ and
    \begin{equation}
    	\min_{k = 0, \ldots, K-1} \| x_{k+1} - x_k \|_{\bb X}^2 
    	\le 2\nu \frac{f(x_0) - \inf_{x \in \bb X}f(x)}{(1 - L \nu) K}.
    \end{equation}
    Moreover, $\{f(x_k)\}_{k \in \N}$ is decreasing.   
    \end{thm}

\begin{pr}
	By definition of $\prox_{\nu h}$, it holds that
	\begin{align} 
		x_{k+1} & \in \argmin_{z\in \R^N} \bigl\{ \tfrac{1}{2\nu} \| z - (x_k - \nu \nabla g(x_{k})) \|_{\bb X}^2 + h(z) \bigr\} \notag\\
		& = \argmin_{z\in \R^N} \bigl\{ \tfrac{1}{2\nu} \| z - x_k \|_{\bb X}^2 + \inner{z - x_k}{\nabla g(x_{k})}  + h(z) \bigr\}.\label{eq:UpFBS}
	\end{align}
	Thus, taking $z = x_k$ in \eqref{eq:UpFBS}, we obtain
	\begin{equation}\label{eq:EstFBSH}
		\tfrac{1}{2\nu} \| x_{k+1} - x_k \|_{\bb X}^2 + \inner{x_{k+1} - x_k}{\nabla g(x_k)} + h(x_{k+1}) \le h(x_k).
	\end{equation}
	On the other hand, applying the Descent \Cref{l: descent lemma} with $x = x_k$ and $h = x_{k+1} - x_k$ gives
	\begin{equation}\label{eq:EstFBSG}
		g(x_{k+1}) \le g(x_k) + \inner{x_{k+1} - x_k}{\nabla g(x_k)} + \frac{L}{2} \| x_{k+1} - x_k \|_{\bb X}^2.
	\end{equation}
	Combining \eqref{eq:EstFBSH} and \eqref{eq:EstFBSG}, we get for $f= g + h$ that
	\begin{align}
		f(x_{k+1}) & \le g(x_k) + \inner{x_{k+1} - x_k}{\nabla g(x_k)} + \tfrac{L}{2} \| x_{k+1} - x_k \|_{\bb X}^2 + h(x_{k+1}) \notag \\
		& \le  f(x_k) - \bigl(\tfrac{1}{2\nu} - \tfrac{L}{2}\bigr) \| x_{k+1} - x_k \|_{\bb X}^2.
	\end{align}
	Since $\nu < 1/L$, it holds that $\tfrac1\nu - L > 0$, and we can from now follow the proof of \Cref{th:grad descent nonconvex} to obtain the claim.
\end{pr}

To show that every accumulation point of FBS is a critical point, we show that the smallest norm of the elements in $\partial f(x_k)$ goes to zero, in analogy to $\nabla f(x_k) \to 0$ for gradient descent. 
We use without proof that for $q \in C^1(\R ^N)$ it holds
\begin{equation}
\partial (q + h)(x) = \nabla q(x) + \partial h(x) \eqqcolon \{ \nabla q(x) +  u \in \R^N : u \in \partial h(x) \}.
\end{equation}
By Fermat's rule $x_{k + 1} \in \prox_{\nu h}(x_{k} - \nu \nabla g(x_{k}))$ if and only if 
\begin{equation}
	0 \in \partial\bigl(h(x_{k+1}) + \tfrac{1}{2\nu} \| x_{k+1} - (x_k - \nu \nabla g(x_k) ) \|^2_{\bb X} \bigr) 
	= \tfrac{1}{\nu} (x_{k + 1} - x_k +  \nu \nabla g(x_{k}) ) + \partial h(x_{k+1}),
\end{equation}
which is equivalent to $\tfrac{1}{\nu} (x_{k} - x_{k+1})
\in \nabla g(x_{k}) + \partial h(x_{k+1})$.
Hence, there exists $\hat u_k \in \partial h(x_{k+1})$ such that $\hat u _k + \nabla g(x_{k}) = \tfrac{1}{\nu} (x_{k} - x_{k+1})$.
%we have $\| x_{k+1} - x_k \|_{\bb X} = \nu \| \partial h(x_{k+1}) + \nabla g(x_{k})  \|_{\bb X}$ and
Consequently, 
\begin{align}
	\min_{u \in \partial f(x_{k+1})} \| u \|_\bb X & \le \| \hat u_k + \nabla g(x_{k+1}) \|_\bb X \notag \\
	& \leq \| \tfrac{1}{\nu} (x_{k} - x_{k+1}) - \nabla g(x_{k}) + \nabla g(x_{k+1})  \|_{\bb X}\notag \\
	& \le \tfrac{1}{\nu} \| x_{k} - x_{k+1} \|_{\bb X}  + L \| x_{k+1} - x_k \|_{\bb X} 
	< \tfrac{2}{\nu} \| x_{k+1} - x_k \|_{\bb X},
\end{align}
which induces a rate of the form \eqref{eq:OptimalityGD}.

\begin{remark}[Convergence of Iterates]
	As for  \Cref{algo:gradient_descent}, FBS does not necessarily converge to a critical point (since  \Cref{algo:gradient_descent} is just a special case of \Cref{algo:FBS}).
	Again, if $f$ satisfies the KL property \cite{AttBolSva2013} and an additional boundedness condition is met, then $\lim_{k\to \infty} x_k$ exists and is a critical point of $f$. 
\end{remark}

As an important special case of  \Cref{algo:FBS}, we derive the FBS updates for \Cref{ex: LASSO}.

\begin{exampleB}[Iterative soft-thresholding algorithm (ISTA)]\label{ex:ista}
    Continuing on \Cref{ex: LASSO}, we split the objective  \eqref{eq:LASSO} as $f = g + h$
    with $g(x) = \frac{1}{2} \| A x - y \|_{\bb Y}^2$ and $h(x) = \alpha \| x \|_1$.
    Then, $\nabla g(x) =  A^\tT (A x -y)$ and $\prox_{h} = S_{\alpha}$ is the soft-shrinkage operator, see \Cref{exampleB:proxofl1norm}.
    Thus, given a stepsize $0 <\nu < 1/\|A\|^2_{\q L(\bb X,\bb X)}$,  \Cref{algo:FBS} becomes
    \begin{align}\label{eq:ISTA}
        x_{k + 1} 
        = S_{\alpha \nu} \left( x_{k} - \nu  A^\tT(A x_{k} - y) \right).
    \end{align}
    The updates \eqref{eq:ISTA} are known as \emphi{iterative soft-thresholding algorithm} (ISTA).
    It is, in particular, important for dictionary-based image reconstruction, which will be considered in the exercises.
    For a generalization of ISTA to Hilbert spaces see \cite{DDM04}.
    The acceleration in \Cref{algo:fFBS} is known as \emph{fast iterative shrinkage algorithm} (FISTA)\index{FISTA} \cite{beck2009fast}.
\end{exampleB}

Similarly as gradient descent, FBS can be accelerated with an extrapolation step.
For the choice $h=0$, \Cref{algo:fFBS} reduces to \Cref{algo:accelerated_gradient_descent}.

\begin{algorithm}[H]
 \KwIn{$g \in C^1(\R^N)$ such that $\nabla g$ is $L$-\textsc{Lipschitz} continuous, a lsc and coercive $h\colon \R^N \to \overline{\R}$ with $\q D(h) \neq \emptyset$, initialization $x_0 \in \R^N$, step size $0 < \nu < 1/L$.}
 Set $z_{0} = x_0$, and $\theta_k = \frac{2}{k+2}$ for $k \ge 0$\;
 \For{\color{black} $k = 0, 1, \ldots$}{
 	$u_k = (1 - \theta_k) x_k + \theta_k z_k $\;
 	$x_{k + 1} \in \prox_{\nu h}\left( u_{k} - \nu \nabla g(u_{k}) \right)$\;
 	$z_{k+1} = x_k + \frac{1}{\theta_k}(x_{k+1} - x_k)$\;
 }
 \caption{Fast FBS}
 \label{algo:fFBS}
\end{algorithm}

Here, we express that rate again in terms of the objective values $\{f(x_k)\}_{k \in \N}$.

\begin{thm}{$1/k^2$ convergence of fast FBS \cite{Ts08}}{proximal gradient accelerated}
	Let $f= g + h$, where $g \in C^1(\R^N)$ is convex with $L$-\textsc{Lipschitz} continuous gradient and $h\colon \R^N \to \overline{\R}$ with $\q D(h) \neq \emptyset$ is convex, coercive and lsc.
%	If $f$ has a minimizer $\hat{x}$, 
	Then, the function values $\{f(x_k)\}_{k \in \N}$  at the iterates generated by \Cref{algo:fFBS} with stepsize $0 < \nu \le \tfrac1L$ satisfy
	\begin{equation}
		f(x_{K}) - \inf_{x \in \bb X}f(x) \leq  \frac{2 L}{(K+1)^2} \min_{x_* \in \argmin(f)} \| x_0 - x_* \|_{\bb X}^2.
	\end{equation}
\end{thm}

\begin{pr}
	
	Firstly, we consider the progress in one step of the algorithm.
	By the Lipschitz continuity of $\nabla g$ and since $\nu \le \frac{1}{L}$, we know by the Descent \Cref{l: descent lemma} that for $k \in \N$ it holds
	\begin{equation} \label{eq:important_1}
		g(x_{k + 1}) \le g(u_k) + \left\langle \nabla g(u_k), x_{k + 1} - u_k \right\rangle + \frac{1}{2\nu} \|x_{k + 1} - u_k\|^2.
	\end{equation}
	Moreover, the variational inequality for proximal operators (\Cref{th:proxprop}\, \circled{2}) implies for all $u \in \bb X$ that 
	\begin{align}
		h(x_{k + 1})
		& \le  h(u) + \tfrac{1}{\nu} \left\langle u_k - \nu \nabla g(u_k) - x_{k + 1}, x_{k + 1} - u \right\rangle \nonumber\\
		& =  h(u) - \left\langle \nabla g(u_k), x_{k + 1} - u \right\rangle  +  \tfrac{1}{\nu} \left\langle u_k - x_{k + 1}, x_{k + 1} - u \right\rangle.
		\label{eq:important_2}
	\end{align}
	Adding the inequalities \eqref{eq:important_1} and \eqref{eq:important_2} yields
	\begin{align}\label{eq:EstFunValueProx}
		f(x_{k + 1})
		&\le f(u) \underbrace{- g(u) +  g(u_k) + \langle \nabla g(u_k), u - u_k \rangle}_{\le 0\text{ (by convexity of $g$)}} \notag\\
		& + \tfrac{1}{2\nu} \|x_{k + 1}-u_k\|^2
		+ \tfrac{1}{\nu} \langle u_k - x_{k + 1}, x_{k + 1} - u\rangle \notag\\
		& \le f(u)+ \tfrac{1}{2\nu} \|x_{k + 1} - u_k\|^2
		+ \tfrac{1}{\nu} \langle u_k - x_{k + 1}, x_{k + 1} - u\rangle.
	\end{align}
	The application of \eqref{eq:EstFunValueProx} with $u \coloneqq \hat x$  and $u \coloneqq x_k$ leads for $\theta_k \in [0,1]$ to
	\begin{align}
		& f(x_{k + 1}) -  f(\hat x) + (1-\theta_k)\left(f(\hat x) - f(x_k) \right) \notag \\
		& \quad = \theta_k \left( f(x_{k + 1}) - f(\hat x) \right) + (1-\theta_k) \left( f(x_{k + 1}) - f(x_k) \right) \notag \\
		& \quad \le \tfrac{1}{2\nu} \| x_{k + 1}-u_k \|^2 + 
		\tfrac{1}{\nu} \bigl\langle u_k-x_{k + 1}, x_{k + 1} - \theta_k \hat x - (1-\theta_k) x_k \bigr\rangle \notag \\
		& \quad = \tfrac{1}{2\nu} \left( 
		\| u_k -  \theta_k \hat x - (1-\theta_k) x_k \|^2  - \| x_{k + 1} -  \theta_k \hat x - (1-\theta_k) x_k\|^2 \right).
%		& = \frac{\theta_k^2}{2\nu} \left( \|z_k - \hat x\|^2 - \|z_{k + 1} - \hat x\|^2 \right).
	\end{align}
	By construction, we have  $u_k - (1-\theta_k) x_k = \theta_k z_k$ and $x_{k+1} - (1 - \theta_{k} )x_k = \theta_k z_{k+1}$.
	This yields
	\begin{equation}
		f(x_{k + 1}) -  f(\hat x) + (1-\theta_k)\left(f(\hat x) - f(x_k) \right) \le \tfrac{\theta_k^2}{2\nu} \left( \|z_k - \hat x\|^2 - \|z_{k + 1} - \hat x\|^2 \right).
	\end{equation}

	Thus, we obtain for a single step that
	\begin{equation}
		\frac{\nu}{\theta_k^2} \left(f(x_{k + 1}) - f(\hat x) \right) + \frac{1}{2} \|z_{k + 1} - \hat x\|^2
		\le \frac{\nu(1-\theta_k)}{\theta_k^2} \left( f(x_k )- f(\hat x) \right) + \frac{1}{2} \| z_k - \hat x\|^2.
	\end{equation}
	Using the relation recursively on the right-hand side and regarding that 
	\begin{equation}
		\frac{1 - \theta_k}{\theta_k^2} = \frac{1 - \frac{2}{k+2} }{ \frac{4}{(k+2)^2}} = 
		\frac{k(k+2)}{4}
		\le \frac{(k+1)^2}{4} = \frac{1}{\theta_{k-1}^2},
	\end{equation}
	we obtain
	\begin{equation}
		\frac{\nu}{\theta_k^2} \left(f(x_{k + 1}) - f(\hat x) \right)
		\le \frac{\nu(1-\theta_0)}{\theta_0^2} \left(f(x_{0}) - f(\hat x) \right)
		+ \frac{1}{2} \| z_{0} - \hat x\|^2
		= \frac{1}{2} \| x_{0} - \hat x\|^2.
	\end{equation}
	This yields the assertion $f(x_{k + 1})- f(\hat x) \le \frac{2}{\nu (k+2)^2}\| x_{0} - \hat x \|^2$.
\end{pr}

Thus, \Cref{th:gradient accelerated} holds by combining \Cref{exampleB:proxofl1norm} and \Cref{th:proximal gradient accelerated}.

\begin{remark}[Generalizations of fast FBS]\label{rem:FBS}
	If $L$ is unknown, again a line search can be incorporated \cite{GS2011,Ne07} into \Cref{algo:fFBS}, see also \Cref{rem:LineSearch}.
	Moreover, we can deploy (inexact) \emph{variable metric strategies}\index{variable metric strategies} \cite{CV2013,PLS08}.
	
    In the nonconvex case, there exist many variants and generalizations of \Cref{algo:fFBS} with various different momentum choices and linesearch methods \cite{CarDucHin2018,LiLin2015} with an overview in \cite{DanDvuGas2022}.
    As for gradient descent, neither of them achieves provably improved rates. 
\end{remark}

\subsubsection{Exercises}
\begin{exerciseB}\label{ex:ProxProp}
	Prove \Cref{th:proxprop}.
\end{exerciseB}
\begin{Solution}
		\begin{enumerate}
		\item We only show it for continuously differentiable as the notion of subgradient is beyond the scope of the course. Substituting $\hat x$ for $x_k$ and $x_{k+1}$ in \eqref{eq: implicit} gives that $\hat x = \prox_{\nu f}(\hat x)$ if and only if $\nabla f(\hat x) = 0$.
		
		\item Since $f$ is convex and $\frac{1}{2 \lambda}\|\cdot\|^2$ is strictly convex, the claim follows by \Cref{th:minima}.
		
		\enquote{$\implies$}:
		Let $x\in \bb X$ and $\hat{x} \in \prox_{\nu f} (x)$.
		
		Then, the definition of the proximal mapping for all $z \in \bb X$ and $t\in(0,1)$ gives
		\begin{align*}
			0&\geq f(\hat x)-f(\hat x+t(z-\hat x))+\frac{1}{2 \nu}(\|\hat x-x\|^2-\|\hat x -x - t(\hat x-z)\|^2)\\
			&=f(\hat x)-f((1-t)\hat x+tz)+\frac{1}{2 \nu}(2t\lsk \hat x-x,\hat x-z\rsk-t^2\|\hat x-z\|^2)
		\end{align*}
		and convexity of $f$ further yields
		\begin{equation*}
			0 \ge t f(\hat{x}) - t f(z) + \frac{t}{\nu} \langle \hat{x} - x, \hat{x}- z \rangle - \frac{t^2}{2\nu}\|y - \hat{x}\|^2.
		\end{equation*}
		Finally, dividing by $t > 0$ and taking $t \searrow 0$ we obtain the variational inequality.	
		
		\enquote{$\impliedby$}:
		Conversely, let the variational inequality hold true. Then we have for all $z \in \bb X$ that
		\begin{align*}
			f(\hat{x}) + \frac{1}{2\nu}\| \hat{x} - x\|^2 
			& \le f(z) - \frac{1}{\nu} \langle \hat{x} - x, \hat{x} - z \rangle + \frac{1}{2\nu} \| \hat{x} - x\|^2 \\
			& \le f(z) + \frac{1}{2 \nu}\|\hat{x} - x\|^2 - \frac{1}{\nu} \langle \hat{x} - x, \hat{x} - z \rangle + \frac{1}{2\nu}\| \hat{x} - z\|^2 \\
			& = f(z) + \frac{1}{2\nu}\|z - x\|^2.
		\end{align*}
		Hence $\hat{x} = \prox_{\nu f}(x)$.
	\end{enumerate}
\end{Solution}

\begin{exerciseB}[Proximal Calculus]\label{ex:proxCalc}
	Let $f \colon \R^N \to \overline{\R}$ be convex and lsc with $\q D(f) \neq \emptyset$ and $C \neq \emptyset \subset \R^N$ be closed and convex.
	Show that we have the following calculation rules for the proximal operator.
	\begin{table}[H]
		\centering
		\begin{tabular}{c|c|c}
			\textbf{property}           & $g(x)$                                                    & $\prox_g(x)$ \\  \hline
			translation                 & $f(x - z)$, $z \in \R^N$                                     & $z + \prox_f(x - z)$ \\
			scaling                     & $f\left(\frac{x}{a}\right)$, $a \ne 0$                    & $a \cdot \prox_{a^{-2} f}\left(\frac{x}{a}\right)$ \\
			reflection                  & $f(- x)$                                                  & $- \prox_f(-x)$ \\ 
		\end{tabular}
	\end{table}
	Conclude that every translation on $\R^N$ (that is, a map $z \mapsto z - a$) is a proximal operator for suitable $f$.
\end{exerciseB}

\begin{Solution}
	\begin{align*}
		\prox_{f(\cdot - z)}(x)
		& = \argmin_{y \in E} \left\{\frac{1}{2} \| x - y \|_2^2 + f(y - z) \right\} \\
		& = \argmin_{y \in E} \left\{\frac{1}{2} \| x - z - (y - z) \|_2^2 + f(y - z) \right\} \\
		& = z + \argmin_{a \coloneqq y - z \in E} \left\{\frac{1}{2} \| x - z - (y - z) \|_2^2 + f(y - z) \right\} \\
		& = z + \argmin_{a \in E} \left\{\frac{1}{2} \| x - z - a \|_2^2 + f(a) \right\}
		= z + \prox_f(x - z).
	\end{align*}
	
	\begin{align*}
		\prox_{f(a^{-1} \cdot)}(x)
		& = \argmin_{y \in E} \left\{ \frac{1}{2} \| x - y \|_2^2 + f\left(\frac{y}{a}\right) \right\}
		= \argmin_{y \in E} \left\{ \frac{1}{2}  \left\| x - a \cdot \frac{y}{a} \right\|_2^2 + f\left(\frac{y}{a}\right)\right\}  \\
		& =a \cdot \argmin_{\frac{y}{a} \in E} \left\{ \frac{1}{2} \left\| x - a \cdot \frac{y}{a} \right\|_2^2 + f\left(\frac{y}{a}\right) \right\}
		= a \argmin_{z \in E} \left\{ \frac{1}{2} \left\| x - a z \right\|_2^2 + f(z) \right\} \\
		& = a \argmin_{z \in E} \left\{ \frac{1}{2a^2} \left\| \frac{x}{a} - z \right\|_2^2+ f(z)\right\}
		= a \prox_{a^{-2} f}\left(\frac{x}{a}\right)
	\end{align*}
	
	\begin{align*}
		\prox_{f(-\cdot)}(x)
		& = \argmin_{y \in E} \left\{ \frac{1}{2} \| x - y \|_2^2 + f(-y) \right\}
		= \argmin_{y \in E} \left\{ \frac{1}{2} \| -x - (- y) \|_2^2 + f(-y) \right\} \\
		& = - \argmin_{z \in E} \left\{ \frac{1}{2} \| -x - z \|_2^2 + f(z) \right\}
		= - \prox_{f}(-x)
	\end{align*}
	
\end{Solution}

\begin{exerciseB}[Convergence of FBS with Convex $h$]\label{exer:fbs convex regularizer}
Consider the setup of \Cref{th:FBS}. Assume that $h$ is additionally convex.
Show that FBS with step size $0 < \nu \le 1/L$ converges with the rate 
\begin{equation}
	\min_{k = 0, \ldots, K-1} \| x_{k+1} - x_{x} \|_{\bb X}^2 
	\le 2\frac{f(x_0) - \inf_{x \in \R^N}f(x)}{\nu K}.
\end{equation}
This matches the convergence rate of gradient descent in \Cref{th:grad descent nonconvex}.
\end{exerciseB}

\begin{Solution}
Similarly to the proof of \Cref{th:proximal gradient accelerated}, the Lipschitz continuity of $\nabla g$ and the step size choice $\nu \le \frac{1}{L}$ together with the descent lemma (\cref{l: descent lemma}) yield for $k \in \N$ that
\begin{equation*}
	g(x_{k + 1}) \le g(x_k) + \left\langle \nabla g(x_k), x_{k + 1} - x_k \right\rangle + \frac{1}{2\nu} \|x_{k + 1} - x_k\|^2.
\end{equation*}
Furthermore, \Cref{th:proxprop}, \circled{2} gives that for $u = x_k$ that 
\begin{align*}
	h(x_{k + 1})
	& \le  h(x_k) + \frac{1}{\nu} \left\langle x_k - \nu \nabla g(x_k) - x_{k + 1}, x_{k + 1} - x_k \right\rangle \nonumber\\
	& =  h(x_k) - \left\langle \nabla g(x_k), x_{k + 1} - u \right\rangle  - \frac{1}{\nu} \| x_k - x_{k + 1} \|_{\bb X}.
\end{align*}
Summing two inequalities yields
\[
f(x_{k+1}) = g(x_{k + 1}) + h(x_{k + 1}) \le g(x_k) + h(x_k) - \frac{1}{2\nu} \|x_{k + 1} - x_k\|^2.
\]
For the rest, we follow the proof of \Cref{th:grad descent nonconvex}. 
\end{Solution}

\begin{exerciseB}[Peak Detection with ISTA]\label{exer:ista}
	Implement the ISTA from \Cref{ex:ista}. 
	For the signal $\hat x=(0,0,0,0,0,2,0,0,5,0,3,0)$ and the matrix $A=I$, generate
	data $y \sim A\hat x + \q N(0,I_{12})$.
	Afterwards, run ISTA for this $y$. 
	Count the number of zero-entries in the reconstruction $x$ depending on the parameter $\alpha$. 
	Redo the same experiment with the Tikhonov functional $J_\alpha = \frac12\|Ax-y\|_2^2+\frac{\alpha}{2}\|x\|_2^2$ from \eqref{eq:tikhonov:funktional} and count again the number of zero-entries for its minimizer.
\end{exerciseB}

\begin{Solution}
	See Jupyter Notebook \url{https://github.com/Sniper2k/inverse_problems_ai_exercises/blob/main/Ex4_3_30.ipynb}.
\end{Solution}

\begin{exerciseB}[Impact of Regularizers]\label{exer:tikh vs. ista}
		\begin{itemize} 
		\item[(i)] Let $N=5$ and $M = 16$.
		For $A \in \R^{M \times N}$ with Gaussian entries $A_{k,j} \sim M^{-1/2} \q N(0,1)$ and $x_* \sim \q N(0,I_N)$, we consider the noise-free inverse problem $Ax=y$.
		Implement Tikhonov regularization (see also \Cref{eq:ImplementTikhonov}) and LASSO regularizations via FISTA (\Cref{algo:fFBS}) for $K=10000$ iterations. 
		Plot the absolute values $|x_\alpha|$ and $\log(|x_\alpha| + 10^{-16})$ of the reconstructions $x_\alpha$ for different regularization parameters $\alpha = 10^{q}$ with $q \in [-5, 5]$. 
		What do you observe? 
		\item[(ii)] Load the $28 \times 28$ image \texttt{10104.png} from the MNIST dataset\footnote{\url{https://github.com/rasbt/mnist-pngs/blob/main/train/4/10104.png}} and convert it to a vector $x_* \in \R^{748}$.
		For $M \in [160, 380]$, consider the inverse problems from (i).
		To obtain a reconstruction, minimize the LASSO and Tikhonov functional with $\alpha = 10^{-1}$ and $\alpha = 10^{-2}$, respectively.
		Visualize the reconstructions for different $M$ (as images).
	\end{itemize}   
\end{exerciseB}

\begin{Solution}
	See Jupyter Notebook \url{https://github.com/Sniper2k/inverse_problems_ai_exercises/blob/main/Ex4_3_31.ipynb}.
\end{Solution}

\subsubsection{Alternating Direction Method of Multipliers}\label{sec:ADMM}

Now, we consider composite functions $f = g  + h(C\,\cdot)$, where $g\colon \R^N \to \overline{\R}$, $h\colon \R^L \to \overline{\R}$ and $C \in \R^{L \times N}$.
This leads to the constrained optimization problem
\begin{equation} \label{problem_pdhg_1}
\argmin_{x \in \R^N, z \in \R^L} \left\{ g(x) + h(z)  \quad {\rm s.t.} \quad C x = z \right\},
\end{equation}
for which we define the \emphi{augmented Lagrangian} for $\gamma > 0$ and a \emph{Lagrange multiplier} $p \in \R^L$ as
\begin{align}\label{eq:AugLagrangian}
	L_\gamma(x,z,p) &\coloneqq g(x) + h(z) + \langle p, Cx - z \rangle + \tfrac{\gamma}{2} \|Cx - z\|_2^2,\\
	&= g(x) + h(z) +  \tfrac{\gamma}{2} \|Cx - z + \tfrac{p}{\gamma} \|_2^2 - \tfrac{1}{2\gamma} \|p\|_2^2.\notag
\end{align}
For $\gamma =0$, this reduces to the well known \emph{Lagrangian}.
If there exists a \emph{saddle point} $(\hat x,\hat z,\hat p)$ of \eqref{eq:AugLagrangian}, namely a point with
\begin{equation}
	L_\gamma(\hat x, \hat z, p) \le L_\gamma(\hat x,\hat z,\hat p) \le L_\gamma(x,z,\hat p),
	\quad \text{for all } x \in \R^N, \ z,p \in \R^L,
\end{equation}
then $ C\hat x = \hat z$ and thus $(\hat x,\hat z)$ minimizes \eqref{problem_pdhg_1}.
Therefore, we can instead search for a saddle point of the augmented Lagrangian \eqref{eq:AugLagrangian} in order to solve the problem \eqref{problem_pdhg_1}.

One approach for finding a saddle point of \eqref{eq:AugLagrangian} is to alternate the minimization with respect to $(x,z)$ and to apply a gradient ascent step with stepsize $\nu>0$ for $p$.
This leads to the \emphii{alternating direction method of multipliers (ADMM)}\index{alternating direction method of multipliers}  \cite{Ga83,GM76} with updates
\begin{equation}\label{eq:uzawa_sep}
\begin{aligned}
x_{k+1}&\in \argmin_{x \in \R^N} L_\gamma(x,z_{k},p_{k})\\
z_{k+1}&\in \argmin_{z \in \R^L} L_\gamma(x_{k+1},z,p_{k})\\
p_{k+1} &=p_{k} + \nu (C x_{k+1} - z_{k+1}).
\end{aligned} 
\end{equation}
By setting $b_{k} \coloneqq p_{k}/\nu$ and $\nu = \gamma$, we obtain a (scaled) ADMM version.

\begin{algorithm}[H]
\KwIn{coercive, lsc $g\colon \R^N\to \overline{\R}$ and $h\colon \R^L \to \overline{\R}$, matrix $C \in \R^{L \times N}$, penalty $\gamma>0$, initializations $z_{0}, b_{0} \in \R^L$}
\For{\color{black}$k = 0, 1, \ldots$}{
$x_{k+1}\in \argmin_{x \in \bb X} \left\{ g(x) +  \frac{\gamma}{2} \|b_{k} + Cx - z_{k}\|_2^2  \right\}$\;
$z_{k+1}\in \prox_{\frac{1}{\gamma} h}(b_{k} + Cx_{k+1})$\;
$b_{k+1}= b_{k} + C x_{k+1} -  z_{k+1}$\;
}
\caption{Alternating Direction Method of Multipliers (ADMM) \label{alg:ADMM}}
\end{algorithm}

Overviews on the ADMM and its applications are given in \cite{BPCPE11,lin2022alternating}, where in particular also the (important) choice of $\gamma >0$ is addressed.

\begin{thm}{Convergence of ADMM \cite{lin2022alternating}}{convegence ADMM}
	Assume that $f\colon \R^N\to \overline{\R}$ with $f = g + h(Cx)$, where $g\colon \R^N\to \overline{\R}$ and $h\colon \R^L\to \overline{\R}$ are coercive, lsc and convex.
	Denote by $(\hat x, \hat z, \hat p)$ the saddle point of \eqref{eq:AugLagrangian} and let $\{(x_k, z_k)\}_{k \in \N}$ be generated by \Cref{alg:ADMM}.
	Then, it holds that
	\begin{equation}
		|g(x_{k}) + h(z_{k}) - g(\hat x) - h(\hat z)| = \q O(1/k), \quad \|Cx_{k} - z_{k}\|_2^2 = \q O(1/k).
	\end{equation}
	Furthermore, for $\tilde x_{k+1} = \frac{1}{k+1} \sum_{j=1}^{k+1} x_{j}$ and $\tilde z_{k+1} = \frac{1}{k+1} \sum_{j=1}^{k+1} z_{j}$ we have
	\begin{equation}
	|g(\tilde x_{k}) + h(\tilde z_{k}) - g(\hat x) - h(\hat z)| = \q O(1/k^2), \quad \|C\tilde x_{k} - \tilde z_{k}\|_2^2 = \q O(1/k^2).
	\end{equation}
\end{thm}

\begin{remark}[Nonconvex Case and Extensions]\label{rem: ADMM}
	Convergence of ADMM for nonconvex $f$ can be established.
	This requires that the Lagrangian $L_\gamma$ satisfies the KL condition and that either $h$ or $g$ is differentiable with $L$-Lipschitz gradient.
	Under these conditions, the convergence is guaranteed for $\gamma > L$ if $g$ is convex and $\gamma > 2L$ if $g$ is nonconvex \cite{lin2022alternating,ThePat2020,WanYinZen2019}.
	Given that we have to require differentiability of one summand, it might seem odd to use ADMM with its implicit update step instead of FBS with the explicit gradient update.
	However, similar as for (stiff) ordinary differential equations (see also \Cref{rem:GradFlow}), implicit updates are sometimes preferable even for the data term $\frac12\|Ax-y\|_2^2$.
	
	In practice, one often uses  iterative methods for solving the first step in \Cref{alg:ADMM}.
	A convergence analysis for this inexact setting is carried out in \cite{He2002,WanYinZen2019}.
	Moreover, we can vary the parameter $\gamma$ in each step \cite{He2000,Kontogiorgis1998}.
	In the convex case, global convergence rates can be found in \cite{hong2017linear,NisLesRec2015} (namely linear convergence for specific parameters) and local rates for quadratic programs are studied in \cite{boley2013local}.
\end{remark}

Even in the convex case, if the first step of \Cref{alg:ADMM} has more than one solution, the generated sequence $\{x_k\}_{k \in \N}$ might diverge as the following example shows.
This is no contradiction to \Cref{th:convegence ADMM} as it only states the convergence of $\|Cx_{k} - z_{k}\|$.

\begin{example}
Let $g(x)=0$ on $\R^{N}$, $h\colon \R \to \overline{\R}$ with $h(0) = 0$ and $h(x)=\infty$ else, and let $C\in\R^{1 \times N}$ be the zero-matrix.
Then, the steps of the ADMM \Cref{alg:ADMM} read as
\begin{equation}
x_{k+1}\in\R^N, \quad z_{k+1}=0, \quad p_{k+1}=p_{k}.
\end{equation}
Choosing $x_{k}=(-1)^k$ leads to a divergent sequence generated by the ADMM.
\end{example}

We now specify \Cref{alg:ADMM} for the sparsity promoting functionals from \Cref{ex: sparsity promoting}.

\begin{example}\label{example:TV_for_inpainting_ADMM}
Let $A\in\R^{M \times N}$ and $B\in\R^{L \times N}$ and consider the functional 
\begin{equation}
	\q J(x) =\tfrac12\|Ax-y\|_2^2+ \alpha \|Bx\|_1.
\end{equation}
In order to apply \Cref{alg:ADMM}, we choose $g(x)=\frac12\|Ax-y\|_2^2$, $h(x)= \alpha \|x\|_1$ and $C = B$.
Then, the first step of \Cref{alg:ADMM} reads as
\begin{align}
x_{k+1}&\in \argmin_{x \in \R^N} \left\{ \tfrac12\|Ax-y\|_2^2 +  \tfrac{\gamma}{2} \|b_{k} + Bx - y_{k}\|_2^2  \right\}.
\end{align}
By Fermat's rule, we obtain that $x_{k+1}$ has to solve the linear system
\begin{equation}\label{eq:OptSys}
	(A^\tT A+ \gamma B^\tT B)x = A^\tT y + \gamma B^\tT(y_{k}-b_{k}).
\end{equation}
Since the system \eqref{eq:OptSys} is symmetric positive definite, we can use the conjugate gradient method (as introduced in any lecture on numerical methods).
The second step reads as
\begin{equation}
z_{k+1}=\prox_{\frac{1}{\gamma} h}(b_{k} + Bx_{k+1})=S_{\frac{\alpha}{\gamma}}(b_{k}+ B x_{k+1}),
\end{equation}
where $S_{\frac1\gamma}$ is the soft-shrinkage from \Cref{exampleB:proxofl1norm}.
Finally, the last step is given by
\begin{equation}
	b_{k+1}= b_{k} + B x_{k+1} -  z_{k+1}.
\end{equation}
To summarize, all steps of \Cref{alg:ADMM} are straight forward to implement.
\end{example}

\begin{remark}[Parallel ADMM]\label{rem:parallel_ADMM}
	With ADMM we can also  tackle minimization problems of the form
	\begin{equation}
		f(x)=g(x)+\sum_{i=1}^n h_i(C_i x).
	\end{equation}
	To this end, we rewrite $f$ as $f(x)=g(x)+H(C x)$ with
	\begin{equation}
		H(z^1,...,z^n) = \sum_{i=1}^n h_i(z^i)\quad\text{and}\quad C=\left(\begin{array}{c}C_1\\\vdots\\C_n\end{array}\right).
	\end{equation}
	Now, we apply the ADMM on the functions $g$ and $H$.
	For
	$x_{k}$, $z_{k}=(z^1_{k},...,z^n_{k})$ and $b_{k}=(b^1_{k},...,b^n_{k})$, this results in the iterations 
	\begin{equation}
		\begin{aligned}
		x_{k+1}&\in\argmin_{x \in \bb X}\Bigl\{ g(x) +  \frac{\gamma}{2}\sum_{i=1}^n \|b^i_{k} + C_i x - z^i_{k}\|_2^2  \Bigr\}\\
		z^i_{k+1}& \in \prox_{\frac{1}{\gamma}h_i}(b^i_{k}+C_i x_{k+1}),\qquad\text{for }i=1,...,n,\\
		b^i_{k+1}&= b^i_{k}+C_i x_{k+1}-z^i_{k+1},\qquad\text{for }i=1,...,n.
	\end{aligned}
	\end{equation}
\end{remark}

The observations in \Cref{rem:parallel_ADMM} can be applied for simplifying the calculations in \Cref{example:TV_for_inpainting_ADMM}.
Still, this does not remove the need of solving a linear system.
Moreover, we can use \Cref{rem:parallel_ADMM} in order to apply the ADMM for minimizing the TV regularized functional from \Cref{ex: sparsity promoting} for the inpainting problem in \Cref{exer:impainting} with its non-differentiable data-fidelity term.

\subsubsection{Exercises}

\begin{exerciseB}\label{exercise:saddle_point_augmented}
Let $f\colon \R^N \to \R$ be given by $f=  g +h(Cx)$, where $g,h$ are convex and lsc and $C\in\R^{L \times N}$.
The associated  \textsc{Lagrangian} and augmented \textsc{Lagrangian} are
\begin{align}
L(x,z,p)=g(x)+h(z)+\lsk p,Cx-z\rsk,\quad L_\gamma(x,z,p)=L(x,z,p)+\frac\gamma2\|Cx-z\|^2,
\end{align}
respectively.
Show that any saddle point $(\hat x,\hat z,\hat p)$ of $L$ is a saddle point of $L_\gamma$.
\end{exerciseB}

\begin{Solution}
If $L(\hat x,\hat z,\hat p)=\infty$, nothing happens. Assume that $L(\hat x,\hat z,\hat p)\in\R$.

Let $(\hat x,\hat z,\hat p)$ be a saddle point of $L$. If it holds $C\hat x\neq \hat z$, we have that
\begin{align*}
L(\hat x,\hat z,\hat p)&=g(\hat x)+h(\hat z)+\lsk \hat \hat p,C\hat x-\hat z\rsk< g(\hat x)+h(\hat z)+\lsk \hat p,C\hat x-\hat z\rsk+\|C\hat x-\hat z\|^2\\
&=g(\hat x)+h(\hat z)+\lsk \hat p+Cx-z,Cx-z\rsk=L(\hat x,\hat z,\hat p+C\hat x-\hat z)
\end{align*}
which contradicts the saddle point property of $(\hat x,\hat z, \hat p)$. Thus, it holds $C\hat x=\hat z$
such that $L(\hat x,\hat z,p)=L_\gamma(\hat x,\hat z,p)$ for all $p$.
In particular, we have
$$
\sup_p L_\gamma(\hat x,\hat z,p)=\sup_p L(\hat x,\hat z,p)=L(\hat x,\hat z,\hat p)=L_\gamma(\hat x,\hat z,\hat p).
$$
Thus, it remains to show that $\inf_{x,z}L_\gamma(x,z,\hat p)=L_\gamma(\hat x,\hat z,\hat p)$.
To show this, we note that
$$
\min_{x,z}\sup_p L_\gamma(x,z,p)\leq \sup_p L_\gamma(\hat x,\hat z,p)=\sup_p L(\hat x,\hat z,p)=L(\hat x,\hat z,\hat p).
$$
Using $L\leq L_\gamma$ and \Cref{lem:saddle_minmax}, we obtain that
$$
L(\hat x,\hat z,\hat p)\geq \min_{x,z}\sup_p L_\gamma(x,z,p)\geq\max_p\inf_{x,z}L_\gamma(x,z,p)\geq \max_p\inf_{x,z}L(x,z,p)=L(\hat x,\hat z,\hat p).
$$
Since the value on the left and the right side the same, we have equality in the above equation.
Thus, we obatin
\begin{equation}
L(\hat x,\hat z,\hat p)=\inf_{x,z}L(x,z,\hat p)\leq\inf_{x,z} L_\gamma(x,z,\hat p)\leq \max_p \inf_{x,z} L_\gamma(x,z,p)=L(\hat x,\hat z,\hat p)=L_\gamma(\hat x,\hat z,\hat p)=L(\hat x,\hat z,\hat p).
\end{equation}
Because we have again equality in the last formula, this yields that $(\hat x,\hat z,\hat p)$ is a saddle point of $L_\gamma$.
\end{Solution}

\begin{exerciseB}[Inpainting]\label{exer:impainting}
	Let us consider inpainting. That is, for a vector $x_* \in \R^{N}$, we measure only entries $x_j$ for $j$ in some index set $\q J \subseteq \{1, \ldots, N\}$.
	This leads to the linear inverse problem $y = D x_*$, where $D$ is diagonal with elements 
	\begin{equation}
	d_{j} = \begin{cases}
		1, & j \in \q J, \\
		0, & \text{otherwise}.
	\end{cases}
	\end{equation}
	For each pixel $j = 1,\ldots, N$, we include $j$ in $\q J$ with probability $0.3$. 
	As ground truth $x_*$, we use the cameraman image, available in the Python package \texttt{skimage} as \texttt{data.camera()}.
	
	\begin{enumerate}
		\item[(i)] Argue why the Tikhonov and LASSO regularization are bad choices for this problem?
		
		\item[(ii)] Recall in \Cref{ex: circulant matrix} we showed that the discrete derivative $(\nabla x)_j = x_{j+1} - x_{j}$ from \Cref{ex: smoothness} can viewed as convolution $\nabla x = x *_N k $ with kernel $k = e_{-1} - e_0$.
		Use this observation to implement the directional derivatives
		\begin{equation}
		(\nabla_1 \, x)_{k,j} = x_{k+1,j} - x_{k,j} 
		\quad \text{and} \quad 
		(\nabla_2 \, x)_{k,j} = x_{k,j+1} - x_{k,j}
		\end{equation}
		and their adjoints with the function \texttt{convolve(img,k, mode = 'constant')} from the Python package \texttt{ndimage}.
		We use \texttt{constant} padding to account for the nonperiodic nature of $x_*$.
		Test your implementation by checking that 
		\begin{equation}
		\langle \nabla_\ell \, u, v \rangle = \langle u, \nabla_\ell^T \, v \rangle, \quad \ell = 1,2
		\end{equation}
		for random $u$ and $v$.  
		
		\item[(iii)] Implement accelerated gradient descent (\cref{algo:accelerated_gradient_descent}) for the variational reconstruction objective with the smoothness prior from \Cref{ex: smoothness}, namely
		\begin{equation}
		J_\alpha(x) = \|Dx - y\|^2_2 + \alpha (\| \nabla_1 x\|^2_2 + \| \nabla_2 x \|^2_2).
		\end{equation}
		Visualize the reconstruction for $\alpha = 0.1$.  
		
		\item[(iv)] Now, we want to use the TV regularization from \Cref{ex: sparsity promoting}, namely
		\begin{equation}\label{eq:InpaintingTV}
		J_\alpha(x) = \|Dx - y\|^2_2 + \alpha (\| \nabla_1 x\|_1 + \| \nabla_2 x \|_1).
		\end{equation}
		Based on \Cref{example:TV_for_inpainting_ADMM} and \Cref{rem:parallel_ADMM}, implement the ADMM with $\gamma = 1$ to minimize \eqref{eq:InpaintingTV} with $\alpha = 10$.
		For the least squares sub-iteration, use the conjugate gradient method (or your implementation of gradient descent) with $10$ iterations and initialization $x_t$ from the previous ADMM iteration. 
		Why are 10 iterations sufficient?
		
		\item[(v)] Plot the absolute errors $| D(x - x_*)|$.
		For inpainting, we want to enforce that the error on the measured pixels is zero. This is done via the functional
		\begin{equation}\label{eq:InpaintingTV2}
		J(x) = \iota_{ \{z \colon Dz = y\}}(x) + \| \nabla_1 x\|_1 + \| \nabla_2 x \|_1,
		\end{equation}
		where 
		\begin{equation}
			\iota_{\q S}(x) = \begin{cases}
			0, & x \in \q S,\\
			\infty, & x \notin \q S.
		\end{cases}
		\end{equation}
		Implement ADMM with $\gamma = 0.01$ for minimizing \eqref{eq:InpaintingTV2} and perform the reconstruction.
		Plot the absolute errors $| D(x - x_*) |$ again and explain the difference.
		
		\item[(vi)] Now, we change the sampling of $\q J$.
		Instead of random  pixels, we exclude a $50 \times 50$ box that contains the cameraman's face.
		Try to reconstruct it with the algorithms from above and explain the outcome.
	\end{enumerate}
\end{exerciseB}

\begin{Solution}
	See Jupyter Notebook \url{https://github.com/Sniper2k/inverse_problems_ai_exercises/blob/main/Ex4_3_38.ipynb}.
\end{Solution}

\subsection{Kurdyka-{\L}ojasiewicz Condition} \label{sec: KL}

As mentioned in the previous sections, the Kurdyka-{\L}ojasiewicz condition plays a crucial role in the convergence analysis of algorithms.
This section provides a brief overview of the relevant literature and adjacent topics without going into too much detail.
Let us start with its formal definition.

\begin{defn}
A lsc function $f \colon \R^{N} \to \overline{\R}$ satisfies the \emphi{Kurdyka-{\L}ojasiewicz (KL) condition} at $\hat x \in \q D(f)$ if there exists a $\eta >0$, a neighborhood $U$ of $\hat x$, and a continuous concave $\varphi\colon [0,\eta) \to \R^+$ such that 
\begin{equation}
\varphi(0) = 0,
\quad \varphi \in C^1((0,\eta)),
\quad \varphi(s) > 0 \text{ for all } s\in(0,\eta),
\end{equation}
and for all $x \in U$ with $f(\hat x) < f(x) < f(\hat x) + \eta$ the Kurdyka-{\L}ojasiewicz inequality holds
\begin{equation}\label{eq: KL inequality}
\varphi'(f(x) - f(\hat x)) \operatorname{dist}(0, \partial f(x)) \ge 1.
\end{equation}
We call $f$ a \emphi{KL function} if the  KL condition is met on $\q D(\partial f) \coloneqq \{ z \in \R^N \colon \partial f(z) \neq \varnothing \}$.
\end{defn}

In \eqref{eq: KL inequality}, $\operatorname{dist}(0, \partial f(x))
\coloneqq \inf_{z \in \partial f(x)} \|z\|$
is the smallest norm of the elements in the subdifferential \eqref{eq: subdifferential}.
If $f$ is continuously differentiable at $x$, then $\operatorname{dist}(0, \partial f(x)) = \| \nabla f(x) \|$. 

The inequality \eqref{eq: KL inequality} controls how fast $\nabla f(x)$ diminishes as $x$ approaches $\hat x$ and the function $\varphi'$ controls this rate.
For $\hat x$ with $0 \notin \partial f(\hat x)$, the KL property is always satisfied \cite[ Lem.\ 2.1]{attouch2010proximal}.
In this case, $\varphi$ can be chosen as
\begin{equation}\label{eq: KL Lojasiewicz}
	\varphi(s) = c s^{1- \theta}, \quad 0 \le \theta < 1,
\end{equation}
see \cite[Lem.~2.1]{li2018calculus}, where we refer to $\theta$ as the \emph{KL exponent}. 
Hence, the difficulty lies in establishing the KL condition for critical points.  

\begin{example}\label{ex: kl strongly convex}
	Let $f\in C^1(\R^N)$ be $\mu$-strongly convex with $L$-Lipschitz continuous gradient.
	By \Cref{th:minima} it has a unique global minimizer $\hat x$ with $\nabla f(\hat x) = 0$.
	Due to \Cref{l: descent lemma} we have for all $x \in \R^N$ that
	\begin{equation}\label{eq:GrowthEst}
	f(x) - f(\hat x) 
	\le \inner{\nabla f(\hat x)}{x - \hat x} + \tfrac{L}{2} \norm{x - \hat x}^2.
	\end{equation}
	Furthermore, by \Cref{ex: strong convexity first order}, we get
	\begin{equation}\label{eq:CoercivityEst}
	\| \nabla f(x) \| 
	= \| \nabla f(x) - \nabla f(\hat x)\|
	\ge \mu \norm{x - \hat x}.
	\end{equation}
	Combining \eqref{eq:GrowthEst} and \eqref{eq:CoercivityEst} gives
	\begin{equation}
	f(x) - f(\hat x) 
	\le \tfrac{L\mu^2}{2\mu^2} \norm{x - \hat x}^2
	\le \tfrac{L}{2\mu^2}\| \nabla f(x) \|^2
	\end{equation}
	or, equivalently,
	\begin{equation}
	\tfrac{\sqrt{L}}{\sqrt 2 \mu} (f(x) - f(\hat x))^{-1/2} \| \nabla f(x) \| \ge 1.
	\end{equation}
	Hence, $f$ satisfies the KL inequality at $\hat x$ for all $\eta> 0$ and $U \subseteq \R^N$ for $\varphi(s) = c s^{1- 1/2}$ with $c = \sqrt{2 L} / \mu$ and $\varphi'(s) = \tfrac{1}{2}c s^{-1/2}$. 
\end{example} 

The importance of the KL property is summarized by \cite[Prop.\ 6]{bolte2016majorization}, which is given in a simplified form below.

\begin{thm}{Local convergence under KL property}{KL convergence}
	Let $f \colon \R^N \to \overline{\R}$ with $\q D(f) \neq \emptyset$ be lsc and bounded from below.
	Assume that  a bounded sequence $\{x_k\}_{k \in \N}$ satisfies 
	\begin{align}\label{eq:SufficientDecrease}
	f(x_{k+1}) 
	& \le f(x_k) - c_1 \| x_{k+1} - x_k \|^2, \\
	\operatorname{dist}(0, \partial f(x_k)) \label{eq:Sharpness}
	& \le c_2 \| x_{k+1} - x_k \|
	\end{align}
	for constants $c_1, c_2 > 0$, and  that $f(x_{k_j}) \to f(x)$ if $x_{k_j} \to x$.
	If $f$ is a KL function, then $\lim_{k \to \infty} x_k = \hat x$.
	Furthermore, if $\varphi$ is of the form $\varphi(s) = c s^{1- \theta}$ with $c>0$ and $0 < \theta < 1$, then there exist $k_0 \in \N$, $c>0$ and $0 \le q < 1$ such that for all $k \ge k_0$ it holds 
	\begin{equation}\label{eq:Rates}
	\| x_k - \hat x \| \le 
	\begin{cases}
		c q^k, & \theta \in [0, 1/2], \\
		c k^{-(1 -\theta) / (2 \theta - 1)}, & \theta \in (1/2, 1).
	\end{cases}
	\end{equation}
\end{thm}
In \Cref{ex: grad descent KL}, we verify that the sequence generated by gradient descent (\Cref{algo:gradient_descent}) satisfies \eqref{eq:SufficientDecrease} and \eqref{eq:Sharpness}.
A similar result for the FBS (\Cref{algo:FBS}) applied to convex functions can be found in \cite{bolte2017error}. 

Note that the convergence rates in \eqref{eq:Rates} together with our observations for strongly convex functions in \Cref{ex: kl strongly convex} lead to the same linear convergence rate as in \Cref{th:grad descent linear}. 
Nevertheless, \Cref{th:KL convergence} applies to  a much broader class of functions.
Describing it requires the definition of o-minimal structures, see also \cite[Def.\ 4.1]{attouch2010proximal}.
A collection $\q O = \{ \q O_n\}_{n \in \N}$ with Boolean algebras $\q O_n$ of subsets of $\R^n$ is an o-minimal structure if
\begin{itemize}
	\item each set in $\q O_1$ is a finite union of intervals and points;
	\item $\q O_n$ contains all the sets $\{x \in \R^n: p(x) = 0\}$ and $\{x \in \R^n: p(x) < 0\}$, where $p$ is a polynomial on $\R^n$;
	\item $U_n \in \q O_n$ and $U_m \in \q O_m$ imply $U_n \times U_m \in \q O_{n+m}$;
	\item $\pi (U_n) = \{ (x_1, \ldots, x_{n-1}) \in \R^n \, \colon \, x \in U_n \} \in \q O_{n-1}$ for all $U_n \in \q O_n$.
\end{itemize}
A set $U \subseteq \R^N$ is \emph{definable} in an o-minimal structure $\q O$ if $U \in \mathcal O_N$.

\begin{thm}{}{definable structure}
There exists o-minimal structures $\q O_{sa}$ and $\q O_{exp}$ such that
\begin{itemize}
	\item[\circled{1}] all semialgebaric sets
	\begin{equation}\label{eq:SemiAlgebraicSets}
		\bigcup_{\ell = 1}^L \{ x \in \R^N \,\colon\, p_{n,k}(x) = 0, \ q_{n,k}(x) < 0, \ k = 1, \ldots, K\}
	\end{equation}
	with polynomials $p_{n,k}$, $q_{n,k}$ are definable in $\q O_{sa}$. 
	\item[\circled{2}] $\{(x,y) \in \R^2: \exp(x) = y\}$ and all the sets from \circled{1} are definable in $\q O_{exp}$.
\end{itemize}
\end{thm}

A detailed discussion regarding $\circled{1}$ is given in \cite[Sec.\ 1 and 2]{van1998tame}.
Item $\circled{2}$ is due to a result by Wilkie \cite{wilkie1996model} and its extension \cite{van1994real}.
A list of other known o-minimal structures can be found in \cite{Kurdyka1998}. 
The functions $F \colon \R^N \to \R^M$ and $f \colon \R^N \to \overline{\R}$ are \emphi{definable} if their graphs $\{ (x, y) \in \R^N \times \R^M \colon F(x) = y\}$ and $\{ (x,t) \in \R^N \times \R \colon f(x) = t\}$ are definable. 
As we see in the following theorem, the set of definable functions is closed under many common operations \cite[App.\ B]{van1996geometric}.

\begin{thm}{Operations with definable functions}{definable operations}
	Let $K \in \N$ and let $f, f_1, \ldots, f_K$ be definable functions in an o-minimal structure $\q O$. Then, the following holds:
	\begin{itemize}
		\item[\circled{1}] $f_1 + f_2$, $\alpha f$, $\alpha \in \R$, are definable in $\q O$;
		\item[\circled{2}] $f_1 \circ f_2$ is definable in $\q O$;
		\item[\circled{3}] $f = (f_1, \ldots, f_K)$ is definable in $\q O$ if and only if $f_1, \ldots, f_k$ are;
		\item[\circled{4}] if $f \in C^1(\R)$, then  $\nabla f$ is definable in $\q O$;
		\item[\circled{5}] if $\q R(f)$ is definable, then the (generalized) inverse $f^{-1}$ is definable;
		\item[\circled{6}] if $U \subseteq \R^N$ is definable in $\q O$, then $\iota_U$ is definable in $\q O$.
		%		\item[\circled{5}] $\min_{y \in \R^M} g(x,y)$, $\max_{y \in \R^M} g(x,y)$  are definable in $\q O$ if $g \colon \R^N \times \R^M \to \overline{\R}$.
	\end{itemize}  
\end{thm}

\begin{example}\label{ex: definable functions}
Most functions appearing in variational regularization are definable.
\begin{itemize}
	\item The functionals $\| \cdot \|_2^2$ and $\| \cdot \|_1$ are piecewise polynomial and thus definable in $\q O_{sa}$.
	Hence, by \Cref{th:definable operations} all convex functionals in the examples from \Cref{sec: variational regularization} are definable in $\q O_{sa}$.
	\item For any $n \in \N$, the functions
	\begin{equation}
	f_n(t) = \begin{cases}
		t^n, & t \ge 0,\\
		\infty & t < 0.
		\end{cases},
	\quad \text{and} \quad 
	f_n^{-1}(t) = \begin{cases}
		t^{-n}, & t \ge 0,\\
		\infty & t < 0.
	\end{cases}, 
	\end{equation}
	are definable in $\q O_{sa}$.
	Then, for $p > 0$ and $p \in \bb Q$, the data fidelities $\| Ax - y \|_2^p$ and the (quasi)norms $\| x \|_p$ are definable in $\q O_{sa}$. For all real $p>0$, they are definable in $\q O_{exp}$, see \Cref{real powers}.  
	\item The log-sum penalty $\sum_{k=1}^N \log(|x_k| + \varepsilon)$ with $\varepsilon \ge 0$ is definable in $\q O_{exp}$.
\end{itemize}  
\end{example}

The concept of o-minimal structures enables us to characterize functions that satisfy the KL condition, see \cite{lojasiewicz1993geometrie} and  \cite[Thm.\ 11]{bolte2007clarke}, namely functions for which  \Cref{th:KL convergence} applies. 

\begin{thm}{}{definability gives KL}
	Let $f\colon \R^N \to \overline{\R}$ be lsc and definable in an o-minimal structure $\q O$ with $\q D(f) \neq \varnothing$.
	Then, it satisfies KL condition at all points $x \in \q D(\partial f)$.
	If $f$ is definable in $\q O_{sa}$, then it admits a KL exponent $\theta \in [0,1) \cap \bb Q$.
\end{thm} 

By combining \Cref{ex: definable functions} and \Cref{th:definability gives KL}, all of the functions that we considered so far have the KL property.
For functions that are definable in $\q O_{exp}$, we cannot guarantee the existence of the KL exponent, see \Cref{ex: no KL exponent}.

\begin{thm}{\cite[Cor.\ 9]{bolte2017error}}{}
Let $f\colon \R^N \to \R$ with $\argmin(f) \neq \emptyset$ be convex and piecewise polynomial with maximal degree $d$.
Then, $f$ is a KL function with exponent $\theta = 1 - 1/((d - 1)^N + 1)$. 
\end{thm}

Consequently, the LASSO functional from \Cref{ex: LASSO} has the KL exponent $\theta = 1/2$.
For this functional, an explicit upper bound for the constant $c$ in \eqref{eq: KL Lojasiewicz} is established in \cite{bolte2017error}. 

\subsubsection{Exercises}

\begin{exerciseB}\label{ex: grad descent KL}
Let $f \in C^1(\R^N)$ with $L$-Lipschitz gradient satisfy the KL condition with exponent $0 < \theta < 1$.
Consider a sequence $\{x_k\}_{k \in \N}$ generated by \Cref{algo:gradient_descent} with step size $\nu \le 1/L$.
Show that the conditions of \Cref{th:KL convergence} are satisfied and conclude that $\{x_k\}_{k \in \N}$ converges. 
\end{exerciseB}
\begin{Solution}
	By construction, we have $x_{k+1} - x_k = \nu \nabla f(x_k)$ giving the second condition of \Cref{th:KL convergence}. The first and the third conditions follow by \Cref{th:grad descent nonconvex}.
\end{Solution}

\begin{exerciseB}
	Verify that the functions in \Cref{ex: definable functions} are definable. 
\end{exerciseB}
\begin{Solution}
	2) The function $f$ is definable, since its graph can be written as 
	\[\{ (t,y) \in \R^2 \, \colon \, t^n - y = 0,\ -t < 0 \} \cup \{ (t,y) \in \R^2 \, \colon \,  t^2 + y^2 = 0 \}.
	\] 
	Since 
	\[
	\q R(f) = [0, +\infty) = \{ y \in \R \, \colon \, -y < 0 \} \cup \{ y \in \R \, \colon \, y = 0 \}
	\] 
	is definable in $\q O_{sa}$, then $f^{-1}$ is definable in $\q O_{sa}$ by \Cref{th:definable operations} \circled{5}. 
	
	Let $p \in \bb Q$, $p > 0$. Then, we can find $n, m \in \N$ such that $p/2 = n/m$. By \Cref{th:definable operations} \circled{2}, $f_{n} \circ f_{-m}$ is definable in $\q O_{sa}$. Since $\| Ax - y \|^p = (\| Ax - y \|^2 )^{p/2} = f_{n} \circ f_{-m}(\| Ax - y \|^2)$ is again definable as superposition. 
	
	For $\| x \|_p$, we notice that $\| x\|_p^p = \sum_{k =1}^{N} |x_k|^p$, where each summand is definable in $\q O_{sa}$ as superposition of $| \cdot |$, ${\,}^p$ and a sum of definable functions, \Cref{th:definable operations} \circled{1}. 
	
	3) Since $(0,+\infty) = \{ t \in \R \, \colon \, -t < 0 \}$ is definable, $\log(t)$ is definable in $\q O_{exp}$ as the inverse of $e^{t}$. The rest follows analogously the the previous part.
\end{Solution}

\begin{exerciseB}\label{real powers}
%Let $\q O_{exp}$ be an o-minimal structure from \cref{th:definable structure}. 
Show that $t \mapsto t^r$ is definable in $\q O_{exp}$ for any $r > 0$ using \Cref{th:definable operations}.  
\end{exerciseB}
\begin{Solution}
	We can present $t^r = \exp( r \log(t))$, which is definable in $\q O_{exp}$ by \Cref{th:definable operations} \circled{2}. 
\end{Solution}

\begin{exerciseB}\label{ex: no KL exponent}
%Let $\q O_{exp}$ be an o-minimal structure from \cref{th:definable structure}. 
Show that $t \mapsto e^{-1/t^2}$ is definable in $\q O_{exp}$ and that it does not admit a KL exponent $0 \le \theta < 1$. 
\end{exerciseB}
\begin{Solution}
In analogy to \Cref{ex: definable functions}, the function
\[
g(t) = \begin{cases}
	-1/t^2 & t \neq 0, \\
	\infty & t = 0
\end{cases}
\]
is definable in $\q O_{sa}$ and $\q O_{exp}$. Then $f(t) = e^{-1/t^2}$ is definable in $\q O_{exp}$. It is a $C^1$ function with derivative $2 t^{-3}e^{-1/t^2}$, $t \neq 0$ that achieves a global minimum at $0$. Assume $f$ admits KL exponent $\theta \in [0,1)$ at $t = 0$. Then, for $t \neq 0 $
\[
1 \le c (1 -\theta) ( e^{-1/t^2} - 0 )^{-\theta} \cdot 2 t^{-3}e^{-1/t^2} = 2 c (1 - \theta) t^{-3} e^{-(1-\theta)/t^2}.
\]
By L'H\^{o}pital's rule, the limit of the right-hand side is zero, giving a contradiction.  
\end{Solution}

\subsection{Outlook}

This brief overview focuses solely on algorithms relevant to our applications.
For a comprehensive treatment of the field, we refer the reader to specialized courses.
Other important algorithms not covered here include
\begin{itemize}
	\item majorize-minimize algorithms based on the construction of majorizing functions as surrogate that are minimized at each iteration instead.
	Similarly to ADMM, the resulting subproblems are designed to be easily solvable.
	Well-established methods in this category are iteratively reweighted Least Squares and the Difference-of-Convex algorithm;
	
	\item (quasi-)Newton methods, which incorporate second order information into the optimization process.
	This usually leads to much faster convergence.
	Due to the increased computational complexity and memory requirement, these are less popular for machine learning and imaging tasks;
	
	\item Douglas--Rachford splitting and the primal dual hybrid gradient algorithm as alternatives to the ADMM for nonsmooth optimization.
	
	\item Penalty and trust-region methods as well as sequential quadratic and semidefinite programming.
	These target constrained optimization problems. 
	Similarly to the second order methods, their iterations are expensive, which is why they are not frequently used in inverse problems. 
\end{itemize}
In principle, one should choose an algorithm for nonsmooth optimization in such a way that no inner loop for computing the involved proximal operators is required.
If this is not possible, then warm-starting the inner iterations with the previous outputs and some controlled inexactness can accelerate the solution process, see \Cref{rem:FBS} and \Cref{rem: ADMM}. 

For several inverse problems, we have to optimize loss functionals $f \colon \bb C^N \to \overline{\R}$.
Since these are nonholomorphic, the complex derivative is not suitable for optimization.
Instead, one treats $f$ as a function of real and imaginary parts, for which the Wirtinger derivatives 
\begin{equation}
	\frac{\partial f}{\partial z} \coloneqq \frac{\partial f}{2 \partial \textrm{Re} z} - \frac{i \partial f}{2 \partial \textrm{Im} z}
	\quad \text{and} \quad 
	\frac{\partial f}{\partial \overline z} \coloneqq \frac{\partial f}{2 \partial \textrm{Re} z} + \frac{i \partial f}{2 \partial \textrm{Im} z}
\end{equation}
are particularly convenient.
All presented iterative methods and their convergence guarantees can be extended to the complex case.

\section{Neural network-based reconstruction}\label{chap:regularization}
The inversion model \eqref{eq:data_fidelity_regularizer} is a promising approach, for which two key challenges remain:
\begin{itemize}
	\item how to \emph{find} a good regularizer $R$;
	\item computational \emph{feasibility} of hyper-parameter selection and evaluation.
\end{itemize}
Here, data-driven reconstruction architectures $\q R_{\theta} \colon \bb Y \to \bb X$ come into play.
They address these challenges by using training data to learn an \textit{optimal} reconstruction operator based on statistical concepts.  
To this end, we parameterize $\q R_{\theta}$ for example by a (deep) neural network or as a variational reconstruction model with a parametric regularizer $R_\theta$.
Then, learning refers to finding the optimal parameters $\widehat{\theta}$.
We will discuss several architectures for $\q R_{\theta}$.
The key inspiration for their design are the observations from the previous sections.

\subsection{Neural Network Architectures}\label{sec:NNArchitectures}
\begin{figure}[tbp]
	\centering
	\includegraphics[width=0.85\textwidth]{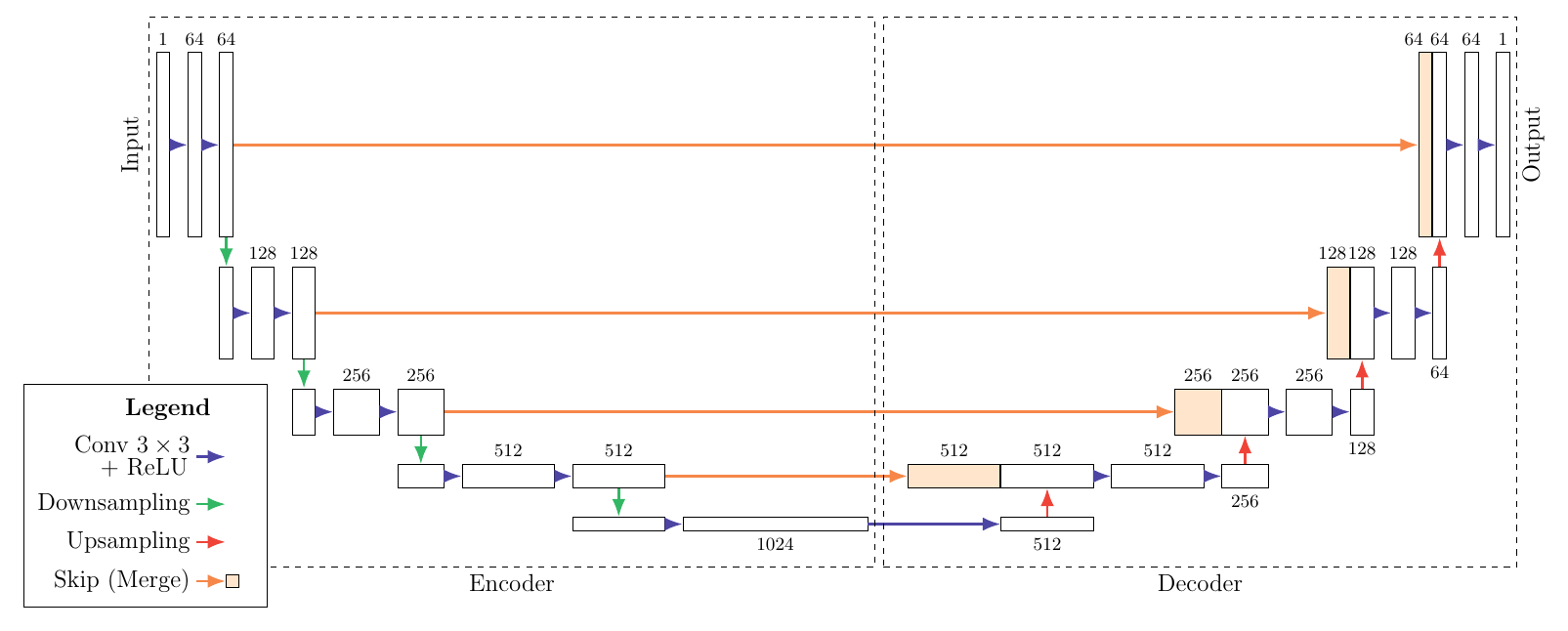}
	\caption{UNet \cite{RonFisBro2015}: The numbers refer to the respective dimension and the arrows indicate the processing direction.
	%The image size decreases due to the valid convolutions.
	}
	\label{fig:UNet}%\vspace{-0.15cm}
\end{figure}

We are mainly interested in processing images $x \in \R^{N_x \times N_y \times C}$, where $N_x$ and $N_y$ refer to the image dimensions and $C$ denotes the channels, see also  \Cref{ex: tomography}.
The channel dimension $C$ is one for gray-scale images and three for colour ones.
In the machine learning literature, such objects are commonly referred to as \emphi{tensors}.
\begin{defn}[Neural Network]\label{def:NN}
	A neural network $\Phi \colon  \R^{N_x \times N_y \times C} \to  \R^{N_x \times N_y \times C}$ is a sequence of $L$ \emph{layers} $h_k$, $k=1,\ldots,L$, and (simple) nonlinear \emph{activation functions}  $\{\sigma_i\}_{i=1}^{L-1}$ in between, namely 
	\begin{equation}
		\Phi (x) = h_L \circ \sigma_{L-1}\cdots \circ \sigma_1\circ h_1 (x).
	\end{equation}
	Here, $\circ$ denotes function composition. 
	Each layer $h_k$ is a parametric function with possibly different input and output dimensions.
	The union of all parameters is denoted by $\theta$.
\end{defn}
\begin{remark}
	In \Cref{def:NN}, we do not consider different input and output dimensions, as this would make the notation more complicated.
	Moreover, for (almost) all discussed methods, input and output dimension of the learnable components coincide.
	This is due to the fact that they act directly on image-like objects instead of the raw data.
\end{remark}
The \emph{architecture} of a neural network refers to the specification of its layers $\{h_k\}_{k=1}^L$ and activations $\{\sigma_k\}_{k=1}^{L-1}$.
In the following, we briefly discuss common components of architectures for imaging tasks, such as the UNet in \Cref{fig:UNet} or the DRUNet in \Cref{fig:DRUNet}.
More detailed expositions are given in any \emph{Machine Learning or Data Science} course.
A short overview of neural networks with visualizations can be found in this  \href{https://stanford.edu/~shervine/teaching/cs-230/cheatsheet-convolutional-neural-networks#layer}{cheatsheet}.\footnote{https://stanford.edu/$\sim$shervine/teaching/cs-230/cheatsheet-convolutional-neural-networks\#layer}

\paragraph{Fully-connected Layers}
Here, $h \colon  \R^{N_x \times N_y \times C} \to  \R^{N_x \times N_y \times K}$ treats the input as a vector and applies an affine transform $h(x) = Wx  + b$ with $W \in \R^{(N_x \times N_y \times K) \times (N_x \times N_y \times C)} $ and $b \in \R^{N_x \times N_y \times K}$.
Regarding imaging applications, fully-connected layers can capture dependency between all pixels, even if they are in the opposite corners of an image.  
Unfortunately, training deep architectures with only fully-connected layers has proven difficult due to the huge number of parameters and vanishing gradients.  
 
\paragraph{Convolution Layers}
These are the backbone of most modern imaging architectures.
A \emphi{convolution layer} extracts local features of the image based on the shift-invariant (discrete) convolution, see \Cref{ex: convolution}.
A convolution layer $h\colon  \R^{N_x \times N_y \times C} \to  \R^{N_x \times N_y \times K}$ with $C$ input channels and $K$ \emphi{filters} is an affine transformation $h = \{h_k\}_{k=1}^K$ with components
\begin{equation}
	h_{k}(x)= \sum_{c = 1}^{C} x_c * f_{k,c} + b_k
\end{equation}
with \emphi{bias terms} $b \in \R^{K}$ and \emph{filters} $f_{k,c} \in \R^{N_x \times N_y}$.
The most common choice is to use zero-padding such that input and output size coincide.
Other popular options are periodic extension (ensures compatibility with the Fourier transform) and the valid mode (we decrease the output size).
Typically, the support $d \times d$ of the filters $f_{k,c}$ is much smaller than $N_x \times N_y$, such as $3 \times 3$ for the UNet in \Cref{fig:UNet}. While $d \times d$ is a small \emphi{field-of-view}, stacking $L$ convolution layers increases it to $1 + (d-1)L$, motivating deep architectures. The learnable parameters $f$ and $b$ total to $(d^2C+1)K$ unknowns.
This is significantly less than for fully-connected layers. 
Still, $(d^2C+1)K$ can be very large.
Thus, the input dimension $C$ is often divided into smaller groups.
For example, with group size $2$, we use only two filters $f_{k,1}$ and $f_{k,2}$, which are applied first to channels $1$ and $2$, then reused for channels $3$ and $4$, and so on.
In the extreme case of group size $1$, this reduces the parameter count to $(d^2+1)K$.

\paragraph{Down- and Upsampling Layers}
In imaging, it is well known that features occur at different scales (an example hierarchy is savanna, zebra, ears, fur, insect), and methods such as the Wavelet transform reflect this in their design.
Many networks, such as the UNet in \Cref{fig:UNet}, are multiscale architectures with down- and upsampling operations.
To get a \emphi{2-fold downsampling} $h\colon  \R^{N_x \times N_y \times C} \to  \R^{\lceil N_x/2 \rceil \times \lceil N_y/2 \rceil \times C}$, we divide the image in small patches of size $2 \times 2 \times C$ (up to boundary handling).
For each channel and patch, we only keep the value in the top left corner.
In most libraries, this form of downsampling is integrated into the convolution layers as \emphi{stride} (which equals 2 here).
The downsampling is a linear operation and its adjoint (sometimes called transposed convolution) can be used for upsampling.
A popular nonlinear downsampling operation is \emphi{max-pooling}, where we instead keep the maximal value in each patch.

\paragraph{Residual/Skip Connections} Empirically, it was observed that with increasing depth $L$, neural networks become hard to train due to vanishing gradients.
To mitigate this issue, many networks contain \emphi{residual/skip connections} are introduced.
Usually, they are of the form $h_k(x) = x + h_j(x)$ for some $j< k$ or $h_k(x)$ merges $x$ and $h_j(x)$ along the channel dimension.
We can interpret these as a \textit{shortcut} in the network.
Taking another look at \Cref{fig:UNet}, we can also motivate them from a different perspective.
The UNet has a multiscale encoder-decoder structure (the image size gets smaller and then larger again), where information is processed on several scales.
The decoder part combines information from the coarser scale with the one from before the downsampling (to avoid loss of information).
Models of the form $\Phi(x) = x + h(x)$ with an input-output skip connection are particularly useful if we expect that $\Phi(x) \approx x$, as for the denoising networks presented in \Cref{sec:Denoising}.

\paragraph{Activation Functions}
The activation functions $\sigma_k \colon \R \to \R$ are usually applied componentwise to each entry of a channel and are frequently also shared across all channels.
Activation functions are typically not parametric and thus do not contribute to the learnable parameters.
Below, we provide a brief overview.
\begin{itemize}
	\item \textbf{Rectified Linear Unit (ReLU):} This is arguably the most popular choice and given by $x \mapsto \max(x,0)$.
	Unfortunately, ReLU is not differentiable.
	\item \textbf{Parametric ReLU}: This generalizes the ReLU as $x \mapsto \max(x,ax)$ with a tunable parameter $a\in \R$.
	The ReLU corresponds to $a=0$.
	\item \textbf{Softplus}: A smooth approximation of ReLU is given by $x \mapsto \log(1 + \exp(x))$.
	\item \textbf{Tanh}: Sometimes, bounded activation functions are preferable and we use the hyperbolic tangent with $x \mapsto (e^{x} - e^{-x}) / (e^{x} + e^{-x})$. 
	\item \textbf{Sigmoid}: Tanh looks similar to the nonnegative sigmoid with $x \mapsto 1 / (e^{-x} + 1)$.
	It is particularly popular in the early work on neural networks.
	\item \textbf{Identity}: For convenience of notation, we also allow the identity  $x \mapsto x$.
\end{itemize}
Sometimes, the choice of $\sigma_k$ can be motivated by theoretical considerations, see \Cref{sec:LearnedReg} and \Cref{sec: unrolling}.
More frequently, however, it is simply justified by trial and error.
Although not relevant towards our applications, we want to mention that there exist non-componentwise $\sigma_k$ such as softmax  (classification) and groupsort (Lipschitz-constrained architectures).

\subsection{Training Problem}\label{sec: training problem}
We exclusively discuss \emph{supervised} learning approaches, where one has pairs of ground-truth images $x$ and their corresponding measurements.
That is, we are given paired training data $(x_i,y_i)\in \bb X \times \bb Y$, $i=1,\ldots,S$, from some (unknown) data distribution.
For a known noise model, such pairs can be easily generated from given $x_i$ via $y_i = A x_i + n_i$.
However, such pairs usually do not reflect reality well and are considered as \emph{inverse crime}.
Based on the samples $(x_i,y_i)$, optimal parameters $\widehat{\theta}$ for $\q R_{\theta}$ are found by empirical risk minimization.
\begin{defn}[Training Problem]
	For training samples $(x_i,y_i)\in \bb X \times \bb Y$, $i=1,\ldots,S$, and a parameteric reconstruction model $\q R_{\theta} \colon \bb Y \to \bb X$, we define the supervised training problem
	\begin{equation}\label{eq:BayesEst}
		\widehat{\theta} \in \argmin_{\theta} 
		\frac{1}{2S}\sum_{i=1}^S
		\bigl \Vert \q R_{\theta}(y_i) - x_i \bigr \Vert_{\bb X}^2.
	\end{equation}
	Although the objective in \eqref{eq:BayesEst} does not explicitly include the forward $A$, it is implicitly incorporated by the assumption that $y_i \approx A x_i$.
\end{defn}
For simplicity, we stick with the squared $\ell^2$-norm as loss in \eqref{eq:BayesEst} and simply note that other choices such as the $\ell^1$-norm or more sophisticated visual perception losses are possible.
In practice, the influence on the obtained results is often small.
After optimizing the parameters $\theta$ based on \eqref{eq:BayesEst}, we expect to get \emph{good} reconstructions from new data $y$ as long as they follow the same distribution.
For out of distribution data $y$, the performance depends highly on the \emph{generalization} of the chosen model.
Generally, the training set $(x_i,y_i)\in \bb X \times \bb Y$, $i=1,\ldots,S$, plays an important role for the performance of the learned $\q R_{\theta}$, and we may observe over-fitting to the training data.
Hence, one commonly monitors the performance of $\q R_{\theta}$ on some validation set during training.
This set should follow the same distribution as the training data, but it must not be part of the samples chosen for \eqref{eq:BayesEst}.
By doing to, we check how well $\q R_{\theta}$ generalizes to these samples.
If everything works as expected, both the training and the validation loss should decrease when minimizing \eqref{eq:BayesEst}.

\begin{remark}
	Alternative training paradigms include \emph{unsupervised learning}, which does not require any groundtruth $x_i$, and \emph{weakly supervised learning}, where we have unpaired data $x_i$ and $y_i$.
	Both paradigms are particularly relevant if obtaining paired training data is prohibitively expensive.
	However, they often do not match the performance of supervised approaches.
	For more details in the context of inverse problems, we refer to \cite{mukherjee2023learned}.
\end{remark}

\subsection{Stochastic Optimization}\label{sec:OptimML}
Given an architecture $\q R_\theta$, we need a method to find a minimizer $\hat{\theta}$ of the loss \eqref{eq:BayesEst}.
In general, this is a highly nonconvex problem, and we can only hope to find local minima. 
For this, we have the large selection of optimization methods from \Cref{sec:Algs} at hand.
In machine learning, most people pursue a (sub)gradient based approach.
We restrict our discussion to \Cref{algo:gradient_descent}, which does not involve any momentum/acceleration steps.
To compute the gradient of \eqref{eq:BayesEst} with respect to $\theta$, we need to differentiate $\q R_\theta$. 
This relies on the chain rule and is commonly referred to as \emphi{backpropagation} and readily implemented in any machine learning package.

\begin{example}[Backpropagation]\label{ex: backprop}
	Here, we use backpropagation to compute the gradient of the training problem \eqref{eq:BayesEst} with a single data pair $(x,y) \in \R^N \times \R^M$ for a 2-layer network $\q R_\theta = \q R(\cdot \,; \theta) = \sigma_2 \circ h_2 \circ \sigma_1 \circ h_1$ consisting of
	\begin{itemize}
		\item fully connected layers $h_k\colon \R^{N_{k-1}} \to \R^{N_k}$, $k = 1,2$,  with $h_k(a; W_k) = W_k a$, $N_0 = M$ and $N_2 = N$; 
		\item and continuously differentiable activation functions $\sigma_1$ and $\sigma_2$.
	\end{itemize}
	To simplify the notation, we define preactivation values $p^{(k)}  \in \R^{N_k}$ and (post)activation values $a^{(k)}   \in \R^{N_k}$ for the input $a^{(0)} \coloneqq y$ as
	\begin{equation}
		a^{(k)}   \coloneqq \sigma_k(p^{(k)}  ), \quad p^{(k)}   \coloneqq h_k(a^{(k-1)} ; W_k), \quad  k=1,2.
	\end{equation}
	These values are computed during the forward pass through the network $\q R_\theta$.%, in order $a_0$, $p_1$, $a_1$, $p_2$, $a_2$.
	
For the least squares loss $\ell(a) = \frac{1}{2} \| a -  x \|^2_\bb X$ and parameters $\theta = \{ W_1, W_2\} \in \R^P$ with $P = NN_1 + MN_1$, the training objective is given by 
	\begin{equation} \label{eq: ex backprop superpos}
		\q L(\theta) = \ell( \q R(y; \theta)) = \ell \circ \sigma_2 \circ h_2 \circ \sigma_1 \circ h_1 (y; \theta).
	\end{equation}
	
	To differentiate \eqref{eq: ex backprop superpos}, we employ the chain rule, which reads for two $C^1$ functions $f\colon \R^L \to \R^N$ and $g\colon \R^M \to \R^L$ as $J_{f \circ g}(x) = J_{f}(g(x)) J_g(x)$.
	Here, $J_{g}(x) \in \R^{L \times M}$ denotes the Jacobian matrix of $g$ at point $x$. 
	With the notation $\nabla_x f \coloneqq (J_{f})^\tT$, where the subscript $\nabla_x$ emphasizes that the derivatives is taken with respect to $x$, the chain rule becomes
	\begin{equation} \label{eq: chain rule grad}
		\nabla_x f \circ g (x) 
		= J_{f \circ g}^\tT(x)
		= (J_{f}(g(x)) J_g(x))^\tT 
		= J_g^\tT(x) J_{f}^{\tT}(g(x)) 
		= \nabla_x g(x) \nabla_x f(g(x)).
	\end{equation}	
Thus, to apply  \eqref{eq: chain rule grad} for our model \eqref{eq: ex backprop superpos}, we need to compute the gradients for each of the components.
Starting with $\ell \colon \R^N \to \R$, we expand it as
	\begin{equation}\label{eq:l2Loss}
		\ell(a) = \frac{1}{2} \sum_{i=1}^N ( a_i -  x_i )^2.
	\end{equation}
	Differentiating \eqref{eq:l2Loss} with respect to $a_i$ gives 
	$\frac{\partial}{\partial a_i} \ell(a) = a_i - x_i$ and $\nabla_{a} \ell(a) = a - x$.
	Secondly, we apply the activation functions $\sigma_k \colon \R^{N_k} \to \R^{N_k}$, $k =1,2$, entrywise as $\sigma_k(p) = (\sigma_k(p_1), \ldots, \sigma_k(p_{N_{k}}))^\tT$.
	Consequently, the partial derivatives are 
	\begin{equation}
		\frac{\partial (\sigma_k)_j(p)}{\partial p_i} = \frac{\partial \sigma_k(p_j)}{\partial p_i} = \begin{cases}
			\sigma_k'(p_i), & i = j, \\
			0, & i \neq j,
		\end{cases}
		\quad j = 1,\ldots, N_k.
	\end{equation}   
	Thus, $\nabla_{p} \sigma_k = \diag(\sigma_k'(p))$.
	Next, we turn to the linear layers $h_k(a;W_k)$. For the derivatives with respect to $a_i$, it holds
	\begin{equation}
		\frac{\partial (h_k)_j(a; W_k)}{\partial a_i} 
		= \frac{\partial}{\partial a_i} \sum_{n=1}^{N_{k-1}} (W_k)_{j,n} a_n
		= (W_k)_{j,i},
		\quad j = 1,\ldots, N_k,
	\end{equation} 
	which yields the gradient $\nabla_a h_k = W_k^\tT$.
	For the derivative with respect to $W_k$, we get
	\begin{equation} 
		\frac{\partial (h_k)_m(a; W_k)}{\partial (W_k)_{j,i}} 
		= \frac{\partial}{\partial (W_k)_{j,i}} \sum_{n=1}^{N_{k-1}} (W_k)_{m,n} a_n
		= \begin{cases}
			a_i, & j = m, \\
			0, & j \neq m,
		\end{cases}
		\quad m = 1,\ldots, N_k.
	\end{equation}
	Therefore, if we vectorize $W_k$, we obtain for any $z \in \R^{N_k}$ that 
	\begin{equation}\label{eq: backprop row comp}
		\nabla_{W_{k}} h_k(a;W_k) z = \begin{bmatrix}
			z_1 a \\
			z_2 a \\ 
			\vdots \\
			z_{N_k} a
		\end{bmatrix} \in \R^{N_k N_{k-1}}.
	\end{equation}
	
	Finally, we use the above computations to build $\nabla_{\theta} \mathcal L$. 
	By \eqref{eq: chain rule grad}, the gradient with respect to $W_2$ is given by
	\begin{equation}\label{eq:DiffW_2}
	\nabla_{W_2} \q L(\theta) = \nabla_{W_2} h_2(a^{(1)} , W_2) \nabla_p \sigma_2 (p^{(2)} ) \nabla_a \ell(a^{(2)}  ),
	\end{equation}
	where $\nabla_p \sigma_2 (p^{(2)}) \nabla_a \ell(a^{(2)}) = \diag( \sigma_2'(p^{(2)})) (a^{(2)} - x)$ and the first factor in \eqref{eq:DiffW_2} is evaluated via \eqref{eq: backprop row comp}.
	Next, we compute the gradient  with respect to $W_1$, which reads
	\begin{align}\label{eq:DiffW_1}
		\nabla_{W_1} \q L(\theta) 
		& = \nabla_{W_1} h_1(a^{(0)}, W_1) \nabla_p \sigma_1 (p^{(1)}) \nabla_a h_2(a^{(1)}, W_2) \nabla_p \sigma_2 (p^{(2)}) \nabla_a \ell(a^{(2)}) \notag \\
		& = \nabla_{W_1} h_1(a^{(0)}, W_1) \diag(\sigma_1'(p^{(1)})) W_2^\tT \diag( \sigma_2'(p^{(2)})) (a^{(2)} - x),
	\end{align}
	which is again evaluated via \eqref{eq: backprop row comp}.
	These formulas can be analogously extended for networks with an arbitrary number of layers.
	Comparing \eqref{eq:DiffW_2} and \eqref{eq:DiffW_1}, we observe that we have to compute the gradients in the reverse order, i.e., $\nabla_{W_2} \q L(\theta)$ and then $\nabla_{W_1} \q L(\theta)$, which is why this process is referred to as backpropagation.
	
	The computational complexity of backpropagation increases both with the number of layers $K$ and the number of parameters. 
	As we will see in \Cref{sec:LearnedReg}, backpropagation is sometimes not the most efficient approach for obtaining gradients (or even computationally unfeasible).
	Then, we use a different custom routine instead.
\end{example}

For the objective $\q L$ in \eqref{eq:BayesEst} with
\begin{equation}
	\q L(\theta) = \sum_{i = 1}^S \q L_i(\theta), \quad \text{and} \quad \q L_i(\theta) = \frac{1}{2S} \bigl \Vert \q R_{\theta}(y_i) - x_i \bigr \Vert_{\bb X}^2,
\end{equation}
 we have to perform the calculations from \Cref{ex: backprop} for each of the training pairs $(x_i,y_i)$, $i = 1, \ldots, S$.
 Since $S$ is usually very large (the data set might not even fit into the memory of our GPU), the full differentiation quickly becomes impractical. 

As a remedy, instead of computing $\nabla \q L$, we only select a small random \emphi{batch} $B_k \subset \{1,\ldots,S\}$ of $m$ indices, and approximate $\nabla \q L$ by the \emphi{stochastic gradient}
\begin{equation}
G_k(\theta) = \frac{S}{m}\sum_{i \in B_k} \nabla \q L_i(\theta) 
\end{equation}
in each step of \Cref{algo:gradient_descent}.
This leads to the stochastic gradient descent \Cref{algo:stochastic_gradient_descent}.

\begin{algorithm}[H]
	\KwIn{$f = \sum_{i=1}^S f_i  \in C^1(\R^N )$, initialization $\theta_0 \in \R^N $, step size $\nu_k > 0$, batch size $m$.}    
	\For{\color{black} $k = 0, 1, \ldots$}{
		sample a random batch $B_k \subset \{ 1,\ldots,S \}$ with $\vert B_k\vert = m$\;
		$\theta_{k + 1} = \theta_k-\nu_k G_k(\theta_k)$;
	}
	\caption{Stochastic Gradient Descent (SGD)}
	\label{algo:stochastic_gradient_descent}
\end{algorithm}
Due to the stochastic nature, we cannot expect that \Cref{algo:stochastic_gradient_descent} converges to a minimizer of $f$ with constant step size $\nu_k = \nu$, which is often referred to as \emphi{learning rate} in machine learning.
Typically, \Cref{algo:stochastic_gradient_descent}  bounces around a critical point.
Thus, convergence is studied in sense of distributions using Markov chain theory, see for instance \cite{yu2021analysis}.
To obtain pointwise convergence, one usually decreases $\nu_k$ successively.
This is known as \emphi{learning rate decay}.
We summarize the convergence of SGD as follows.  
\begin{thm}{Convergence of SGD}{convSGD}
Let $f = \sum_{i=1}^S f_i$, where each $f_i \in C^1(\R^N)$ has Lipschitz continuous gradient. Suppose that the step sizes $\nu_k$ satisfy
\begin{equation}\label{eq: stoc grad step size}
\sum_{k=0}^\infty \nu_k = +\infty 
\quad \text{and} \quad
\sum_{k=0}^\infty \nu_k^2 < +\infty,
\end{equation}
and that the indices in $B_k$ are sampled independently and uniformly with replacement. 
Then, the sequence $\{x_k\}_{k \in \N}$ generated by \Cref{algo:stochastic_gradient_descent} satisfies $\lim_{k \to \infty} \Vert \nabla f(x_k) \Vert_{\bb X} = 0$ almost surely, in which case we have 
\begin{equation}
	\min_{k = 0, \ldots, K-1} \| \nabla f(x_k) \|_{\bb X}^2 
	= \q O\biggl( \Bigl( \sum_{k=0}^{K-1} \nu_k \Bigr)^{-1} \biggr), \quad \text{as} \quad K \to \infty.
\end{equation}
If $\lim_{k \to \infty} x_k = x^*$, then $\nabla f(x^*) = 0$.
\end{thm}

\Cref{th:convSGD} is a compilation of the statements in \cite{khaled2023better, liu2022almost, orabona2020almost}.
In \cite{khaled2023better}, a non-asymptotic convergence result for SGD under rather general conditions is provided.
This result is refined in \cite{liu2022almost}, where stronger asymptotic (namely $k \to \infty$) rates are obtained.
Finally, \cite{Bertsekas2000} establishes that with probability one the so-called last-iterate convergence $\lim_{k \to \infty} \Vert \nabla f(x_k) \Vert_{\bb X} = 0$ applies. They proof was later refined and condensed in \cite{orabona2020almost}. 

The conditions in \eqref{eq: stoc grad step size} date back to results by Robbins and Siegmund \cite{robbins1971convergence}. In particular, taking $\nu_k = (k+1)^{-1 + \varepsilon}$ with $0 < \varepsilon < 1/2$ satisfies \eqref{eq: stoc grad step size} and results in the convergence rate $\q O(k^{-\varepsilon})$.
This is slower than the $\q O(k^{-1})$ rate of gradient descent.
Yet, the computational cost of a single iteration of SGD is much lower, leading to overall faster training.   

\begin{remark}\label{rem: ADAM}
	In practice, gradient-based methods with acceleration, momentum, and gradient scaling are preferable over the vanilla SGD in \Cref{algo:stochastic_gradient_descent}.
	The most popular one in recent years is ADAM \cite{kingma2014adam}.
\end{remark}

\begin{remark}
	When deploying random sampling with replacement for a training set of size $S$, a given sample \emph{will not} be used for $S$ consecutive iterations with probability $(1 -\frac1S)^S \approx 36,7\%$ for large $S$.
	Thus, an alternative is random reshuffling \cite{mishchenko2020random, gurbuzbalaban2021random}, where all samples are observed once by looping over 
	a permuted index set $\{1, \ldots, S\}$.
	For each loop (commonly referred to as \emph{epoch}), a different (random) permutation is used.
	This prevents order-induced bad convergence and/or over-fitting. 
\end{remark}

\subsection{Denoising and Learned Post-Processing}\label{sec:Denoising}
\begin{figure}[htb]
	\centering
	\includegraphics[width=0.99\textwidth]{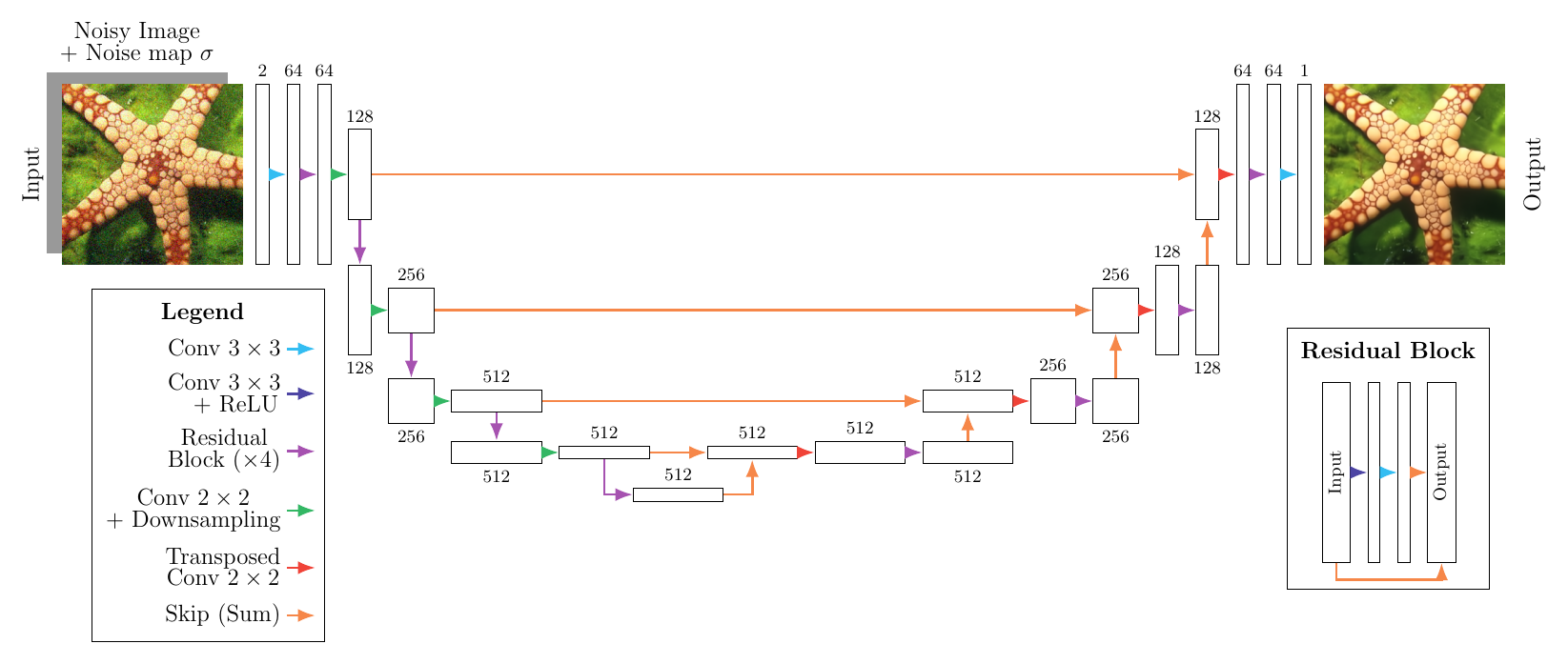}
	\caption{DRUNet denoiser \cite{ZhaLiZuo2022}.
		%: ``SConv''  and ``TConv'' represent  convolution with downsampling and a transposed convolution, respectively.
		%Original figure from \cite{ZhaLiZuo2022}.
	The \emph{starfish} image is part of BSDS500.}
	\label{fig:DRUNet}%\vspace{-0.15cm}
\end{figure}

Now, we briefly discuss two data-driven reconstruction methods.
A short introduction to the training of neural networks using PyTorch is given in their  \href{https://docs.pytorch.org/tutorials/beginner/basics/intro.html }{tutorial}.\footnote{https://docs.pytorch.org/tutorials/beginner/basics/intro.html } 
\paragraph{Denoising}
First, we briefly discuss the supervised learning of a Gaussian denoiser, where the underlying (inverse) problem is $y = x + \eta$ with $\eta \sim \q N(0,\sigma^2)$ for some $\sigma \in [0,50]$.
This task is not ill-posed in the sense of Hadamard (it is actually the most well-posed task possible), and the difficulty arises solely from the noise.
As reconstruction map $\q R_\theta$, commonly referred to as \emphi{denoiser}, we deploy the DRUNet architecture depicted in \Cref{fig:DRUNet}.
Note that $\q R_\theta$ takes $\sigma$ as an input parameter (a so-called non-blind setting).

\begin{figure}[htb]
	\centering
	\includegraphics[width=0.99\textwidth]{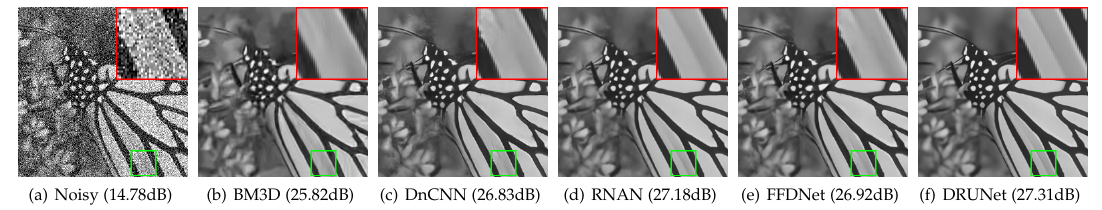}
	\caption{Denoising of the \emph{monarch} image for $\sigma = 50$ for some competing methods from the literature.
	Discussing those methods is beyond the scope of this lecture.}
	\label{fig:Monarch}%\vspace{-0.15cm}
\end{figure}

Common training sets for natural gray-scale images are BSDS500 \cite{arbelaez_contour_2011}, div2k \cite{timofte2017ntire} or Flickr2k, and we create the respective data pairs $(x,y)$ artificially (on the fly) by adding Gaussian noise to them.
For computational efficiency, we do not process full images in \eqref{eq:BayesEst}.
Instead, we randomly extract $128 \times 128$ patches from each training image when loading a batch.
We discuss the implementation for a simplified setting in Exercise \cref{ex:Denoiser}.
An example reconstruction is given in \Cref{fig:Monarch}, where the reconstruction quality is assessed in terms of the \emph{peak signal-to-noise ratio} \emphii{PSNR}
\begin{equation}\label{eq: PSNR}
	\mathrm{PSNR}(x_1,x_2) = 10 \log_{10}\biggl(\frac{255^2 N_xN_y}{\Vert x_1- x_2 \Vert^2_2}\biggr),
\end{equation}
which is measured in dB.
Higher values correspond to {\it better} results.

\paragraph{Learned Post-Processing}
As motivation, we consider inverse problems for which $A^*A$ is a convolution.
This class obviously includes denoising and deconvolution (\Cref{ex: convolution}), but one can show that it also includes MRI and CT (\Cref{ex: tomography}).
Since the Moore--Penrose inverse satisfies $A^\dagger = (A^*A)^{-1} A^*$ where $(A^*A)^{-1}$ corresponds to a filtering operation, we refer to $A^\dagger$ as \emphi{filtered backprojection}.
In  \Cref{ex: svd convolution}, we have shown that we can use the Fourier transform for an efficient evaluation.
The actual implementation is covered in \Cref{ex:DeepInv} .
A more sophisticated reconstruction approach is to solve the associated $L_1$-regularized (\Cref{ex:ista}) or TV-regularized (\Cref{example:TV_for_inpainting_ADMM}) problem, which requires repeated convolutions and point-wise nonlinearities.
In summary, CNN architectures should offer a promising solution approach.

\begin{figure}[htb]
	\centering
	\includegraphics[trim = 0mm 50mm 0mm 40mm,clip=true,width=0.99\textwidth]{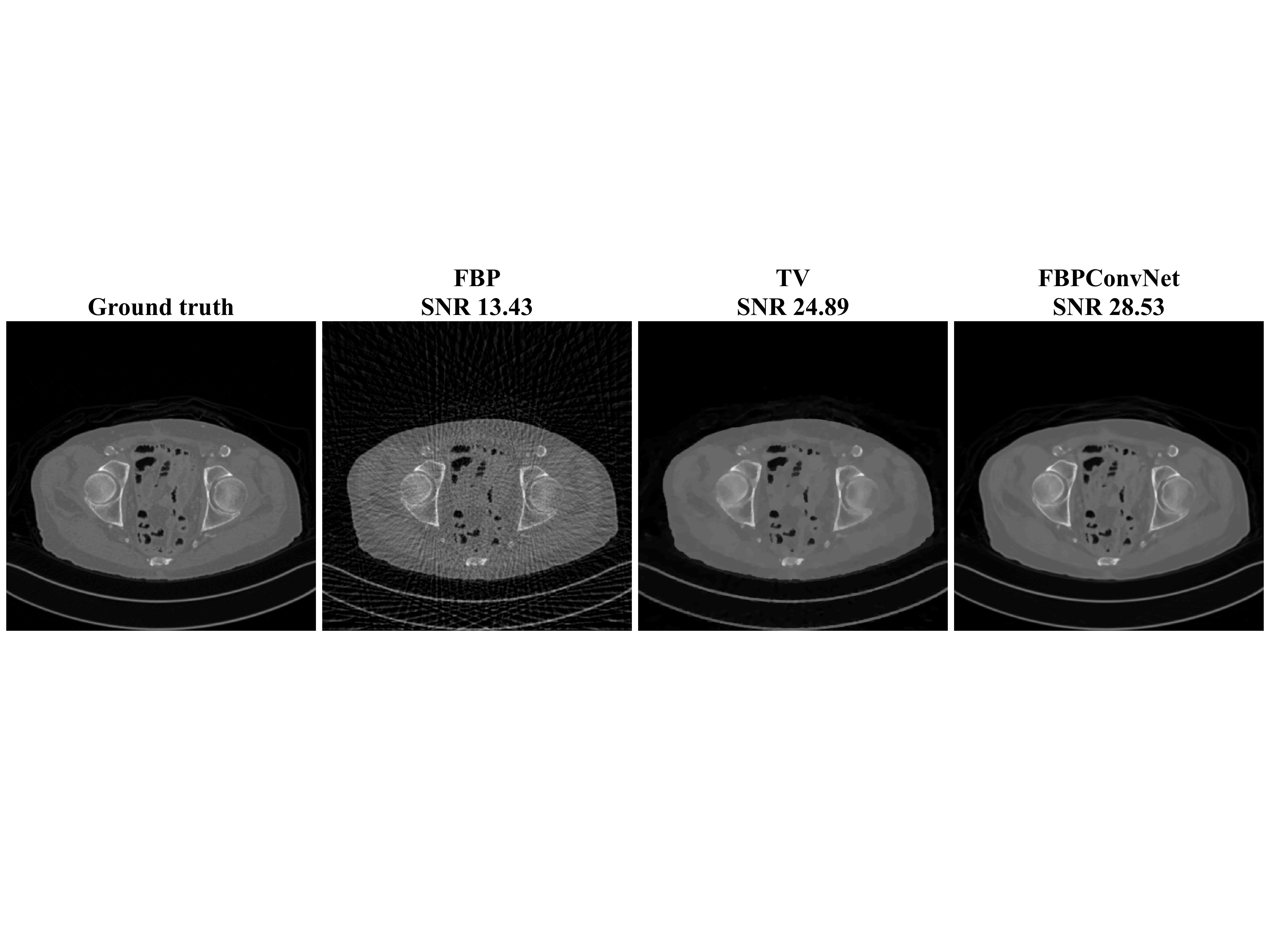}
	\caption{Postprocessing of the FBP for CT \cite{JinMcCFro2017}.
	The network $\Phi$ significantly enhances the quality and removes the (typical) stripe artifacts from the FBP.}
	\label{fig:PostProcess}
\end{figure}

Now, we want to design CNN-based reconstruction architecture $\q R_\theta \colon \bb Y \to \bb X$ for the  Radon transform with 143 views (sparse view setting) as forward map $A$, see also  \Cref{ex: tomography}.
As proposed in \cite{JinMcCFro2017}, we first apply the generalized inverse $A^\dagger$ (one could also apply a simple spectral regularization $\tilde{\q R}_\alpha \colon \bb Y \to \bb X$), which is then followed by a UNet $\Phi_\theta \colon \bb X \to \bb X$.
Then, we get the overall architecture $\q R_\theta = \Phi_\theta \circ A^\dagger$.
While we could also train a network $\Phi_\theta \colon \bb Y \to \bb X$ directly on the data, applying $A^\dagger$ first simplifies the learning and encapsulates the forward map $A$ into the reconstruction architecture.
In particular, the $\Phi_\theta$ would need to encode the coordinate change, which is avoided when $A^\dagger$ is used as input layer for the network instead.
After having fixed the architecture, we need to specify our supervised loss \eqref{eq:BayesEst}.
Given ground-truth images $x_i$, we can simulate training pairs $(x_i,y_i)$ via $y_i = A x_i + n_i$.
In \Cref{fig:PostProcess} we show reconstruction results from \cite{JinMcCFro2017}, for which the training set consisted of 475 CT images from a database by the Mayo clinic.

\begin{figure}[htb]
	\centering
	\includegraphics[width=0.485\textwidth]{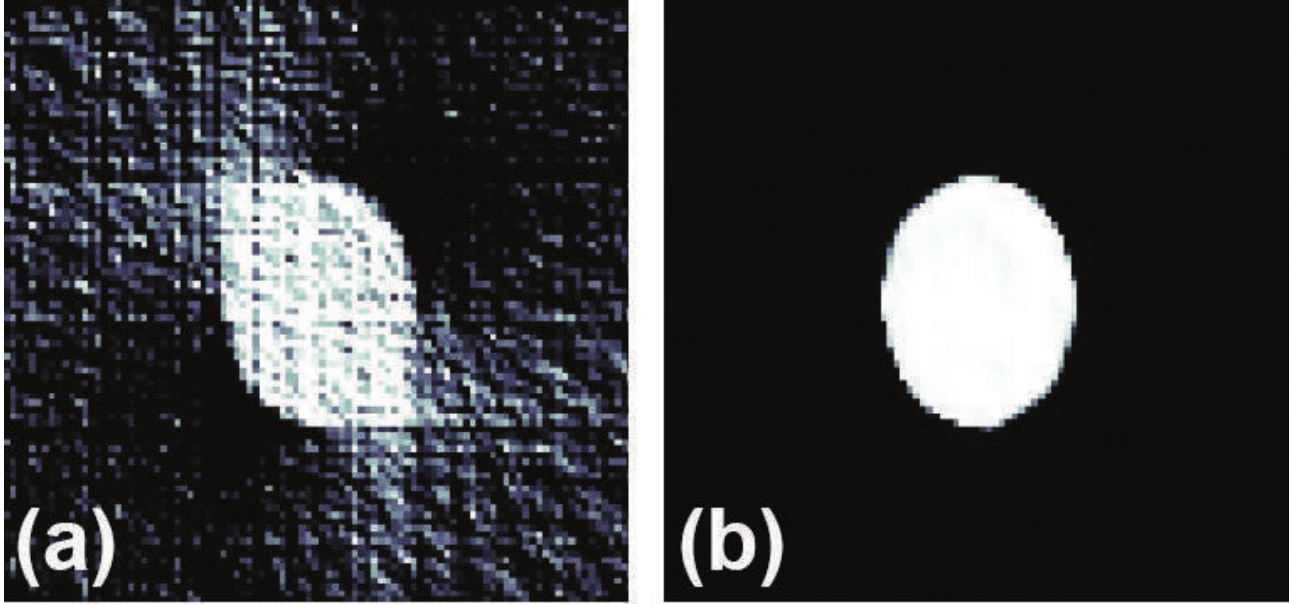}
	\includegraphics[width=0.49\textwidth]{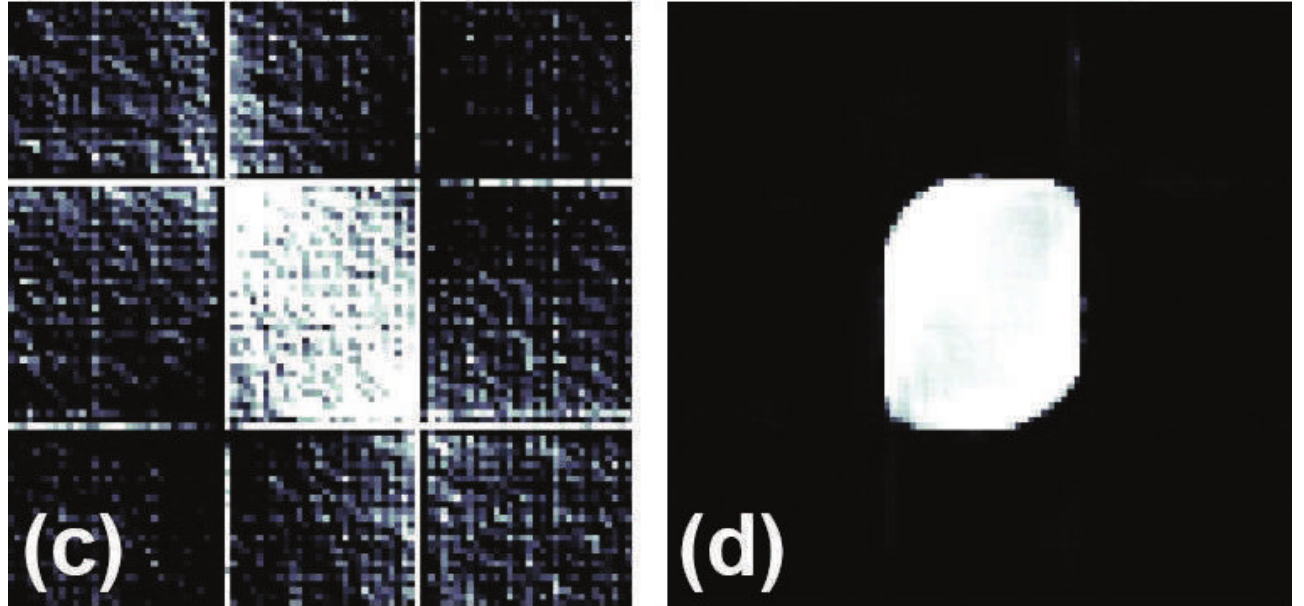}
	\caption{Ground-truth $x_i$ for (a) and (b) is an ellipsoid, and a square for (c) and (d).
	The reconstructions (a) and (c) correspond to $A^\dagger$, and (b) and (d) to post-processed solutions.}
	\label{fig:PostProcessFail}
\end{figure}

As this example shows, learned post-processing methods can be very powerful with functionality beyond simple memorization of the training dataset, see also \cite{GenGuhMac2022}.
Still, the quality of the data is crucial for every learning process.
An important question to consider is: What are the consequences when certain structures or pathologies are absent during training but appear in test images?
For this, we investigate a limited-angle CT setting, where  $A$ is the Radon transform with a restricted range of angles in $[0,\pi/2]$.
Following a recent SIAM News article\footnote{ \href{https://www.siam.org/publications/siam-news/articles/learning-in-image-reconstruction-a-cautionary-tale/}{\it Learning in Image Reconstruction: A Cautionary Tale}
Code for the experiment is available at \url{https://github.com/TimRoith/LearningInImageReconstruction-SIAM-News}}, we consider ellipses randomly sampled on the domain $[0,1]^2$ as simpler images $x_i$.
Otherwise, the reconstruction and training pipeline remains as for \Cref{fig:PostProcess}.
As we can see in \Cref{fig:PostProcessFail}, the method works well if we try to reconstruct an ellipsoid image even  if this specific one was not part of the training set (in distribution data).
If we instead apply the method to a square (out of distribution data), the method \emph{hallucinates} the missing information to match an ellipsoid.
In other words, we have an overfitting to the training set.
Hence, we want to explore more \emph{robust} approaches to solve inverse problems in the remaining part of the lecture.
\subsubsection{Exercises}
\begin{exerciseB}\label{ex:DeepInv}
	Familiarize yourself with the \href{https://deepinv.github.io/deepinv/}{DeepInverse}\footnote{
		\url{https://deepinv.github.io/deepinv/}} library \cite{tachella2025deepinv}.
	Try to set up different forward operators (CT, MRI, Blur), apply them to an example image, and depict the pseudoinverse.
	Then, implement more sophisticated reconstruction methods such as Tikhonov and total-variation regularization.
	Do these improve the results?
	
	You may use the following Jupyter Notebook as a basis for this exercise, \url{https://github.com/Sniper2k/inverse_problems_ai_exercise_templates/blob/main/Ex5_4_1_template.ipynb}.
\end{exerciseB}

\begin{Solution}
	See Jupyter Notebook \url{https://github.com/Sniper2k/inverse_problems_ai_exercises/blob/main/Ex5_4_1.ipynb}.
\end{Solution}

\begin{exerciseB}\label{ex:Denoiser}
	Train a UNet denoiser with 3 scales using DeepInverse \cite{tachella2025deepinv}.
	As training data, take the CBSD68 dataset and separate it into 60 training and 8 validation color images. 
	For testing, we use the set3c dataset containing 3 images. 
	During training, we extract random $128 \times 128$ patches  from the training images.
	For the validation and the test set, we use $128 \times 128$ center patches.
	For all datasets, Gaussian noise is generated on the fly with $\sigma = 0.05$.
	Can you achieve a reconstruction PSNR of 30 on the test set?
	
	Now, change the physics to the inpainting problem from \cref{exer:impainting}.
	That is, with probability $0.3$ the pixel is observed.
	Train the UNet again as reconstruction network.
	How does the inpainting PSNR compare to the denoising PSNR?
	
	You may use the following Jupyter Notebook as a basis for this exercise, \url{https://github.com/Sniper2k/inverse_problems_ai_exercise_templates/blob/main/Ex5_4_2_template.ipynb}.
\end{exerciseB}

\begin{Solution}
	See Jupyter Notebook \url{https://github.com/Sniper2k/inverse_problems_ai_exercises/blob/main/Ex5_4_2.ipynb}.
\end{Solution}

\begin{exerciseB}\label{ex:PostprocessingUNet}
	Train a post-processing network for the described CT setup.
	You should use the predefined CT operator from DeepInverse \cite{tachella2025deepinv}, and generate the training and test sets through artificial simulation from the groundtruth Mayo clinic images.
	What happens if you deploy the post-processing network to other images, e.g., the Shepp-Logan phantom.\footnote{\url{https://www.google.com/search?q=Shepp-Logan+phantom}}
	
	You may use the following Jupyter Notebook as a basis for this exercise, \url{https://github.com/Sniper2k/inverse_problems_ai_exercise_templates/blob/main/Ex5_4_3_template.ipynb}.
\end{exerciseB}

\begin{Solution}
	See Jupyter Notebook \url{https://github.com/Sniper2k/inverse_problems_ai_exercises/blob/main/Ex5_4_3.ipynb}.
\end{Solution}

\section{Learned regularizers} \label{sec:LearnedReg}

The variational reconstruction model
\begin{equation}\label{eq:data_fidelity_regularizer2}
	\argmin_{x\in\R^N} \Vert Ax - y \Vert_2^2 + \alpha R(x),\quad \alpha>0
\end{equation}
from \Cref{sec: variational regularization} depends heavily on the choice of $R$ and $\alpha$.
Now, we discuss the search for $R$ that leads to good reconstructions.
For this, we rely on the machine learning principles that we developed in \Cref{chap:regularization}.
For any parametric architecture $R_\theta$, \eqref{eq:data_fidelity_regularizer2} induces a regularization operator $\q R_{\theta} \colon \bb \R^M \to \R^N$, where we include $\alpha$ into the learnable parameters $\theta$.
A simple instance with just one learnable parameter would be to find the optimal regularization strength $\alpha$ for the non-parameteric TV regularizer $\Vert \nabla x \Vert_1$ from \Cref{ex: sparsity promoting}.
More generally, we can try to learn a better sparsifying transform $B$ to replace the gradient $\nabla$.
Taking this as a starting point, we aim to design more regularization architectures $R_\theta$.

To this end, given paired training data $\{(x_i,y_i)\}_{i=1}^S$, the parameters $\theta \in \R^P$ of $R_{\theta}$ can be learned with the supervised loss \eqref{eq:BayesEst}.
This leads to the \emphi{bilevel learning problem}
\begin{align}
	&\hat \theta \in \argmin_{\theta} \biggl\{\mathcal L(\theta) = \frac{1}{2S}\sum_{i=1}^S
	\Vert \hat{x}(\theta,y_i) -  x_i \Vert^2_{\bb X}\biggr\}\label{eq:bilevel-prob}\\
	&\hat{x}(\theta; y)\in \argmin_{x}\bigl\{\Phi_{y}(x,\theta) = \Vert A x -  y \Vert_2^2 + R_\theta(x) \bigr\}.\label{eq:lower-level}
\end{align}
Given parameters $\theta$, we compute $\hat{x}(\theta; y)$ using one of the algorithms from \Cref{sec:Algs}.
Then, we minimize \eqref{eq:bilevel-prob} using stochastic optimization as discussed in \Cref{sec:OptimML}.
The gradient of the training loss $\q L (\theta)$ can be computed by applying the chain rule \eqref{eq: ex backprop superpos}, which leads to
\begin{equation}\label{eq:hypergrad_chain}
	\nabla_\theta \q L (\theta) = \frac{1}{S} \sum_{i=1}^S \nabla_{\theta}\hat{x}(\theta;y_i) \bigl(\hat{x}(\theta; y_i) -  x_i \bigr) .
\end{equation}
The formula \eqref{eq:hypergrad_chain} involves the (stochastic) gradient of $\hat{x}(\theta; y)$ with respect to $\theta$, which is also known as \emphi{hypergradient} in bilevel  optimization. 
A comparison  of different approaches for its estimation can be found in \cite{JiYanLia2021,LiaXioFet2018,ZucSac2022} and a complete overview on learning regularization functionals including implementations is given in \cite{hertrich2025learning}.

\subsection{Estimating the Hypergradient}\label{sec:EstHyper}
For simplicity, we restrict the notation to a single data point $(x,y)$ with associated loss summand $\q L_y$.
Now, we discuss three different approaches for estimating $\nabla_{\theta}\hat{x}(\theta;y)$.

\subsubsection{Implicit Differentiation Approach}
This approach is built around the optimality condition of the lower-level problem \eqref{eq:lower-level}, which reads as
\begin{equation}\label{eq:ift1}
	\nabla_{x}\Phi_y(\hat{x}(\theta), \theta)=0.
\end{equation}
The solutions to the implicit equation \eqref{eq:ift1} are characterized by the following theorem.
\begin{thm}{Implicit function theorem (IFT)}{ImplicitFunction}   
	Let $\Phi_y \in C^2(\R^N \times \R^P)$ and $\hat{x} \in \R^N$, $\bar \theta \in \R^P$ with $\nabla_{x}\Phi_y(\hat{x}, \bar \theta)=0$.
	Assume that
	\begin{equation}\label{eq:IFTHess}
		\nabla_{x,x}^2\Phi_y\bigl(\hat{x}(\bar \theta),\bar \theta\bigr) = 2 A^T A +  \nabla_{x,x}^2 R_{\bar \theta}(\hat{x}(\bar \theta))\in \R^{N \times N}
	\end{equation}
	is invertible.
	Then, there exists an open neighborhood $U$ of $\bar \theta \in \R^P$ and $\hat x \colon U \to \R^N$ such that $\nabla_{x}\Phi_y(\hat{x}(\theta), \theta)=0$ for all $\theta \in U$.
	If for each $\theta \in U$ the Hessian $\nabla_{x,x}^2\Phi_y(\hat{x}(\theta),\theta)$ is invertible, then it holds that
	\begin{equation}\label{eq:JacS}
		\nabla_\theta \hat x (\theta) = - \nabla_{x,\theta}^2\Phi_y(\hat{x}(\theta),\theta)  \bigl(\nabla_{x,x}^2\Phi_y(\hat{x}(\theta),\theta)\bigr)^{-1}.
	\end{equation}
\end{thm}
\begin{pr}
	We only comment on \eqref{eq:JacS} and refer to introductory analysis courses for a complete proof.
	Differentiating \eqref{eq:ift1} on both sides with respect to $\theta$ leads to
	\begin{equation}\label{eq:ift2}
		0= \nabla_\theta \hat x (\theta) \nabla_{x,x}^2\Phi_y(\hat{x}(\theta),\theta) +\nabla^2_{x,\theta}\Phi_y(\hat{x}(\theta),\theta),
	\end{equation}
	from which the claim readily follows.
\end{pr}
\begin{example}
	Let $\Phi \colon \R \times \R \to \R$ with $\Phi(x, \theta) = \frac13 x^3 + \theta^2x - x$.
	It holds that $\nabla_x \Phi(x,\theta) = x^2 + \theta^2 - 1$ and $\nabla^2_{x,x} \Phi(x, \theta) = 2x$.
	Hence, we cannot apply the IFT for the points $(0,\pm 1)$.
	For all other points $(\hat x(\theta), \theta)$ satisfying $\nabla_x \Phi(\hat x(\theta),\theta)= 0$, we compute $\nabla_{x,\theta}^2\Phi(x,\theta) = 2 \theta$ and $(\nabla^2_{x,x} \Phi(x,\theta))^{-1} = \tfrac{1}{2x}$.
	Thus, by \Cref{th:ImplicitFunction} we have $\nabla_\theta \hat x(\theta) = - \theta /  \hat x(\theta)$. 
	On the other hand, the implicit equation $x^2 + \theta^2 = 1$ describes the unit circle and its local solutions are $\hat x(\theta) = \pm \sqrt{1 - \theta^2}$, depending on the sign of $\hat x$.
	As expected, we obtain the derivative from before, namely
	\begin{equation}
		\nabla_\theta \hat x(\theta) = \pm \frac{1}{2\sqrt{1 - \theta^2}} \cdot (-2 \theta) = - \frac{\theta}{\hat x(\theta)}.
	\end{equation}
	As we often do not know the solution function $\theta \mapsto \hat x(\theta)$, the relation \eqref{eq:JacS} is very useful.
\end{example}
The matrix \eqref{eq:IFTHess} is invertible if $R_\theta$ is strictly convex or more generally if $\hat{x}$ is a nondegenerate local minimum of $x \mapsto \Phi_y(x, \theta)$, see also \Cref{th:second_order_sufficient}.
Inserting \eqref{eq:JacS} into \eqref{eq:hypergrad_chain} leads to
\begin{equation}\label{eq:formula_grad}
	\nabla_\theta \q L_y (\theta) = -\nabla_{x,\theta}^2\Phi_{y}(\hat{x}(\theta;y),\theta) \bigl(\nabla_{x,x}^2\Phi_{y}(\hat{x}(\theta;y),\theta)\bigr)^{-1} \bigl(\hat{x}(\theta; y) -  x \bigr).
\end{equation}
In practice, we do not compute the inverse in \eqref{eq:formula_grad} explicitly and use an iterative linear solver (such as CG or MINRES) instead (also called backward solve).
Its computational complexity and memory requirements are independent of the chosen algorithm for computing $\hat{x}(\theta;y)$.
In particular, the complexity does not depend on the number of steps in the forward solve.
However, the backward solve requires several Hessian vector products for each CG/MINRES step, which can become costly for architectures $R_\theta$ with many parameters in case that the products are computed using automatic differentiation.

\begin{remark}[Convex Regularizer]\label{eq:RemLip}
	For a $C^2(\R^N)$ and convex regularizer $R$, we can use \Cref{th:ImplicitFunction} to estimate the Lipschitz constant of $y \mapsto \hat x(y)$. 
	For an injective forward $A$, we have that the objective \eqref{eq:lower-level} is strongly convex.
	Thus, a unique reconstruction exists by \Cref{th:minima}.
	Moreover, the conditions of \Cref{th:ImplicitFunction} are satisfied, and
	for given $y_1, y_2 \in U$ we get
	\begin{align}
		\Vert \hat x(y_1) - \hat x(y_2) \Vert_2 \leq& \sup_{y \in U} \Vert \nabla_y \hat x (y) \Vert_{2} \Vert y_1- y_2\Vert_2
		 = \bigl\Vert A (\nabla_{x,x}^2 \Phi(x,y))^{-1} \bigr\Vert_{2} \Vert y_1- y_2\Vert_2 \label{eq:RegularizerHessian} \\
		\leq & \bigl\Vert A (A^T A)^{-1} \bigr\Vert_{2} \Vert y_1- y_2\Vert_2 
		= \| A^\dagger\|_2 \Vert y_1- y_2\Vert_2.
	\end{align}
	The detailed derivation of the above bounds is left as \Cref{ex: stongly convex via Banach}.
	This means that $\hat x \colon U \subset \R^M \to \R^N$ is Lipschitz continuous with Lipschitz constant at most $\Vert A^{\dagger}\Vert_2$.
	This is a much stronger property than the mere continuity established in \Cref{th:stability variational}.
	Moreover, since our estimation ignores the stabilizing effect of the regularizer $R$ in \eqref{eq:RegularizerHessian}, we expect even better stability in practice.
\end{remark}

\subsubsection{Unrolling Approach}
To compute a root $\hat{x}(\theta; y)$ of the implicit equation \eqref{eq:ift1}, we deploy one of the iterative algorithms from \Cref{sec:Algs}.
These algorithms take the general form $\hat{x}(\theta; y) \approx x^{(T)}$ with
\begin{equation}\label{eq:IterativeUpdate}
	x^{(k+1)}=\Gamma(x^{(k)};\theta), \text{ for }k=0,1,\dots,T-1,
\end{equation}
where $T$ corresponds to the (fixed) number of iterations and $x^{(0)}$ is the initialization.
Here, we only discuss the gradient descent step $\Gamma(x;\theta) = x - \nu \nabla_{x}\Phi_{y}(x,\theta)$ from \Cref{algo:gradient_descent}, where the step size $\nu$ is included into the learnable parameters $\theta$.
Since our aim is to minimize \eqref{eq:lower-level}, we have to choose $T$ large enough such that the scheme \eqref{eq:IterativeUpdate} has converged to a stationary point, see also \Cref{th:grad descent nonconvex}.
Different from the IFT approach, we then compute  each $\nabla_{\theta}\hat{x}(\theta; y)$ by applying the chain rule to the iterative scheme \eqref{eq:IterativeUpdate}, which leads to
\begin{align}
	\nabla_\theta \hat x(\theta;y)
	& \approx \nabla_\theta x^{(T)} = \nabla_\theta\Gamma(x^{(T-1)}(\theta);\theta) \nonumber \\
	& = \nabla_\theta \Gamma(x;\theta)\big|_{x = x^{(T-1)}} + (\nabla_\theta x^{(T-1)}) \nabla_x \Gamma(x;\theta)\big|_{x = x^{(T-1)}}. \label{eq: backprop unrolling derivation}
\end{align}
By setting
\begin{equation}
B_{k+1} = \nabla_{x} \Gamma(x;\theta)\big|_{x = x^{(k)}}, \qquad C_{k+1} = \nabla_{\theta} \Gamma(x;\theta)\big|_{x = x^{(k)}}
\end{equation}
and computing the gradients of $x^{(T-1)}, \ldots, x^{(1)}$ analogously, we arrive at the estimate
\begin{align}\label{eq:GradUnrolling}
	&\nabla_{\theta}\hat{x}(\theta;y) = \sum_{k=1}^{T} C_kB_{k+1}\ldots B_{T},
\end{align}
which is often referred to as \emphi{unrolling} \cite{MonLiEld2021}.
The hypergradient \eqref{eq:GradUnrolling} corresponds to the exact numerical gradient.
In some sense, we can interpret unrolling as \emph{learning to optimize} \cite{CheCheChe2022} the problem \eqref{eq:lower-level}, since by learning $\theta$ we tune the optimization method for optimal performance.
The expression \eqref{eq:GradUnrolling} can be evaluated using backpropagation, see \Cref{ex: backprop}, namely from right to left and starting with the summand $T$.
This is readily implemented in any automatic differentiation package.
As it requires the iterates $x^{(k)}$ in reverse order, we cannot  evaluate \eqref{eq:GradUnrolling} during the forward pass and need to store $\{x^{(k)}\}_{k=1}^{T}$ instead.
Hence, the memory requirement grows linearly with $T$, which restricts the number of steps $T$ that we can take.
This contrasts our observation that $T$ should be large enough to find a minimizer of \eqref{eq:lower-level}.
We can reduce the memory requirements by checkpointing, i.e., by only storing some of the iterates $\{x^{(k)}\}_{k=1}^{T}$. The remaining ones are recomputed on the fly with \eqref{eq:IterativeUpdate} and,
as the trade-off, computational cost increases.
Still, we do not need to solve a linear system as for the IFT approach.

Below, we provide a general convergence result for the iterative scheme \eqref{eq:IterativeUpdate}.
\begin{thm}{Banach Fixed Point Theorem}{BanachFix}
	If $\Gamma(x;\theta) \in C^1(\R^N)$ and $\|  \nabla_{x} \Gamma(x;\theta) \|_2 \leq \rho < 1$ for all $x \in \R^N$, then $\Gamma(\cdot;\theta)$ admits a unique fixed point $x_*$ and the iterations \eqref{eq:IterativeUpdate} satisfy
	\begin{equation}\label{eq:BanachRate}
		\Vert x^{(k)}- x_* \Vert_2 \leq \rho^k \Vert x^{(0)} - x_* \Vert_2.
	\end{equation}
\end{thm}
\begin{pr}
	We only provide a brief sketch. The assumption on $\Gamma$ implies that $\Gamma$ is $\rho$-Lipschitz continuous in $x$.
	Then, the existence proof builds upon the inequality
	\begin{equation}\label{eq:contract}
		\Vert x^{(k)}- x^{(k-1)} \Vert_2 = \bigl\Vert \Gamma(x^{(k-1)};\theta)-\Gamma(x^{(k-2)};\theta) \bigr\Vert_2 \leq \rho \Vert x^{(k-1)} -  x^{(k-2)} \Vert_2,
	\end{equation}
	with which we can show that $\{x^{(k)}\}_{k \in \N}$ is a Cauchy sequence.
	The uniqueness follows from the contraction property.
	To establish \eqref{eq:BanachRate}, we can use an estimate similar to \eqref{eq:contract}.
\end{pr}
In the context of gradient descent, $\| \nabla_{x} \Gamma(x;\theta) \|_2 < 1$ can only hold if the underlying $\Phi_y$ is strongly convex.
Thus, \Cref{th:BanachFix} is strongly connected with \Cref{th:grad descent linear}, and we discuss this in \Cref{ex: stongly convex via Banach}.

In view of \Cref{th:BanachFix}, a common trick to stabilize the training in the context of unrolling is to replace \eqref{eq:bilevel-prob} with the \emphi{regularized bilevel problem}
\begin{equation}\label{eq:RegBilevelProblem}
	\argmin_{\theta} \mathcal L(\theta) + \lambda  \frac{1}{S} \sum_{i=1}^S \|  \nabla_{x} \Gamma(\hat{x}(\theta;y_i);\theta) \|_2, \qquad \lambda > 0.
\end{equation}
This encourages contractivity of $\Gamma(\cdot;\theta)$ during training.
While we evaluate the Jacobian at $\hat{x}(\theta;y_i)$, other choices for the evaluation points are also possible.
To estimate the spectral norm $\| \nabla_{x} \Gamma(x;\theta) \|_{2}$, we have two options.
One approach is to use power iterations (discussed later in \eqref{eq:PowerIters}).
Alternatively, we can exploit the inequality $\| \cdot \|_{2} \leq \| \cdot \|_F$ and instead estimate the Frobenius norm.
The latter can be computed efficiently using the Hutchinson trace estimator
\begin{equation}
	\| \nabla_{x} \Gamma(x;\theta) \|^2_F = \mathbb{E}_{v \sim \mathcal{N}(0,I_N)} \|  \nabla_{x} \Gamma(x;\theta) v \|^2_2.
\end{equation}
Since we process multiple samples $x_i$ in each batch, sampling a single random vector $v$ per batch is often sufficient.
This approach requires only a single Jacobian-vector product, which is considerably more efficient than power iterations.
However, the Frobenius norm may significantly overestimate the spectral norm, leading to overly conservative regularization.

Another stabilization technique is to perform the iterative updates \eqref{eq:IterativeUpdate} with a random number of steps $T$ and a random initialization $x^{(0)}$ during training.
This adds stochasticity that prevents the model from overfitting to a specific unrolling depth and path.

\subsubsection{Truncated Backpropagation Approach}
A memory-efficient intermediate regime between unrolling and the IFT approach is \emphi{truncated backpropagation} \cite{ShaCheHat2019}. 
Instead of applying the chain rule to the whole scheme \eqref{eq:IterativeUpdate}, we only do so for the last $K$ steps, leading to the hypergradient
\begin{equation}\label{eq:TruncGrad}
	\nabla^K_{\theta}\hat{x}(\theta,y) \approx \sum_{k=T-K+1}^{T} C_kB_{k+1}\ldots B_{T}.
\end{equation}
By doing so, we ignore the early iterates for the gradient estimation, which means that we do not need to store them.
For the special case of \emphi{Jacobian-free backpropagation} (JFB) \cite{FunHeaLi2022} with $K=1$, we do not require any Jacobian matrices $B_k$.

First, we discuss the relation of $\nabla^K_\theta \q L_y (\theta) =  \nabla^K_{\theta}\hat{x}(\theta,y)(\hat{x}(\theta; y) -  x)$ with $\nabla_\theta \q L_y (\theta)$ based on the hypergradient \eqref{eq:GradUnrolling} from unrolling. 
\begin{thm}{Vanishing Bias}{exp convergence}
	Let $\Phi_{y}\in C^2(\R^N \times \R^P)$ have $L$-Lipschitz gradient and be  $\alpha$-strongly convex.
	Then, for $\nu < \frac{1}{L}$, it holds that
	\begin{align} 
		\norm{\nabla^K_\theta \q L_y (\theta) - \nabla_\theta \q L_y (\theta)}_2 \leq \frac{(1 -\nu \alpha)^K }{\nu \alpha}  M_C \norm{\hat{x}(\theta; y) -  x}_2.    \label{eq:convex bound}
	\end{align}
	where  $M_C = \max_{k\in\{0,\dots,T-K\} } \| C_k \|_{2}$.	
\end{thm}
\begin{pr}
	By definition of $\nabla^K_{\theta}\hat{x}(\theta,y)$, we obtain for $e_K \coloneq \nabla^K_\theta \q L_y (\theta) - \nabla_\theta \q L_y(\theta)$ that
	\begin{align}
		e_K = \biggl( \sum_{k=0}^{T-K}  C_{k} B_{k+1} \cdots B_{T-K}\biggr) B_{T-K+1} \cdots B_{T} \bigl(\hat{x}(\theta; y) -  x\bigr).
	\end{align}
	Therefore, if $\Phi_y $ is $\alpha$-strongly convex, then $\Vert B_k \Vert_2 \leq 1 - \nu \alpha$ by \Cref{ex: stongly convex via Banach}, and we estimate \eqref{eq:EstLocalStronglyConvex} as
	\begin{align}\label{eq:EstLocalStronglyConvex}
		\norm{e_K}_2 &\leq \biggl \Vert \sum_{k=0}^{T-K}  C_{k} B_{k+1} \cdots B_{T-K} \biggr \Vert_2 \bigl \Vert B_{T-K+1} \cdots B_{T}  \bigl(\hat{x}(\theta; y) -  x\bigr)\bigr \Vert_2\notag\\
		&\leq (1 - \nu \alpha)^{K} \norm{\hat{x}(\theta; y) -  x}_2  \max_{k\in\{0,\dots,T-K\} } \norm{C_k}_2 \sum_{k=0}^{T-K} (1 - \nu \alpha)^{k}.
	\end{align}
	Now, the estimate~\eqref{eq:convex bound} follows from
	$\sum_{k=0}^{T-K} (1 - \nu \alpha)^{k}  \leq  \sum_{k=0}^{\infty} (1 - \nu \alpha)^{k}  = \frac{1}{\nu \alpha}$.
\end{pr}

%Since $0 \leq 1 -\nu \alpha < 1$, \Cref{th:exp convergence} implies that if $\hat{x}^{(T)}$ close to a strict local minimum of \eqref{eq:lower-level}, then the error of using the approximation $\nabla^K_{\theta}\hat{x}(\theta,y)$ decays exponentially in $K$. 

We observe that the error of using the approximation $\nabla^K_{\theta}\hat{x}(\theta,y)$ decays exponentially in $K$.  
\begin{remark}
	A similar bound to \eqref{eq:convex bound} can be obtained for nonconvex $\Phi_y$, if the function is $\alpha$-strongly convex around the iterates $x^{(T)}, \ldots$, $x^{(T-K+1)}$, i.e, they are in proximity of a nondegenerate local minimum.
\end{remark}

Next, we show that $\nabla^K_\theta \q L(\theta)$ is actually a sufficient descent direction.
This is a much stronger property than the error decay established in \Cref{th:exp convergence}.

\begin{thm}{Sufficient Descent}{bound of many terms}
	Let $\Phi_y \in C^2(\R^N \times \R^P)$ be $\alpha$-strongly convex with $L$-Lipschitz gradient.
	Assume that $C_k$, $k= 0,\ldots,T$, has full column rank and that   $x \mapsto \nabla_x \Gamma(x;\theta)$ and  $x \mapsto \nabla_\theta \Gamma(x;\theta)$ are $L_{\Gamma}$-Lipschitz continuous.  
	Then, for all $K\geq 1$, with $T$ large enough and $\nu < \frac1L$ small enough, there exists a $c > 0$ such that $\nabla^K_\theta \q L_y(\theta)^T \nabla_\theta \q L_y(\theta) \geq  c \norm{\hat{x}(\theta; y) -  x}_2^2$.
\end{thm}
\begin{pr}
	We only prove the result for $K=1$.
	Since $\nabla^1_\theta \q L_y (\theta)= C_T  (\hat{x}(\theta; y) -  x)$, we obtain for $d_1 = 	\nabla^1_\theta \q L_y (\theta)^T \nabla_\theta \q L_y (\theta)$ with \eqref{eq:GradUnrolling} that
	\begin{align}
	 d_1= \norm{\nabla^1_\theta \q L_y (\theta)}_2^2 + \nabla^1_\theta \q L_y (\theta)^T \biggl(  \sum_{k=0}^{T-1}  C_{k} B_{k+1}  \cdots B_{T-1} \biggr)  B_{T} \bigl(\hat{x}(\theta; y) -  x\bigr).\label{eq:inner_prod}
	\end{align}
	Now, we lower bound the terms in this sum. 
	This is a lengthy computation with the final result given in \eqref{eq:FinalEstAlign}.
	First, we get that
	\begin{equation}\label{eq:LowerBoundStep1}
		\!\! \nabla^1_\theta \q L_y(\theta)^T   C_{k} B_{k+1} \cdots B_{T}   (\hat{x}(\theta; y) -  x)	\geq  \nabla^1_\theta \q L_y(\theta)^T  C_{T} B_{k+1} \cdots B_{T} 
		 (\hat{x}(\theta; y) -  x) - \Delta_1,
	\end{equation}
	where $\Delta_1$ is obtained by bounding    
	\begin{align}
		& \| \nabla_\theta^1 \q L_y(\theta)^T (C_T - C_{k}) B_{k+1} \cdots B_{T}   (\hat{x}(\theta; y) -  x) \| \nonumber \\
		& \quad \le \norm{\nabla^1_\theta \q L_y(\theta)}_2 \|C_T - C_{k} \| \norm{ B_{k+1} \cdots B_{T}  (\hat{x}(\theta; y) -  x) }_2 \nonumber \\
		& \quad \le L_{\Gamma} \norm{\nabla^1_\theta \q L_y(\theta)}_2 \norm{x_{T-1} - x_{k-1}}_2 \norm{ B_{k+1} \cdots B_{T}  (\hat{x}(\theta; y) -  x) }_2 \eqqcolon \Delta_1 . \label{eq:Error1}
	\end{align}
	In the last step of \eqref{eq:Error1}, we used the Lipschitzness of $\nabla_{\theta} \Gamma$ to obtain the bound 
	\begin{equation}
		\|C_T - C_{k} \|_2 = \| \nabla_{\theta} \Gamma(x^{(T)} ;\theta) - \nabla_{\theta} \Gamma(x^{(k)};\theta) \|_2 \le L_{\Gamma} \| x_T - x_k \|_2.
	\end{equation}
	Then, we estimate the first term in \eqref{eq:LowerBoundStep1} using a telescopic sum, which leads to
	\begin{align}
		& \nabla^1_\theta \q L_y(\theta)^T  C_{T} B_{k+1} \cdots B_{T} 
		 (\hat{x}(\theta; y) -  x) \nonumber \\
		& \quad =  \nabla^1_\theta \q L_y(\theta)^T  C_{T} \left[ B_{T}^{T-k}  + \sum_{j = k+1}^{T-1} B_{k+1} \dots B_{j-1} (B_j - B_T) B_{T}^{T-j} \right] (\hat{x}(\theta; y) -  x)  \nonumber \\
		& \quad \geq \nabla^1_\theta \q L_y(\theta)^T  C_{T} B_{T}^{T-k} 
		 (\hat{x}(\theta; y) -  x)   - \Delta_2, \label{eq:LowerBoundStep2}
	\end{align}
	where
	\begin{equation}\label{eq:Error2}
		\Delta_2 \!=\!  L_{\Gamma} \norm{\nabla^1_\theta \q L_y(\theta)^T  C_{T} }_2 \norm{ \hat{x}(\theta; y) -  x}_2 \!\!\sum_{j=k+1}^{T-1}\norm{x_{T-1} - x_{j-1}} \norm{ B_{k+1} \cdots B_{j-1}} \norm{ B_{T} }^{T-j}.
	\end{equation}
	As last step, we estimate the first term in \eqref{eq:LowerBoundStep2} as
	\begin{equation}
		\nabla^1_\theta \q L_y(\theta)^T  C_{T} B_{T}^{T-k} 
		 (\hat{x}(\theta; y) -  x) \geq (1 - \nu \alpha )^{T-k} \Vert \nabla^1_\theta \q L_y(\theta) \Vert_2^2  - \Delta_3,
	\end{equation}
	where with $\nabla^1_\theta \q L_y(\theta) = C_T (\hat{x}(\theta; y) -  x)$ we have
	\begin{equation}\label{eq:Error3}
		\Delta_3 =  \norm{\nabla^1_\theta \q L_y(\theta)^T  C_{T} }_2 \norm{ \hat{x}(\theta; y) -  x}_2 \Vert B_T - (1 - \nu \alpha) I_N \Vert^{T-k}.
	\end{equation}
	It remains to bound the three terms $\Delta_1$, $\Delta_2$ and $\Delta_3$. 
	For this, we use the improved version of \Cref{th:grad descent linear} from \Cref{ex: stongly convex via Banach} that yields 
	\begin{equation}\label{eq:ImprovedEst}
		\|B_{k} \| \le 1 - \nu \alpha
		\quad \text{and} \quad 
		\| x_{k} - \hat x(\theta;y) \|_2 \le (1 - \nu \alpha)^{k} \| x_0 - \hat x(\theta;y) \|.
	\end{equation}
	Now, for $\Delta_1$ as in \eqref{eq:Error1}, we obtain with \eqref{eq:ImprovedEst} that
	\begin{align}
		\!\! \Delta_1 
		\! & \! \leq L_{\Gamma} \norm{\nabla^1_\theta \q L_y(\theta)}_2\norm{ x_{T-1} - x_{k-1}}_2 \norm{\hat{x}(\theta; y) -  x }_2 (1 - \nu\alpha)^{T-k}  \notag \\
		& \! \leq L_{\Gamma} \norm{\nabla^1_\theta \q L_y(\theta)}_2 \big( \norm{x_{T-1} - \hat{x}(\theta; y)} + \norm{ x_{k-1} - \hat{x}(\theta; y) }_2 \big) \norm{\hat{x}(\theta; y) -  x }_2 (1 - \nu\alpha)^{T-k}  \notag\\
		& \! \leq 2L_{\Gamma}\norm{\nabla^1_\theta \q L_y(\theta)}_2  \norm{\hat{x}(\theta; y) -  x_0 }_2 \norm{\hat{x}(\theta; y) -  x }_2 (1 - \nu\alpha)^{T-1}  \notag\\
		& \! \leq 2 L_{\Gamma} e^{-\nu \alpha (T-1)} \norm{C_T}_2 \norm{\hat{x}(\theta; y) -  x }_2^2 \norm{ \hat{x}(\theta; y) -  x_0}_2,
	\end{align}
	where we used $(1-x) \le e^{-x}$ in the last step.
	Similarly, we obtain for $\Delta_2$ as in \eqref{eq:Error2} that
	\begin{align}
		\Delta_2 & 
		\leq \norm{ C_{T} }^2_2 \norm{ \hat{x}(\theta; y) -  x}^2_2   \sum_{j=k+1}^{T-1} 2 L_{\Gamma} \norm{ \hat{x}(\theta; y) -  x_0}_2 (1 - \nu \alpha)^{(j-1) + (j - k-1) + (T-j)} \notag\\
		& \le 2 L_{\Gamma} \norm{ C_{T} }^2_2 \norm{ \hat{x}(\theta; y) -  x}^2_2 \norm{ \hat{x}(\theta; y) -  x_0}_2 e^{-\nu \alpha (T-1)}\ \sum_{j=0}^{T-k-2} e^{-\nu \alpha j}\notag\\
		& \leq 2 L_{\Gamma} \frac{e^{- \nu \alpha(T-1)}}{1 - e^{-\nu \alpha}} \norm{C_{T} }^2_2 \norm{ \hat{x}(\theta; y) -  x}^2_2 \norm{ \hat{x}(\theta; y) -  x_0}_2.
	\end{align}
	Turning to $\Delta_3$ as in \eqref{eq:Error3}, we bound the eigenvalues of $B_T - (1-\nu \alpha) I_N$ using the relation \eqref{eq: B eigenvalue bound} from \Cref{ex: stongly convex via Banach} as
	\begin{equation}
		\nu \alpha - \nu L =  1 - \nu L - (1-\nu \alpha) \le \lambda_n(B_T - (1-\nu \alpha) I_N) \le 1-\nu \alpha- (1-\nu \alpha) \le 0. 
	\end{equation}
	Thus, we finally obtain that
	\begin{equation}
		\Delta_3 \leq  \norm{C_{T} }_2^2 \norm{ \hat{x}(\theta; y) -  x}_2^2 ( \nu L - \nu \alpha)^{T-k}.
	\end{equation}
	In summary, we get that 
	\begin{align}\label{eq:FinalEstAlign}
		&\nabla^1_\theta \q L_y(\theta)^T   C_{k} B_{k+1} \cdots B_{T}   (\hat{x}(\theta; y) -  x) \geq (1 - \nu \alpha )^{T-k} \Vert \nabla^1_\theta \q L_y(\theta) \Vert_2^2  - \Delta, \notag\\
		&\Delta \coloneqq \Delta_1 + \Delta_2 +\Delta_3 = \norm{ \hat{x}(\theta; y) -  x }_2^2 \ \q O\Bigl( \frac{e^{-\nu \alpha(T-1)}}{1- e^{-\nu \alpha}} +   \left(  \nu (L-\alpha)  \right)^{T-k}  \Bigr).
	\end{align}
	Therefore, we have the estimate
	\begin{align}
		d_1
		\geq  \norm{\nabla^1_\theta \q L_y (\theta)}_2^2 -   \norm{ \hat{x}(\theta; y) -  x}_2^2 \ \q O \biggl(  \sum_{k=0}^{T-1}  \frac{e^{-\nu \alpha(T-1)}}{1- e^{-\nu \alpha}} +   \left(  \nu (L-\alpha)  \right)^{T-k}\biggr).
	\end{align}
	Because 
	\begin{align}
		\sum_{k=0}^{T-1}\left(  \nu (L-\alpha)  \right)^{T-k} 
		&= \sum_{k=1}^{T}\left(  \nu (L-\alpha)  \right)^k \leq  \frac{\nu (L-\alpha)}{1-\nu (L-\alpha)}
		&\text{for $\nu \leq \frac1L$}
	\end{align}
	and  since $C_T$ is non-singular by assumption, we obtain that
	\begin{align}\label{eq:AlignmentLowerBound}
	d_1 & \geq  \sigma^2_{\min}(C_T) \norm{ \hat{x}(\theta; y) -  x}_2^2
		-   \norm{\hat{x}(\theta; y) -  x}_2^2 \ \q O \Bigl(   \frac{ T e^{-\nu\alpha(T-1)}}{1- e^{-\nu\alpha}}  +  \frac{\nu (L-\alpha)}{1-\nu (L-\alpha)} \Bigr)\notag\\
		 & \geq c  \norm{\hat{x}(\theta; y) -  x}_2^2 
	\end{align}
	for some $c > 0$ if $T$ is large enough and $\nu$ is small enough.
\end{pr}
\begin{remark}
	In practice, we want to choose $T$ large anyways, so that the requirement is no real restriction.
	If $R_\theta$ is convex, then the denoising objective $\Phi_y(x;\theta)= \frac12 \Vert x - y \Vert_2^2 + \alpha R_\theta(x)$ is 1-strongly convex in $x$.
	Then, the optimal stepsize is $\nu= 2/(L+ 1)$ with corresponding contraction factor $1 - \frac{2\mu}{L+\mu} < 1 -\frac{\mu}{L}$, see \Cref{ex: stongly convex via Banach}.
	If the contraction factor is small, the first summand in \eqref{eq:AlignmentLowerBound} vanishes quickly as $T$ increases.
	However, we need to choose $\nu$ smaller than optimal to decrease the second summand.
	This seems artificial and might be  an artifact of the deployed proof techniques.
	Still, \Cref{th:bound of many terms} indicates that truncated backpropagation should work particularly well for well-conditioned problems where $\mu \approx L$.
	Since checking the assumptions of \Cref{th:exp convergence} and \Cref{th:bound of many terms} is infeasible in practice, these should merely serve as motivation.
\end{remark}

\begin{example}
Let us consider the optimization of the training loss $\q L_y$ with plain gradient descent using the constant learning rate $\tilde \nu >0$ (not to be confused with the step size $\nu$ in $\Gamma(x,\theta)$), namely
\begin{equation}
\theta^{(t+1)} = \theta^{(t)} - \tilde \nu \nabla^K_\theta \q L_y (\theta^{(t)}).
\end{equation}
Assume that $\q L_y$ has $\tilde L$-Lipschitz gradient. 
Then, the Descent \Cref{l: descent lemma} gives
\begin{equation}\label{eq:DescentEstimate}
\q L_y(\theta^{(t+1)}) \le \q L_y(\theta^{(t)}) - \tilde \nu \bigl\langle \nabla^K_\theta \q L_y(\theta^{(t)}),\nabla_\theta \q L_y(\theta^{(t)}) \bigr \rangle
+ \frac{\tilde L \tilde \nu^{2}}{2} \| \nabla^K_\theta \q L_y(\theta^{(t)}) \|_2^2.
\end{equation}
Under the assumptions of \Cref{th:bound of many terms}, we have $\nabla^K_\theta \q L_y(\theta)^T \nabla_\theta \q L_y(\theta) \geq  c \norm{\hat{x}(\theta; y) -  x}_2^2$ and similarly to \eqref{eq:convex bound} the bound 
$
\| \nabla^K_\theta \q L_y(\theta) \|_2 \le \frac{M_C(\theta^{(t)})}{\nu \alpha} \| \hat{x}(\theta^{(t)}; y) -  x \|_2
$
applies. Inserting these into \eqref{eq:DescentEstimate} yields
\begin{equation}
\q L_y(\theta^{(t+1)}) \le \q L_y(\theta^{(t)}) - \tilde \nu \Bigr(c - \tilde \nu \frac{\tilde L M_C^2(\theta^{(t)}) }{2 \nu^2 \alpha^2} \Bigl) \underbrace{ \| \hat{x}(\theta^{(t)}; y) -  x \|_2^2 }_{= 2 \q L_y(\theta^{(t)})}.
\end{equation}
Thus, for $\tilde \nu$ small enough (depending also on $t$), we get linear convergence of the loss $\q L_y$.
\end{example}

A similar analysis as the one in \Cref{th:exp convergence} and \Cref{th:bound of many terms} for the relation between the JFB and IFT approach was recently carried out in \cite{BolPauVai2023}.
However, one can in general not expect that JFB gets arbitrarily close to the IFT gradient.
Here, we only establish a simple relation between truncated backpropagation and the IFT approach.
For this, we recall that $\hat x(\theta;y)$ is a fixed point of the iteration \eqref{eq:IterativeUpdate}, namely that
\begin{equation}
\hat x(\theta;y) = \Gamma(\hat x(\theta;y), \theta).
\end{equation}
Hence, similarly to \eqref{eq: backprop unrolling derivation}, its gradient can be computed recursively as 
\begin{align}\label{eq:DerivativeFixedPoint}
	\nabla_\theta \hat x(\theta;y)
	= \nabla_\theta \Gamma(\hat x(\theta;y);\theta) 
	= \nabla_\theta \Gamma(x;\theta)\big|_{x =\hat x(\theta;y)} + \nabla_\theta \hat x(\theta;y) \nabla_x \Gamma(x;\theta)\big|_{x= \hat x(\theta;y)}.
\end{align}
Rearranging  \eqref{eq:DerivativeFixedPoint} and setting $B_{\infty} = \nabla_{x} \Gamma(x;\theta)\big|_{x =\hat x(\theta;y)}$ and $C_{\infty} = \nabla_{\theta} \Gamma(x;\theta)\big|_{x =\hat x(\theta;y)}$ gives
\begin{equation}
\nabla_\theta \hat x(\theta;y) = C_{\infty} (I - B_{\infty})^{-1}.
\end{equation}
On the other hand, repeating the chain rule $K$ times leads to 
\begin{equation}
\nabla_\theta \hat x(\theta;y) = C_\infty \sum_{k=0}^{K-1} B_\infty^{k} + B_\infty^{K} \nabla_\theta \hat x(\theta;y).
\end{equation}
\begin{thm}{Neumann Series}{neumann}
	 If $\Vert B_{\infty} \Vert_2 < 1$, the Neumann series $\sum_{k=0}^{\infty} B_\infty^k$ is convergent with limit $(I-B_\infty)^{-1}$.
\end{thm}
\begin{pr}
	Using the geometric series, we can show that $S_K = \sum_{k=0}^{K} B_\infty^k$, $K \ge 0,$ is a Cauchy sequence, and, thus, converges absolutely to some $S_\infty$.
	Moreover, it holds that $\lim_{k \to \infty}  B_\infty^k = 0$.
	Thus,
	\begin{equation}
		(I-B_\infty) S_\infty = \lim_{K \to \infty} (I-B_\infty) S_K = \lim_{K \to \infty} \sum_{k=0}^{K} B_\infty^k -  B_\infty^{k+1} =  \lim_{K \to \infty} I - B_\infty^{K+1} = I.
	\end{equation}
	Similarly, we can show that $S_\infty$ is a left inverse.
\end{pr}
Thus, close to the equilibrium $\hat{x}(\theta, y)$, we can approximate \eqref{eq:TruncGrad} with the Neumann series
\begin{equation}\label{eq:NeumannGrad}
	\nabla_\theta \hat x(\theta;y)
	=  C_\infty \sum_{k=0}^{\infty} B_\infty^{k} \approx  C_\infty \sum_{k=0}^{K-1} B_\infty^{k}.
\end{equation}
The formula \eqref{eq:NeumannGrad} is known as Neumann backpropagation \cite{LiaXioFet2018}.
The key difference from the truncated backpropagation formula \eqref{eq:TruncGrad} is that \eqref{eq:NeumannGrad} evaluates the matrices $B_t$ and $C_t$ only at $\hat{x}(\theta, y_i)$, rather than using the iterates from the last $K$ steps of \eqref{eq:IterativeUpdate}. Hence, we don't need to store the iterates $x_t$.

To understand the behavior of Neumann backpropagation, we consider the case of unrolled gradient descent \eqref{eq:GradUnrolling}.
Here, the matrix $B_{\infty} = I - \nu \nabla^2_{x,x} \Phi_y(\hat{x}(\theta, y);\theta)$ determines the contraction properties.
If $\| B_{\infty} \|_2 < 1$, then the Neumann series in \eqref{eq:NeumannGrad} converges to the implicit function theorem (IFT) gradient \eqref{eq:formula_grad}. 
This condition holds provided that the step size $\nu$ is small enough and that the Hessian $\nabla^2_{x,x} \Phi_y(\hat{x}(\theta, y);\theta)$ is positive definite.
The latter is guaranteed when $\hat{x}(\theta, y)$ is a local minimizer of $x \mapsto \Phi_y(x;\theta)$. 
The finite series is only accurate approximation for large enough $K$.
In practice, even $K=1$ often works surprisingly well.
However, moderately increasing $K$ to values like 3–5 can improve gradient accuracy without substantial computational overhead.

\subsubsection{Exercises}

\begin{exerciseB}\label{ex: solution lipschitz continuity}
	Work out the details of \Cref{eq:RemLip} using the implicit function theorem.
\end{exerciseB}

\begin{Solution}
	Let us compute
	\begin{equation}
		\nabla_{x} \Phi(x,y) = 2 A^\tT(Ax - y) + \nabla_{x} R(x),
		\quad \nabla_{x,x}^2 \Phi(x,y) = 2 A^\tT A + \nabla^2_{x,x} R(x)
		\quad \text{and} \quad \nabla_{x,y}^2 \Phi(x,y) = 2 A.
	\end{equation}
	Since $A$ is injective, $A^\tT A$ is inverible positive definite matrix. Therefore, $\nabla_{x,x}^2 \Phi(x,y)$ is independent of $y$ and invertible for all $x\in \R^N$, $y \in \R^M$. Then, conditions of \Cref{th:ImplicitFunction} are satisfied and \eqref{eq:JacS} implies $\nabla_y \hat x(y) = - A (\nabla_{x,x}^2 \Phi(x,y))^{-1}$.
	
	Therefore, for given $y_1, y_2 \in U$ we get
	\begin{align}
		\Vert \hat x(y_1) - \hat x(y_2) \Vert_2 \leq& \sup_{y \in U} \Vert \nabla_y \hat x (y) \Vert_{2} \Vert y_1- y_2\Vert_2
		= \bigl\Vert A (\nabla_{x,x}^2 \Phi(x,y))^{-1} \bigr\Vert_{2} \Vert y_1- y_2\Vert_2 
	\end{align}
	The matrix $\nabla_{x,x}^2 \Phi(x,y) - A^\tT A = \nabla^2_{x,x} R(x)$ is symmetric positive semidefinite since $R$ is convex.
	Moreover, $A^\tT A$ and $\nabla_{x,x}^2 \Phi(x,y)$ symmetric positive definite and so are their inverses.
	Therefore, 
	\begin{equation}
		( \nabla_{x,x}^2 \Phi(x,y)^{-1} - (A^\tT A)^{-1}
		= \nabla_{x,x}^2 \Phi(x,y)^{-1} ( \nabla_{x,x}^2 \Phi(x,y) - A^\tT A) (A^\tT A)^{-1} 
	\end{equation}
	is symmetric positive semidefinite as a product of symmetric positive semidefinite matrices.
	Consequently, 
	for symmetric positive semidefinite matrix $A ( A^\tT A)^{-2} A^\tT$ with the eigenvector $z$, $\|z\|_2 =1$, corresponding to the largest eigenvalue, we get
	\begin{align}
		\| ( A^\tT A)^{-1} A^\tT \|_2^2
		&= \| A ( A^\tT A)^{-2} A^\tT \|_2
		= z^\tT A ( A^\tT A)^{-2} A^\tT z 
		\le z^\tT A (\nabla_{x,x}^2 \Phi(x,y))^{-2} A^\tT z \notag \\
		& \le \| A (\nabla_{x,x}^2 \Phi(x,y))^{-2} A^\tT\|_2
		= \| A (\nabla_{x,x}^2 \Phi(x,y))^{-1} \|_2^2.
	\end{align}
	In view of \Cref{ex:ReverseStatementInverse} (i), we obtain
	\begin{equation}
		\Vert \hat x(y_1) - \hat x(y_2) \Vert_2 \le \| A ( A^\tT A)^{-1} \|_2 \| y_1 - y_2 \|_2 
		= \| A^\dagger \|_2 \| y_1 - y_2 \|_2,
	\end{equation}
	which means that $\hat x \colon U \subset \R^M \to \R^N$ is Lipschitz continuous with Lipschitz constant at most $\Vert A^{\dagger}\Vert_2$.
	
\end{Solution}

\begin{exerciseB}
	Prove \eqref{eq:GradUnrolling} using the chain rule (assuming that all functions are sufficiently smooth).
\end{exerciseB}

\begin{exerciseB}\label{ex: stongly convex via Banach}
	Consider a $\mu$-strongly convex function $\Phi_y(\cdot, \theta) \colon \R^N \to \R$ with $L$-Lipschitz continuous gradient, and let $x \in \R^N$ be arbitrary.
	Show that the eigenvalues $\lambda_n$, $n = 1,\ldots, N$, of $\Gamma(x; \theta) = x - \nu \nabla_x \Phi_y(x,\theta)$ with $0 < \nu < 2/L $ satisfy
	\begin{equation}\label{eq: B eigenvalue bound}
	1 - \nu L \le \lambda_n \le 1 - \nu \mu 
	\quad \text{and} \quad \|\nabla_x \Gamma(x; \theta)\|_{2} \le \max\{1 - \nu \mu, \nu L -1\}<1.
	\end{equation}
	Based on \Cref{th:BanachFix}, establish \Cref{th:grad descent linear} without the squares on both sides of \eqref{eq:ConvSpeedStronglyConvex}.
	Further, infer that the optimal stepsize is $\nu = \frac{2}{L+\mu}$.
\end{exerciseB}

\begin{Solution}
	%%% Solution
	The gradient of it is a symmetric matrix $\nabla_x \Gamma(x; \theta) = I - \nu \nabla_{x,x} \Phi_y(x,\theta)$. Strong convexity of $\Phi_y(\cdot, \theta)$ and Lipschitz continuity of $\nabla_x \Phi_y(\cdot, \theta)$ can be expressed in the terms of Hessian as
	\[
	\mu \|z\|_2^2 \le \inner{\nabla_{x,x} \Phi_y(x,\theta) z}{z} \le L \|z\|_2^2, \quad \text{for all} \ z \in \R^N.
	\]
	Then, the eigenvalues of $\nabla_x \Gamma(x; \theta)$ are bounded as
	\[
	(1 - \nu L) \|z\|_2^2 \le \inner{\nabla_x \Gamma(x; \theta) z}{z} \le (1 - \nu \mu) \|z \|_2^2. 
	\]  
	If $\nu \le 1/L$, then $\nu L -1 < 0$ yielding $\|\nabla_x \Gamma(x; \theta)\|_{2} \le 1 - \nu \mu$ and \eqref{eq: B eigenvalue bound}. If $1/L < \nu < 2/L$, then $\|\nabla_x \Gamma(x; \theta)\|_{2} \le \max\{1 - \nu \mu, |1- \nu L| \}$, which is exactly \eqref{eq: B eigenvalue bound}. 
	Banach fixed point theorem then recovers the results of \Cref{th:grad descent linear}. Furthermore, the rate is optimal, if elements in the maximum are equal, which gives $2 = \nu(\mu + L)$ and the stated step size.
\end{Solution}

\subsection{Regularizer Architectures}\label{sec: reg architecture}
To conclude this section, we discuss a few regularizer architectures from the literature.
\subsubsection{Field-of-Experts}
\begin{figure}[tb]
	\centering
	\includegraphics[width=\linewidth]{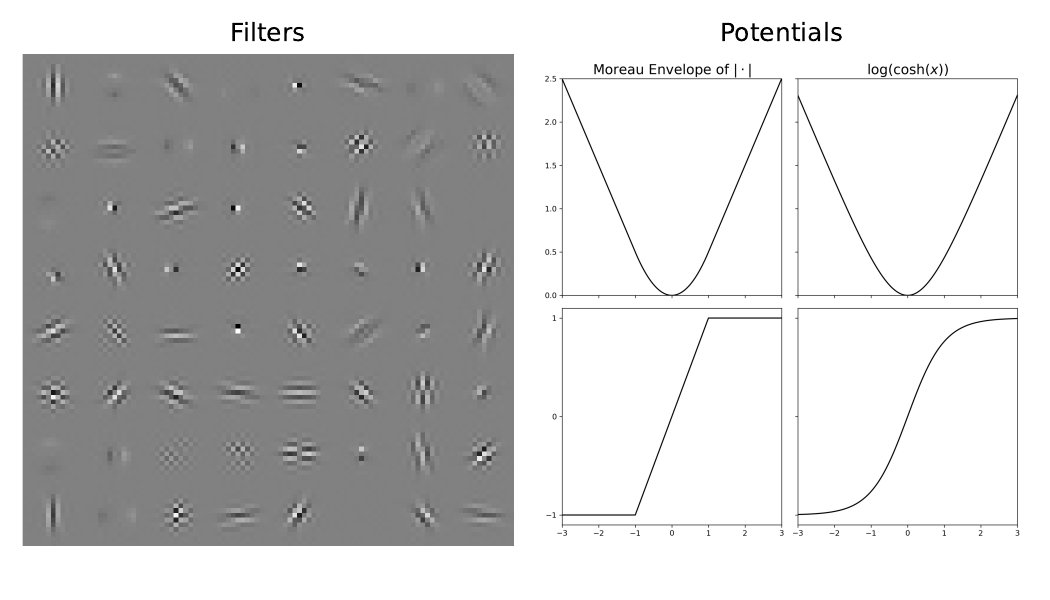}
	\caption{Left: The 64 convolution kernels of a learned FoE.
		Black corresponds to negative and white to positive values.
		Right: The two discussed potentials $\psi^1$ and their derivative $\varphi^1=(\psi^1)'$.
		Both approximate the absolute value.}
	\label{fig:RR}
\end{figure}

The \emphi{field-of-experts} regularizer (FoE) \cite{GouNeuBoh2022,RotBla2009} takes the form
\begin{equation}\label{eq:RegBase}
	R(x) = \sum_{j=1}^{c} \langle \1_{N}, \psi_j(W_{j} x)\rangle.
\end{equation}
%In \eqref{eq:RegBase}, 
For each of the $c$ summands, the potentials $\psi_j\colon \R \to \smash{\R_+}$ are applied componentwise and the $W_{j} \colon \R^N \to \R^N$ are convolution operators.
Hence, we can interpret $R$ as a spatial penalization of multiple filter responses.
As a special case, we obtain for $c = 2$, $W_1 = \nabla_1$, $W_2 = \nabla_2$, and $\psi_1= \psi_2 = \vert \cdot \vert$ the total variation regularization from \Cref{exer:impainting}.
If all $\psi_j$ are differentiable, the same holds for $R$ and we get that
\begin{equation}
	\label{eq:gradmodel}
	\nabla R(x) = \sum_{j=1}^c W_j^T \psi_j'(W_j x).
\end{equation}
In \Cref{fig:RR}, we visualize a possible specification of the FoE model \eqref{eq:RegBase}, which involves learnable $W_j$ and specifically chosen $\psi_j$.
The two depicted $\psi$ are approximations of $ \vert \cdot \vert$.

To facilitate the learning of the stacked convolution operator $W= (W_j)_{j=1}^c$, we parameterize it as a multi-convolution.
More precisely, we decompose $W$ as $W = V^3 V^2 V^1$ with 3 zero-padded (stacked) convolutions with kernels of size $k_s\times k_s$, $k_s$ odd, and an increasing number of output channels.
Recall that $W$ defined this way is indeed a (stacked) convolution by \Cref{ex: circulant matrix} (iii).
For the intermediate convolutions $V^i$, we use $4$, $8$ and $64$ channels, respectively, resulting in a total of $c=64$ kernels.
Later, we will use 3 kernels of size $k_s = 5$, which allows to explore the effective field of view ($13 \times 13$) more efficiently.
Moreover, the convolution kernels of $V^1$ are constrained to have zero mean, leading to $W \1_d = 0$.
This ensures that we do not penalize constant shifts in the images.
To avoid scaling ambiguities between the $W_j$ and $\psi_j$, we enforce $\Vert (W_j)_{j=1}^c \Vert_2= 1$ via the normalization $(W_j)_{j=1}^c /b_\infty$ based on the power iterations
\begin{equation}\label{eq:PowerIters}
	b_{k+1} = \frac{ W^T W b_k}{ \| W^T W b_k \|_2} = \frac{\sum_{j=1}^c W_j^T W_j b_k}{\bigl \Vert \sum_{j=1}^c W_j^T W_j b_k\bigr \Vert_2}, 
\end{equation}
for some $b_0\in \R^{N}$, $\|b_0\|_{2} = 1$, preferably random to avoid bad starting points. See \Cref{Ex:PowerIters} for details on the power method.
These can be embedded into the computational graph for finding the lower-level solution based on Jacobian-free backpropagation.
Empirically, we observe that the kernel $K$ corresponding to the learned $\sum_{j=1}^c W_j^T W_j$ is often close to a shifted $2$D discrete Kronecker impulse.
The shift occurs since $K$ must have mean zero.

Regarding $\psi_j$, we impose ${\psi_j}(t) = \smash{1/\alpha_j^2 \psi^\beta(\alpha_j t)}$, where the learnable $\alpha_j \in \R$ adapt a shared potential $\psi^\beta$.
In practice, we use either $\psi^\beta( x) = \beta^{-1}\log(\cosh(\beta x))$ or the Huber function
\begin{equation}\label{eq:Huber}
	\psi^\beta(x) = \begin{cases}
		\vert x \vert  - \frac{\beta^{-1}}{2}, & \vert x \vert > \beta^{-1},\\
		\frac{\beta}{2}x^2, & \vert x \vert \leq \beta^{-1}.
	\end{cases}
\end{equation}
This simplification facilitates the interpretability (see \Cref{fig:RR}) and reduces the number of learnable parameters drastically.
The pre-factor $1/\alpha_j^2$ ensures that the (weak) derivative $\psi_j''$ remains independent of $\alpha_j$.
This is key for convexity investigations based on \eqref{eq:ConvexHess}.
\begin{lemma}\label{lem:CvxCRR}
	If $\psi^\beta \in C^2(\R)$ and $(\psi^\beta )'' \geq 0$, then $R$ is convex.
	Moreover, if $(\psi^\beta )'' \geq -1$ and $\Vert W \Vert_2=1$, then $\inner{\nabla^2 R(x) z}{z} \ge - \|z\|^2_2$ for all $z \in \bb R^N$.
\end{lemma}
\begin{pr}
	It holds $\nabla^2 R(x) = \sum_{j=1}^c W_j^T \mathrm{diag}(\psi_j''(W_jx))W_j$.
	From this, we directly infer convexity.
	Moreover, it holds that
	\begin{equation}
		\inner{\nabla^2 R(x) z}{z} =  \sum_{j=1}^c \inner{\mathrm{diag}(\psi_j''(W_jx))W_jz}{W_jz} \geq -   z^T \sum_{j=1}^c W_j^T W_j z  \geq - \Vert z \Vert_2^2.
		\end{equation}
		This concludes the proof.
\end{pr}
If $R$ is convex, the same holds for \eqref{eq:data_fidelity_regularizer2} and the global convergence theory from \Cref{sec:Algs} applies.
For the case of denoising with $A = I
$, \eqref{eq:data_fidelity_regularizer2}  remains convex  even if $\inner{\nabla^2 R(x) z}{z} \ge -\|z\|^2_2$ for all $z \in \bb R^N$ since the data fidelity $\Vert x - y \Vert_2^2$ is 1-strongly convex.
For this setting, we instead use the potentials $\tilde \psi^\beta(x) = \psi^\beta(x) - \psi^1(x)$ with $\beta\geq 0$ and one of the $\psi$ from above.

\subsubsection{ICNN}
Several approaches to learning $R$ leverage \emphi{Input Convex Neural Network}s (ICNNs) $R\colon \R ^N \to \R$  \cite{AXK2016}. 
These take the form $R( x)=z_\ell(x)$, where $z_\ell(x)$ is defined recursively as 
\begin{equation}\label{eq:feedforward}
	 z_1(x) = \psi_0(V_0 x+ b_0),\quad z_{i+1}(x) = \psi_i(W_i z_i(x)+V_i x+ b_i), \quad i=1,\dots,\ell-1,
\end{equation}
where 
$\theta=\{(W_i)_{i=1}^{\ell-1},(V_i)_{i=0}^{\ell-1},( b_i)_{i=0}^{\ell-1}\}$ are learnable parameters, and $\psi_i$ are componentwise non-linear activation functions.
\begin{lemma}\label{l: ICNN}
	Let $\psi_0 \colon \R \to \R$ be convex, $\psi_i \colon \R \to \R$, $i \geq 1$, be convex and non-decreasing, and $W_i$ have non-negative entries.
	Then, $R \colon \R^N \to \R$ as defined by \eqref{eq:feedforward} is a convex function.
\end{lemma}
\begin{pr}
	The proof uses two well-known properties of convex functions.
	First, it holds that $\sum_{i=1}^K a_i f_i(x)$ is convex if all $f_i$ are convex and $a_i \ge 0$.
	Second, it holds that $f_1 \circ f_2$ is convex if $f_1\colon \R \to \R$ is convex and non-decreasing, and $f_2 \colon \R \to \R$ is convex.
	Both properties can be shown directly using the definition.
	Similarly, we verify that each component of $z_1(x)$ is convex.
	Now, we proceed inductively and assume that each component of $z_i(x)$ is convex.
	Then, it holds that by the first property that each component of $W_i z_i(x)+V_i x+ b_i$ is convex.
	Since $\psi_i$ is convex and non-decreasing, the second property readily implies that every component of $z_{i+1}(x)$ is convex.
\end{pr}
\begin{example}
	If we set $\ell = 2$, $b_0 = b_1 = 0$, $W_1 = \1_{Nc}^T$, $V_1 = 0$, $\psi_1(x) = x$, and if $V_0$ is a (stacked) convolution matrix, then the ICNN regularizer \eqref{eq:feedforward} reduces to FoE \eqref{eq:RegBase}.
	More precisely, we get
	\begin{equation}
		z_1(x) = \psi_0(V_0 x)  
		\quad \text{and} \quad
		z_2(x) = \1_{Nc}^T z_1(x) = \sum_{j=1}^c \inner{\1_N}{\psi_0( (V_0)_j x)}.
	\end{equation}
	Hence, we can see the ICNN as a generalization of the convex FoE regularizer.
\end{example}
During training, we enforce the non-negativity of $W_i$ by projecting each entry onto the non-negative numbers after every update of \Cref{algo:stochastic_gradient_descent} (or its generalizations).
As for $\psi_i$, common choices are the softplus, ReLU and leaky ReLU activations.
To utilize the full capcity of the ICNN model, we can in principle also learn the activations $\psi_i$ as done in \cite{DGB2022}.

\begin{example}
	In \cite{MukDit2021}, the authors consider an  ICNN architecture with convolution layers $W_l$, $V_l$, $l=1,\ldots,\ell-2$, and global average pooling for $W_{\ell-2}$. 
	The functions $\psi_i$,  $l=0,\ldots,\ell-2$, are chosen to be the leaky ReLU with $a=0.2$, and $\psi_{\ell-1} (x) = x$.
	For their sparse-view CT setup, the authors recommended $\ell=10$ layers, and convolution kernels of size $5\times5$ with stride 1 and 48 output channels. 
\end{example}

\begin{example}
	The authors of \cite{fang2023s} utilize gradients of ICNNs for learning a non-convex $R$.
	More precisely, they implicitly learn $R$ via its proximal operator, and we will study this approach in \Cref{sec: PnP}.
\end{example}

\subsubsection{Least Squares Residual}
The role of the regularizer $R$ is at least twofold.
First, it should stabilize the reconstruction and reduce the noise; second, the $R$-minimizing solutions should be in a feasible set $\mathcal M$.
One way of characterizing $\mathcal M$ is as the fixed points of a map $U \colon \R^N \to \R^N$, acting as a nonlinear projection onto $\mathcal M$.
Then, penalizing the residuals $x - U(x)$ yields a regularizer
\begin{equation}\label{eq:LSR}
	R(x) = \frac12 \Vert x - U(x) \Vert^2, 
\end{equation}
which we refer to as the \emphi{Least Squares Residual Regularizer} (LSRR).
The gradient of the regularizer \eqref{eq:LSR} is given by 
\begin{equation}\label{eq:GradLSRR}
	\nabla R(x) =  (I_N - J_U (x))^T (x - U(x)),
\end{equation}
where the matrix-vector product with $J_U^T(x)$ can be computed using backpropagation if the activation functions in $U$ are smooth.
Typically, $U$ is chosen as an encoder-decoder architecture such as the UNet \cite{RonFisBro2015} or DRUNet \cite{ZhaLiZuo2022}, which we introduced in \Cref{chap:regularization}.
Since we cannot estimate the Lipschitz constant of $\nabla R$ efficiently, a line search is inevitable when minimizing \eqref{eq:lower-level}.
Commonly, the architecture $U$ is initialized as a pretrained image denoising or restoration network \cite{ZouLiuWoh2023}. 
The architecture \eqref{eq:LSR} has been used for sparse-view CT in \cite{LiSch2020}, and in combination with unrolling for MRI reconstruction in \cite{aggarwal2018modl}.

\begin{remark}
	In the context of bilevel learning, the regularizer \eqref{eq:LSR} is usually pretrained based on the gradient-step denoising framework.
	More precisely, given a noise level $\sigma$, we train the reconstruction model $\mathcal R_\theta(y) = y - \sigma^2 \nabla R(y)$.
	Then, we fit the hyperparameter $\alpha$ in \eqref{eq:data_fidelity_regularizer2} and start the bilevel training from this configuration.
\end{remark}

\subsection{Learning a Convex Field-of-Experts}\label{eq:LearnFoE}
Now, we discuss the learning of the FoE architecture \eqref{eq:RegBase}, where we deploy ${\psi_j}(t) = \smash{1/\alpha_j^2 \psi^\beta(\alpha_j t)}$ with the convex Huber potential \eqref{eq:Huber} and $\alpha_j>0$.
Regarding $W= (W_j)_{j=1}^c$, we stick to the described configuration.
According to \Cref{lem:CvxCRR}, the resulting regularizer $R$ is convex, and we thus obtain the convex reconstruction problem
\begin{equation}\label{eq:RecProblemCRR}
	\argmin_{x\in\R^N} \frac12\Vert Ax - y \Vert_2^2 + \lambda \sum_{j=1}^{c} \langle \1_{N}, \psi_j(W_{j} x)\rangle.
\end{equation}

The Huber potential $\psi^\beta$ is piecewise smooth with (weak) derivatives
\begin{equation}
	(\psi^\beta)^\prime (x) = \begin{cases}
		\sign(x), & \vert x \vert > \beta^{-1},\\
		\beta x, & \vert x \vert \leq \beta^{-1}.
	\end{cases}
	\quad \text{and} \quad
	(\psi^\beta)''(x) = \begin{cases}
		0, & \vert x \vert > \beta^{-1},\\
		\beta, & \vert x \vert \leq \beta^{-1},
	\end{cases}
\end{equation}
so that also the Lipschitz constant of $\psi_j'$ is $L_{\psi_j'} = \beta$.
Hence, $R$ is differentiable with gradient
\begin{equation}
	\nabla R(x) = \sum_{j=1}^c W_j^T \psi_j'(W_j x).
\end{equation}
Now, to compute the Lipschitz constant of $\nabla R$, we bound the (weak) Hessian of $R$ for all $x,z \in \R^N$ as
\begin{align}
z^T \nabla^2 R(x) z
& = \sum_{j=1}^c z^T W_j^T \diag( \psi_j''(W_j x)) W_j z 
\leq \beta \sum_{j=1}^c z^T W_j^T W_j z  \le \beta \|z\|_2^2,
\end{align}
where we used that $\|W\|_2 = 1$ and $\psi_j'' \geq 0$.
Hence, $\nabla R$ is $\beta$-Lipschitz continuous, and we can minimize \eqref{eq:RecProblemCRR} efficiently using \cref{algo:accelerated_gradient_descent} with step size $\nu = (\Vert A \Vert^2_2 + \alpha\beta)^{-1}$.
In particular, we do not need to deploy a line search.

For the actual learning of $R_\theta$, we want to use a simple inverse problem for efficiency reasons.
The arguably simplest one is denoising with $A = I$.
To this end, let $\{x_i\}_{i=1}^S$ be a set of clean $(64\times64)$ patches that are extracted from the BSDS500 data set.
Each of these patches $x_i$ is corrupted with Gaussian noise as $ y_i= x_i + \sigma n_i$, where $n_i\sim \mathcal N(0, I)$ and $\sigma =25$.
Note that we do not need to corrupt the patches before training and simply do so when sampling a clean batch.
Now, we have defined all ingredients of the bilevel problem \eqref{eq:bilevel-prob}.
We set $\lambda =1$ in \eqref{eq:RecProblemCRR}, such that the learnable parameters are $W= (W_j)_{j=1}^c$, $\beta$ and $\{\alpha_j\}_{j=1}^c$.

\begin{figure}[tbp]
	\centering\includegraphics[width=\linewidth]{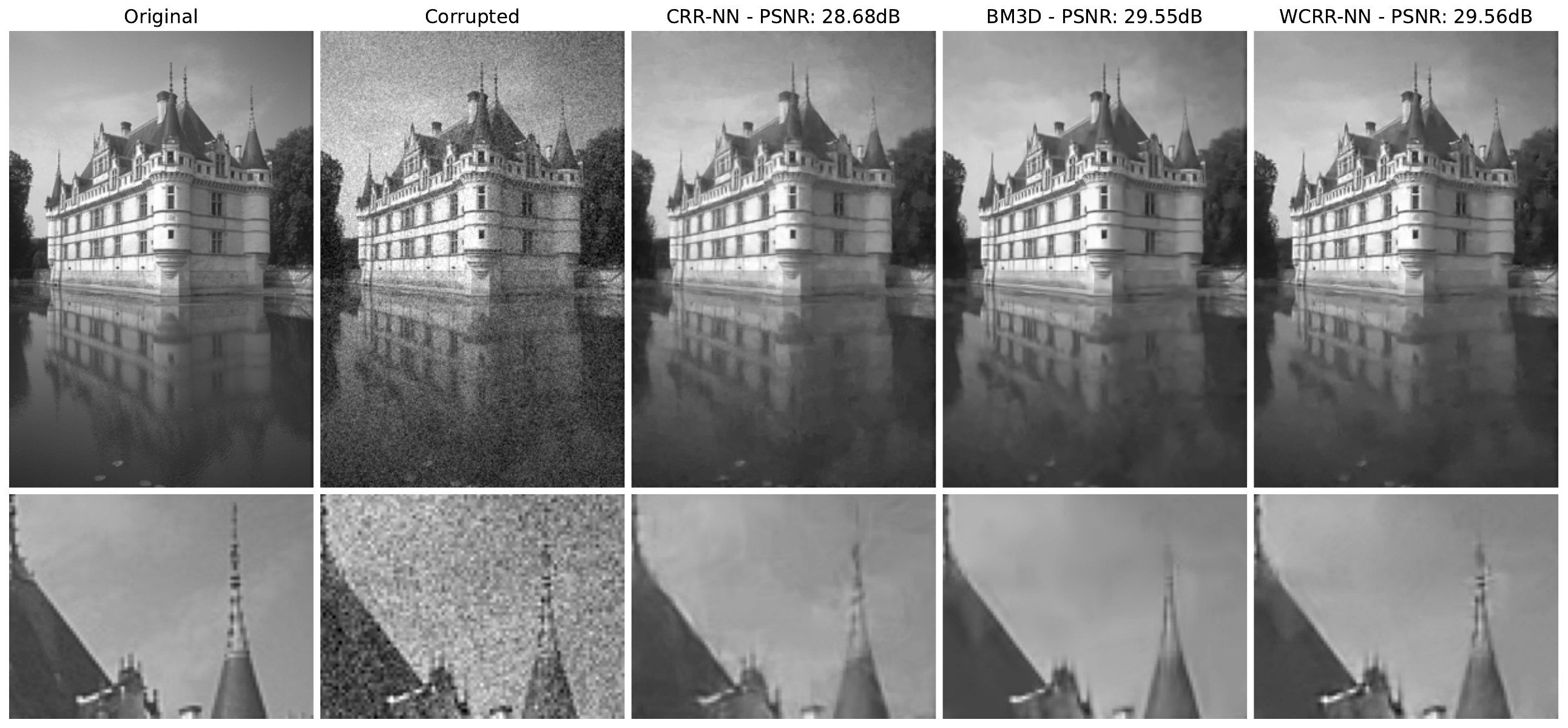}
	\caption{Denoising of the ``castle" image from the BSD68 test set for noise level $\sigma=25$.}
	\label{fig:DenoiseCompare}
\end{figure}

To keep the optimization of the bilevel objective \eqref{eq:bilevel-prob} as simple as possible, we use Jacobian-free backpropagation to estimate the hypergradient.
For the upper-level problem, we use the Adam optimizer with a batch size of 32 patches and with a learning rate of $0.01$.
This rate is decayed by a factor of $0.99$ after every epoch (namely after using all the training images once), and we train for 200 epochs in total.
We return the model with the best reconstruction PSNR on a predefined validation set.
A visual reconstruction result is given in \Cref{fig:DenoiseCompare}, where CRR-NN corresponds to our learned FoE, and WCRR-NN is the discussed weakly convex extension.
As we can see, this relaxation indeed improves the performance.
The BM3D denoiser is a classical non-learned baseline, which has no interpretation as variational energy minimization.
We should stress that denoising is only the training task and our main interest is to apply the learned $R$ to inverse problems.

So far, we learned $R$ for denoising with $\sigma=25$.
To solve a generic inverse problem based on \eqref{eq:RecProblemCRR}, in addition to the regularization strength $\lambda$, we propose to tune $\beta$.
This scaling parameter determines where the transition between the quadratic and linear behavior of the Huber potential $\psi^\beta$ occurs, which should be adapted to the problem.
The search for these two parameters can be carried out with a grid search on some (small) validation set.
In \Cref{ex:CRR}, we will investigate the generalization capabilities of a learned FoE regularizer based on a CT task.
While a universally trained $R$ usually works reasonable for any inverse problem, operator-specific training can increase performance, especially when structured nullspace components, such as CT undersampling artifacts, must be removed.

\subsubsection{Exercises}

\begin{exerciseB}\label{Ex:PowerIters}
	Assume that the two largest eigenvalues of $W$ satisfy $\vert \lambda_1 \vert> \vert\lambda_2\vert$.
	Prove that the power iterations \eqref{eq:PowerIters} converge and that they find $\Vert W \Vert_2 = \vert \lambda_1 \vert$. 
	How does the convergence speed depend on the ratio $\rho =\vert \lambda_2 \vert / \vert\lambda_1\vert$?
\end{exerciseB}

\begin{exerciseB}\label{ex:CRR}
	Train the convex FoE model as described in \Cref{eq:LearnFoE}.
	Then, deploy it to the CT task that we investigated for the learned post-processing in \Cref{ex:PostprocessingUNet} and compare the results.
	You may use the following Jupyter Notebook as a basis for this exercise, \url{https://github.com/Sniper2k/inverse_problems_ai_exercise_templates/blob/main/Ex6_3_2_template.ipynb}.
\end{exerciseB}

\begin{Solution}
See Jupyter Notebook \url{https://github.com/Sniper2k/inverse_problems_ai_exercises/blob/main/Ex6_3_2.ipynb}.
\end{Solution}

\section{Plug-and-Play algorithms} \label{sec: PnP}

Now, we discuss the \emphi{Plug-and-Play} (PnP) framework for solving the variational problem
\begin{equation}\label{eq:pnp_mot_opt_problem}
	\argmin_{x\in \R^N} \tfrac{1}{2\alpha}\|Ax - y\|^2_\bb Y + R(x).
\end{equation}
Historically, it was first proposed for the ADMM \cref{alg:ADMM}, see \cite{Venkatakrishnan_2013}.
To this end, we write \eqref{eq:pnp_mot_opt_problem} as
\begin{equation}
	\argmin_{x,z \in \R^N, x = z} \tfrac{1}{2 \alpha}\|Ax - y\|^2_\bb Y + R(z).
\end{equation}
Then, the augmented Lagrangian with the rescaled variable $b$ is given by
\begin{equation}
	L_\gamma(x,z,b) = \tfrac{1}{2\alpha}\|Ax - y\|^2_\bb Y + R(z) + \tfrac{\gamma}{2} \| x- z + b\|_2^2 - \tfrac{\gamma}{2} \| b \|_2^2,
\end{equation}
and the updates of the scaled ADMM (\Cref{alg:ADMM}) with $C = I$ and $\gamma >0$ read as 
\begin{equation}\label{eq:ADMMupdates}
	\begin{aligned}
		x_{k+1} & \in \argmin_{x \in \R^N} \bigl\{ \tfrac{1}{2\alpha}\|Ax - y\|^2_\bb Y +  \tfrac{\gamma}{2} \|x + b_{k} - z_{k}\|_2^2  \bigr\} \\
		z_{k+1} & \in \prox_{\tfrac{1}{\gamma} R}(b_{k} + x_{k+1}) = \argmin_{z \in \bb X} \bigl\{ R(z) + \tfrac{\gamma}{2} \|b_{k} + x_{k+1} - z\|_2^2  \bigr\} \\
		b_{k+1} &= b_{k} + x_{k+1} -  z_{k+1}.
	\end{aligned}
\end{equation}
Due to the modular structure of \eqref{eq:ADMMupdates}, we can replace some of the components in the updates by learned methods.
This is the basis of the PnP paradigm.
In particular, we note that the purpose of \smash{$\prox_{\frac{1}{\gamma} R}$} is to remove artifacts from $b_{k} + x_{k+1}$.
Hence, it appears reasonable to replace it by a denoiser $f_\theta \colon \R^N \to \R^N$, yielding the following \emphi{PnP-ADMM} algorithm.      

\begin{algorithm}[H]
	\KwIn{denoiser $f_\theta \colon \R^N \to \R^N$, initializations $z_{0}\in \R^N$ and $b_{0} \in \R^N$}
	\For{\color{black}$k = 0, 1, \ldots$}{
		$x_{k+1}\in \argmin_{x \in \R^N} \left\{ \frac{1}{2\alpha}\|Ax - y\|^2_\bb Y +  \frac{\gamma}{2} \|b_{k} + x - z_{k}\|_2^2  \right\}$\;
		$z_{k+1} =  f_{\theta}(b_{k} + x_{k+1})$\;
		$b_{k+1}= b_{k} + x_{k+1} -  z_{k+1}$\;
	}
	\caption{Plug-and-Play ADMM (PnP-ADMM) \label{alg:PnPADMM}}
\end{algorithm}
A similar strategy can be applied to other optimization algorithms involving proximal mappings, such as the FBS (\Cref{algo:FBS}).
Yet, FBS involves the step size $\nu$, and we need to enforce $\nu \leq 1 / \Vert A \Vert^2$ to ensure convergence, see \Cref{th:FBS}.
This often limits its performance, making ADMM preferable in practice \cite{hurault2022proximal}.
If we cannot ensure convergence, it remains unclear how many iterations we should perform in \Cref{alg:PnPADMM}.
Moreover, one typically sets $\gamma = 1$ so that $f_\theta \approx \prox_{R}$.

For an overview of PnP methods, we refer to the review \cite{kamilov2023plug}. 
In practice, Plug-and-Play algorithms often showcase improved reconstruction performance  compared to traditional methods for imaging applications such as tomography \cite{Sreehari_2016, eksioglu2016decoupled, He_2019} and diffraction imaging \cite{Sun_2019}.
However, due to the learned nature of the denoiser in the update $z_{k+1} = f_{\theta}(b_{k} + x_{k+1})$, it is hard to theoretically proof convergence of  \Cref{alg:PnPADMM} to a stationary point, ideally even to a critical point of the energy \eqref{eq:pnp_mot_opt_problem}.
For the latter, we expect that $f_{\theta}$ has to be the proximal operator of some regularizer $R$.
To ensure this, two prominent strategies have been proposed, and we discuss them in this section.    
These strategies can be extended for Bregman Proximal Gradient algorithm \cite{hurault2023convergent}, which replaces the Euclidean norm in the proximal mapping of FBS (\Cref{algo:FBS}) with Bregman divergence. The use of Bregman divergences is preferred for nonsmooth data fidelities, which, for instance, appear for Poisson noise that is commonly used for modeling counting noise of detectors in imaging applications. 

\subsection{Denoising with the Gradient of an ICNN} \label{subsec:ProximalPNP}

Our motivation for replacing \smash{$\prox_{\frac{1}{\gamma} R}$} with a denoiser $f_\theta$ was that both act similarly.
To design architectures $f_\theta$ that are the proximal mapping of a (implicitly defined) regularizer $R_\theta$, we require the characterization of proximal mappings for convex functions by Moreau \cite{moreau1965proximite} and its extension to nonconvex functions \cite[Thm.\ 4a]{GriNik2020}. 

\begin{thm}{Characterization of proximal mappings}{prox characterization}	
	Let $U \neq \emptyset \subset \R^N$ be connected and $f \colon U \to \R^N$.
	Then, the following are equivalent.
	\begin{enumerate}
		\item There exists $R \colon \R^N \to \overline{\R}$ such that $f(x) \in \prox_{R}(x)$ for all $x \in U$.
		\item There exists a convex lsc $\varphi \colon \R^N \to \overline{\R}$ such that $f(x) \in \partial \varphi(x)$ for all $x \in U$.
	\end{enumerate}
	For any such pair $R$ and $\varphi$, there exists $K \in \R$ with
	\begin{equation}\label{eq: prox-reg relation}
		R(f(x))  = \inner{f(x)}{x} - \frac{1}{2}\|f(x)\|^2_2 - \varphi(x) + K \quad \text{for all} \quad x \in U.
	\end{equation}
	Moreover, if $U= \R^N$, then  $R$ is convex if and only if $f$ is 1-Lipschitz.
\end{thm}

\begin{remark}
	Here, $\partial \varphi(x)$ denotes the subdifferential of $\varphi$ at $x$ as defined in \eqref{eq: subdifferential}.
	Since we only work with $\varphi \in C^1(\R^N)$, the condition \circled{2} becomes $f = \nabla \varphi$.
	
	As proximal mapping is invariant to constant shifts, i.e., $\prox_{R} = \prox_{R-K}$ for all $K \in \R$, we can choose $K=0$, when defining $R$ through $\varphi$ and $f$.  
\end{remark}
 
\begin{pr}
	First, we establish the equivalent characterizations of proximity operators.
	
	$\Rightarrow$: We define
	\begin{equation}
			\rho(x) = \begin{cases}
			\langle f(x), x\rangle - \frac12 \Vert f(x) \Vert_2^2 - R(f(x)), \quad & \text{if } x \in U,\\
			\infty, \quad & \text{else.}
		\end{cases}
	\end{equation}
	Since $f(x)$ minimizes $z \mapsto \frac12 \Vert z -x \Vert_2^2 + R(z)$, we obtain for every $x,z \in U$ that
	\begin{align}
		\rho(z) - \rho(x) & = \langle f(z), z\rangle - \tfrac12 \Vert f(z) \Vert_2^2 - R(f(z)) - \langle f(x), x\rangle + \tfrac12 \Vert f(x) \Vert_2^2 + R(f(x)) \notag \\
		& = \langle f(z), z - x\rangle
		+ \tfrac{1}{2} \| f(x) - x \|_2^2 + R(f(x)) - \tfrac{1}{2} \| f(z) - x \|_2^2 - R(f(z)) \notag \\
		& \le  \langle f(z), z - x\rangle. \label{eq: prox rho upper bound} 
	\end{align}
%	Hence, $f(z) \in \partial \rho(z)$ for all $z \in U$.
	The function $\varphi(x) \coloneqq \sup_{z \in U} \rho(z) + \langle f(z), x-z \rangle \ge \rho(x)$ is convex and lsc since it is the pointwise supremum of affine functions.
	Moreover, by \eqref{eq: prox rho upper bound} it also holds that $\varphi(x) \le \rho(x)$ for all $x \in U$.
	Thus, we have that $\varphi(x) = \rho(x)$ for all $x \in U$ and from \eqref{eq: prox rho upper bound} we obtain $f(z) \in \partial \varphi(z)$ as desired.
	
	$\Leftarrow$: Define $h \colon U \to \overline{\R}$ as $h(x) = \langle f(x), x\rangle-\varphi(x)$. 
	Since $\varphi$ is convex and $f(x) \in \partial \varphi(x)$ for all $x \in U$, we obtain by \eqref{eq: strong convexity} that $\varphi(x) \geq \varphi(x^\prime) + \langle f(x^\prime),x -x^\prime \rangle$ for all $x, x^\prime \in U$.
	This leads to
	\begin{equation}\label{eq:EstExistR}
		h(x^\prime) - h(x) = \langle f(x^\prime), x^\prime \rangle -\varphi(x^\prime) - \langle f(x), x \rangle + \varphi(x) \geq \langle f(x^\prime) - f(x), x \rangle.
	\end{equation}
	For any $z \in \q R (f)$, there exists $x \in U$ with $z = f(x)$.
	If there exist $x_1 \neq x_2 \in U$ with $f(x_1) = f(x_2) = z$, then $h(x_1) - h(x_2) \geq 0$ and by interchanging $x_1$ and $x_2$ we get $h(x_1) = h(x_2)$.
	Hence, we can define $\theta \colon \R^N \to \bar \R$ via
	\begin{equation}
		\theta(z) = \begin{cases}
			h(f^{-1}(z)) \quad & \text{if } z \in \q R (f)\\
			\infty \quad & \text{else.}
		\end{cases}
	\end{equation}
	Using \eqref{eq:EstExistR} and the definition of $\theta$, for all $z^\prime \in \q R (f)$ we get that $\theta(z^\prime) - \theta(f(x)) = h(f^{-1}(z^\prime)) - h(x) \geq \langle z^\prime -f(x), x \rangle$, namely that
	\begin{equation}\label{eq:EstTheta}
		\theta(z^\prime) - \langle z^\prime, x \rangle \geq \theta(f(x))  - \langle f(x), x \rangle.
	\end{equation}
	Now, we set $R(z) = \theta(z) - \frac12 \Vert z \Vert_2^2$.
	Inserting this into \eqref{eq:EstTheta} leads to
	\begin{equation}
		R(z^\prime) + \tfrac12 \Vert z^\prime \Vert_2^2 - \langle z^\prime, x \rangle + \tfrac12 \Vert x \Vert_2^2 \geq R(f(x))  + \tfrac12 \Vert f(x) \Vert_2^2 - \langle f(x), x \rangle + \tfrac12 \Vert x \Vert_2^2 
	\end{equation}
	for all $z^\prime \in \q R (f)$.
	In particular, we obtain that $f(x) \in \prox_R(x)$.
	
	Now, it remains to establish the remaining claims.
	Let $R$ and $\tilde \varphi$ be an admissible pair.
	Since $U$ is connected, we get from part one that $\varphi$ and $\tilde \varphi$ differ by a constant.
	Finally, the original result by Moreau \cite[Cor.\ 10.c]{moreau1965proximite} states that $f \colon \R^N \to \R^N$ is a proximity operator of some convex $R$ if and only if \circled{2} holds and $f$ is nonexpansive.
	This concludes the proof.
\end{pr}

A direct consequence of \Cref{th:prox characterization}	is the following.
If we take a convex network $\varphi_\theta \in  C^1(\R^N)$ and set $f_\theta = \nabla_x \varphi_\theta$, then the resulting architecture is a proximal mapping.
To this end, we can deploy the ICNN architecture discussed in \Cref{sec: reg architecture}.
Recall that an ICNN with $\ell$ layers is given by $\varphi_\theta(x) = h_{\ell}(x)$, where the $z_i$ are defined recursively via
\begin{equation}\label{eq:ICNN_Layers}
	h_1(x) = \psi_0(V_0 x+ b_0),\quad h_{i+1}(x) = \psi_i(W_i h_i(x)+V_i x+ b_i), \text{ for }i=1,\dots,\ell-1,
\end{equation}
with parameters $\theta=\{(W_i)_{i=1}^{\ell-1},(V_i)_{i=0}^{\ell-1},( b_i)_{i=0}^{\ell-1}\}$.
By \Cref{l: ICNN}, $\varphi_\theta$ is convex if the $W_i$ have non-negative entries and the $\psi_i \colon \R \to \R$ are convex and non-decreasing.
To ensure differentiability, we require $\psi_i \in C^{1}(\R)$.
This leads to the following corollary.

\begin{cor}\label{cor: implicit reg}
	Let the assumptions of \Cref{l: ICNN} hold.
	If additionally $\psi_i \in C^{1}(\R)$, then for $\varphi_\theta(x) = h_{\ell}(x)$ there exists $R_\theta \colon \R^N \to \overline{\R}$ such that $\nabla_x \varphi_\theta(x) \in \prox_{R_\theta}(x)$ for all $x \in \R^N$. 
\end{cor} 

\Cref{cor: implicit reg} guarantees the existence of an underlying regularizer $R_\theta$.
As next step, we want to guarantee that the variational reconstruction problem \eqref{eq:pnp_mot_opt_problem} with $R=R_\theta$ admits a global minimum. 
According to  \Cref{cor:existence variational}, this holds provided that $R_\theta \colon \R^N \to \overline \R$ is (lower semi-) continuous and coercive. 
The following technical corollary to \Cref{th:prox characterization} is helpful in establishing the continuity of $R_\theta$.

\begin{cor}[{\cite[Cor.\ 6b]{GriNik2020}}]\label{cor: prox characterization continuity}
	Consider the setup of \Cref{th:prox characterization} and let $U \subset \q R(f)$ be open. Then, $R \in C^{k+1}(U)$ if and only if the restriction $\tilde f$ of $f$ to $f^{-1}(U) \coloneqq \{ x \in \R^N \colon f(x) \in U \}$ is injective and $\tilde f^{-1} \in C^{k}(U)$.
\end{cor}

Consequently, we have to modify $f_\theta$ in such a way that $\q R(f_\theta) = \R^N$ and that a continuous inverse of $f_\theta$ exists.
First, we provide intuition on the scalar case.
If $\varphi\colon \R \to \R$ is continiously differentiable and strongly convex, then its derivative $f = \varphi^\prime$ is continuous, and strictly increasing.
In particular, we get $\q R(f) = \R$ and that a continuous inverse $f^{-1}$ exists.
Thus, including a strong convexity term to the ICNN $\varphi_\theta$ appears promising.
\begin{thm}{Continuously differentiable regularizer}{implicit reg cont}
	Under the assumptions of \Cref{cor: implicit reg}, any ICNN $\varphi_\theta(x) = h_{\ell}(x) + \frac{\mu}{2}\|x\|_2^2$ with $0 < \mu < 1$ induces a regularizer $R_\theta \in C^1(\R^N)$.
\end{thm}
\begin{pr}
	To apply \Cref{cor: prox characterization continuity}, we need to show $\q R(f_\theta) = \R^N$ and $f_\theta^{-1} \in C^0(\R^N)$.
	We first establish  $\q R(f_\theta) = \R^N$. 
	Let $v \in \R^N$ and consider the strongly convex and continuous function
	\begin{equation}\label{eq:ObjectiveExplicitR}
		\rho_v(x) = h_{\ell}(x) + \tfrac{\mu}{2} \|x\|_2^2 - \inner{v}{x}.
	\end{equation}
	By \Cref{cor:uniqueness variational}, the function $\rho_v(x)$ admits a unique global minimizer $\hat x$ and by Fermat's rule, $\hat x$ satisfies $v = \nabla_x h_{\ell}(\hat x) + \mu x = f_\theta(\hat x)$. 
	Thus, $f_\theta^{-1}(v) = \argmin_x \rho_v(x) = \hat x$.
	Furthermore, by \Cref{rem:ContinuityOptim}, it holds that $f_\theta^{-1}\in C^0(\R^N)$.
\end{pr}

\begin{remark}
In principle, we can also evaluate the (implicitly defined) regularizer $R_\theta$.
Using \eqref{eq: prox-reg relation}, we get that 
\begin{equation}\label{eq: implicit reg computation}
	R_\theta(x) = \inner{f_\theta^{-1}(x)}{x} - \tfrac{1}{2}\| x \|^2_2 - \varphi_\theta(f_\theta^{-1}(x)),
\end{equation}
where $f_\theta^{-1}(x)$ is computed by minimizing $\rho_x$ as defined in \eqref{eq:ObjectiveExplicitR}.
In the context of \Cref{alg:PnPADMM}, we have $f_\theta^{-1}(z_{k+1}) = b_k + x_{k+1}$ by construction.
Hence, we can efficiently track $R(z_k)$ and the convergence behavior of \Cref{alg:PnPADMM} without computing an inverse.
\end{remark}

To show coercivity of $R_\theta$, we first show that $h_{\ell}$ and $f_\theta$ are Lipschitz continuous if
\begin{equation}\label{eq: scaled softmax}
	\psi_i(t) = \tfrac{1}{\beta} \log(1 + \exp(\beta t)), \qquad \beta > 0.
\end{equation}
In fact, any other smooth and Lipschitz activation with slow growth at infinity also works.

\begin{lemma}\label{l: implicit reg Lipschitz}
	If the assumptions of \Cref{th:implicit reg cont} hold and the $\psi_i$ are chosen as \eqref{eq: scaled softmax}, then both $h_{\ell}$ and $f_\theta$ are Lipschitz continuous.  
\end{lemma} 
\begin{pr}
	We only show that $h_{\ell}$ is Lipschitz.
	The derivation for $f_\theta$ is more technical and can be found in \cite[Lem.\ C.10]{fang2023s}.  
	First, we observe that \eqref{eq: scaled softmax} is 1-Lipschitz since
	\begin{equation}
	| \psi'(t) | = \frac{1}{\beta} \cdot \frac{\beta \exp(\beta t)}{1 + \exp(\beta t)} 
	\le 1 \quad \text{for all} \quad t \in \R.
	\end{equation}
	Then, $h_1$ is Lipschitz continuous since for all $x,u \in \R^N$ it holds
	\begin{align}
		\| h_1(x) - h_1(u) \|_2^2
		& = \sum_{j = 1}^{N_1} | \psi( (V_0 x)_j  + (b_0)_j ) - \psi( (V_0 u)_j  + (b_0)_j ) |^2 \notag\\
		&\le \sum_{j = 1}^{N_1} | (V_0 x)_j  - (V_0 u)_j  |^2 
		= \| V_0 (x-u)\|_2^2 \le \| V_0 \|_{2}^2 \| x - u \|_2^2.
	\end{align}
	By induction, we obtain the Lipschitz continuity of $h_{\ell}$ as a composition of Lipschitz continuous functions.
\end{pr}

Finally, we can establish the coercivity of $R_\theta$. 

\begin{thm}{Coercive implicit regularizer}{implicit reg coercive}
	Under the assumptions of \Cref{l: implicit reg Lipschitz}, any ICNN $\varphi_\theta(x) = h_{\ell}(x) + \frac{\mu}{2}\|x\|_2^2$ with $0 < \mu < 1$ induces a coercive regularizer $R_\theta \in C^1(\R^N)$.
\end{thm}
\begin{pr}
	Using \eqref{eq: prox-reg relation}, we rewrite $R_{\theta}(f_{\theta}(x))$ as
	\begin{align}\label{eq:implcitReg}
		R_{\theta}(f_{\theta}(x))
		& = \inner{x}{f_\theta(x)} - \tfrac{1}{2}\|f_\theta(x)\|^2_2 - \varphi_\theta(x) \notag\\
		&=
		\langle x, \nabla h_{\ell}(x) + \mu x \rangle
		- \tfrac{1}{2} \| \nabla h_{\ell}(x) + \mu x \|_2^2
		- 
		h_{\ell}(x) - \tfrac{\mu}{2} \| x \|_2^2 \notag \\
		&=
		(1-\mu) \langle x, \nabla h_{\ell}(x) \rangle
		+ \tfrac{\mu(1-\mu)}{2}\| x \|_2^2
		- \tfrac{1}{2} \| \nabla h_{\ell}(x) \|_2^2
		- h_{\ell}(x). 
	\end{align}
	Using the Cauchy-Schwarz inequality, we further estimate \eqref{eq:implcitReg} as
	\begin{equation}\label{eq:RthetaCS}
		R_{\theta}(f_{\theta}(x)) \ge (\mu -1) \| x \|_2 \| \nabla h_{\ell}(x) \| + \tfrac{\mu(1-\mu)}{2}\| x \|_2^2
		- \tfrac{1}{2} \| \nabla h_{\ell}(x) \|_2^2
		- |h_{\ell}(x)|.
	\end{equation}
	Since $h_{\ell}$ is $L$-Lipschitz, we have that  $\|\nabla h_{\ell}(x) \|_2 \le L$ and $|h_{\ell}(x) | \le L \| x \|_2 + |h_{\ell}(0)|$.
	Inserting this into \eqref{eq:RthetaCS} leads to
	\begin{equation}
			R_{\theta}(f_{\theta}(x)) \ge \| x \|_2 \bigl ( \tfrac{\mu(1-\mu)}{2}\| x \|_2 + (\mu-2) L \bigr) - \tfrac{L^2}{2} - |h_{\ell}(0)|.
	\end{equation}
	Thus, for $x$ satisfying
	\begin{equation}\label{eq: implicit coercivity bounds}
		\| x \|_2 \ge M \coloneqq \max\biggl\{\frac{4 (2-\mu) L}{\mu(1-\mu)}, \Bigl(\frac{4 L^2 + 8 |h_{\ell}(0)|}{\mu(1-\mu)} \Bigr)^{1/2} \biggr\} > 0
	\end{equation}
	it holds that
	\begin{equation}\label{eq:LowerBoundReg}
		R_{\theta}(f_{\theta}(x)) 
		\ge \tfrac{\mu(1-\mu)}{4}\| x \|_2^2 - \tfrac{L^2}{2} - |h_{\ell}(0)| \ge \tfrac{\mu(1-\mu)}{8}\| x \|_2^2.
	\end{equation}
	
	Next, we want to insert $f_\theta^{-1}(x)$ into \eqref{eq:LowerBoundReg}, for which we have to show $\| f_\theta^{-1}(x) \|_2 \ge M$.
	Moreover, to establish coercivity of $R_\theta$, we need to lower bound $	\| f_\theta^{-1}(x) \|_2 $ in terms of $\| x \|_2$.
	Since  $f_\theta$ is $L_\nabla$-Lipschitz continuous, we get for $v = f_\theta^{-1}(x)$ and $u = f_\theta^{-1}(0)$ that
	\begin{equation}
		L_\nabla^{-1} \| x - 0 \|_2 \le \| f_\theta^{-1}(x) - f_\theta^{-1}(0)\|_2 \le \| f_\theta^{-1}(x) \| + \| f_\theta^{-1}(0)\|_2. 
	\end{equation}
	Thus, for $\|x\|_2 \ge 2 L_\nabla \max\{\| f_\theta^{-1}(0)\|_2 , M\}$ we have that
	\begin{equation}
	\| f_\theta^{-1}(x) \|_2 \ge \tfrac{1}{2 L_\nabla} \| x \|_2 \ge M
	\quad \text{and} \quad 
	R_{\theta}(x) 
	\ge \tfrac{\mu(1-\mu)}{8}\| f_\theta^{-1}(x) \|_2^2 
	\ge \tfrac{\mu(1-\mu)}{32 L_\nabla^2}\| x \|_2^2. 
	\end{equation}
	Hence,  $R_{\theta}$ is coercive.
\end{pr}

\Cref{th:implicit reg cont,th:implicit reg coercive} ensure that the conditions of \Cref{cor:existence variational} are satisfied, namely that a minimizer of \eqref{eq:pnp_mot_opt_problem} exists.
Next, we show that \Cref{alg:PnPADMM} converges based on \cite[App.\ C]{fang2023s}. 

\begin{thm}{Convergence of PnP-ADMM}{PnP-ADMM}
	Let the assumptions of \Cref{th:implicit reg coercive} hold and $\{ (x_{k}, z_k, b_k) \}_{k \in \bb N}$ be generated by \Cref{alg:PnPADMM} with $f_\theta = \nabla_x \varphi_\theta(x)$, $\gamma=1$ and $\alpha > \| A \|_2^2$.
	Then, there exists a constant $c>0$ such that the augmented Lagrangian satisfies
	\begin{equation}
		L_\gamma(x_{k+1}, z_{k+1}, b_{k+1}) \le L_\gamma(x_{k}, z_{k}, b_{k}) - c \| x_{k} - z_{k} + b_{k}\|_2^2.
	\end{equation}
	Moreover, the residual satisfies
	\begin{equation}
		\min_{k = 0,\ldots, K-1} \| x_k - z_k + b_{k}\|_2^2 = \q O(1/K).
	\end{equation}
	and the iterates converge to a limit point $(x^*, z^*, b^*)$ such that
	\begin{equation}
		x^* = z^*, \quad \nabla_x R_\theta (z^*) = b^* = - \tfrac{1}{\alpha} A^\tT (A x^* - y).
	\end{equation}
	in particular, we get that $x^*$ is a critical point of \eqref{eq:pnp_mot_opt_problem} with $R = R_\theta$. 
\end{thm} 
\begin{pr}
	We show that the mentioned convergence conditions in  \Cref{rem: ADMM} are satisfied.
	First, we observe that the data fidelity term is $C^1(\R^N)$ with $\|A\|_2^2/\alpha$-Lipschitz gradient.
	Second, we need to verify that $L_\gamma$ is a KL function.
	For this, we can rely on \Cref{th:definable structure} and in particular \Cref{th:definable operations}.
	Below, we directly refer to it's 6 items \circled{1} to \circled{6}.
	
	Since both $x \mapsto \exp(x)$ and the interval $(0,\infty)$ are definable,  \circled{5} implies that $x \mapsto \log(x)$ is definable.
	Thus, by combining \circled{1} and \circled{2}, the activation functions $\psi_i$ as defined in \eqref{eq: scaled softmax} are definable.
	Consequently, each layer $h_k(x)$ as defined in \eqref{eq:ICNN_Layers} is definable by  \circled{1} and \circled{2}.
	Thus, also $\varphi_\theta(x) = h_\ell(x) + \mu \| x \|^2$ is definable by \circled{1}.
	Since $\varphi_\theta \in C^1(\R^N)$ is definable, its gradient $f_\theta$ is definable by \circled{4}.
	Moreover, we showed that $\q R(f_\theta) = \R^N$ in \Cref{th:implicit reg cont}, which is a definable set.
	Hence, \circled{5} implies that $f^{-1}_\theta$ is definable.
	Finally, by \eqref{eq: implicit reg computation}, we now get that $R_\theta$ is definable and, consequently, $L_\gamma$ as well.    
\end{pr}

Note that the underlying regularizer $R_\theta$ for $f_\theta$ is not necessarily convex.
As we saw in \Cref{th:prox characterization}, convexity corresponds to $f_\theta$ being 1-Lischitz, which is a hard to check in practice. 
Moreover, we have a restriction on the regularization strength, namely $\alpha > \| A \|_2^2$.
If the regularizing effect is too strong, then $f_\theta$ should be trained for a smaller noise level.

\subsection{Gradient-Step Denoiser} 

Alternatively, we can choose $f_\theta$ as the gradient step denoiser $f_\theta = x - \nabla g_\theta$, where $g_\theta$ is a differentiable neural network \cite{hurault2022proximal}. 

The motivation for the gradient step denoiser architecture comes from the optimal transport problem, see for instance \cite{peyre2019computational}. Given two probability measures $\mu$ and $\nu$ in $\R^N$, the goal is to find a transport map $f \colon \R^N \to \R^N$ such that $\nu = \mu \circ f^{-1}$ and $f$ minimizes the transport cost \mbox{$\int_{\R^N} \frac{1}{2} \| x - f(x)\|_2^2 \diff \mu (x)$}.  
Brenier established in \cite[Prop.~3.1]{brenier1991polar} that 
such $f$ exist and is the gradient of a potential function $\varphi = \frac{1}{2} \|\cdot\|_2^2 - g$ if $\varphi$ is convex and  $g$ is a solution of underlying dual problem. 
%to optimal transport. 
In our setting, $\mu$ and $\nu$ represent noisy and, respectively, noiseless data distributions. Thus, we can view training of gradient step denoiser $f_\theta$ as learning of the optimal transport map between the two distributions. 

Again, we can ensure that $f_\theta$ is a proximal operator of some $R$ via \Cref{th:prox characterization}.
Essentially, the obtained gradient-step denoiser then simply corresponds to a different parameterization of convex functions.

\begin{thm}{Gradient-step denoiser}{PR-DRUNet properties}
	Let $g_\theta \in C^{2}(\R^N)$ and $\nabla g_\theta$ be $L$-Lipschitz with $L < 1$.
	Then, it holds that
	\begin{enumerate}
		\item the potential $\varphi_\theta = \frac{1}{2} \| x\|_2^2 - g_\theta(x)$ with $f_\theta = \nabla \varphi_\theta$ is $(1-L)$-strongly convex;
		\item the denoiser $f_\theta$ satisfies $f_\theta(x) = \prox_{R_\theta}(x)$, where $R_\theta \colon \R^N \to \R$ is defined by
		\begin{equation}\label{eq:EvalReg}
		R_\theta(x) \coloneqq g_\theta(f^{-1}_{\theta}(x)) - \tfrac{1}{2}\|f^{-1}_{\theta}(x) - x\|^2_2;
		\end{equation}
		\item the regularizer satisfies $R_\theta(x) \ge g_\theta(x)$ for all $x \in \R^N$;
		\item moreover, $R_\theta \in C^{1}(\R^N)$ with $\frac{L}{1-L}$-Lipschitz gradient.
	\end{enumerate} 	
\end{thm}

\begin{pr}
	\circled{1}\,: Since $g_\theta \in C^2(\R^N)$, we can apply the second-order condition \eqref{eq: convexity grad}.
	Using the results of \Cref{ex:Lipschitz Hessian}, we get for all $x,z \in \R^N$ that
	\begin{equation}
	\inner{\nabla^2 \varphi_\theta(x) z}{z} 
	= \inner{I_N z}{z} - \inner{\nabla^2 g_\theta(x) z}{z}
	\ge \|z\|_2^2 - L \|z\|_2^2 = (1-L)\|z\|_2^2.  
	\end{equation}
	\circled{2}\,: The existence of a regularizer $R_\theta$ satisfying \eqref{eq: prox-reg relation} follows from \Cref{th:prox characterization}. For any $v \in \R^N$, the function $\rho_v(x) = \varphi_\theta (x) - \inner{v}{x}$ is strongly convex and continuous. Hence, by \Cref{cor:uniqueness variational}, it admits a unique global minimizer $\hat x \in \R^N$, which, by Fermat's rule, satisfies $v = \nabla_x \varphi_\theta(\hat x) = f_\theta(\hat x)$ and thus $f_\theta^{-1}(v) = \argmin _x \rho_v(x)$.
	By inserting $f_\theta^{-1}$ into  \eqref{eq: prox-reg relation}, we then get that
	\begin{align}
		R_\theta(x) & = \inner{f_\theta^{-1}(x)}{x} - \tfrac{1}{2}\| x \|^2_2 - \varphi_\theta(f_\theta^{-1}(x)) \notag \\
		& = - \tfrac{1}{2} \left( \| x \|^2_2 -  2 \inner{f_\theta^{-1}(x)}{x} + \| f_\theta^{-1}(x) \|_2^2 \right) + g_\theta(f_\theta^{-1}(x)) \notag \\
		& = g_\theta(f_\theta^{-1}(x))	- \tfrac{1}{2} \| f_\theta^{-1}(x) -x \|_2^2.
	\end{align}
	\circled{3}\,: From \circled{2}, we obtain that
	\begin{equation}
	R_\theta(f_{\theta}(x)) = g_\theta(x) - \tfrac{1}{2} \| x - f_\theta (x) \|_2^2.
	\end{equation}
	Since $f_\theta(x) = \prox_{R_\theta}(x)$, it holds that
	\begin{equation}
	g_\theta(x) 
	= R_\theta(f_{\theta}(x)) + \tfrac{1}{2} \| x - f_\theta (x) \|_2^2 
	= \min_{v \in \R^N} R_\theta(v) + \tfrac{1}{2} \| x - v \|_2^2 
	\le R_\theta(x).
	\end{equation}
	\circled{4}\,:
	By \Cref{rem:ContinuityOptim}, it holds that $f_\theta^{-1}\in C^0(\R^N)$ such that \Cref{cor: prox characterization continuity} further implies $R_\theta \in C^1(\R^N)$.
	Since it holds that $f_\theta(x) = \argmin_{z \in \R^N} R_\theta(z) + \frac{1}{2}\|z - x\|_2^2$, Fermat's rule yields 
	\begin{equation}\label{eq:PnPFermat}
	0 = \nabla R_\theta(f_\theta(x)) + f_\theta(x) - x =  \nabla R_\theta(f_\theta(x)) - \nabla g_\theta(x).
	\end{equation}
	Due to the invertibility of $f_\theta$ , for every $u,v \in \R^N$ there exist $x,z  \in \R^N$ with $u = f_\theta(x)$ and $v = f_\theta(z)$, such that \eqref{eq:PnPFermat} and the $L$-Lipschitz continuity of $\nabla g_\theta$ yield
	\begin{equation}\label{eq:PnP1}
	\|\nabla R_\theta(u) - \nabla R_\theta(v) \|_2
	= \| \nabla g_\theta(x) - \nabla g_\theta(z) \|_2
	\le L \| x - z \|_2.
	\end{equation}
	Furthermore, by the strong convexity of $\varphi_\theta$ derived in \circled{1} and \Cref{ex: strong convexity first order}, we get
	\begin{equation}\label{eq:PnP2}
	(1-L)\| x - z \|_2
	\le \| \nabla \varphi_\theta(x) - \nabla \varphi_\theta(z) \|_2
	= \| f_\theta(x) -f_\theta(z) \|_2
	= \| u - v\|_2 .
	\end{equation}
	Combining \eqref{eq:PnP1} and \eqref{eq:PnP2} concludes the proof.
\end{pr}

In \cite{hurault2022proximal}, the authors take $g_\theta(x) = \frac{1}{2}\| x - \Phi_\theta(x) \|_2^2$ (which is also discussed as Least Squares Residual Regularizer in \Cref{sec: reg architecture}) with the DRUNet architecture $\Phi_\theta$ discussed in \Cref{chap:regularization}.
To make $\Phi_\theta$ smooth, the ReLU activation in the original DRUNet architecture are replaced by the softmax function \eqref{eq: scaled softmax}. 
This leads to the smooth denoiser
\begin{equation}\label{eq: grad step denoiser}
f_\theta = x - \nabla g_\theta 
= x - (I - J_{\Phi_\theta}(x))^\tT (x - \Phi_\theta(x))
= \Phi_\theta(x) + J_{\Phi_\theta}(x)^\tT (x - \Phi_\theta(x)),
\end{equation}
namely to applying $\Phi_\theta$ with a linear correction.
Moreover, we have $g_\theta \ge 0$ by construction. 

A coercive regularizer $R_\theta$ can be guaranteed by adjusting $g_\theta$.
Using the Lipschitz constant $L<1$ of $\nabla g_\theta$ and some $0 < \varepsilon < 1$, we define the new function
\begin{equation}\label{eq: adjusted PR-DRUNet}
	\tilde g_\theta(x) \coloneq g_\theta(x) + \frac{\varepsilon(1-L)}{2}\|x\|_2^2
	\quad \text{with} \quad
	\nabla \tilde g_\theta(x) = \nabla g_\theta(x) + \varepsilon(1-L) x.
\end{equation}
Its gradient is $\tilde L$-Lipschitz with $\tilde L = (1 - \varepsilon) L + \varepsilon < 1$, and $\tilde g_\theta$ is coercive since $\tilde g_\theta(x) \ge \frac{\varepsilon (1- L)}{2}\|x\|_2^2$.
By \Cref{th:PR-DRUNet properties}\,\circled{3}, we thus get coercivity of $R_\theta$. 
Consequently, we can establish the convergence of the PnP-ADMM \Cref{alg:PnPADMM} analogously to \Cref{th:PnP-ADMM}.

\begin{thm}{Convergence of PnP-ADMM 2}{PnP-ADMM PR-DRUNet}
	Let $g_\theta \in C^2(\R^N)$ be as in \eqref{eq: adjusted PR-DRUNet} and $\nabla g_\theta$ be $L$-Lipschitz continuous with $L < 1$.
	Then, the results of \Cref{th:PnP-ADMM} apply to the sequence generated by PnP-ADMM with $\gamma=1$ and $\alpha  > \| A \|^2_2$.
\end{thm}

Both \Cref{th:PnP-ADMM} and \Cref{th:PnP-ADMM PR-DRUNet} require $\alpha  > \| A \|^2_2$.
This is due to the condition on $\gamma$, see also \Cref{rem: ADMM}.

We can get rid of this dependency if we swap the order in which we perform the minimization with respect to the variables $x$ and $z$ in \cref{alg:PnPADMM}.  

\begin{algorithm}[H]
	\KwData{denoiser $f_\theta \colon \R^N \to \R^N$, initializations $x_{0}\in \R^N$ and $b_{0} \in \R^N$}
	\For{\color{black}$k = 0, 1, \ldots$}{
		$z_{k+1} = f_{\theta}(b_{k} + x_{k})$\;
		$x_{k+1}\in \argmin_{x \in \R^N} \left\{ \frac{1}{2\alpha}\|Ax - y\|^2_\bb Y +  \frac{\gamma}{2} \|b_{k} + x - z_{k+1}\|_2^2  \right\}$\;
		$b_{k+1}= b_{k} + x_{k+1} -  z_{k+1}$\;
	}
	\caption{Plug-and-Play ADMM (PnP-ADMM-alt) \label{alg:PnPADMM2}}
\end{algorithm}

This shifts the differentiability requirement from the data fidelity to the regularizer $R_\theta$ when establishing the convergence of ADMM.
In particular, we can rely on \Cref{alg:PnPADMM2} instead of \Cref{alg:PnPADMM} if we want to replace $\Vert Ax -y \Vert^2_2$ in \eqref{eq:pnp_mot_opt_problem} by a non differentiable data term.

\begin{thm}{Convergence of PnP-ADMM-alt}{PnP-ADMM PR-DRUNet 2}
	Let $g_\theta \in C^2(\R^N)$ be as in \eqref{eq: adjusted PR-DRUNet} and $\nabla g_\theta$ be $L$-Lipschitz continuous with $L < 1/2$.
	Then, for $\gamma = 1$ and any $\alpha > 0$, the convergence guarantees from \Cref{th:PnP-ADMM} apply to the iterates of PnP-ADMM-alt.
\end{thm}

\begin{pr}
	We provide an intuitive proof for $L < 1/3$. 
	If we change the roles of the data fidely and regularization term, then $g = R_\theta$ in the ADMM \Cref{alg:ADMM}, which is nonconvex.
	Thus, we require $1> 2 L_{\nabla R_\theta}$ according to \Cref{rem: ADMM}.
	Since the Lipschitz constant $L_{\nabla R_\theta}$ of $\nabla R_\theta$ satisfies $L_{\nabla R_\theta} = \frac{L}{1-L}$ by \Cref{th:PR-DRUNet properties} \circled{4}, it holds that
	\begin{equation}
	1 > \frac{2L}{1-L}, 
	\quad \text{which is equivalent to} \quad 1 > 3 L.
	\end{equation}
	The weaker bound $L < 1/2$ can be obtained by investigating the convergence of ADMM for nonconvex functions in more detail, see \cite{ThePat2020}.
\end{pr}

To ensure that the more restrictive Lipschitz bound in  \Cref{th:PnP-ADMM PR-DRUNet 2} is satisfied, we replace $f_\theta$ by the averaged network
\begin{equation}\label{eq: grad step denoiser constant control}
\tilde{f}_\theta = \beta f_\theta + (1-\beta) I_N = I_N - \beta \nabla g_\theta 
\quad \text{with} \quad 
\beta \leq \tfrac12.
\end{equation}
Since this approach reduces the influence of the initial denoiser $f_\theta$ (by staying closer to the identity), it helps us sidestep the restrictive lower bound on the regularization strength $\alpha$ that was required in \Cref{th:PnP-ADMM PR-DRUNet}.

To ensure convergence the applicability of \Cref{th:PnP-ADMM PR-DRUNet} and \Cref{th:PnP-ADMM PR-DRUNet 2}, it now remains to ensure that $\nabla g_\theta$ is $L$-Lipschitz with $L<1$ after the training.
Since the $L$-Lipschitz continuity of $\nabla g$ is equivalent to $\| J_{\nabla g}(x) \|_2 \le L$ for all $x \in \R^N$, we regularize the training problem with $\max \{\| J_{\nabla g}(x) \|_2, 1 - \xi\}$ with $\xi > 0$.
We have already discussed a similar reqularization in \Cref{sec:EstHyper}, and defer the details to the next section.

\subsection{Training a Gradient-Step Denoiser for PnP}\label{sec:trainingPnP}

In this section, we train a gradient-step denoiser
\begin{equation}
f_\theta(x) = x - \nabla g_\theta(x)
\quad \text{with} \quad
g_\theta(x) = \frac{1}{2}\| x - \Phi_\theta(x) \|_2^2
\end{equation}
for Algorithm \ref{alg:PnPADMM2}.

As a denoiser architecture $\Phi_\theta$, we use the DRUNet \cite{ZhaLiZuo2022} from \Cref{sec:Denoising} with a few changes.
Namely, the ReLU activation is replaced by the softmax \eqref{eq: scaled softmax}, the number of residual blocks is two (instead of four), and  the number of convolution channels is reduced by half, i.e., to 32/64/128/256 channels for each scale.
We refer to this model as the "4m" model since it has $\approx 4$ million trainable parameters.  

The training consists of three stages.
For each one, we train on the BSDS500 dataset \cite{arbelaez_contour_2011} consisting of color images.
This set is split into 200 training, 100 validation and 200 test images\footnote{Available at \url{https://github.com/BIDS/BSDS500/tree/master}.}.
Out of 200 training and 200 test images in BSDS500, we select the first 40 images as our validation set and the remaining 360 ones are used as our training set.
From the 100 validation images, we take the 68 ones corresponding to the BSD68 dataset as our test set.

In each epoch, patches of size $128 \times 128$ are randomly cropped from the training images, and then randomly flipped and rotated by 90 degrees.
Afterwards, we generate the corresponding noisy patches $z_i = x_i + \q N(0, \sigma_i^2 I)$ with $\sigma_i$ sampled uniformly at random from $[0,0.1]$.
For the validation and test images, we crop a centered $321 \times 321$ patch, and add noise with $\sigma = 0.03$ on the fly.
These non-deterministic sets lead to a standard deviation of around $0.0025$ in PNSR \eqref{eq: PSNR} across $50$ evaluation runs.
Thus, the random effect can be neglected.   

In the first training stage, we train the DRUNet $\Phi_\theta$ directly by minimizing the loss
\begin{equation}
	\frac{1}{2S} \sum_{i=1}^S \| \Phi_\theta(z_i) - x_i \|_2^2 + \lambda \|\theta\|^2_2
\end{equation}
using ADAM with batch size 16 for 10000 epochs. Here $\lambda=10^{-8}$ imposes a small weight decay. 
The initial learning rate of $2 \cdot 10^{-4}$ is decreased by a factor of $0.99975$ every epoch. 

\begin{remark}[Number of epochs]
	In the end, we want to perform roughly 200000 steps of ADAM.
	Given that the batch size is 16 and the training set contains 360 images, we perform $22.5$ iterations per epoch.
	Hence, we set the number of epochs to 10000.      
\end{remark}

\begin{remark}[Learning rate and its decay]
	Typically, one selects the initial learning rate $\nu_0$ as large as possible such that the training does not diverge, i.e., the PSNR on the training set drops down to small values.
	Usually, divergence occurs immediately during the early training epochs, so that we can tune the $\nu_0$ efficiently.
	During the training, we don't want to decrease the learning rate too aggressively, as this will prevent the model from learning.
	As a rule of thumb, we decrease the initial $\nu_0$ by a factor of $\approx 0.1$ over the course of training. 
	In our case, after $K = 10000$ epochs, we have the learning rate $\nu_K = 0.99975^K \nu_0 \approx 0.082 \nu_0$. 
\end{remark}

We select the model with the highest PSNR on the validation set. 
For reference, the noisy test images $z_i$ have a PSNR of 30.458 dB against the groundtruth, and we achieved a reconstruction PSNR of 37.790 dB with the 4m model.
This is of course up to the randomness in training.
For comparison, the pretrained DRUNet in DeepInverse with 16 million parameters gives a PSNR of 37.981 dB, which is slightly higher.

Our actual goal is to train a gradient-step denoiser $f_\theta  = \nabla \varphi_\theta$.
Therefore, the second stage is to train $f_\theta$ by minimizing the loss
\begin{equation}
	\frac{1}{2S} \sum_{i=1}^S \| f_\theta(z_i) - x_i \|_2^2 + \lambda \|\theta\|^2_2.
\end{equation}
In principle, we could skip the first part of the training and start from here.
However, since $f_\theta$ features the Jacobian of $\varphi_\theta$, computing derivatives of $f_\theta$ is computationally heavy, and pretraining $\varphi_\theta$ significantly reduces the number of required epochs.
Thus, we employ ADAM with batch size 16 for only 2000 epochs.
This time, the starting learning rate is $10^{-4}$, and the decay factor is $0.999$.

Already after 10 epochs, the PSNR on the validation set is 37.95 dB.
 After 100 epochs, we reach  38.0 dB and the remaining 1900 epochs only improve the PSNR to 38.05dB.   
With $f_\theta$, we achieve a test PSNR of 37.855 dB for the 4m model.
For comparison, the pretrained GSDRUNet from \cite{hurault2022gradient} with 16 million parameters achieves a PSNR of 37.849 dB. 

The third training stage enforces the Lipschitz constraint on $\nabla g$, which is required in \Cref{th:PR-DRUNet properties} and \Cref{th:PnP-ADMM PR-DRUNet}.
The Lipschitz constant of $\nabla g_\theta$ being less than one is equivalent to having
\begin{equation}
\sup_{x \in \R^N} \| J_{\nabla g_\theta} (x) \|_{2} \le 1.
\end{equation}
As computing the supremum is infeasible, we restrict ourselves to the points in the dataset.
This leads to the penalized training loss
\begin{equation}\label{eq:PenalizedLoss}
	\argmin_{\theta} \frac{1}{2S} \sum_{i=1}^S \left[ \| \Phi_\theta(z_i) - x_i \|_2^2 + \lambda \|\theta\|^2_2
	 + \mu \max\{ \| J_{\nabla g_\theta} (w_i) \|_{2}, 1 - \xi \} \right].   
\end{equation}
In \eqref{eq:PenalizedLoss}, we choose $\xi = 0.1$ and $\mu = 5 \cdot 10^{-4}$ to control the strength of the penalty.
In every epoch, new interpoaltions  $w_i = \gamma_i z_i + (1 - \gamma_i) x_i$  between $z_i$ and $x_i$ are drawn at random with uniformly distributed $\gamma_i \in [0,1]$.
With this, we aim to control the Lipschitzness beyond the training images $\{x_i\}_{i=1}^S$.  

To optimize the penalty, we employ use Jacobian-free propagation \eqref{eq:TruncGrad}. 
Recall that the eigenvector $v_1$ of $J_{\nabla g_\theta} (w_i)$ with the largest magnitude is the maximizer of $\| J_{\nabla g_\theta} (w_i) v\|_2$ with constraint $\|v\|_2 = 1$. The maximizer can be found through the power iterations
\begin{equation}
v^{k+1} = \frac{J_{\nabla g_\theta} (w_i) v^k}{\| J_{\nabla g_\theta} (w_i) v^k\|_2},
\end{equation}
see also \eqref{eq:PowerIters}.
We approximate $v_1$ by 50 power iterations, and use $\| J_{\nabla g_\theta} (w_i) v^{50}\|_2$ as an approximation of $\| J_{\nabla g_\theta} (w_i) \|_{2}$.
The required gradient $\nabla_\theta \| J_{\nabla g_\theta} (w_i) \|_{2}$ is estimated via Jacobian-free backpropagation.

We reduce the batch size to 4 (due to memory limitations), and use a fixed learning rate of $10^{-4}$ with decay rate $0.997$ for 20 epochs. 
We select the model with the best PSNR on the validation set among those epochs for which
$\max_{i = 1, \ldots, S} \| J_{\nabla g_\theta} (x_i) \|_{2} < 1$. 
Due to the additional penalty, we expect a lower test performance.
Indeed, the resulting test PSNR is only 37.347 dB for the 4m model.
We also select a model for reduced denoising effect $\beta = 1/2$ \eqref{eq: grad step denoiser constant control} with PSNR 34.169 for \Cref{alg:PnPADMM2} and its theoretical requirements in \Cref{th:PnP-ADMM PR-DRUNet 2}.

The complete training results are summarized in \Cref{tab: denoiser training}.

\begin{table}[h!]
	\resizebox{\linewidth}{!}{
		\begin{tabular}{l @{\qquad}  c c c c}
			Model & PSNR of $D_\theta$ (first) & PSNR of $f_\theta$ (second) & PSNR of $f_\theta$ (third) & PSNR of $f_\theta$ (third, $\beta = 1/2$) \\
			\midrule
			4m model & 37.790 & 37.855 & 37.347 & 34.169 \\
			Pretrained & 37.981 &  37.849 & --- & --- \\
			\bottomrule
	\end{tabular}}
	\caption{Denoising performance for $\sigma = 0.03$ after the first, second and third training stage.}
	\label{tab: denoiser training}
\end{table}

Now, we proceed to deploying the PnP algorithms.
As underlying inverse problem, we consider a deconvolution problem combined with 2-fold downsampling.
The measurements are corrupted by Gaussian noise, namely
\begin{equation}\label{eq:InvProblemPnP}
y = Ax = D [x * h] + \q N(0, \sigma_*^2 I).
\end{equation}
In \eqref{eq:InvProblemPnP}, $D\colon \R^{3 \times 320 \times 320} \to \R^{3 \times 160 \times 160} $ is a linear 2-fold downsampling, $h \in \R^{1 \times 160 \times 160}$ is a convolution kernel (applied channel-wise), and $\sigma_* = 0.03$.
Here, we have $\|A\|_2 =1$ since the kernel $h$ has nonnegative entries that sum up to one.

\begin{figure}
	\begin{center}
		\subfloat[Ground truth. \label{fig: pnp example gt}]{
			\includegraphics[width= 128pt, height = 128pt]{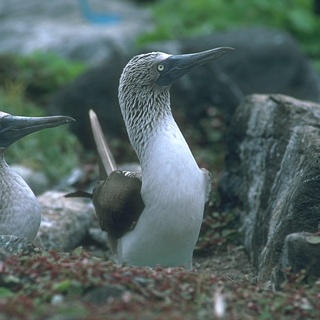}% 
		} \hfill
		\subfloat[Kernel. \label{fig: pnp example kernel}]{
		\includegraphics[width= 128pt, height = 128pt]{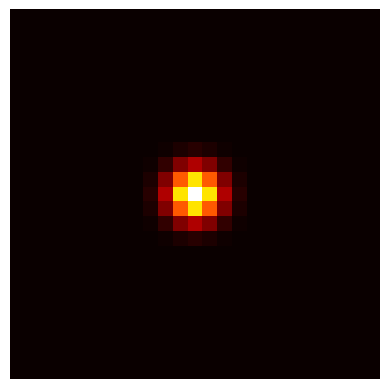}%
		}
		\hfill
		\subfloat[Measurements. \label{fig: pnp example meas}]{
			\includegraphics[width= 64pt, height = 64pt]{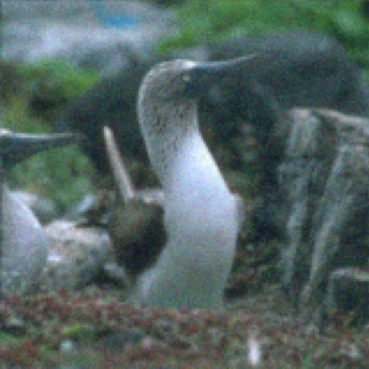}%
		}
		\caption{Setup of the inverse problem.\label{fig: pnp example}}
	\end{center}
\end{figure}

To make dimensions even, we crop $320 \times 320$ center patches instead of $321 \times 321$ ones for the validation and test sets.
We use image \texttt{0011.png}\footnote{\url{https://github.com/clausmichele/CBSD68-dataset/blob/master/CBSD68/original_png/0011.png}} as an example.
It is shown in \Cref{fig: pnp example} together with the kernel $h$ and the resulting measurements $y \in \R^{3 \times 160 \times 160}$.

\begin{table}[b]
	\resizebox{\linewidth}{!}{
		\begin{tabular}{l @{\qquad}  c c c c}
			Method & Regularization $\alpha$ & Denoiser strength $\sigma$ & PSNR & Avg. Runtime (sec) \\
			\midrule
			\texttt{ADMM-TV} & $10^{-2.5}$ & --- & 25.839 & 28.8 \\
			\texttt{PnP-GS-Pre} & 1.1 $\|A\|_2^2$ & 0.0125 & 26.168 & 59.2 \\
			\midrule
			\texttt{PnP-GS-4m}  & 1.1 $\|A\|_2^2$ & 0.025 & 27.248 & 15.2 \\
			\texttt{PnP-GS-Lip-4m} & 1.1 $\|A\|_2^2$ & 0.03 & 26.385 & 16.7 \\
			\texttt{PnP-GS-LHalf-4m} & $\|A\|_2^2$& 0.0275 & 26.642 & 19.6 \\
			\bottomrule
	\end{tabular}}
	\caption{Performance of PnP-ADMM-alt with trained gradient-step denoisers.}
	\label{tab: PnP performance}
\end{table}
 
As a baseline method, we report the performance for the TV regularization discussed in \cref{example:TV_for_inpainting_ADMM}.
The regularization parameter $\alpha$ is chosen to maximize the PSNR on the validation set among the values $\{ 10^{-3 + 0.25 k} : k = 0, \ldots, 8 \}$.
We compare the result with several versions of the PnP-ADMM (Algorithm \ref{alg:PnPADMM}) and PnP-ADMM-alt (Algorithm \ref{alg:PnPADMM2}).
Namely, we denote by \texttt{PnP-GS-4m} and \texttt{PnP-GS-Lip-4m} the result of Algorithm \ref{alg:PnPADMM} with the gradient-step denoiser after the second and third training stage, respectively.
We set $\alpha = 1.1 \|A\|_2^2 = 1.1$ to ensure that the conditions of \Cref{th:PnP-ADMM PR-DRUNet} hold.
As another reference, we include PnP-ADMM-alt with the pretrained GSDRUNet from \cite{hurault2022gradient}, and refer to it as  \texttt{PnP-GS-Pre}.

The denoising strength $\sigma$, which is an input argument of the DRUNet, is chosen with a grid search on $\{ 0.01 + 0.0025 k : k = 0, \ldots, 10 \}$ to maximize the validation PSNR.

Furthermore, we also include PnP-ADMM-alt as
\texttt{PnP-GS-LHalf-4m}, which incorporates rescaling  \eqref{eq: grad step denoiser constant control} with $\beta = 1/2$. 
By \Cref{th:PnP-ADMM PR-DRUNet 2}, we can choose $\alpha >0$ and thus, in addition to grid search for denoising strength $\sigma$ as for other \texttt{PnP-GS} methods, 
we also select $\alpha \in \{ 10^{-1 +0.5k} : k =0, \ldots, 4 \}$ based on validation PSNR. 

For all methods, we perform $1000$ ADMM iterations.
The algorithm terminates early if the relative residual $\| x^{k} -x^{k-1} \|_2 / (\| x^k \|_2 + 10^{-6})$ is below the convergence threshold $10^{-5}$. 

\begin{figure}[!h]
	\begin{center}
		\includegraphics[width=\linewidth]{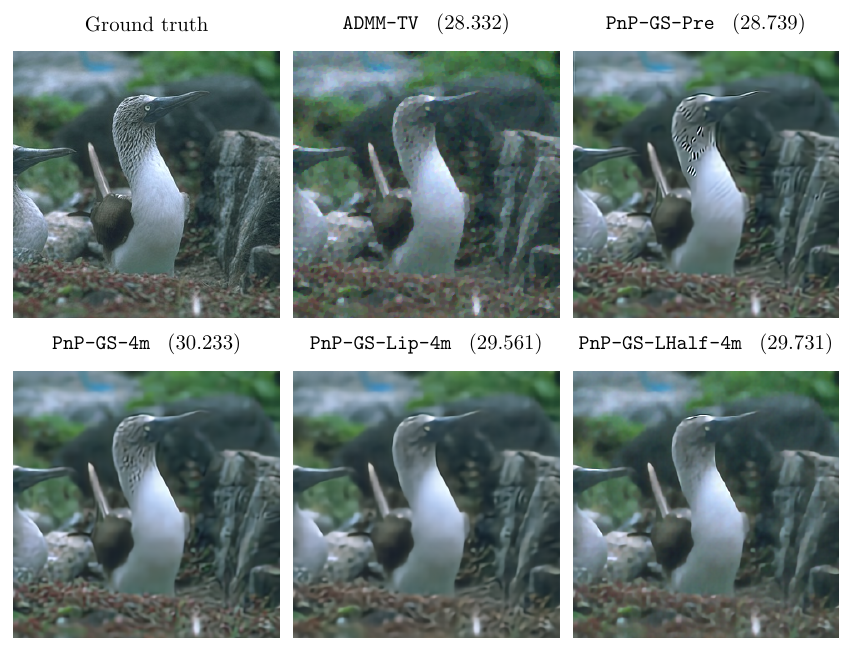}
		\caption{Reconstructions with multiple PnP methods.
			The PSNR value is reported in brackets. \label{fig: pnp results}}
	\end{center}
\end{figure}

After the parameter selection, we evaluate the methods' performance on the test set. The summary of the optimal parameters and the resulting performance metrics can be found in \Cref{tab: PnP performance}.
A visual reconstruction example is given in \Cref{fig: pnp results}.

We observe that PnP methods outperform TV both in terms of PSNR value and runtime (for 4m model).
The runtime comparison is not entirely fair, as it does not include the training time.
Among the PnP methods, \texttt{PnP-GS-4m} shows the best performance.

\subsubsection{Exercises}

\begin{exerciseB}
	In this exercise we implement PnP-ADMM in DeepInverse \cite{tachella2025deepinv}. As the inverse problem, we consider the same setup as in \Cref{sec:trainingPnP} and solve downsampled convolution described in \eqref{eq:InvProblemPnP}. The only change is reduction of image sizes to $64 \times 64$. As a validation dataset, 4 random images are sampled from CBSD68 dataset. As a test dataset, we use set3c dataset. For both datasets, we crop $64 \times 64$ center patches.	
	
	Firstly, via DeepInverse functionality implement ADMM for TV-regularized \cref{ex: sparsity promoting} as a baseline method. Secondly, implement PnP-ADMM and PnP-ADMM-alt (Algorithms \ref{alg:PnPADMM} and \ref{alg:PnPADMM2}). As a denoiser, use a pretrained GSDRUNet provided by DeepInverse and also try your denoiser from \cref{ex:Denoiser} with controlled denoising strength $D_t(x) = (1-t) x + t D(x)$, $t\in[0,1]$. 
	
	All of the above algorithms require parameter fitting: regularization $\alpha$ for TV and denoising strength $\sigma$ and $t$ for GSDRUNet and your denoiser, respectively. Use the validation dataset to select parameters yielding optimal performance. Then, compute PSNR of the tuned methods on the test dataset. Do learned methods outperform TV regularization? 
	 
	You may use the following Jupyter Notebook as a basis for this exercise, \url{https://github.com/Sniper2k/inverse_problems_ai_exercise_templates/blob/main/Ex7_3_1_template.ipynb}.
\end{exerciseB}

\begin{Solution}
	See Jupyter Notebook \url{https://github.com/Sniper2k/inverse_problems_ai_exercises/blob/main/Ex7_3_1.ipynb}.
\end{Solution}

\section{Unrolling techniques} \label{sec: unrolling}
In the seminal work of Gregor and LeCun \cite{gregor2010learning}, a technique called algorithm \emphi{unrolling} (or \emph{unfolding}) was developed.
It draws inspiration from iterative algorithms such as the ISTA in  \eqref{eq:ISTA} to build neural networks.
Following this work, unrolling has been applied for many signal and image processing tasks.
Examples include sparse recovery (see \Cref{ex: sparsity promoting}) with unrolled ADMM \cite{yang2018admm}, and blind deconvolution \cite{li2019algorithm} (see \Cref{ex: convolution} with \emph{both} unknown signal $x$ and kernel $g$).
For an overview, we refer to \cite{MonLiEld2021}.

For unrolling, we interpret the update step of an optimization algorithm as a neural network layer and $K$ updates of an algorithm form a (deep) neural network where the parameters are shared between layers. We saw an application of unrolling in \Cref{sec:LearnedReg} for purely learning regularizers, but this idea can be taken further. Some (or even all) components of the update, such as the regularization operators and the step size, can then be learned.  
The trained network can then be interpreted as an algorithm with parameters that are optimized to the data and task at hand, overcoming the lack of interpretability in conventional networks.

Compared to learned regularizers and plug-and-play methods, we do not only change the regularization, but can also deal with inaccuracies in the forward model $A$.
This allows to account for features that are missed during the modeling of the problem.

\subsection{Learned ISTA}

\begin{figure}[tbp]
	\centering
	\begin{tabular}{cc}
		\hspace{-2mm}
		\subfloat[Fixed point of ISTA.]{
			\includegraphics[width=0.3\linewidth]{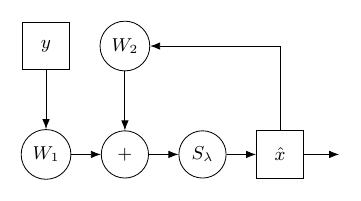}
			\vspace{0.3cm}
			\label{fig:rnn}
		}
		&
		\hspace{-4mm}
		\subfloat[Unrolled ISTA.
		]{
			\includegraphics[width=0.615\linewidth]{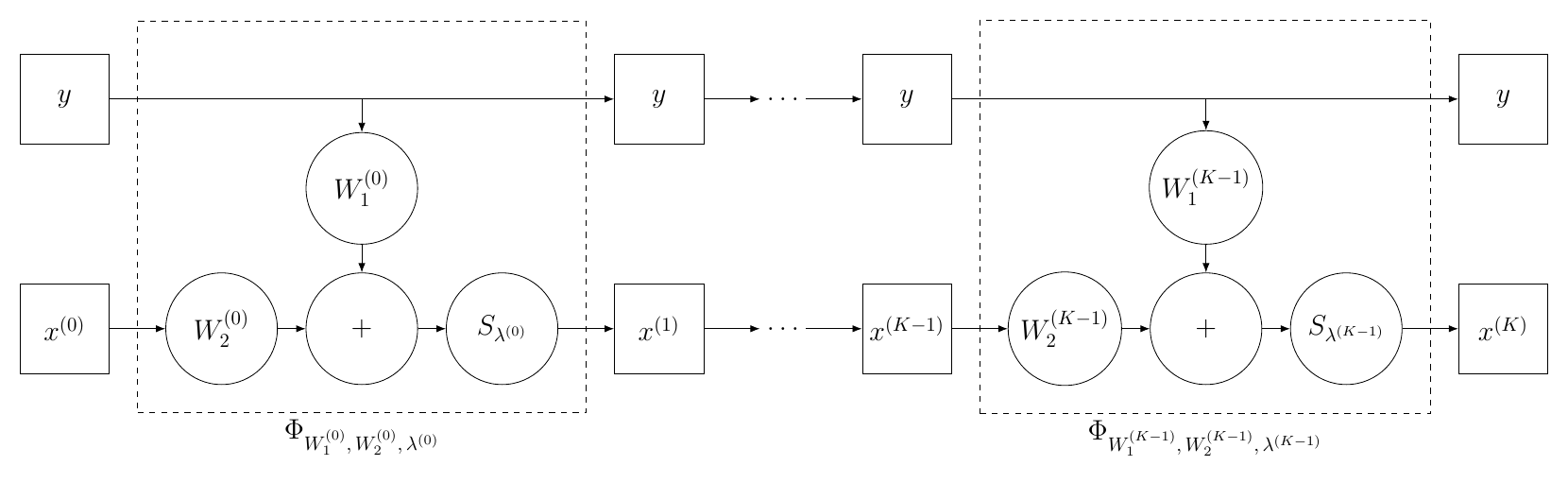}
			\label{fig:lista}
		}
	\end{tabular}
	\caption{Diagrams of ISTA and LISTA.}
	\vspace{-1em}
	\label{fig:rnn_lista}
\end{figure}

We study unrolling for the sparse recovery problem from \Cref{ex: LASSO}.
There, our goal is to solve
\begin{equation}\label{eq: sparse recovery}
	y = Ax^* + \varepsilon,
\end{equation}
where $x^* \in \R^N$ is the $s$-sparse ground truth vector, i.e., $\|x^*\|_0 \le s$, $A \in \R^{M \times N}$ is the measurement matrix with typically $M < N$, and $\varepsilon \in \R^M$ is noise.
In  \Cref{ex: LASSO}, we  proposed to recover $x$ from \eqref{eq: sparse recovery} by minimizing the LASSO functional
\begin{equation}\label{eq:Lasso}
	\min_{x} \tfrac{1}{2}\| Ax - y\|_2^2 + \alpha \|x\|_1, \qquad \alpha>0.
\end{equation}
In \Cref{ex:ista}, we derived for this purpose the ISTA
algorithm
\begin{equation}\label{eq:ISTA2}
x^{(k+1)} 
= S_{\alpha \nu} \bigl( (I - \nu  A^\tT A) x^{(k)} + \nu  A^\tT y \bigr),
\end{equation}
where $S_{\lambda}$ is the soft-thresholding function and $\nu>0$ is a step size.
As discussed in our motivation, we can interpret \eqref{eq:ISTA2} as a recurrent neural network with weights
\begin{equation}
W_1 = \nu A^\tT \in \R^{N \times M}, \quad W_2 = I - \nu A^\tT A \in \R^{N \times N}, 
\quad \text{and} \quad \lambda = \nu \alpha,
\end{equation}
see Figure~\ref{fig:rnn}. 
More precisely, for $K$ iterations of ISTA, we stack the layer $\Phi_{W_1, W_2, \lambda} (x,y)= (S_{\lambda}(W_1 y + W_2 x),y)$ $K$-times, leading to the neural network
\begin{equation}\label{eq: LISTA shared weights}
(x^{(k)},y) = \underbrace{\Phi_{W_1, W_2, \lambda} \circ \ldots \circ \Phi_{W_1, W_2, \lambda}}_{K\ \text{times}}(x^{(0)},y).
\end{equation}
As discussed in \Cref{sec: training problem}, the parameters $\theta = (W_1,W_2,\lambda)$ can be learned from training data $(x_i,y_i)\in \bb X \times \bb Y$, $i=1,\ldots,S$, based on the supervised learning problem 
\begin{equation}
	\hat \theta \in \argmin_{\theta} \biggl\{\mathcal L(\theta) = \frac{1}{2S}\sum_{i=1}^S
	\bigl\Vert x^{(K)}(\theta,y_i,x_0) -  x_i \bigr\Vert^2_{\bb X}\biggr\}.
\end{equation}
We provide some training details after investigating the properties of unrolled ISTA.

To make \eqref{eq: LISTA shared weights} more flexible, one can assign different parameters  $\theta^{(k)} = (W_1^{(k)}, W_2^{(k)}, \lambda^{(k)})$ to each layer.
The resulting iterations   
\begin{equation}\label{eq:gen_ista}
	x^{(k+1)} = S_{\lambda^{(k)}}\bigl( W^{(k)}_1 y + W^{(k)}_2 x^{(k)}\bigr), \quad k = 0,1,\ldots,K-1,
\end{equation}
are known as \emph{learned iterative soft thresholding algorithm} (\emphi{LISTA})  \cite{gregor2010learning, sprechmann2015learning}.
Empirically, it has significantly accelerated the recovery even for unobserved data as long as it follows the training distribution.
This leads to the question whether this behavior can be proven theoretically for sparse $x \in \R^N$.
Below, we give a brief overview.

First, LISTA can be viewed as inexact FBS algorithm and the convergence of the latter was investigated in \cite{giryes2018tradeoffs}.
By selecting an appropriate factorization of $A^\tT A$, it is possible to accelerate ISTA \cite{moreau2017understanding}, and extensions to other structured sets of parameters $\theta$ are possible.
The work \cite{chen2018theoretical} established conditions on the parameters $\theta$ if LISTA is expected to generalize to  all sparse $x \in \R^N$.
Moreover, the existence of $\theta$  leading to linear convergence was established \cite{chen2018theoretical,liu2019alista}.
In the following, we explore the last two results in more detail.

\subsubsection{Necessary Conditions for Convergence of LISTA}

In \Cref{sec: classical algorithms}, we established that for any $y \in \bb R^M$, (F)ISTA converges to a minimizer of the LASSO functional \eqref{eq:Lasso}.
Despite the empirical success of LISTA, it is unclear whether good performance can be guaranteed for all $x \in \R^N$.
Here, we derive necessary conditions for uniform convergence on the $s$-sparse set
\begin{equation}\label{eq: lista vector assumptions}
	\q X(B,s) \coloneqq \bigl\{ x \in \R^N  :  \| x \|_\infty \leq B, \| x\|_0 \leq s \bigr\}, \qquad s \geq 2.
\end{equation}
Given the index set $\q S \subseteq \{1, \ldots, N\}$, let $D_{\q S,:} \in \R^{|\q S| \times K}$ consist of the rows of  $D \in \R^{N \times K}$ with indices in $\q S$.
Analogously, we use the notation $x_{\q S}$ for vectors, $D_{:,\q S}$ for column extraction, and $D_{\q S, \q S}$ for the extraction of both rows and columns.
When $\q S = \{i\}$, we simply write $D_{i,:}$.
Next, we recall some reformulations of the soft thresholding operator.

\begin{lemma}\label{l: soft thresholding}
	For $S_\lambda$, $\lambda > 0$, as defined in \eqref{eq: soft thresholding} and $x \in \R$, we can rewrite $z = S_\lambda(x)$ as
	\begin{equation*}
	z = \sign(x) \max\{|x| - \lambda, 0\}  = x - \lambda \rho \quad \text{with} \quad \rho
	= \begin{cases}
		 \sign(z), & |x| > \lambda \Leftrightarrow z \neq 0, \\
		 (|x| / \lambda) \sign(x), & |x| \le \lambda \Leftrightarrow z = 0, 
	\end{cases}
	\end{equation*}
	where in any case $\rho \in [-1,1]$.
\end{lemma}

Now, we are ready to investigate LISTA.

\begin{thm}{Thresholding parameters}{lista necessary 1}
	Let $\{x^{(k)}\}_{k \in \N}$ denote the LISTA iterations \eqref{eq:gen_ista} with $\theta = \{W_1^{(k)}, W_2^{(k)}, \lambda^{(k)}\}_{k = 0}^\infty$, exact data $y = Az$, and initialization $x^{(0)} = 0$.
	If $\{\|W_2^{(k)}\|_2\}_{k = 0}^\infty$ is bounded, then the uniform convergence $\lim_{k \to \infty} \sup_{z \in \q X(B,s)} \Vert x^{(k)}(\theta,Az,0) - z\Vert_2 = 0$ implies that $\lim_{k \to \infty}\lambda^{(k)} = 0$.
\end{thm}

\begin{pr}
	We abbreviate $x^{(k)}(z) \coloneq x^{(k)}(\theta, Az, 0)$.
	As $\lim_{k \to \infty} \sup_{z \in \q X(B,s)} \Vert x^{(k)}(z) - z\Vert_2 = 0$, there is $k_1 \in \bb N$ with $\Vert x^{(k)}(z)-z\Vert_\infty < B/2$ for all $k \ge k_1$ and $z \in \q X(B,s)$.
	For $x^* =  \frac{B}{2} e_1 \in \q X(B/2,s)$, we get that
	\begin{equation}\label{eq: lista xik nonzero}
	x_1^{(k)}(x^*)= x^*_1 + x_1^{(k)}(x^*) - x^*_1\ge x^*_1 - |x_1^{(k)} (x^*)- x^*_1| > B/2 - B/2 = 0,
	\end{equation}
	which implies $x_1^{(k)}(x^*) >0$ for all $k \ge k_1$.
	Thus, \Cref{l: soft thresholding} with iteration \eqref{eq:gen_ista} yields
	\begin{equation}\label{eq:UpdateStep}
	x^{(k+1)}_{1}(x^*) = (W^{(k)}_2)_{1,:} \, x^{(k)}(x^*) + (W^{(k)}_1)_{1,:} \, A x^* - \lambda^{(k)}, \qquad k \geq k_1.
	\end{equation}
	Now, we rewrite \eqref{eq:UpdateStep} as 
	\begin{align}\label{eq:UpdateStep2}
			[ (I - W^{(k)}_2 - W^{(k)}_1 A ) x^* ]_1= - \lambda^{(k)} + \xi^{(k)}_1(x^*), \qquad k \geq k_1,
	\end{align}
	where
	\begin{equation}\label{eq:FormXi}
		\xi^{(k)}(z) = W^{(k)}_2 (x^{(k)}(z) -	z) - (x^{(k+1)}(z) - z) \in \R^N.
	\end{equation}
	With the same argumentation, we obtain for $2 x^* \in \q X(B,s)$ that
	\begin{equation}\label{eq:UpdateDoubled}
		[ (I - W^{(k)}_2 - W^{(k)}_1 A ) 2 x^* ]_1= - \lambda^{(k)} + \xi^{(k)}_1(2x^*), \qquad k \geq k_1.
	\end{equation}
	Subtracting \eqref{eq:UpdateDoubled} from two times \eqref{eq:UpdateStep2} gives $\lambda^{(k)} = 2 \xi^{(k)}_1(x^*) - \xi^{(k)}_1(2x^*)$ for $k \geq k_1$.
	Since the weights $\{\|W_2^{(k)}\|_2\}_{k \in \N}$ are bounded and $\lim_{k \to \infty} \sup_{z \in \q X(B,s)} \Vert x^{(k)}(z) - z\Vert_2 = 0$, we can find for every $\eta > 0$ some $k_2 \geq k_1$ such that $\| \xi^{(k)}(z)\|_2 \le \eta$ for all $k \ge k_2$ and $z \in \q X(B,s)$.
	Thus, we obtain for $k \ge k_2$ that
	\begin{equation}
	\lambda^{(k)} = |\xi^{(k)}_1(x^*) - \xi^{(k)}_1(2x^*)|
	\le 2 \|\xi^{(k)}(x^*)\|_2 + \|\xi^{(k)}(2x^*)\|_2
	\le 3\eta,
	\end{equation}
	which implies that $\lim_{k \to \infty}\lambda^{(k)} = 0$. 
\end{pr}

At first glance, \Cref{th:lista necessary 1} may seem counter intuitive:
For ISTA,  which is a special case of LISTA with fixed and untrained weights, we have chosen the constant sequence $\lambda^{(k)} = \nu \alpha$.
We can resolve this paradox as follows.
The minimizer of the LASSO functional \eqref{eq:Lasso} is in general not the solution of the inverse problem \eqref{eq: sparse recovery} due to the additional regularization term.
To obtain a $\|\cdot\|_1$-minimizing solution of  \eqref{eq: sparse recovery} (which may not be unique), we would need to choose vanishing regularization strength $\alpha^{(k)} \to 0$ as discussed in \Cref{sec: variational regularization}.
Indeed, a common strategy in sparse recovery is to start with an initial $\alpha^{(0)}$ and then decrease it as, e.g., $\alpha^{(k+1)} = 0.1 \alpha^{(k)}$.
Each time,  \eqref{eq:Lasso} is minimized again with the new $\alpha^{(k)}$ and $x^{(0)}$ being the previous reconstruction.
This choice matches the behavior established in \Cref{th:lista necessary 1}.

\begin{thm}{Weight coupling}{lista necessary 2}
	Consider the setup of \Cref{th:lista necessary 1}.
	Then, the weights $\{ W_1^{(k)}, W_2^{(k)} \}_{k = 0}^\infty$ have to satisfy
	\begin{equation}
	\lim_{k \to \infty} W_2^{(k)} - (I - W_1^{(k)} A) = 0.
	\end{equation}
\end{thm}

\begin{pr}
	Let $\q S \subset \{1,\ldots,N\}$ be a sparsity pattern with $\vert \q S \vert \le 2$, and $x^* \in \q{X}(B,2)$ with $x_i^* \neq 0$ if and only if $i \in \q S$.
	By \Cref{l: soft thresholding}, there exists $\rho^{(k)}(x^*) \in [-1,1]^N$ such that it holds
	\begin{equation}\label{eq:UpdateRepresntation}
		x^{(k+1)}_{\q S} (x^*)
		= (W^{(k)}_2)_{\q S,:} \, x^{(k)}(x^*) + (W^{(k)}_1)_{\q S,:} \, A x^* - \lambda^{(k)} \rho^{(k)}_{\q S}(x^*).
	\end{equation}
	Again, we rewrite \eqref{eq:UpdateRepresntation} as 
	\begin{align}\label{eq:ReformulationUpdate}
	[ (I - W^{(k)}_2 - W^{(k)}_1 A ) x^* ]_{\q S} = - \lambda^{(k)}\rho^{(k)}_{\q S}(x^*) + \xi^{(k)}_{\q S}(x^*),
	\end{align}
	where $\xi^{(k)}(z) \in \R^N$  is defined as in \eqref{eq:FormXi}.
	Inserting $W^{(k)}_2 \, x^* = (W^{(k)}_2)_{:,\q S} x^*_{\q S}$ and $W^{(k)}_1 A x^* = ( W^{(k)}_1 A)_{:, \q S} x^*_{\q S}$ into \eqref{eq:ReformulationUpdate} yields
	\begin{equation}\label{eq:ReformLISTA}
	\bigl( I - (W^{(k)}_2)_{\q S,\q S} - (W^{(k)}_1 A)_{\q S, \q S} \bigr) x^*_{\q S} =- \lambda^{(k)}\rho^{(k)}_{\q S}(x^*) + \xi^{(k)}_{\q S}(x^*).
	\end{equation}
	As shown in \Cref{th:lista necessary 1}, for every $\eta > 0$ there exists $k_1 \in \N$ such that $\| \xi^{(k)}(z)\|_2 \le \eta$ for all $k \ge k_1$ and $z \in \q X(B,s)$.
	Thus, taking the norm on both sides of \eqref{eq:ReformLISTA} gives
	\begin{equation}
	\bigl\| \bigl( I - (W^{(k)}_2)_{\q S,\q S} - (W^{(k)}_1 A)_{\q S, \q S} \bigr) x^*_{\q S} \bigr\|_2 
	\le \|\xi^{(k)}_{\q S}(x^*)\|_2 + \lambda^{(k)}\Vert \rho^{(k)}_{\q S}(x^*)\|_2
	\le \eta + \lambda^{(k)} \sqrt{2}.
	\end{equation}
	By \Cref{th:lista necessary 1}, $\lim_{k \to \infty} \lambda^{(k)} = 0$ and there exists $k_2 \geq k_1$ with $\lambda^{(k)} \le \eta$ for $k \geq k_2$.
	For all $k \ge k_2$, we thus obtain
	\begin{equation}\label{eq:IneqNorm}
	\bigl\| \bigl( I - (W^{(k)}_2)_{\q S,\q S} - (W^{(k)}_1 A)_{\q S, \q S} \bigr) x^*_{\q S} \bigr\|_2 
	\le (1 + \sqrt{2}) \eta.
	\end{equation}
	Using \eqref{eq:IneqNorm}, we bound $\| I - (W^{(k)}_2)_{\q S,\q S} - (W^{(k)}_1 A)_{\q S, \q S} \| _2$ as
	\begin{align}
		\bigl\| I - (W^{(k)}_2)_{\q S,\q S} - (W^{(k)}_1 A)_{\q S, \q S} \bigr \|_2
		&= \sup_{ \substack{\| x\|_2 = 1,\\
		\supp(x) \subseteq \q S} } \bigl\| \bigl( I - (W^{(k)}_2)_{\q S,\q S} - (W^{(k)}_1 A)_{\q S, \q S} \bigr) x_{\q S} \bigr\|_2 \notag\\
		& \le B^{-1} (1 + \sqrt{2}) \eta.  
	\end{align}
	Consequently, $\lim_{k \to \infty} I - (W^{(k)}_2)_{\q S,\q S} - (W^{(k)}_1 A)_{\q S, \q S} = 0$ for all sparsity patterns $\q S$ with $\vert \q S \vert = 2$.
	Hence, it holds that $\lim_{k \to \infty} [I-W^{(k)}_2-W^{(k)}_1 A] = 0$.
\end{pr}

\Cref{th:lista necessary 2} implies that the weights have to asymptotically satisfy $W_2^{(k)} = I - W^{(k)}_1 A$.
Renaming $W_1^{(k)}$ as $W^{(k)}$ and substituting $W_2^{(k)} = I - W^{(k)} A$ into \eqref{eq:gen_ista} gives
\begin{equation}\label{eq: lista cp}
	x^{(k+1)} 
	= S_{\lambda^{(k)}} \bigl( x^{(k)} - W^{(k)} (A x^{(k)} - y) \bigr),\quad k = 0,1,\cdots,K-1.
\end{equation}
Here, we can interpret the argument of $S_{\lambda^{(k)}}$ as a gradient descent step with respect to a changing inner product (variable metric methods).
The iterations \eqref{eq: lista cp} are referred to as \emphi{LISTA-CP} shortly for \emph{LISTA with coupled parameters}. 

\subsubsection{Linear Convergence of LISTA}

Following \cite{chen2018theoretical}, we show that LISTA-CP achieves linear convergence if the weights $W$ are chosen appropriately.
For this, we require the notion of coherence.
Throughout, we assume that the columns of $A$ are normalized, namely that
\begin{equation}\label{eq: LISTA normalization}
\| A_{:,j} \|_2 = 1, \quad \text{for all} \quad  j = 1, \ldots, N.  
\end{equation}

\begin{defn}
		
	The \emphi{generalized mutual coherence} $\mu \colon \R^{M \times N} \to \R$ is defined as
	\begin{equation} \label{eq:mutual_coher}
		\mu(A) = \inf_{\substack{W \in \R^{N \times M} \\ W_{i,:}  A_{:,i} = 1,1 \leq i \leq N} } \Bigl\{ \max_{1 \leq i \neq j \leq N }| W_{i,:} A_{:,j} |\Bigr\}.
	\end{equation}
\end{defn}

We first show that $\mu(A)$ is well-defined.

\begin{lemma}\label{lemma:mucoher}
	The infimum in \eqref{eq:mutual_coher} is attained for some $W \in \R^{N \times M}$. 
\end{lemma}
\begin{pr}
	 The problem \eqref{eq:mutual_coher} is feasible since $A^\tT$ satisfies the constraints due to \eqref{eq: LISTA normalization}.
	 Moreover,  it holds that $\mu(A) \geq 0$, and we can rewrite \eqref{eq:mutual_coher} as the linear program
	 \begin{equation}\label{eq:LinearProgram}
	\begin{aligned}
		\min_{W, t} \quad & t \\
		\text{subject to} \quad & W_{i,:}A_{:,i} = 1, \quad \text{for all } i = 1,\ldots, N, \\
		& -t \leq W_{i,:}A_{:,j} \leq t, \quad \text{for all } i \neq j.
	\end{aligned}
	\end{equation}
	 Thus, by the fundamental theorem of linear programming, see for example \cite[Prop.\ 3.3]{bertsimas1997introduction}, there exists a minimizer of \eqref{eq:mutual_coher}. 
\end{pr}

We denote the minimizers of \eqref{eq:mutual_coher} by 
\begin{equation}
	\q W(A) \coloneqq \bigl\{ W \in \R^{N \times M} : W \ \text{minimizes \eqref{eq:mutual_coher}} \bigr\} \neq \emptyset.
\end{equation}

Now, we establish that there exist parameters $\{ W^{(k)}, \lambda^{(k)} \}_{k = 0}^{K-1}$ so that \eqref{eq: lista cp} achieves uniform linear convergence on the set $\q X(B,s)$. 

\begin{thm}{Convergence of LISTA-CP}{convergence lista}
	Let $A$ satisfy \eqref{eq: LISTA normalization} as well as $(2s -1) \mu(A) <1$.
	Given a sequence $\{\gamma^{(k)}\}_{k \in N}$ with  $0 < \gamma^{(k)} \leq \frac{2}{ (2s - 1)\mu(A) + 1}$ and  a matrix $W \in \q W(A)$, we define the weights $W^{(k)} = \gamma^{(k)} W$.
	Moreover, for a noise level $\sigma>0$ and $C_W = \max_{i,j} |W_{i,j}|$, we define $\lambda^{(k)}$ recursively, depending on the current worst-case error of the iterate $x^{(k)}$ generated by \eqref{eq: lista cp}, as
	\begin{equation}\label{eq: lista weights analytical}
	\lambda^{(k)} = \gamma^{(k)}\bigg[  \mu(A) \sup_{\substack{ z \in \q X(B,s) \\ \varepsilon \in \R^M, \|\varepsilon\|_1 \le \sigma}} \| x^{(k)}(\theta, A z + \varepsilon, 0) - z \|_1  + C_W \sigma \bigg].
	\end{equation}
	Then,  
	for $x^* \in \q X(B,s)$ and $\varepsilon \in \R^M$ with $\| \varepsilon \|_1 \le \sigma$, we have for the iterates $x^{(k)} \coloneq x^{(k)}(\theta, A x^* + \varepsilon, 0)$ that
	\begin{equation}\label{eq:RestrictedSupport}
		\supp(x^{(k)}) \subseteq \supp(x^*)
	\end{equation} 
	and
	\begin{equation}\label{eq:ALISTAConvEst}
		\| x^{(K)} - x^* \|_1 \leq s B \exp\biggl(- \sum_{k = 0}^{K-1}  c^{(k)}\biggr)  + 2 \sigma s C_W \sum_{k = 0}^{K-1} \gamma^{(k)} \exp\biggl(- \sum_{l = k+1}^{K-1} c^{(l)}\biggr),
	\end{equation}
	where the $c^{(k)}>0$ depend only on $s$,  $\gamma^{(k)}$ and $\mu(A)$. 
\end{thm}

\begin{pr}
	Let $x^* \in \q X(B,s)$ with $\q S = \supp(x^*)$ and $\varepsilon \in \R^M$ with $\|\varepsilon\|_1 \le \sigma$.
	Given the noisy observation $y = A x^* + \varepsilon$, we write $x^{(k)} \coloneq x^{(k)}(\theta, y , 0)$ for brevity.
	
	\paragraph{Support property.}
	We prove \eqref{eq:RestrictedSupport} by induction.
	Since $x^{(0)} = 0$, \eqref{eq:RestrictedSupport} holds for $k=0$.
	Now, assume that $\supp(x^{(k)}) \subseteq \q S$.
	For $i \notin \q S$, $x^{(k+1)}_i$ is given according to \eqref{eq: lista cp} as
	\begin{align}\label{eq:UpdateSoft}
		x^{(k+1)}_i \!= \!S_{\lambda^{(k)}} \bigl(x^{(k)}_i - W^{(k)}_{i,:} (A x^{(k)} - y) \bigr)\!= \!S_{\lambda^{(k)}} \Bigl(  -\! \sum_{j \in \q S} W^{(k)}_{i,:} A_{:,j} (x^{(k)}_j - x^*_j) + W^{(k)}_{i,:} \varepsilon \Bigr),
	\end{align}
	where we used the induction hypothesis $x^{(k)}_i = 0$ and $y = \sum_{j \in \q S} A_{:,j}  x^*_j  + \varepsilon$.
	Since $W^{(k)} = \gamma^{(k)} W$ and $W \in \q W(A)$, we further  have that
	\begin{align}\label{eq:EstSupport}
	\biggl| - \sum_{j \in \q S} W^{(k)}_{i,:} A_{:,j} (x^{(k)}_j - x^*_j) + W^{(k)}_{i,:} \varepsilon \biggr| 
	& \le \sum_{j \in \q S} \gamma^{(k)} |W_{i,:} A_{:,j}| |x^{(k)}_j - x^*_j| + \gamma^{(k)} |W_{i,:} \varepsilon| \notag\\
	& \le \gamma^{(k)}  \mu(A) \sum_{j \in \q S} |x^{(k)}_j - x^*_j| + \gamma^{(k)} \Vert W_{i,:}\Vert_\infty \Vert \varepsilon\Vert_1 \notag\\
	& \le \gamma^{(k)} \mu(A) \| x^{(k)} - x^*\|_1 + C_W \gamma^{(k)} \|\varepsilon\|_1 \le \lambda^{(k)}.
	\end{align}
	Thus, the update \eqref{eq:UpdateSoft} becomes $x^{(k+1)}_i = 0$ for all $i \notin \q S$, which implies $\supp(x^{(k+1)}) \subseteq \q S$.
	
	\paragraph{Error bound for $x^*$ and $\varepsilon$.}
	
	By Step 1, it suffices to consider $i \in \q S$.
	Using \Cref{l: soft thresholding}, we find $\rho^{(k)} \in [-1,1]^N$ such that
	\begin{align}
		x^{(k+1)}_i & = x^{(k)}_i - \sum_{j \in \q S} W^{(k)}_{i,:} A_{:,j} (x^{(k)}_j - x^*_j) + W^{(k)}_{i,:} \varepsilon - \lambda^{(k)} \rho^{(k)}_i \notag\\
		& = x^{(k)}_i - \gamma^{(k)} \sum_{j \in \q S} W_{i,:} A_{:,j} (x^{(k)}_j - x^*_j) + \gamma^{(k)} W_{i,:} \varepsilon - \lambda^{(k)} \rho^{(k)}_i .
	\end{align}
	Since $W_{i,:} A_{:,i} = 1$ by construction, we further get
	\begin{align}\label{eq:DiffIterates}
		x^{(k+1)}_i - x^*_i  = (1- \gamma^{(k)}) (x^{(k)}_i - x^*_i) - \gamma^{(k)} \!\!\!\!\sum_{j \in \q S, j \neq i}\!\!\! W_{i,:} A_{:,j} \bigl(x^{(k)}_j - x^*_j\bigr) + \gamma^{(k)} W_{i,:} \varepsilon - \lambda^{(k)} \rho^{(k)}_i.
	\end{align}
	Similarly as for \eqref{eq:EstSupport}, taking the absolute value on both sides of \eqref{eq:DiffIterates} leads to
	\begin{align}\label{eq:EstSingleEntery}
		|x^{(k+1)}_i \!- x^*_i|  \le |1- \gamma^{(k)}| | x^{(k)}_i \!- x^*_i | + \gamma^{(k)} \mu(A) \!\!\!\!\sum_{j \in \q S, j \neq i}\!\! |x^{(k)}_j \!- x^*_j| 
		+ \gamma^{(k)} C_W \sigma + \lambda^{(k)}.
	\end{align}
 	Summing over $i \in \q S$ and using that $\vert \q S \vert \leq s$ and $\supp(x^{(k)}) \subseteq \q S$ for all $k \in \N$, we arrive at
 	\begin{equation}\label{eq:EstSingleEntery2}
 		\| x^{(k+1)} - x^* \|_1   \le \bigl(|1- \gamma^{(k)}| + (s - 1)\gamma^{(k)}  \mu(A) \bigr) \| x^{(k)} - x^* \|_1 + s \gamma^{(k)} C_W \sigma + s \lambda^{(k)}.
 	\end{equation}
 For brevity, we denote the supremum over $x^* \in \q X(B,s)$ and $ \varepsilon \in \R^M, \|\varepsilon\|_1 \le \sigma$ with $\sup_{x^*, \varepsilon}$.
 Since $\lambda^{(k)}$ is chosen according to \eqref{eq: lista weights analytical}, we can further estimate \eqref{eq:EstSingleEntery2} as
 	\begin{equation}\label{eq: lista cp proof}
 		\|x^{(k+1)}-x^*\|_1
 		\le \bigl(|1- \gamma^{(k)}| + (2s - 1) \gamma^{(k)} \mu(A)\bigr) \sup_{x^*, \varepsilon} \|x^{(k)}-x^*\|_1 + 2s \gamma^{(k)} C_W \sigma.
 	\end{equation}
 	By assumption, it holds $(2s - 1)  \mu(A) < 1$, which for $\gamma^{(k)} \leq 1$ leads to
 	\begin{equation}
 		|1- \gamma^{(k)}| + (2s - 1) \gamma^{(k)} \mu(A) 
 		<1- \gamma^{(k)} + \gamma^{(k)} = 1.
 	\end{equation}
	For $1 < \gamma^{(k)} < \frac{2}{ \mu(A) (2s - 1) + 1}$, we also obtain
 	\begin{equation}
 	|1- \gamma^{(k)}| + \gamma^{(k)} \mu(A) (2 s - 1) = \gamma^{(k)} ( \mu(A) (2 s - 1) + 1)  - 1 < 2 -1 = 1.
 	\end{equation}
 	Hence, we can define
 	\begin{equation}\label{eq:DefCk}
 		c^{(k)} \coloneqq - \log\bigl(|1- \gamma^{(k)} | + \gamma^{(k)}  \mu(A) (2 s - 1) \bigr) > 0.
 	\end{equation}
 	With the constant $c^{(k)}$ and by taking the supremum over $x^* \in \q X(B,s)$ and $\varepsilon \in \R^M$ with $\|\varepsilon\|_1 \le \sigma$, we can rewrite \eqref{eq: lista cp proof} more compactly as
 	\begin{equation}\label{eq:EstimateSingle}
 		 \sup_{x^*, \varepsilon} \|x^{(k+1)}-x^*\|_1
 		\le \exp(-c^{(k)}) \sup_{x^*, \varepsilon}  \|x^{(k)}-x^*\|_1 + 2\sigma s C_W  \gamma^{(k)} .
 	\end{equation}
 	Now, applying \eqref{eq:EstimateSingle} successively for $K-1$ times and combining $x^{(0)} = 0$ with $\|x^*\|_1 \le |\q S| \| x^*\|_\infty \le s B$ for all $x^* \in \q X(B,s)$ gives the claim
 	\begin{align}
 		\|x^{(K)}-x^*\|_1 
 		& \le \prod_{k = 0}^{K-1}\exp(-c^{(k)}) \sup_{x^*, \varepsilon}  \|x^{(0)}-x^*\|_1 + 2 s \sigma C_W \sum_{k = 0}^{K-1} \gamma^{(k)} \prod_{l = k+1}^{K-1} \exp(-c^{(l)}) \nonumber \\
 		& = \exp\biggl(- \sum_{k = 0}^{K-1} c^{(k)}\biggr) s B + 2 s \sigma C_W \sum_{k = 0}^{K-1} \gamma^{(k)} \exp\biggl(- \sum_{l = k+1}^{K-1} c^{(l)}\biggr). \label{eq:EstCompleteIter}
 	\end{align}
 	This concludes the proof.
\end{pr}

Now, we discuss the implications of  \Cref{th:convergence lista}.
First, we observe that linear convergence can only occur if $A$ is mutually coherent, namely if $A$ satisfies $(2s -1) \mu(A) <1$, see also the lower bounds in \cite{liu2019alista}.
Moreover, to recover non-sparse solutions, we need higher mutual coherence of $A$, which can be seen as a form of regularity requirement.

For the noiseless case $\varepsilon = 0$, the second summand in \eqref{eq:ALISTAConvEst} vanishes, resulting in
\begin{equation}\label{eq:LinearConvLista}
\| x^{(k)} - x^* \|_1 \le s B \exp\biggl(- \sum_{k = 0}^{K-1} c^{(k)}\biggr).
\end{equation}
For $\gamma^{(k)} = 1$, we obtain from the defining equation  \eqref{eq:DefCk} of $c^{(k)}$ that $c^{(k)} = c > 0$, and \eqref{eq:LinearConvLista} turns into a linear convergence rate.
However, we cannot compute $ \lambda^{(k)} $ explicitly as \eqref{eq:ALISTAConvEst} includes a maximum. 
In particular, \Cref{th:convergence lista} only guarantees that an admissible schedule for linear convergence exists.
During learning, we can then try to approximate this schedule based on the data.

In the general noisy case with $\sigma>0$, this argument does not necessarily hold anymore.   
In particular, we only obtain a bound of the form 
\begin{align}
	\| x^{(K)} - x^* \|_1 \leq& s B \exp (- K c)  + 2 \sigma s C_W \sum_{k = 0}^{K-1} \exp\bigl(- (K-k-1) c \bigr)\notag \\
	\leq & s B \exp (- K c)  + \sigma \frac{2s C_W}{1 - \exp(-c)} ,
\end{align}
which means that we can only expect linear convergence behavior until we reach the magnitude of the noise, see also the discussion in \cite{verdun2024fast}. 
The key equation capturing this behavior in the proof of  \Cref{th:convergence lista} is \eqref{eq:EstimateSingle}, and we would need to have $\lim_{k \to \infty} \gamma^{(k)} = 0$ in order to suppress the second summand. 
Thus, to improve the convergence behavior under the presence of noise, we can try to also learn a non-constant $\gamma^{(k)}$.

In summary, the main implication of \Cref{th:convergence lista} is that we can analytically select $W^{(k)} = \gamma^{(k)} W$ for any $W \in \q W(A)$.
Consequently, it remains to learn $\gamma^{(k)}$ and $\lambda^{(k)}$.
Thus, we consider \emph{analytic LISTA} (\emphi{ALISTA}) given by
\begin{equation}\label{eq: alista}
	x^{(k+1)} = S_{\lambda^{(k)}} \bigl( x^{(k)} - \gamma^{(k)} W (A x^{(k)} - y) \bigr), \quad k = 0,1,\ldots,K-1,
\end{equation}
with learnable parameters $\theta = \{ \gamma^{(k)}, \lambda^{(k)} \}_{k =0}^{K-1}$.
The empirical performance of ALISTA matches LISTA and LISTA-CP, while involving only $2K$ parameters.
Nevertheless, we need to find $W \in \q W(A)$ numerically before the training. 
Miminimizing \eqref{eq:mutual_coher} directly is challenging due to $\max_{i \neq j} |W_{i,:} A_{:,j}|$.
A feasible alternative is to minimize the majorization 
\begin{equation}\label{eq:MajorizationMU}
\mu(A) \le \inf_{\substack{W \in \R^{N \times M} \\ W_{i,:}  A_{:,i} = 1,1 \leq i \leq N} } \Bigl\{ \Bigl( \sum_{1 \leq i \neq j \leq N }| W_{i,:} A_{:,j} |^2 \Bigr)^{1/2} \Bigr\}.
\end{equation}
With the convex and closed constraint set $\q C(A) \coloneqq \{ W \in \R^{N \times M}  :  W_{i,:}  A_{:,i} = 1,1 \leq i \leq N \}$, we rewrite the right hand side of \eqref{eq:MajorizationMU} as
\begin{equation}
W \in \argmin_{W \in \q C(A)} \Bigl( \sum_{1 \leq i \neq j \leq N} | W_{i,:} A_{:,j} |^2 \Bigr)^{1/2} 
= \argmin_{W \in \R^{N \times M}} \sum_{i,j =1}^{N} | W_{i,:} A_{:,j} |^2 + \iota_{\q C(A)}(W). 
\end{equation}
This problem can be minimized with FBS algorithm, see \Cref{exampleB:proxofl1norm}.
Alternatively, we can solve the (possibly huge) linear program \eqref{eq:LinearProgram}.

\subsubsection{Training Strategy}

For training any of the versions of LISTA, it is benefical to perform the training in stages, where we increase the number of steps until we reached the desired number $K$.
By $\theta^{(k)}$ we denote the parameters required for $k$ iterations of LISTA.
Then, the training strategy is the following. 
For $k=0, \ldots, K-1$, we 
\begin{enumerate}
	\item initially only train $\theta^{(k)} \backslash \theta^{(k-1)}$, i.e., the parameters corresponding to $k$th update.
	The parameters of the first $k-1$ updates stay fixed;
	\item then train all parameters $\theta^{(k)}$ with a decreasing learning rate;
	\item proceed to stage $k+1$.   
\end{enumerate}
This procedure empirically improves the results, especially for a large depth $K$.

\subsubsection{Exercises}

\begin{exerciseB}
	Implement and train LISTA, LISTA-CP and ALISTA. Compare their performance to ISTA and FISTA.
	
	You may use the following Jupyter Notebook as a basis for this exercise, \url{https://github.com/Sniper2k/inverse_problems_ai_exercise_templates/blob/main/Ex8_1_7_template.ipynb}.
\end{exerciseB}

\begin{Solution}
See Jupyter Notebook \url{https://github.com/Sniper2k/inverse_problems_ai_exercises/blob/main/Ex8_1_7.ipynb}.
\end{Solution}

\subsection{Outlook for learning-based methods}

To summarize our introduction to learned methods for inverse problems, we look back at the discussed techniques: denoising and learned post-processing in \Cref{chap:regularization}, learned regularizers in \Cref{sec:LearnedReg}, plug-and-play methods in \Cref{sec: PnP}, and unrolling in \Cref{sec: unrolling}. 

We started with learned post-processing that acts as a corrector on top of a baseline reconstruction, removing the artifacts associated with algorithm and the noise.
This is, in certain sense a "sharp denoising": 
we get the result and apply the denoiser once to obtain the final reconstruction.
 In contrast, PnP methods use "soft denoising": we gradually denoise the reconstructed object, steering it towards a more regular solution.

We also observed that by taking the denoiser as a gradient of a (strongly) convex $C^2$-function, we can interpret it as an implicitly defined regularizer that can be computed numerically.
On the other hand, direct bi-level learning of an explicit regularizer allows for a more controlled choice of the underlying architecture. 

The last technique, namely unrolling, rather focuses on improving an optimization algorithm of choice. Due to introduction of separate parameters for each update, we essentially use a different regularizer/denoiser each time.
While we lose the interpretability as actual optimization of a fixed loss, it avoids the complexity of bi-level learning associated with gradient computation.
 Furthermore, the full power of unrolling lies in the ability to learn the actual optimization algorithm.

In this book, we left out Bayesian approaches such as diffusion \cite{lu2026mathematical} and flow matching \cite{wald2025flow}.
These rely on the generation of samples for simple reference distribution and transforming them to the data distribution using the optimal transport techniques.
In the context of inverse problems, they can be used as denoisers \cite{martin2025pnp}.
Another technique prevalent in inverse problems is measurement denoising, see for instance \cite{wu2024blind}.
Its roots stem from the fact that ground truth objects or high-quality reconstructions are not available in many settings.
Yet, relatively clean and noisy measurements are accessible, allowing to train a denoiser. As the measurements strongly depend on experimental metadata, the challenge is to incorporate it in the models architecture.    

From the perspective of theoretical results that we discussed, augumenting the data can lead to large errors, see \eqref{eq: SVD error scaling}. On the other hand, if denoiser acts as a projection to the range of the forward map, the error rates for spectral methods improve according to \Cref{th:spektral:apriori}. In analogy to the learned post-processing and PnP, we can also include gradual measurement denoising through ADMM framework \cite{Nguyen.2026}.

In all of the above settings, we require supervised datasets relying on known objects or clean measurements.
If only noisy measurement are available, it is still possible to train a denoiser using semi-supervised learning techniques \cite{tachella2026self}, however this generally comes with a few caveats.
For an overview of denoising methods, see the extensive overview \cite{elad2023image}.

%%%%%%%%%%%%%%%%%%%%%%%%%%%%%%% Index %%%%%%%%%%%%%%%%%%%%%%%%%%%%%%%%%%%%%%%%%%
%%%%%%%%%%%%%%%%%%%%%%%%%%%%%%%%%%%%%%%%%%%%%%%%%%%%%%%%%%%%%%%%%%%%%%%%%%%%%%%%
\addcontentsline{toc}{section}{Index}
\printindex

\clearpage
%% References
\bibliographystyle{alpha}
\bibliography{Bibliography}
            
   \clearpage
            
\end{document}